\documentclass[twoside,12pt]{article}
\usepackage{amsmath, amssymb, epsfig}

\newtheorem{theorem}{Theorem}[section]
\newtheorem{proposition}[theorem]{Proposition}
\newtheorem{lemma}[theorem]{Lemma}
\newtheorem{corollary}[theorem]{Corollary}

\newtheorem{definition}[theorem]{Definition}
\newtheorem{remark}[theorem]{Remark}
\newenvironment{proof}{\par\smallskip\noindent{\bf Proof.}\/}{\quad\par\medskip}

\newcommand{\RR}{{\rm\bf R}}

\newcommand{\BB}{ {\rm\bf B}}
\newcommand{\CC}{{\rm\bf C}}

\newcommand{\NN}{{\rm\bf N}}
\newcommand{\ZZ}{{\rm\bf Z}}

\newcommand{\qed}{\rule{0cm}{0cm} \hfill \rule{2mm}{2mm}\rule{0.1cm}{0cm}}
\newcommand{\str}{\mathop{\,\rm str}\nolimits}
\newcommand{\so}{\mathop{\,\rm\bf so}\nolimits}

\newcommand{\tr}{\mathop{\,\rm tr}\nolimits}

\newcommand{\RC}{\mathop{\rm RC}\nolimits}

\newcommand{\Str}{\mathop{\, \rm Str}\nolimits}
\newcommand{\End}{\mathop{\, \rm End}\nolimits}
\newcommand{\Id}{\mathop{\, \rm Id}\nolimits}
\newcommand{\REG}{\mathop{\,\rm reg}}
\newcommand{\RES}{\mathop{\,\rm res}}

\newcommand{\END}{\mathop{\,\rm END}\nolimits}
\newcommand{\Cl}{\mathop{\,\rm Cl}\nolimits}

\newcommand{\llangle}{\langle\!\langle}
\newcommand{\rrangle}{\rangle\!\rangle}

\newcommand{\sgn}{\mathop{\,\rm sgn}\nolimits}
\newcommand{\Diff}{\mathop{\,\rm Diff}\nolimits}

\newcommand{\ind}{\mathop{\,\rm ind}\nolimits}

\newcommand{\Ch}{\mathop{\,\rm Ch}\nolimits}

\newcommand{\cdotinfty}[1]{\stackrel{\mbox{\,\Huge .}}{C}\!{}^{\infty}_{#1}}
\newcommand{\delx}{\frac{\partial}{\partial x}}
\newcommand{\dely}[1]{\frac{\partial}{\partial y^{#1}}}
\newcommand{\delz}[1]{\frac{\partial}{\partial z^{#1}}}
\newcommand{\delyy}{\frac{\partial}{\partial {\bf y}}}
\newcommand{\delzz}{\frac{\partial}{\partial {\bf z}}}
\newcommand{\phiv}{{}^{\phi}{\cal V}}
\newcommand{\phiomega}{{}^{\phi}{\Omega}}
\newcommand{\bomega}{{}^{b}{\Omega}}
\newcommand{\CB}{{\rm CB}}
\newcommand{\KD}{{\rm KD}}
\newcommand{\omfib}{\Omega_{\rm fib}}
\newcommand{\phitx}{{}^{\phi}{TX}}

\newcommand{\psisus}{\Psi_{\rm sus}}
\newcommand{\phitsx}{{}^{\phi}{T^*X}}

\newcommand{\phidir}{ {\sf D}^{\phi}}
\newcommand{\phidiff}{{\mbox{Diff}_{\phi}}\,}
\newcommand{\phindx}{{}^{\phi}\! N\partial X}
\newcommand{\phinsdx}{{}^{\phi}N^*\partial X}
\newcommand{\dx}{d\, x}
\newcommand{\phicl}{{\bf c}_{\phi}}
\newcommand{\psiphi}{\Psi_{\phi}}

\newcommand{\fahbf}{{\rm hbf}}
\newcommand{\fahlf}{{\rm lbf}}
\newcommand{\fahrf}{{\rm rbf}}
\newcommand{\faphibf}{\phi{\rm bf}}
\newcommand{\fatff}{{\rm tff}}
\newcommand{\fatf}{{\rm tf}}
\newcommand{\fatb}{{\rm tb}}
\newcommand{\fabf}{{\rm bf}}
\newcommand{\faff}{{\rm ff}}
\newcommand{\farf}{{\rm rf}}
\newcommand{\falf}{{\rm lf}}
\newcommand{\nnabla}{\mbox{\boldmath$\nabla$}}

\newcommand{\dd}{{\,\it d}}
\newcommand{\Dir}{{\sf D}}

\newcommand{\spec}{\mathop{\rm spec}\nolimits}

\newcommand{\cl}{{\bf c}}

\newcommand{\dvol}{\mathop{\it dvol}\nolimits}

\newcommand{\Norm}[1]{\left\|#1\right\|}
\newcommand{\Betrag}[1]{\left|#1\right|}

\newcommand{\sowie}{\mathop{\,\stackrel{\wedge}{=}\,}\nolimits}

\newcommand{\ext}{\mathop{\mbox{\boldmath $\varepsilon$}}\nolimits}
\newcommand{\eps}{\mathop{\mbox{\boldmath $\epsilon$}}\nolimits}
\newcommand{\Int}{\mathop{\mbox{\boldmath $\iota$}}\nolimits}

\renewcommand{\Re}{\mathop{\,\rm Re}\nolimits}
\renewcommand{\Im}{\mathop{\,\rm Im}\nolimits}

\newcommand{\Const}{\mathop{\rm C}\nolimits}

\newcommand{\zwzw}[4]{\begin{pmatrix} #1 & #2 \\ #3 & #4
\end{pmatrix}}

\begin{document}
\pagenumbering{roman}
\title{Index- and Spectral Theory for Manifolds with Generalized Fibred Cusps
\footnote{ This is the author's doctoral thesis submitted at the University of Bonn.
Written with support
from the DFG through funds from the SFB 256 and a fellowship
from the ``Graduiertenkolleg
f\"ur Mathematik der Universit\"at Bonn''. Financial support was
also obtained from the DAAD through a
``HSP III''-scholarship.}}
\author{Boris Vaillant \\
 Department of Mathematics, University of Bonn}
\maketitle
\par
\bigskip
\begin{abstract}{\small
Generalizing work of W. M\"uller
we investigate the spectral theory for the Dirac
operator $\Dir$ on a  noncompact manifold $X$ with generalized fibred cusps
$$ C(M)=M\times [A,\infty[_r, \qquad g=\dd r^2+\phi^*g_Y+e^{-2cr}g_Z,$$
at infinity. Here
$\phi:M^{h+v}\rightarrow Y^h$ is a compact fibre bundle with fibre $Z$ and
a distinguished horizontal space $HM$.  The metric $g_Z$ is a metric in the fibres
and $g_Y$ is a metric on the base of the fibration.
We also assume that the kernel of the vertical Dirac operator at infinity
forms a vector bundle over $Y$.\par
Using the ``$\phi$-calculus'' developed by R. Mazzeo and R. Melrose
  we explicitly construct the meromorphic
continuation of the resolvent $G(\lambda)$ of $\Dir$ for small
spectral parameter as a special ``conormal distribution''. From this
 we deduce a description of the generalized
eigensections and of the spectral measure of $\Dir$.
\par
Complementing  this, we perform an explicit construction of the
heat kernel $[\exp(-t\Dir^2)]$ for finite and small times $t$,
corresponding to large
spectral parameter $\lambda$. Using a generalization of
Getzler's technique, due to R. Melrose, we can describe the singular terms in the heat
kernel expansion for small times in the interior of the manifold
as well as at spatial infinity. This then allows us to  prove
an index formula for $\Dir$,
$${\rm ind}_-(\Dir)= \frac{1}{(2\pi
i)^{n/2}}\int_X\widehat{A}(R)\Ch(F^{E/S})+\frac{1}{(2\pi
i)^{(h+1)/2}}\int_Y\widehat{A}(R^Y)\widehat{\eta}(\Dir^{V})+\frac{1}{2}
\eta(\Dir_Y),$$  which calculates the extended
$L^2$-index of $\Dir$, in terms of the usual local expression, the family eta invariant
for the family of vertical Dirac operators at infinity and the eta invariant
for the horizontal ``Dirac'' operator at infinity.}
\end{abstract}
\par
\medskip
\newpage
\tableofcontents
\newpage
\pagenumbering{arabic}

\section*{Introduction}
\addcontentsline{toc}{section}{\numberline{}Introduction}
\markboth{INTRODUCTION}{}
In his book \cite{MuellerBook} W. M\"uller, using methods from
harmonic analysis, described the spectral decomposition  and index
theory of the Dirac operator over a locally symmetric space $$
\Gamma \setminus G/ K $$ of rank one. This is a manifold equipped
with a natural riemannian metric. Near infinity it has the
structure of a ``fibred cusp'' $$ C(M)=M\times [A,\infty[_r, $$
where $\phi:M\rightarrow Y$ is a compact fibre bundle with fibre
$Z$ and a distinguished horizontal space $HM$. The metric on the
cusp is given by
\begin{equation}\label{cuspmetric}
g=\dd r^2+\phi^*g_Y+e^{-2cr}g_Z,
\end{equation} with $g_Z$ a
metric in the fibres and $g_Y$ a metric on the base of the
fibration.
\par
\smallskip
In this work we investigate the spectral theory of the Dirac
operator, denoted by $\Dir^d$ below, over a general complete
manifold $X$, noncompact only in generalized cusps as above. We
use the methods developed in \cite{Meaomwc}, \cite{Melaps} and
\cite{MelMazmwfb}. Following their philosophy  we compactify $X$
to a manifold with fibred boundary $\partial X=M$ by introducing
$x=e^{-cr}$ as a defining function for the boundary. The metric on
$X$ will then be denoted by $g_d$. This  is  a degenerate (and
singular) metric near the boundary, where it has the form
$$g_d=\left(\frac{\dd x}{x}\right)^2+\phi^*g_Y+x^2 g_Z,$$ (here we
have dropped the factors $c^2$ and $c^{-2}$). However, this metric
is a {\em true} metric on the fibres of a vector bundle ${}^dTX$
whose space of sections ${}^d{\cal V}(X)$ can be identified with
 the space of sections of $X$ into $TX$, possibly singular at $\partial
X$, but of finite $g_d$-length.
\par
In terms of local fibre coordinates ${\bf y},{\bf z}$ in $M$ and a
 product decomposition of $X$ near $\partial X$, vector fields
 in ${}^d{\cal V}(X)$ are locally spanned over $C^{\infty}(X)$ by $$x\delx, \quad
\delyy=(\dely{1},\ldots, \dely{h}), \quad
\frac{1}{x}\delzz=(\frac{1}{x}\delz{1},\ldots,
\frac{1}{x}\delz{v}).$$ It is easily seen from this representation
that ${}^d{\cal V}(X)$ is {\em not} a Lie algebra, such that
there is no good notion of $d$-differential operators etc.. This
prompts the introduction of the {\em nondegenerate} conformal
metric
\begin{equation}\label{phiintro}
g_{\phi}=\frac{1}{x^2}g_d=\left(\frac{\dd
x}{x^2}\right)^2+\frac{1}{x^2}\phi^*g_Y+g_Z.
\end{equation}
The  associated space of vector fields of finite
$g_{\phi}$-length, $${}^{\phi}{\cal V}(X):=x\cdot {}^{d}{\cal
V}(X)\subset\Gamma(X,TX),$$ consists of true vector fields on $X$
and can be seen to be the space of sections of the vector bundle
$$\phitx=x\,{}^dTX.$$
 It is easy to verify that it forms a Lie
algebra and we can define the notions of $\phi$-differential and
$\phi$-pseudodifferential operators. The
 analysis of these operators has been begun in \cite{MelMazmwfb}, using the
 general framework set out in \cite{Meaomwc}.
 \par\smallskip
Now the  spinor Dirac operator $\Dir^d$, associated to the metric
$g_d$, is a selfadjoint operator on $L^2(X,E,\dvol_d)$. In
Appendix \ref{conformal} it is shown that
\begin{equation}\label{introconf}
\Dir^d=\cl_d\circ\nabla^{E,d}= \frac{1}{x}\Dir^{\phi}
+\frac{n-1}{2}\cl_{\phi}(\frac{\dd x}{x^2}),
\end{equation}
where $\Dir^{\phi}\in \Diff_{\phi}^1(X,E)$ is the spinor Dirac
operator associated to $g_{\phi}$. It follows from this -- and is
actually true for  {\em any} Dirac operator $\Dir^d$ on a Clifford
bundle associated with a metric $g_d$ as above -- that $x\Dir^d$
is an elliptic $\phi$-differential operator.
 This will allow us to apply the analysis developed
for $\phi$-operators in the first part of the text to the spectral
analysis of $\Dir^d$.
\par
\smallskip
Before describing the results obtained in this work, let us
investigate the structure of the problem somewhat further. First,
the restriction of the operator $x\Dir^d$ to the boundary in the
usual sense of differential operators is an elliptic vertical
family of first order differential operators $\Dir^{\phi,V}$. It
is easy to see that $\Dir^d$ has compact resolvent, and therefore
pure point spectrum, if this family is invertible. The other
extreme case, where one allows arbitrary behavior of the null
spaces of $\Dir^{\phi,V}$, turns out to be too difficult in
general. In this work we  make the additional {\em assumption},
fulfilled for instance in the case of the signature operator, that
the null spaces  of $\Dir^{\phi,V}$ form a vector bundle $${\cal
K}\rightarrow Y,\qquad \Gamma(Y,{\cal K})={\rm null}
(\Dir^{\phi,V}),$$ over the base $Y$ of the boundary fibration. We
will denote the projection onto this null space by $\Pi_{\circ}$
and write $B=Y\times\RR_{+,x}$, which can be thought of as the
cylindrical base of the model cusp (compactified at ``infinity''
$x=0$ ). Then we can introduce the operator
\begin{equation}\label{introibd}
I_{b}(\Dir^d)="\Pi_{\circ}\Dir^d\Pi_{\circ}|_{\partial
X,{}^bTX}"=\Dir_Y+\cl_d(\frac{\dd x}{x})x\frac{\partial}{\partial
x }: C^{\infty}(B,{\cal K})\rightarrow C^{\infty}(B,{\cal K}).
\end{equation}
This is described  more precisely in Section \ref{normalopsbf}. It
turns out that -- just as in the $b$-case, or the case of
asymptotically cylindrical ends -- the continuous spectrum of
$\Dir^d$ is governed by this operator, with bands of continuous
spectrum starting at the eigenvalues of $\Dir_Y$ and going out to
infinity. Especially, the operator $\Dir^d$ is Fredholm, if and
only if $\Dir_Y$ is invertible.
\par\medskip
This result will be obtained as a corollary of the construction of
the analytic continuation of the resolvent $$
G_{\Dir^d}(\lambda)=(\Dir^d-\lambda)^{-1}$$ of $\Dir^d$ described
in Chapter 3. Using the concepts developed in \cite{MelMazmwfb},
\cite{Meaomwc}, \cite{Melaps} etc., we explicitly construct the
resolvent as an integral kernel on the space $X^2_{\phi}$. This is
a compactification of the interior of $X^2$, which is  adapted to
the specific asymptotics of the problem at infinity. The main
result, summarizing this procedure and stated as Theorem
\ref{mainresult1}, gives a detailed account of the asymptotics of
the resolvent. It is then possible to use the methods known from
the $b$-calculus to derive the usual applications of this
construction, such as a description of the spectral measure, $$\dd
E_{\Dir^d}(r)=\lim_{\varepsilon\rightarrow 0}
          \frac{1}{2\pi i}[ G_{\Dir^d}(r+i\varepsilon)-
          G_{\Dir^d}(r-i\varepsilon)]\dd r,$$
 of $\Dir^d$ in terms of its generalized eigenfunctions.  This is done in Section
 \ref{secconsequences}.
\par\smallskip
The general framework for the constructions in Chapter 3 will be
recalled in Chapter 2. Here the space $X^2_{\phi}$  is described
as a blow up of $X^2$ and the notion of $\phi$-pseudodifferential
operators is introduced as a special class of distributions on
$X^2_{\phi}$. Building on the work in \cite{MelMazmwfb}, we then
analyze the mapping- and composition properties of such operators.
\par
\smallskip
Chapter 1 contains an introduction to the geometry of manifolds
with fibred boundaries endowed with a metric of type
(\ref{phiintro}) as described above. Here we also introduce the
concept of coding the asymptotic behavior of sections of a vector
bundle into the structure of a new vector bundle. An instance of
this has already been presented in the definition of the bundles
${}^dTX$, $\phitx$, and the general concept will be used
repeatedly throughout the text. The more detailed results about
the behavior near the boundary of the connections and curvatures
associated to $g_{\phi}$ and $g_d$, obtained in the later parts of
this Chapter, will however not be used seriously until Chapter 5.
\par
\medskip
The second  important topic in this work is the proof of an index
theorem for $\Dir^d$. As a corollary to the description of the
spectral measure, we can describe the large time behavior of the
heat operator \begin{equation}\label{introheatkernel}
\exp(-t\Dir^2)=\int_{\RR}e^{-tr^2}\dd E_{\Dir^d}(r)
\end{equation}
and prove that the extended $L^2$-index of $\Dir^d$ (is finite
dimensional and) can be calculated as the large time limit of the
regularized heat supertrace: $${\rm ind}_-(\Dir^d)=
\lim_{T\rightarrow \infty}{\rm
reg}_{z=0}\Str(x^z[\exp(-T(\Dir^d)^2)]).$$ This will allow us to
derive a formula for the index using a variant of the
McKean-Singer argument, i.e. by calculating the small time limit
of the regularized heat supertrace and the difference between its
large- and small-time limits. However, as can be seen from
(\ref{introheatkernel}), the heat operator at small times involves
the resolvent at large spectral parameter. Our construction in
Chapter 3 does not give us sufficient control of the resolvent in
that area.
\par\smallskip
 Thus, in
Chapter 4, we present a direct construction of the heat kernel for
finite and small times. Again, this is done by introducing an
appropriate compactification of the kernel space and a
corresponding ``heat calculus'', within which the heat kernel can
be constructed. Theorem \ref{mainresult2} gives a certain amount
of control of the asymptotics of the heat kernel at small times
and near the boundary.
\par\smallskip
Chapter 5 is then devoted to a refinement of that result,
providing more information about the Clifford degree of the
coefficients in the asymptotic expansion. In the philosophy of
Patodi, Getzler, etc. this will allow us to eliminate the
divergent terms in the small time expansion(s) of the local heat
supertrace. This then turns out to be sufficient to prove the
index theorem. Denoting by $\eta(\Dir_Y)$ the classical
eta-invariant associated to the operator $\Dir_Y$ and by
$\widehat{\eta}(\Dir^{\phi,V},E)$ the family eta invariant
associated to the vertical family $\Dir^{\phi,V}$ (see Appendix
\ref{vertfam}) we prove the following formula for the extended
$L^2$-index of $\Dir$ in Theorem \ref{indextheorem}: $$ {\rm
ind}_-(\Dir)= \frac{1}{(2\pi
i)^{n/2}}\int_X\widehat{A}(R^d)\Ch(F^{E/S,d})+\frac{1}{(2\pi
i)^{(h+1)/2}}\int_Y\widehat{A}(R^Y)\widehat{\eta}(\Dir^{\phi,V},E)+\frac{1}{2}
\eta(\Dir_Y).$$ This formula generalizes the corresponding formula
obtained by W. M\"uller in \cite{MuellerBook}. It reduces to the
usual APS (or b-) index formula in the case that there is no
fibre, i.e. $Z=\{{\rm pt}\}$.
\par
 In Section
\ref{signatureop} we  specialize our results to the case of the
signature operator ${\sf S}^d=\dd_d +\dd_d^*$. It is a well known
fact in the $b$-case that the extended $L^2$-index of this
operator equals the $L^2$-signature. This is shown to be also true
in our case. We compare our formulas with the adiabatic limit
formula for the eta invariant in \cite{BC} and \cite{Daithesis}.
More serious applications are however left for future examination.
Finally some auxiliary results and background material have been
collected in the Appendix.
\par
\smallskip
It is our declared purpose in this work, to give a thorough
introduction to a set of methods for the spectral analysis of
geometric operators associated to  $d$-metrics. Especially, our
proof of the index theorem  is far from being ``minimal''. Remark
\ref{shorterproof} at the end of the text sketches several
shortcuts for the proof.\par
 Some remarks about notation: For
simplicity we will assume throughout that $M$ is {\em connected},
i.e. there is only one cusp. However all the results presented
here generalize directly to the case of several cusps. Whenever
there is a grading we use {\em graded} calculation rules. Thus
commutators $[A,B]$  and tensor products $A\otimes B$ of graded
objects are always graded. Notations like $\str$ or ${\rm sdim}$
are understood.
\par
\bigskip
\noindent{\bf Acknowledgements}
\par
\smallskip
First, special thanks are due to both my advisors, Werner M\"uller
and Richard Melrose. The work of Werner M\"uller  is at the basis
of this thesis. He suggested that I try to generalize his results
in \cite{MuellerBook} and supported my research for many years. I
thank Richard Melrose for inviting me to MIT in 97/98 and taking
the time to teach me his approach to the spectral analysis of
operators in noncompact settings. Most of the methods used in this
thesis are based on his ideas and I would not have got very far
without his support.
\par
Many more people have contributed to various parts of this work.
Especially, I wish to thank my fellow students in Bonn and
Cambridge, Sang Chin, Robert Lauter, Paul Loya, Sergiu Moroianu,
Gorm Salomonsen, Oliver Sick and Jared Wunsch, who helped me learn
the techniques used here and with whom I discussed many of the
problems in detail. Some of their work (like \cite{RobertSergiu1},
\cite{RobertSergiu2} and \cite{Paul}) is closely related to mine.

\newpage
\section{The Geometry of $\phi$-Manifolds}
\par
In this Chapter, we introduce the main geometrical notions
associated to manifolds with a fibred boundary. The concepts and
notation introduced in the first Section, as well as the notion of
an exact $\phi$-metric introduced at the beginning of Section
\ref{exactphisection}, are fundamental to all the further parts of
this work. On the other hand,  most of the more detailed
calculations in Sections \ref{exactphisection} to
\ref{boundarygeom} will only be used in later Chapters. They can
be skipped at first reading and be referred back to as needed.
\par
\subsection{$\phi$-Manifolds}\label{secphimanifolds}
\bigskip
Let $X$ be a compact manifold of dimension $n$ with connected and fibred boundary
$\partial X=M$, $$\phi:\partial X\stackrel{Z}{\longrightarrow}Y.$$
Denote by $v$ the dimension of the fibre $Z$ and by $h$ the
dimension of the base $Y$ of this fibration. We call X a {\em
$\phi$-manifold}. Fix a defining function $x$ for the boundary
$\partial X$.\par Given these data we can define the Lie algebra
of $\phi$-vector fields on $X$ by $${}^{\phi}{\cal V}(X):=\{ V\in
\Gamma(TX)\,|\quad V|_{\partial X}\quad \mbox{tangent to $Z$
and}\quad V x = O(x^2)\}.$$ The space ${}^{\phi}{\cal V}(X)$ can
be regarded as the space of sections of a vector bundle over $X$: $$
\phiv(X) \equiv \Gamma (\phitx), $$ where we call $\phitx$ the
$\phi$-tangent bundle to $X$. The dual bundle, $\phitsx$, and all
other associated tensor bundles will be distinguished by the
symbol ${}^{\phi}$. Being isomorphic over the interior of $X$, the
bundles $\phitx$ and $TX$ are abstractly isomorphic over all of
$X$. However, there is no {\em natural} isomorphism between them,
and it will be shown that $\phitx$ carries some extra structure
that $TX$ does not have.
\par
 It will be convenient to recall
that the fact that $\Gamma(\phitx)$ is a Lie algebra is equivalent
to the fact that exterior derivation maps $$\dd
:\Gamma(\phitsx)\rightarrow\Gamma(\Lambda^2\,\phitsx).$$ Now, the
inclusion ${}^{\phi}{\cal V}(X)\hookrightarrow \Gamma(TX)$
 gives a natural
map $\phitx\rightarrow TX$, which is not an isomorphism over the
boundary, since its image consists of the vectors tangent to the
fibres at the boundary. Restricted to the boundary we get a short
exact sequence
\begin{equation}\label{shortphindx}
0\rightarrow \phindx \rightarrow \phitx|_{\partial X}\rightarrow
  V\partial X\rightarrow 0,
\end{equation}
where $V\partial X$ is the ``vertical'' bundle of vectors tangent
to the fibres in the boundary. Also, this sequence is used to {\em
define} the vector bundle $\phindx$ over $\partial X$. The dual
sequence
\begin{equation}\label{shortvsx}
0\rightarrow V^*\partial X\rightarrow \phitsx|_{\partial
X}\rightarrow \phinsdx \rightarrow 0,
\end{equation}
then defines the vector bundle $V^*\partial X$ over $\partial X$.
\par
The existence of these subbundles over the boundary allows us to
make sense of the vector bundles $[\phindx]$ and
$[{}^{\phi}V^*\partial X]$ over $X$ (!) by defining their spaces
of sections as
\begin{eqnarray} \Gamma(X,[\phindx])&:=& \{ W\in
\Gamma(X,\phitx)|\quad W|_{\partial X,\phitx}\in \phindx
\}\label{[phindx]}\\
 \Gamma(X,[{}^{\phi}V^*\partial X])&:=& \{ \xi\in
\Gamma(X,\phitsx)|\quad\xi|_{\partial X,\phitx}\in V^*\partial X
\}.\label{[vsx]}
\end{eqnarray}
Here the notation is intended to emphasize that restriction to the
boundary is understood in the sense of elements in $\phitx$ or
$\phitsx$. Thus, the definition of these bundles really depend on
the {\em pairs} $(\phindx,\phitx)$ and $(V^*\partial X,\phitsx)$.
We will include this into the notation more explicitly, whenever
there is a possible ambiguity. Now the obvious inclusions of the
spaces of sections give natural maps
 $$ [\phindx]\longrightarrow \phitx,\qquad [{}^{\phi }V^*\partial
X]\longrightarrow \phitsx, $$ which, combined with restriction to
the boundary in that sense, give
 $$ [\phindx]\stackrel{|_{\partial
X,\phitx}}{\longrightarrow}\phindx,\qquad [{}^{\phi }V^*\partial
X] \stackrel{|_{\partial X, \phitsx}}{\longrightarrow} V^*\partial
X.$$  Note that the bundle $[\phindx]$ is well known as the
``scattering bundle'' over X (compare
\cite{Melsastftloaes},\cite{Melgeomscat}). Recalling that the
``$b$-bundle'' ${}^bTX$ is defined as the
 vector bundle whose space of sections consists of sections of $TX$,
 which are tangent to the boundary (\cite{Melaps}), we have
 $$[\phindx]={}^{\rm sc}TX=x{}^bTX.$$
However, in our framework, sections in $[\phindx]$ should rather
be thought of as sections in $\phitx$ which are ``asymptotically
horizontal''. In contrast to that, at this stage, we cannot make
invariant sense of the notion of ``asymptotically vertical''
sections in $\phitx$!
\par\bigskip
\noindent
{\bf Local Description}
\par\medskip
Let us emphasize that -- once the boundary defining function $x$
has been fixed -- the constructions above are completely natural.
It will be useful however, to have coordinate representations.
 So assume that we have chosen a collar neighborhood
$C(\partial X)\cong
\partial X \times [0,\varepsilon[_x$ of the boundary and a local
trivialization of $\phi$, so that locally near the boundary $X$
looks like $$ U_X\cong U_Y\times Z \times
[0,\varepsilon[\stackrel{\varphi}{\longrightarrow} U_Y.$$
 Then we may choose coordinates
\begin{itemize} \item $x :$ a boundary defining function of
$Y\times \{0\}$ in $Y\times [0,\varepsilon[$, lifted via $\varphi$
to $U_X$,
\item  ${\bf y}=(y_1, \ldots, y_h)$ : local coordinates on
$U_Y\subset Y$, lifted to $U_X$ via $\varphi$,
\item
${\bf z}=(z_1, \ldots, z_v)$: local coordinates on $Z$ extended to
$U_X$.
\end{itemize}
In these coordinates sections in $\phitx$
 are locally spanned by $$ x^2\delx, \quad
x\delyy=(x\dely{1},\ldots, x\dely{h}), \quad   \delzz=(\delz{1},
\ldots, \delz{v}),$$ and the fact that ${}^{\phi}{\cal V}(X)$ is
indeed a Lie algebra can be read off from this representation.
Also elements of $\Gamma([\phindx])$ are locally spanned by $$
x^2\delx, \quad x\delyy=(x\dely{1},\ldots, x\dely{h}), \quad
x\delzz=(x\delz{1}, \ldots, x\delz{v}),$$ and sections of
$[{}^\phi V\partial X]$ locally look like $$ x^3\delx, \quad
x^2\delyy=(x^2\dely{1},\ldots, x^2\dely{h}), \quad
\delzz=(\delz{1}, \ldots, \delz{v}).$$
\par
A  change of product structure and/or coordinates compatible
with $x$ can be regarded as a family $\Psi=(\Psi_x:\partial
X\rightarrow \partial X)_{x\in [0,1]}$ of diffeomorphisms on
$\partial X$ such that $\Psi_0$ is a fibre bundle diffeomorphism
of $\phi:\partial X\rightarrow Y$ fixing the base. Then
 $$\Psi^*\frac{\dd y_j}{x}=\frac{\dd
y_j}{x}+O(x^0)\langle \dd {\bf y}, \dd {\bf z}\rangle,\quad
\mbox{and}\quad \Psi_*\delz{j}= O(x^0)\langle
\delzz,x\delyy\rangle,$$ which exemplifies the fact that we cannot
identify ``vertical'' vectors within $\phitx$, for now!
\par\bigskip
Thus, as an additional piece of structure, {\em choose} a
splitting of the Lie algebra sequence (\ref{shortphindx},
\ref{shortvsx}), i.e. a decomposition
\begin{equation}\label{choosevertical}
\phitx|_{\partial X} =\phindx \oplus V\partial X, \qquad
\phitsx|_{\partial X}=\phinsdx\oplus V^*\partial X. \end{equation}
For instance, such a decomposition in the sense of Lie algebras
could be obtained  by  choosing a product structure near the
boundary as above.\par
 The structure (\ref{choosevertical}) will be {\em fixed
once and for all}, the corresponding projections will be denoted
by ${\sf v}$ and ${\sf n}$. Among other things, this allows us to
introduce the maps $$\phi_*:\phitx|_{\partial X}\rightarrow
{}^{\rm sc}TB|_Y,\qquad \phi^*:{}^{\rm sc}T^*B\rightarrow
\phitsx|_{\partial X},$$ {\em over the boundary}, with $B=Y\times
[0,\infty[_x $. Also, we can now define the bundles
$[{}^{\phi}V\partial X]$, the space of ``asymptotically vertical''
vectors in $\phitx$, and $[\phinsdx]$ in analogy with our
definition in (\ref{[phindx]}, \ref{[vsx]}). \begin{remark}\rm  The
reader should be aware of the following {\bf pitfalls}. First
$$[\phindx]^*={}^{\rm sc}T^*X=\frac{1}{x}[{}^{\phi}V^*\partial X
],\quad\mbox{thus in general}\quad[\phindx]^*\neq[\phinsdx] !!$$
Also, readers should convince themselves that in general
$$[\phindx]\otimes[\phindx]\neq
[\phindx\otimes\phindx]_{\phitx\otimes\phitx}  !!$$ The first
bundle really encodes some second order behavior at the boundary
of sections in $\phitx\otimes\phitx$, while the second bundle
clearly only carries first order information.
\end{remark}
\par
\smallskip
The following Lemma restates the fact that we have chosen a
decomposition compatible with the Lie algebra structure of
$\phitx$:
\begin{lemma}
\begin{enumerate}
\item $[\Gamma([{}^{\phi}V\partial X]),\Gamma([{}^{\phi}V\partial X])]\subset
      \Gamma([{}^{\phi}V\partial X])$
\item $[\Gamma([\phindx]),\Gamma([\phindx])]\subset
       O(x)\Gamma([\phindx])$
\item $[\Gamma([\phindx]),\Gamma(\phitx)]\subset
       \Gamma([\phindx])$. \qed
\end{enumerate}
\end{lemma}
\par\medskip
Hoping that the reader has by now got used to the above way of
coding the boundary (or asymptotic) behavior of sections into some
bundle structure, we proceed with introducing one further
modification which is going to be used in the next Section:
\par
The vector bundle ${}^{\rm e}TX$ over $X$ is defined (compare
\cite{mazzthesis}) as the bundle whose space of sections is
$$\Gamma({}^{\rm e}TX)\equiv \{ W\in \Gamma(TX) | \quad
W|_{\partial X} \in V\partial X \}.$$ Again, this is a Lie algebra
(thus exterior derivation is well defined on $\Lambda{}^{\rm
e}T^*X$ ) and we have a map $\phitx\rightarrow {}^{\rm e}TX$.
Thus, the sequence analogous to (\ref{shortphindx}) is $$
0\rightarrow  {}^{\rm e}N\partial X\rightarrow {}^{\rm e}TX
|_{\partial X}\rightarrow
  V\partial X\rightarrow 0,$$
and it has a splitting map $V\partial X\rightarrow {}^{\rm e}TX$
induced by  our splitting (\ref{choosevertical}). In the above set
of local coordinates, the bundle ${}^{\rm e}TX$ is spanned by $$
x\delx, \quad x\delyy=(x\dely{1},\ldots, x\dely{h}), \quad
\delzz=(\delz{1},\ldots, \delz{v}).$$  We leave it to the reader
to define the bundles $[{}^{\rm e}N\partial X]$, $[{}^{\rm
e}V^*\partial X]$, $[{}^{\rm e}N^*\partial X]$, and $ [{}^{\rm
e}V\partial X]$ over X and to write down their local descriptions.

\subsection{Exact $\phi$-Metrics}\label{exactphisection}
A $\phi$-metric is a fibre metric $g\in\Gamma(S^2\phitsx)$ on
$\phitx$. We call it a {\em compatible} $\phi$-metric when the
decomposition at the boundary $\phitx|_{\partial X}=\phindx\oplus
V\partial X$ is orthogonal. Another way to state this property is
to demand that $g$ has a decomposition $$ g\in
\Gamma(S^2[\phinsdx])+\Gamma(S^2[{}^{\phi}V^*\partial X]).$$\par
An {\em exact $\phi$-metric} $g_{\phi}$ on $\phitx$ is a
compatible $\phi$-metric with a more refined decomposition of the
form
\begin{eqnarray}\label{exactphimetric}
\lefteqn{\qquad\qquad\qquad\qquad\qquad g_{\phi}=
\frac{\dx^2}{x^4} + h + g_Z, \qquad\mbox{with}}
\\
 h&\in& \Gamma(S^2[{}^{\rm e }N^*\partial X])\subset \Gamma(S^2[\phinsdx]),
 \quad h|_{\partial X,{}^{\rm e}TX}=\frac{1}{x^2}\phi^*g_B, \quad
 g_B\in \Gamma(S^2T^*B|_Y)\nonumber
 \\ g_Z&\in&\Gamma(S^2T^*X)\subset \Gamma(S^2[{}^{\phi}V^*\partial
 X]),
\nonumber
\end{eqnarray}
where  the notation $h|_{\partial X,{}^{\rm e}TX}$ means that we
regard $h$ as an element in $S^2\,{}^{\rm e}T^*X$ via the map
$[{}^{\rm e}N^*\partial X]\rightarrow{}^{\rm e}T^*X$ when
restricting to the boundary. We sometimes write
$h=[\frac{1}{x^2}\phi^*g_B]$ to symbolize this construction. We
also make use of the notation $g_{B,b}=\frac{\dd x^2}{x^2}+g_B$
for the underlying $b$-metric on $B$.  Note that {\em the
decomposition (\ref{exactphimetric}) is by no means unique}. For
instance, only the $S^2V^*\partial X$-part of the restriction
$g_Z|_{\partial X}$ makes invariant sense!
\par
\bigskip
\noindent{\bf Product Metrics}
\par
\medskip
The standard example of an exact $\phi$-metric is
a metric of type (\ref{phiintro}) as described in the Introduction.
There we assumed that there is a fixed choice of product structure,
$$C(M)=M\times [0,\varepsilon[_x \stackrel{\varphi}{\longrightarrow} Y\times
[0,\varepsilon[_x,$$
in a collar neighborhood of the boundary $\partial X$ and a choice of
horizontal space $HM$. Especially, the fibration $\phi$
and the horizontal space of the boundary  are extended to $C(M)$. Then a metric $g_Y$
on $Y$ can be lifted to a metric $\varphi^*g_Y$ on $HM$ over $C(M)$.
A metric, which over $C(M)$ is of the form
\begin{equation}\label{productref}
g_{\phi}=\left(\frac{\dd
x}{x^2}\right)^2+\frac{1}{x^2}\varphi^*g_Y+g_Z,
\end{equation}
where  $g_Z$ is a metric in the fibres, will then be called a
{\em product $\phi$-metric} on X.
It is easily seen that a local basis
$$ V_1,\ldots, V_v, \in VM \quad\mbox{and}\quad H_1,\ldots H_h \in HM,$$
which is orthonormal w.r.t. $\phi^*g_Y+g_Z$, extends to a local
basis
$$ V_1,\ldots, V_v, \quad xH_1,\ldots xH_h, \quad x^2\frac{\partial}{\partial x}
={\sf X}_{\phi}\quad  \in \phitx
 $$
near $\partial X$, which is orthonormal w.r.t. $g_{\phi}$ and
$\nabla^{\phi}$-parallel along $\frac{\partial }{\partial x}$.
Using this, the connection and curvature of such a product metric
can be calculated rather explicitly along the lines described
in Appendix \ref{geomtens}.
\par
The notion of an exact $\phi$-metric is  the natural
generalization
of  metrics of this type to the
``$\phi$-setting''. In fact, an exact $\phi$-metric can be thought
of as being asymptotically of the form (\ref{productref}).
In this general setting, the techniques described above are not
available to us and the geometry at the boundary is somewhat more
intricate.
However, we will sometimes refer to the product case to simplify
our formulas.
\par
\medskip
The following Lemma lists the essential properties of (the inverse of) an
exact $\phi$-metric at the boundary:
\begin{lemma}[{\boldmath Inverse of an Exact $\phi$-Metric at $\partial X$}]
\label{gbasics}\quad \par \noindent
 Let $\alpha_i\in\Gamma(T^*X)$,
$\beta_j \in \Gamma(T^*B|_Y)$ and choose extensions
$[\frac{1}{x}\phi^*\beta_i]\in\Gamma([{}^{\rm e }N^*\partial X])$
\begin{enumerate}
\item $g_{\phi}^{-1}(\frac{\dd x}{x^2},\frac{\dd
x}{x^2})=1+\phi^*O(x^2)$
\item $g_{\phi}^{-1}(\frac{\dd x}{x^2},[\frac{1}{x}\phi^*\beta_1])
=\phi^*O(x)$
\item $g_{\phi}^{-1}(\frac{\dd x}{x^2},\alpha)= O(x^2)$
\item $g_{\phi}^{-1}([\frac{1}{x}\phi^*\beta_1],[\frac{1}{x}\phi^*\beta_2])
=g_{B,b}^{-1}(\beta_1,\beta_2)=\phi^*O(x^0)$
\item $g_{\phi}^{-1}([\frac{1}{x}\phi^*\beta_1],\alpha)
=O(x) $.
\item $g_{\phi}^{-1}(\alpha_1,\alpha_2)=O(x^0)$
\end{enumerate}
Here the notation $\phi^*O(x^l)$ means that the term is $O(x^l)$
but the coefficient of $x^l$ is of the form $\phi^*f$.\qed
\end{lemma}
\par\smallskip
We will use this result to analyze the boundary behavior of the
Levi-Civita connection
for an exact $\phi$-metric.
First, we need to introduce some more notation.
Denote by ${\sf X}_{\phi}$ the vector field
      $g_{\phi}^{-1}(\frac{\dd x}{x^2})\in \Gamma(x^2TX)$.
\begin{lemma}\label{thehatmap}
 The map $A\mapsto \widehat{A}:=A-\frac{\dd x}{x^2}(A){\sf X}_{\phi} $
 maps $TX\rightarrow {}^bTX$ (and $xTX\rightarrow {}^{\rm sc}TX$) and
 has the following properties for $A,B\in TX$:
 \begin{enumerate}
 \item $\dd x(\widehat{A})=O(x^2)$
 \item $\widehat{\widehat{A}}=\widehat{A}+O(x^0){\sf X}_{\phi}$
 \item $g_{\phi}(\widehat{x A},\widehat{x
 B})=\phi^*g_B(A,B)|_{\partial X}+O(x)$.
 \item $\frac{1}{x^2}\phi^*g_B({\sf X}_{\phi},A)=\dd x (A)O(x^0)$,
where $O(x^0)$ is independent of $A$.\qed
\end{enumerate}
\end{lemma}
\par
The central geometric result in our setting is now stated in
\begin{proposition}
[{\boldmath Connection and Curvature of an Exact $\phi$-Metric at
$\partial X$}]\label{exactconnection} {\quad}
\par\smallskip
\noindent
The Levi-Civita connection for an exact $\phi$-metric is a {\em true}
connection, i.e. $$\nabla^{\phi}:\Gamma(\phitx)\longrightarrow
\Gamma(T^*X\otimes \phitx), \quad \mbox{and}\quad R^{\phi} \in \Gamma(\Lambda^2T^*X,
\End(\phitx)).$$
 Also, for $T, T_1, T_2\in \Gamma({}^bT X)$, $N\in
\Gamma(TX)$, $B\in\Gamma(\phitx)$, $\omega\in \Gamma({}^{\rm
sc}T^* B|_Y)$, $L\in \Gamma({}^{\rm sc}T B )$:
\begin{enumerate}
\item $\nabla^{\phi}_N\frac{\dd x}{x^2}(B)|_{\partial X}=
   -\frac{1}{x^2}\phi^*g_B(\widehat{x N},B)|_{\partial B}$
\item $\nabla^{\phi}_N{\sf X}_{\phi}|_{\partial X,\phitx}=
   -\widehat{x N}$
\item $\nabla^{\phi}_T:\Gamma([\phindx])\rightarrow\Gamma([\phindx])$, and
   $\nabla^{\phi}_T:\Gamma([{}^{\phi}V \partial X])\rightarrow\Gamma([{}^{\phi}V
   \partial X])$
\item $\nabla^{\phi}_T|_{\partial X}={\sf v}\circ \nabla^{\phi}_T\circ{\sf v}+
\phi^*\nabla^{B,\rm sc}_T\circ {\sf n}$ on
$\Gamma(\phitx|_{\partial X})$
\item $R^{\phi}(T_1,T_2)|_{\partial X} ={\sf v}\circ R^{\phi}(T_1,T_2)\circ{\sf v}\oplus
\phi^*R^{B,sc}(T_1,T_2)$ in $\End(\phitx|_{\partial X})$.
\end{enumerate}
 Especially
$\nabla^{\phi}_T \phi^* L|_{\partial X}=0$ in $\phitx$ and
$\nabla^{\phi}_T \phi^*\omega|_{\partial X} = 0$ in
 $\phitsx$ for vertical $T$.
\end{proposition}
\begin{proof}
The full proof of the mapping property of $\nabla^{\phi}$ is
tedious and we do not give it here. To see how to proceed, let
us prove (a) in detail. First note that the expression on the RHS
of (a) is well defined and independent of the ambiguities in the
decomposition of $g_{\phi}$. To prove the equality, let $N, B$ be
as above. Then, writing again ${\sf X}_{\phi}=g_{\phi}^{-1}
(\frac{\dd x}{x^2})$ the Koszul formula yields $$
2\nabla^{\phi}_N\frac{\dd x}{x^2}(B)=(L_{{\sf
X}_{\phi}}g_{\phi})(N,B)+\dd \frac{\dd x}{x^2}(N,B)=(L_{{\sf
X}_{\phi}}g_{\phi})(N,B) .$$
 The Lie derivative of the metric can
now be calculated using the decomposition of the metric and
Cartan's formula:
\begin{eqnarray*}
L_{{\sf X}_{\phi}}\left( \frac{\dd x}{x^2} \right)^2(N,B)&=&
\left( \dd \Int(X_{\phi})\frac{\dd x}{x^2}\otimes \frac{\dd x
}{x^2}+\frac{\dd x}{x^2}\otimes \dd \Int(X_{\phi})\frac{\dd
x}{x^2}\right)(N,B)\\
&=&  \dd(1+\phi^*O(x^2))(N)\frac{\dd
x}{x^2}(B)+\dd(1+\phi^*O(x^2))(B)\frac{\dd x}{x^2}(N) =O(x)\\
 L_{{\sf X}_{\phi}}g_Z(N,B)&=& {\sf X}_{\phi}g_Z(N,B)-g_Z([{\sf X}_{\phi},N],B)-
 g_Z(N,[{\sf X}_{\phi},B])=O(x),
 \end{eqnarray*}
 where we have used Lemma \ref{gbasics} (a) and (b), especially that
 ${\sf X}_{\phi}\in \Gamma(x^2TX)$.\par
 The calculation of the contribution of the term $h$ in the
 decomposition of $g_{\phi}$ is a little bit more involved:
 First, write
$$ (L_{{\sf X}_{\phi}} h)(N,B)=\frac{1}{x^2}(L_{{\sf
X}_{\phi}}x^2h)(N,B)-2(1+\phi^*O(x^2))h(xN,B) .$$ We can assume
that $x^2h$ is (a sum of elements) of the form
$[\phi^*\alpha]\otimes[\phi^*\beta]\in [x{\,}^{\rm e}N^*\partial
X]^2$. Note that for elements in $[x{\,}^{\rm e}N^*\partial X]$ of
this type we have $ \dd [\phi^*\alpha], \dd [\phi^*\beta] \in
\Lambda^2 [x{\,}^{\rm e}N^*\partial X].$ Now
\begin{eqnarray*}
\frac{1}{x^2}(L_{{\sf X}_{\phi}}x^2h)(N,B)&=&
\frac{1}{x^2}(L_{{\sf
X}_{\phi}}[\phi^*\alpha]\otimes[\phi^*\beta])(N,B)\\ &=&
\frac{1}{x^2}\left[ \dd [\phi^*\alpha]({\sf
X}_{\phi},N)[\phi^*\beta](B)+(N[\phi^*\alpha]({\sf
X}_{\phi}))[\phi^*\beta](B)\right.\\ & &\left. \qquad\qquad + \dd
[\phi^*\beta]({\sf
X}_{\phi},B)[\phi^*\alpha](N)+(B[\phi^*\beta]({\sf
X}_{\phi}))[\phi^*\alpha](N )\right],
\end{eqnarray*}
and the first and the third summand here are $O(x)$ leaving us
with
\begin{eqnarray*}
\frac{1}{x^2}(L_{{\sf
X}_{\phi}}x^2h)(N,B)&=&\frac{1}{x^2}(N[\phi^*\alpha]({\sf
X}_{\phi}))[\phi^*\beta](B)+\frac{1}{x^2}(B[\phi^*\beta]({\sf
X}_{\phi}))[\phi^*\alpha](N )+O(x)\\ &=& 2(Nx)[\phi^*\alpha]
(\frac{1}{x^2}{\sf X}_{\phi})[\phi^*\beta](\frac{1}{x}B) +
\frac{1}{x^2}[\phi^*\alpha](N)B\phi^*O(x^2)+O(x)\\ &=& 2 (Nx) h
(\frac{1}{x^2}{\sf X}_{\phi}, x B)+ O(x).
\end{eqnarray*}
Summing the two results, the equation is proved. \par The proof of
(c) uses the same method. For $T\in\Gamma({}^bTX)$, $A\in
\Gamma([\phindx])$ and $V\in\Gamma([V\partial X])$  one can check
that $$
2\langle\nabla^{\phi}_TV,A\rangle_{\phi}=(L_{V}g_{\phi})(T,A),$$
and direct calculation, using the decomposition of $g_{\phi}$ as
above, shows that the RHS is of order $O(x)$. Parts (d) and (e)
follow easily. \qed
\end{proof}

\subsection{$\phi$-Differential Operators and the Dirac Operator $\Dir^d$}
Since $\Gamma(\phitx)$ is a Lie algebra of vector fields, the associated space of
$\phi$-differential operators on $X$, $\phidiff(X)$, can be
defined in the usual way. Given vector bundles $E, F$ this
definition can be extended to $\phi$-differential operators
between these bundles, for example by using connections. Of
interest to us will be the case of a {\em $\phi$-Clifford bundle
}, which is a hermitian vector bundle $E\rightarrow X$ with
hermitian metric $h^E$, connection $\nabla^{E,\phi}$, and
$(\nabla^{E,\phi},\nabla^{\phi})$-parallel Clifford action $$
\cl_{\phi}:\Gamma(\phitsx)\rightarrow\Gamma(\End(E)),$$ or,
equivalently $\phicl\in \Gamma(\phitx\otimes\End(E))$. The
notational conventions used w.r.t. this Clifford action have been
listed in Appendix \ref{cliffordconv}. The associated Dirac
operator, given by $$\Dir^{\phi}=\phicl\circ\nabla^{E,\phi},$$ is a
$\phi$-differential operator with coefficients in $E$. The
restriction of this operator to the boundary $\partial X$ is a
vertical family of Dirac operators which we sometimes denote by
$\Dir^{\phi,V}$.
\par
\medskip
Our real interest lies in ``degenerate'' metrics of the type
$g_d=x^2g_{\phi}$, where $g_{\phi}$ is an exact $\phi$-metric,
and the associated Dirac operator $\Dir^d$. We will
use the notation ${}^d$ to label geometric objects associated to
this metric. Thus $$g_d=x^2g_{\phi},\quad{}^dTX=x^{-1}\phitx,\quad
{}^dT^*X=x\phitsx,\quad {\sf X}_d=g_d^{-1}(\frac{\dd x}{x})$$
 and so on. Note
that the sections of ${}^dTX$ may be singular as sections in $TX$
and do {\em not} form a Lie algebra, which is why we have to
consider the bundle $\phitx$ in the first place.
\par
The Levi Civita connection
$$ \nabla^d:\Gamma({}^dTX)\rightarrow \Gamma(T^*X\otimes {}^dTX),
\quad \nabla^d:\Gamma({}^dT^*X)\rightarrow \Gamma(T^*X\otimes
{}^dT^*X),$$ associated to the metric $g_d$, is a true connection
and can be described
in terms of the connection $\nabla^{\phi}$ as in Appendix
\ref{conformal}.
Thus for
vectors $N\in TX$, $W\in\Gamma(\phitx)$, and a differential form
$\alpha\in\Gamma(\Lambda\phitsx)$,
 \begin{eqnarray}
\nabla^d_N \frac{1}{x}W &=&\frac{1}{x}\nabla^{\phi}_N
W+\frac{1}{x}G(N)W,\, G(N)W :=\frac{\dd x }{x^2}(W) x
N-g_{\phi}(W,x N){\sf X}_{\phi}\label{nabladphivector}\\
 \nabla^d_N\ x^{\sf N} \alpha &=& x^{\sf N}\nabla^{\phi}_N\alpha
 +x G(N)\alpha,\,
 G(N)\alpha:=\ext_{\phi}(x N)\Int_{\phi} (\frac{\dd x}{x^2})\alpha-
\ext(\frac{\dd x}{x^2})\Int(x N) \alpha,\label{nabladphiform}
\end{eqnarray}
where ${\sf N}$ is the counting operator.
It will be
useful to remember that for $N\in TX$
$$G(N)-G(\widehat{N})=O(x^2)\End(\phitx)$$ by Lemma
\ref{thehatmap}. The map $G(N)$ is of order $O(x)$, if $N$ is vertical
at the boundary. In general $G(N)$ maps
\begin{equation}\label{gnice}
G(N):\Gamma([{}^{\phi}V\partial X])\longrightarrow
O(x)\Gamma(\phitx),\qquad G(N): \Gamma([\phindx]) \longrightarrow
\Gamma([\phindx]).
\end{equation}
\par
We can now introduce our basic objects of study.
Given $g_d$ as above, we consider a {\em $d$-Clifford bundle
}, which is a hermitian vector bundle $E\rightarrow X$ with
parallel
hermitian metric $h^E$, connection
$\nabla^{E,d}:\Gamma(E)\rightarrow \Gamma(T^*X\otimes E) $, and
$(\nabla^{E, d},\nabla^{d})$-parallel Clifford action $$
\cl_{d}:\Gamma({}^dT^*X)\rightarrow\Gamma(\End(E)).$$
The associated Dirac operator is
$$\Dir^d=\cl_d\circ \nabla^{E,d}.$$
It is shown in Appendix \ref{conformal} -- or can be checked by hand --
 that the definitions
\begin{eqnarray}\nabla^{E,\phi}_N \xi&:=&
\nabla^{E,d}_N\xi+ \frac{1}{2}\cl_{\phi}(\widehat{x
N})\cl_{\phi} (\frac{\dd x}{x^2})\xi, \label{nablaedephi}\\
\cl_{\phi}(\alpha)&=&x\cl_d(\alpha),\quad \cl_{\phi}(A)
=\frac{1}{x}\cl_{d}(A),\quad\mbox{for }
A\in\phitx,\alpha\in\phitsx,\nonumber
\end{eqnarray}
endow the bundle $E$ with the structure of a $\phi$-Clifford
bundle. The relationship between the Dirac operators $\Dir^d$
and $\Dir^{\phi}$ is given by
$$\Dir^d=\cl_d\circ\nabla^{E,d}= \frac{1}{x}\Dir^{\phi}
+\frac{n-1}{2}\cl_{\phi}(\frac{\dd x}{x^2})
=x^{-(n-1)/2}\frac{1}{x}\Dir^{\phi} x^{(n-1)/2}.$$
 This will allow
us to apply the analysis developed for $\phi$-operators in the
next Chapters to the operator $\Dir^d$.
\par\smallskip
It should be emphasized that the $d$-Clifford structure on $E$ and the
Dirac operator $\Dir^d$ are regarded as our {\em basic objects}.
For example, as explained in
Appendix \ref{conformal},  if $\Dir^d$ is the signature operator $\dd_d+\dd_d^*$
on $E=\Lambda\,{}^dT^*X$, then $\Dir^{\phi}=x\Dir^d$
{\em differs} from the conformally transformed signature operator
by an endomorphism.
Whenever we want to  calculate explicit formulas involving the Clifford
bundle we will therefore be careful to use the $d$-structure.
\par
\medskip
 Let us end this Section with  some useful remarks about the
 connection $\nabla^d$ and about
the relationship between the curvatures $R^{\phi}$ and $R^d$:
\begin{lemma}[{\boldmath Connection and Curvature for $g_d$ at $\partial X$}]
\label{cab}
 Let $T\in \Gamma({}^bTX)$. Then
\begin{enumerate}
\item $\nabla^d_T\frac{\dd x}{x}(A)|_{\partial X}=0$ for $A\in\Gamma({}^dTX)$
\item $\nabla^d_T|_{\partial X}={\sf v}\circ\nabla^d_T\circ{\sf v}+\phi^*\nabla_T^{B,b}
\circ{\sf n}$ on $\Gamma({}^dTX|_{\partial X})$.
\item $\frac{1}{x}\nabla^d_A\frac{\dd x}{x}(\frac{1}{x}B)|_{\partial X}=g_Z(A,
B)$ for $A, B\in \Gamma(\phitx)$
\item $\nabla^d_{\cdot}{\sf X}_d|_{\partial X}
    = {\sf v}\cdot $, the vertical
projection ${}^dTX|_{\partial X}\rightarrow {}^dV^*\partial X$
\item The tensor $C=R^{{}^dT^*X}-R^{\phitsx}\in\Gamma(\Lambda^2T^*X\otimes
\End(\phitsx))$ is given on {\em tangent} vectors $T_1, T_2$ by
$$C(T_1,T_2)|_{\partial X}= \ext_{\phi}(\widehat{x
T_1})\Int(\widehat{x T_2}) - \ext_{\phi}(\widehat{x
T_2})\Int(\widehat{x T_1})$$
\end{enumerate}
\end{lemma}
\begin{proof}
To prove (a) we will use  Lemma \ref{exactconnection} (a). Thus,
for $B\in \Gamma(\phitx)$ we have
 \begin{eqnarray*}
 \nabla^d_T\frac{\dd x}{x}(\frac{1}{x}B)&=&\nabla^{\phi}_T\frac{\dd x}{x^2}(B)
 +G(\widehat{T})\frac{\dd x}{x^2}(B/x)\\
 &=&-h(\widehat{x T},B)+g_{\phi}(\widehat{xT},B)-\frac{\dd x}{x^2}(\widehat{x T})
 \frac{\dd x}{x^2}(B)+ O(x )\\
 &=& g_Z(\widehat{x B},T)+O(x)=O(x).
\end{eqnarray*}
To prove (c), we just have to observe that the $O(x)$-term in (a)
is $O(x^2)$ if $A, B$ are in $\Gamma(\phitx)$. This should be
clear from the proof of Lemma \ref{exactconnection} (a). Part (d)
is just a reformulation of (c).
\par
To prove (e), let $A,B\in \Gamma({}^bTX)$. Then first
\begin{eqnarray*}
R^{{}^dT^*X}(A,B)-R^{\phitsx}(A,B)&=&[\frac{1}{x}\nabla^d_Ax,\frac{1}{x}\nabla^d_Bx]-
          \frac{1}{x}\nabla^d_{[A,B]}x-[\nabla^{\phi}_A,\nabla^{\phi}_B]
          +\nabla^{\phi}_{[A,B]}\nonumber \\
     &=&[\nabla^{\phi}_A,G(B)]-[\nabla^{\phi}_B,G(A)]-
     G([A,B])+[G(A),G(B)]
\end{eqnarray*}
The different terms in this expression are
\begin{eqnarray}
[\nabla^{\phi}_A,G(B)]&=&\frac{\dd
x}{x}(A)G(B)+G(\nabla_A^{\phi}B)-\ext_{\phi}(x B)\Int(\widehat{x A
})+\ext_{\phi}(\widehat{x A})\Int(x B)\label{nablaagb} \\
-[\nabla^{\phi}_B,G(A)]&=&-\frac{\dd
x}{x}(B)G(A)-G(\nabla_B^{\phi}A)+\ext_{\phi}(x A)\Int(\widehat{x B
})-\ext_{\phi}(\widehat{x B})\Int(x A) \label{nablabga}\\
 \, [G(A),G(B)]&=& [ G(\widehat{A}),G(\widehat{B}) ]=\ldots
=\ext_{\phi}(\widehat{x B}) \Int(\widehat{x
A})-\ext_{\phi}(\widehat{x A})\Int(\widehat{ x B}).\nonumber
\end{eqnarray}
Now the first summand in (\ref{nablaagb}) and the last two
summands in (\ref{nablabga}) give $$ \frac{\dd
x}{x}(A)G(\widehat{B})+\ext_{\phi}(x A)\Int(\widehat{x B
})-\ext_{\phi}(\widehat{x B})\Int(x A)= \ext_{\phi}(\widehat{x
A})\Int(\widehat{x B })-\ext_{\phi}(\widehat{x B})\Int(\widehat{x
A})$$ and we get the same contribution fom the first summand in
(\ref{nablabga}) and the last two summands in (\ref{nablaagb}).
Putting everything together gives the result.
\par
 Note that we have made extensive use of the fact that $A$ and $B$
 are tangent to the boundary, since we needed every summand to be regular in
 $\End(\phitsx)$!
\qed
\end{proof}
\par\smallskip
The expression for $C(T_1, T_2)$ only contains horizontal tensors.
Therefore it maps $$C(T_1, T_2):[{}^{\phi}V^*\partial
X]\longrightarrow x\phitsx,\qquad C(T_1, T_2):
[{}^{\phi}N^*\partial X]\longrightarrow [{}^{\phi}N^*\partial X],
$$ in general and even vanishes at the boundary whenever $T_1$ or
$T_2$ are vertical. This means that
 $C|_{\partial X}\in \phi^* \Gamma(\Lambda^2 T^*Y\otimes \End(\phitsx))$.
 Using Lemma \ref{exactconnection}(e) we get the following
\begin{corollary} \label{curvboundary}
$R^d(T_1,T_2)|_{\partial X}  = {\sf v}\circ R^{\phi}(T_1,
T_2)\circ{\sf v} \oplus \phi^*R^{B,b}(T_1,T_2)$ for $T_1, T_2\in
{}^bTX$.\qed
\end{corollary}

\subsection{Geometry of the Boundary}\label{boundarygeom}
The essentials of the geometry of fibre bundles $M=\partial
X\rightarrow Y$ are recalled in Appendix \ref{geomtens}. In this
Section we want to show how the geometry of $(X, g_{\phi})$
relates to the geometry of the boundary $\partial X$. To do this,
a metric on $\partial X$ has to be defined.
\par\smallskip
 The first step is an extension-Lemma:
Fix a normal vector field $\nu\in \Gamma(TX)$. Then any section
$A\in\Gamma(\partial X,\phitx|_{\partial X})$ can be extended to
the interior by parallel transport along $\nu$. More concretely,
in a small neighbourhood $U(\partial X)$ of the boundary, we write
$\overline{A}$ for the extension of $A$ into $\phitx$  such that
$\frac{1}{x}\overline{A}$ is parallel w.r.t. $\nabla^{d}$ along
$\nu$. This choice of extension may seem a little odd at this
stage, but will prove to be useful in Proposition
\ref{rdatthedel} and Chapter 5. In any case it
has the following nice properties
\begin{lemma}[{\boldmath Special Extension}]
\label{extension}  Let $\nu$, $\nu_1$, $\nu_2$ be normal vector
fields in $TX$ and index the associated extensions accordingly.
Let $N\in \Gamma(TX)$ and $V\in\Gamma(\partial X,V\partial X)$,
$A\in\Gamma(\partial X,\phindx)$.
\begin{enumerate}
\item $\overline{A}\in \Gamma(U(\partial X),[\phindx])$ and
$\overline{V}\in \Gamma(U(\partial X),[{}^{\phi}V\partial X])$
\item $\nabla^{d}_N\frac{1}{x}\overline{A}\in
\Gamma(U(\partial X),[{}^{\rm d}N\partial X])$ and
$\nabla^{d}_N\frac{1}{x}\overline{V}\in \Gamma(U(\partial X),
[{}^{\rm d }V\partial X])$.
\item $\nabla^{\phi}_N\overline{A}\in \Gamma(U(\partial X),[\phindx])$
and $\nabla^{\phi}_N\overline{V}\in \Gamma(U(\partial
X),\,[{}^{\phi}V\partial X])$.
\item  $\frac{1}{x}(\overline{A}^1-\overline{A}^2)\in
\Gamma(U(\partial X),[\phindx])$ and
$\frac{1}{x}(\overline{V}^1-\overline{V}^2)\in \Gamma(U(\partial
X),[{}^{\phi}V\partial X])$
\item $\frac{1}{x}\overline{A}|_{\partial X}\in
{}^bTX|_{\partial X}$ is independent of the choice of $\nu$.
\end{enumerate}
\end{lemma}
\begin{proof}
(a) is clear by construction, (b) and (c) follow from Lemma
\ref{exactconnection}(c) and (\ref{gnice}). Consider (d): Taking
two extensions $\overline{V}^1, \overline{V}^2\in \Gamma(\phitx)$
of $V\in \Gamma(V\partial X)$ it is clear that $\overline{V}^1-
\overline{V}^2\in O(x)\Gamma(\phitx)$. Thus $$\frac{1}{x}(
\overline{V}^1-\overline{V}^2)= \frac{1}{Nx}\nabla_N^{\phi}(
\overline{V}^1-\overline{V}^2)\in\Gamma([{}^{\phi}V\partial X])$$
by (c).\qed
\end{proof}\par
\smallskip
As a consequence, this construction yields  a canonical choice of
horizontal space in ${}^bTX|_{\partial X}$ as the span $$
{}^bH\partial X = \langle \frac{1}{x}\overline{A}|_{\partial
X,{{}^bTX}}\,|\quad A\in \Gamma(\partial X,\phindx) \rangle. $$
The induced horizontal and vertical projections in $T\partial
X=TM$ and ${}^bTX|_{\partial X}$ will be denoted by ${\sf h}$ and
${\sf v}$. Using these and the boundary defining function $x$ we
can now define the maps
\begin{eqnarray*}
\delta_{h}&:& \phitx|_{\partial X}\longrightarrow
{}^bTX|_{\partial X}\qquad A\longmapsto [\frac{1}{x}\overline{{\sf
n}A}+{\sf v}A]|_{\partial X,{}^bTX}
\\ \delta_{h}^{-1}&:&{}^bTX|_{\partial X}\longrightarrow \phitx|_{\partial X} \qquad
A\longmapsto x[{\sf h}A]_b|_{\partial X,\phitx}+[{\sf
v}A]_b|_{\partial X,TX}
\\ & &\qquad\qquad\qquad\qquad\qquad\qquad\qquad\qquad=x[A]_b|_{\partial X,\phitx}+[{\sf
v}A]_b|_{\partial X,TX},
\end{eqnarray*}
where $[A]_b$ is an extension of $A\in  {}^bTX|_{\partial X}$ to
${}^bTX$.
 The mapping
properties of $\delta_h$ and the justification of the
``inverse''-notation are given in
\begin{lemma}\label{deltah}
\begin{enumerate}
\item $\delta_{h}\in \Gamma(\partial X, \phitsx\otimes {}^bTX|_{\partial X})$
\item $\delta_{h}\circ \delta_{h}^{-1}=\Id_{{}^bTX}$, and \quad
   $\delta_{h}^{-1}\circ \delta_{h}=\Id_{\phitx}$
\item $\delta_h^{-1}(TM)=\langle\frac{\dd x}{x^2}\rangle^{\perp}\subset
\phitx$, and  \quad $\delta_h^{*}(T^*M)={}^eT^*X$.\qed
\end{enumerate}
\end{lemma}
We can now simply define the metric on the fibres of
${}^bTX|_{M}\rightarrow M$ as
$$g_{M,b}=(\delta_{h}^{-1})^{*}(g_{\phi}).$$ Using the
identification of $TM$ with the kernel of $\frac{\dd x}{x}$ in
${}^bTX|_{\partial X}$ we can also restrict the metric $g_{M,b}$
to the metric $g_M$ on $M$. However, since it involves no extra
work, we prefer to stick with the bundle ${}^bTX|_M$. Thus, the
Koszul formula defines the ``Levi-Civita''-connection
$\nabla^{M,b}$ associated to this metric. To be able to use this
definition in calculations, we need a more detailed description of
$g_{M,b}$. Parts (b) and (c) of the following Lemma will therefore
be of particular importance:
\begin{lemma}[\boldmath Metric on the Boundary]\label{metbound}
\quad \par \noindent
\begin{enumerate}
\item $g_{M,b}^{-1}=\delta_{h}(g_{\phi}^{-1})$
\item $g_{M,b}(V,T)=g_{\phi}(\overline{V},T)$ for $V\in VM$, $T\in {}^bTX$
\item $g_{M,b}(T_1,T_2)=g_{\phi}({\sf v}T_1,{\sf v}T_2)+g_d(T_1,T_2)$ for $T_1, T_2
\in {}^bTX$
\item $\frac{1}{x}\overline{A}|_{\partial X}$ for $A\in \Gamma(\phindx)$
 is horizontal w.r.t. $g_{M,b}$.
\end{enumerate}
\end{lemma}
\begin{proof}
Part (a) is straightforward.  To prove (b), we first show that
$g_{\phi}(\overline{V},T)$ is well defined by showing that it is
well defined for the different terms in the decomposition of
$g_{\phi}$. To do this note first that $$ N\frac{\dd x
}{x^2}(\overline{V})= \nabla^{d}_N\frac{\dd
x}{x}(\frac{1}{x}\overline{V})+\frac{\dd
x}{x}(\nabla^{d}_N\frac{1}{x}\overline{V})=O(x),$$ i.e. the
expression $\frac{\dd x}{x^2}(\overline{V})$ is of order $O(x^2)$.
Thus
 $$\frac{\dd x}{x^2}\otimes \frac{\dd
x}{x^2}(\overline{V},T)=O(x^0)\dd x (T)=O(x).$$ To see that the
expression $h(\overline{V},T)$ only depends on $V,T|_{\partial
X}$, note that the extension $\overline{V}$ pairs with any element
$\varphi\in \Gamma([{}^{\rm e}N^*\partial X])$ to $O(x)$ and there
is a natural map $x[{}^{\rm e}N^*\partial X] \rightarrow T^*X$.
But the pairing of $T$ with any element in $\Gamma(\partial
X,T^*X)$ clearly only depends on $T|_{\partial X}$. The
corresponding statement for $g_Z$ is clear.
\par
 To prove the equality in (b) note that
it is trivially true for vertical $T$. On the other hand $$
\langle \overline{V}, {\sf h}[T]\rangle_{\phi}= \langle
\overline{V}, \frac{1}{x}\overline{A} \rangle_{\phi}= \frac{1}{Nx}
N \langle \frac{1}{x}\overline{V},
\frac{1}{x}\overline{A}\rangle_{d}=
\frac{1}{Nx}\left[\langle\nabla^{d}_N\frac{1}{x}\overline{A},
\frac{1}{x}\overline{V}\rangle_{d}+\langle\frac{1}{x}\overline{A},
\nabla^{d}_N\frac{1}{x}\overline{V}\rangle_{d}\right]$$ is $0$ by
definition of the extension. This also proves (c) and (d). \qed
\end{proof}
\par
As an aside, note that part (c) of the Lemma implies that the
metric $g_{M,b}$ (or $g_M$) is really of the type
(\ref{specialmetric}) required in Appendix \ref{geomtens}, since
$$g_{M,b}({\sf h}T_1,{\sf h}T_2)=g_{d}({\sf h}T_1,{\sf h}T_2)=
\phi^*g_{B,b}(T_1,T_2).$$ The following Lemma lists some ways to
calculate the horizontal and vertical projections
\begin{lemma}\label{halpha}
Let $\alpha\in \Gamma({}^bT^*X)$, $N\in\Gamma(TX)$ and $T\in
\Gamma({}^bTX)=\Gamma([{}^dN^*\partial X])$. Then $x\alpha$ is a section in
$[{}^dV^*\partial X]$ and
\begin{enumerate}
\item
$\nabla_N^d(x\alpha)(T)|_{\partial X}=
\nabla_N^{\phi}\alpha(xT)|_{\partial X}=\dd x(N)\alpha({\sf
h}T)|_{\partial X}.$
\item $x\nabla_N^{\phi}T|_{\partial X,{}^bTX}=-\dd x(N){\sf h}T$
\item $x\nabla_N^{d}T|_{\partial X,{}^bTX}=\dd x(N){\sf v}T$
\end{enumerate}
\end{lemma}
\begin{proof}
We use the method developed in Section \ref{exactphisection}. For
simplicity, assume that $\alpha\in\Gamma(T^*X)$. Also write
$A:=(g_{\phi})^{-1}(\alpha)$. Then
\begin{eqnarray*}\nabla^d_N(x\alpha)(T)&=&
\nabla^{\phi}_N\alpha(xT)+O(x)\\ &=&(L_Ag_{\phi})(N,xT)+\dd
\alpha(N,xT)+O(x)=(L_Ag_{\phi})(N,xT)+O(x).
\end{eqnarray*} The
different parts of the Lie derivative of the metric can be
calculated as in Lemma \ref{exactconnection}:
\begin{eqnarray*}
L_A(\frac{\dd x}{x^2}\otimes\frac{\dd x}{x^2})(N,xT)&=&
O(x^{-1})\dd x\otimes\dd x(N,xT) +O(x)=O(x)
\\L_Ag_Z(N,xT)&=& O(x)\\
L_Ah(xN,T)&=& (Nx)h(A,T)=(Nx)(\alpha(T)-g_Z(A,T))=(Nx)\alpha({\sf
h}T)
\end{eqnarray*}
This proves the claim for $\alpha\in\Gamma(T^*X)$. It is easy to
see that it also holds for forms of the type $\alpha=f\frac{\dd x
}{x}$. Parts (b) and (c) are just reformulations of (a).\qed
\end{proof}
\par \smallskip
 We now want to describe
 $\nabla^{M,b}$ and its associated tensors in terms of $g_d$ and
 $\nabla^d$. We  start with  the following
\begin{definition}\label{Sphi}
\begin{enumerate}
\item The curvature $\Omega_{\phi}\in
  C^{\infty}(\partial X,\Lambda^2\,\phitsx \otimes V\partial X)$
is defined by $$\Omega_{\phi}(A,B):={\sf
v}[\frac{1}{x}\overline{{\sf n}A},\frac{1}{x}\overline{{\sf
n}B}]|_{\partial X}=\Omega_M(\frac{1}{x}\overline{{\sf n
}A},\frac{1}{x}\overline{{\sf n}B}),\quad\mbox{ for $A, B\in
\phitx$}$$
\item The tensor $S_{\phi}\in C^{\infty}(\partial X,V^*\partial X\otimes \phinsdx
\otimes V\partial X)$ is defined by $$ S_{\phi}(V)A:={\sf v
}\frac{1}{x}\nabla^{\phi}_{\overline{V}}\overline{A}|_{\partial X
},\quad\mbox{ for $V\in V\partial X$, $A\in \phindx$}.$$ It could
be called the $\phi$-second fundamental form.
\item A variant of (b) is $S_{d}\in C^{\infty}(\partial X,V^*\partial X\otimes
{}^dN^*\partial X \otimes V\partial X)$, defined by $$
S_{d}(V)\frac{1}{x}A:={\sf v
}\nabla^{d}_{\overline{V}}\frac{1}{x}\overline{A}|_{\partial X
}\stackrel{(\ref{nabladphivector})}{=}S_{\phi}(V)A+\frac{\dd
x}{x^2}(A)V,\quad\mbox{ for $V\in V\partial X$, $A\in \phindx$}.$$

\item $B_{\phi}\in C^{\infty}(\partial X,\Lambda^2 V^*\partial X
\otimes \phindx)$ is defined by $$ B_{\phi}(V,W) :={\sf
n}\frac{1}{x}\left( [\overline{V},\overline{W}]- \overline{[V,W]}
\right),\quad\mbox{ for $V,W\in V\partial X$}.$$
\end{enumerate}
\end{definition}
\par
\smallskip
The reader should verify that $B_{\phi}$ is indeed a well defined
tensor: Since $[\overline{V},\overline{W}]$ and $\overline{[V,W]}$
have the same boundary value as elements in $\phitx$, their
difference
\begin{equation}\label{diffvw}
 [\overline{V}, \overline{W}]-\overline{[V,W]} \mbox{\quad lies in \quad}
 x\phitx.
\end{equation}
The nonvanishing of the tensor
$B_{\phi}$ is linked the fact that $g_{\phi}$ might not be a
product $\phi$-metric near the boundary! In the product case (\ref{diffvw})
vanishes near the boundary and therefore $B_{\phi}=0$.
\par
For a general exact $\phi$-metric we have the following analogue
of Lemma \ref{omegaform}:
\begin{lemma}\label{lemnablab}
Let $T_0\in\Gamma(TM), T_1, T_2\in\Gamma(M,{}^bTX|_M)$. Then
\begin{eqnarray*}
&&\langle\nabla^{M,b}_{T_0}T_1,T_2\rangle_{M,b}=\langle\nabla_{T_0}^{\phi}
\overline{{\sf v}T_1},\overline{{\sf v}T_2}\rangle_{\phi}+
\langle\nabla_{T_0}^{d} \overline{{\sf h}T_1},\overline{{\sf
h}T_2}\rangle_{d}
\\ &&\qquad\qquad +\frac{1}{2}\left[
\langle \Omega_{\phi}(\overline{{\sf h}T_0},\overline{{\sf
h}T_1}),\overline{{\sf v}T_2}\rangle_{\phi} -\langle
\Omega_{\phi}(\overline{{\sf h}T_1},\overline{{\sf
h}T_2}),\overline{{\sf v}T_0}\rangle_{\phi} -\langle
\Omega_{\phi}(\overline{{\sf h}T_2},\overline{{\sf
h}T_0}),\overline{{\sf v}T_1}\rangle_{\phi}\right]
\\ &&\qquad\qquad+\langle S_{\phi}(\overline{{\sf v}T_0})\overline{{\sf
h}T_1},\overline{{\sf v}T_2}\rangle_{\phi} -\langle
S_{\phi}(\overline{{\sf v}T_0})\overline{{\sf
h}T_2},\overline{{\sf v}T_1}\rangle_{\phi}
\\ &&\qquad\qquad +\frac{1}{2}\left[
\langle B_{\phi}(\overline{{\sf v}T_1},\overline{{\sf
v}T_2}),x\overline{{\sf h}T_0}\rangle_{\phi} -\langle
B_{\phi}(\overline{{\sf v}T_2},\overline{{\sf
v}T_0}),x\overline{{\sf h}T_1}\rangle_{\phi} -\langle
B_{\phi}(\overline{{\sf v}T_0},\overline{{\sf
v}T_1}),x\overline{{\sf h}T_2}\rangle_{\phi}\right]
\end{eqnarray*}
\end{lemma}
\begin{proof}
This is an application of the Koszul  formula using  Lemma
\ref{metbound}(c). We only give an example of how the calculation
differs from the one in Lemma \ref{omegaform}:
\begin{eqnarray*}
2\langle S_{\phi}(V)A,W\rangle_{\phi} &=&
-\overline{W}\langle\overline{V},\frac{1}{x}\overline{A}\rangle_{\phi}
+\overline{V}\langle\frac{1}{x}\overline{A},\overline{W}\rangle_{\phi}
+\frac{1}{x}\overline{A}\langle\overline{W},\overline{V}\rangle_{\phi}\\
&
&\qquad-\langle\overline{V},[\frac{1}{x}\overline{A},\overline{W}]\rangle_{\phi}
+\langle
\frac{1}{x}\overline{A},[\overline{W},\overline{V}]\rangle_{\phi}
+\langle\overline{W},[\overline{V},\frac{1}{x}\overline{A}]\rangle_{\phi}
\\
 &=&
-\overline{W}\langle\overline{V},\frac{1}{x}\overline{A}\rangle_{M,b}
+\overline{V}\langle\frac{1}{x}\overline{A},\overline{W}\rangle_{M,b}
+\frac{1}{x}\overline{A}\langle\overline{W},\overline{V}\rangle_{M,b}\\
&
&-\langle\overline{V},[\frac{1}{x}\overline{A},\overline{W}]\rangle_{M,b}
+\langle
\frac{1}{x}\overline{A},[\overline{W},\overline{V}]\rangle_{M,b}
+\langle\overline{W},[\overline{V},\frac{1}{x}\overline{A}]\rangle_{M,b}
+\langle
\frac{1}{x}\overline{A},[\overline{W},\overline{V}]\rangle_{\phi}
\\
&=&2\langle S_{M}(V)\frac{1}{x}\overline{A},W\rangle_{M,b}
+\langle
\frac{1}{x}\overline{A},[\overline{W},\overline{V}]\rangle_{\phi}
\\ &=&2\langle S_{M}(V)\frac{1}{x}\overline{A},W\rangle_{M,b}-\langle
B_{\phi}(V,W),A\rangle_{\phi},
\end{eqnarray*}
and so on. \qed
\end{proof}
\par\smallskip\noindent
As a special case of this Lemma, note that
\begin{equation}\label{minorrole}
\langle\nabla^{M,b}_{T_0}{\sf v }T_1,{\sf v
}T_2\rangle_{M,b}=\langle\nabla_{T_0}^{\phi} \overline{{\sf
v}T_1},\overline{{\sf v}T_2}\rangle_{\phi}+\frac{1}{2}\langle
B_{\phi}(\overline{{\sf v}T_1},\overline{{\sf
v}T_2}),x\overline{{\sf h}T_0}\rangle_{\phi}.
\end{equation}
This will play a minor role in the calculation of
(\ref{minorrole2}).
\par\medskip
We can now calculate parts of the curvature $R^d$ at the boundary
more precisely.
\begin{proposition}\label{rdatthedel}
Let $N\in \Gamma(TX)$ and $W,A,B\in \Gamma(\phitx)$. Then at the boundary
\begin{enumerate}
\item $\langle R^d(W,N)A,B\rangle_{\phi}=
           \langle R^d(W,N)\frac{1}{x}A,\frac{1}{x}B\rangle_{d}$
           \par\smallskip
$ \qquad\qquad = \langle R^d(W,N){\sf v}A,{\sf v}B\rangle_{\phi}
       +\dd x(N)[ \langle S_d(W)\frac{1}{x}{\sf n}B,{\sf
       v}A\rangle_{\phi}-\langle S_d(W)\frac{1}{x}{\sf n}A,{\sf
       v}B\rangle_{\phi}]$
\item $\langle R^d(W,N){\sf n}A,{\sf n}B\rangle_{\phi}$ is
$O(x)$ and  the corresponding coefficient depends on the
choice of extension of  ${\sf n}A$, ${\sf n}B$
into the interior. Still, one has
$$ \langle (\nabla^d_NR^d)(W,N){\sf n} A,{\sf n}B\rangle_{\phi}
=(\dd x(N))^2 \langle \Omega_{\phi}(A,B), W\rangle_{\phi}.
$$
\end{enumerate}
\end{proposition}
\begin{proof}
First note that the LHS in (a) is $O(x)$, whenever
$W\in\Gamma([\phindx])$. We thus  can assume right away that
$W\in\Gamma([{}^{\phi}V\partial X])$.
Then, using the special extension, w.r.t. the normal vector field
$N$, we get
\begin{eqnarray*}
\langle R^d(W,N)\frac{1}{x}{\sf n}A,\frac{1}{x}{\sf v}B\rangle_d &=& \langle
R^d(\overline{W},N)\frac{1}{x}\overline{{\sf n}A},
\frac{1}{x}\overline{{\sf v}B}\rangle_d = -N \langle
\nabla^{d}_{\overline{W}}\frac{1}{x}\overline{{\sf n}A},\frac{1}{x}
\overline{{\sf v}B}
\rangle_{d}\\ &=& -\dd x(N)\langle
\nabla^{d}_{\overline{W}}\frac{1}{x}\overline{{\sf n}A},\overline{{\sf v}B}
\rangle_{\phi}=-\dd x(N)\langle
S_d(\overline{W})\frac{1}{x}\overline{{\sf n}A},\overline{{\sf
v}B}\rangle_{\phi}.
\end{eqnarray*}
Furthermore, the LHS,
\begin{equation}\label{rdproof}
\langle R^d(W,N)\frac{1}{x}A,\frac{1}{x}B\rangle_d = \langle
R^d(A,N)\frac{1}{x}B, \frac{1}{x}W\rangle_d -\langle
R^d(B,N)\frac{1}{x}A, \frac{1}{x}W\rangle_d,
\end{equation}
is of order $O(x)$, whenever $A, B$ are in $\Gamma([\phindx])$.
This finishes the proof of (a). To prove (b), assume that $A, B$
are in $\Gamma([\phindx])$. Then
\begin{eqnarray*}
\langle (\nabla_N^dR^d)(W,N)\frac{1}{x}A,\frac{1}{x}B\rangle_d &=&
N\langle R^d(\overline{W},N)\frac{1}{x}\overline{A},\frac{1}{x}
\overline{B}\rangle_d
\\ &=&
(Nx)\left[ \langle
R^d(\frac{1}{x}\overline{A},N)\frac{1}{x}\overline{B},
\frac{1}{x}\overline{W}\rangle_d -\langle
R^d(\frac{1}{x}\overline{B},N)\frac{1}{x}\overline{A},
\frac{1}{x}\overline{W}\rangle_d \right]
\\ &=& (Nx)N \left[
\langle
\nabla^{d}_{\frac{1}{x}\overline{A}}\frac{1}{x}\overline{B},
\frac{1}{x}\overline{W}\rangle_{d}- \langle
\nabla^{d}_{\frac{1}{x}\overline{B}}\frac{1}{x}\overline{A},
\frac{1}{x}\overline{W}\rangle_{d} \right]
\\ &&\qquad\qquad\qquad +(Nx)\langle
\nabla^{d}_{[\frac{1}{x}\overline{A},N]}\frac{1}{x}\overline{B}
-\nabla^{d}_{[\frac{1}{x}\overline{B},N]}\frac{1}{x}\overline{A},
\frac{1}{x}\overline{W}\rangle_{d}
\\ &=& (Nx)^2 \langle
[\frac{1}{x}\overline{A},\frac{1}{x}\overline{B}],\overline{W}
\rangle_{\phi} =(Nx)^2 \langle \Omega_{\phi} (A,B), W
\rangle_{\phi}.
\end{eqnarray*}
This finishes the proof.\qed
\end{proof}
\par
\smallskip
Using Definition
\ref{Sphi} and (the proof of) Lemma \ref{lemnablab}, the RHS in
Proposition \ref{rdatthedel}(a) can be further decomposed
into contributions from the intrinsic and the extrinsic geometry
of the boundary fibration. For instance, the term
$$ \langle R^d(W,\nu){\sf v}A,{\sf v}B\rangle_{\phi}|_{\partial X}=
 \langle R^{\phi}(W,\nu){\sf v}A,{\sf v}B\rangle_{\phi}|_{\partial
 X},$$
 where $\nu=\frac{1}{x}{\sf X}_d$, should  be considered as
 a tensor describing some aspect of the extrinsic geometry, since
 it
 vanishes whenever $g_{\phi}$ is a $\phi$-product metric.

\newpage
\section{The $\phi$-Calculus}
\par
As before, $X$ denotes a manifold with fibred boundary. By the
Schwartz kernel theorem, continuous operators from
$\cdotinfty{}(X)$ to $C^{-\infty}(X)$ are represented by elements
of $C^{-\infty}(X^2)$. In this Chapter we will give an overview of
the ideas developed in \cite{MelMazmwfb}  and \cite{Meaomwc},
\cite{Melaps} for the analysis of such operators whose kernels are
distributions conormal to the diagonal and to the boundaries in a
certain sense. The composition formula for these $\phi$-
pseudodifferential operators proved in Section
\ref{compositionformula} is an extension of the corresponding
result in \cite{MelMazmwfb}.
\subsection{The $b$- and the $\phi$-Blow up}\label{bandphisection}
\par
Denote by $\beta_{L}, \beta_{R}$ the left and right projections
from $X^2$ to $X$. Given two sets of local coordinates on $X$ near
the boundary $x, {\bf y}, {\bf z}$ and $x', {\bf y}', {\bf z}'$,
as in Section \ref{secphimanifolds} we get coordinates on $X^2$ by
lifting the first set from the left and the second set from the
right to $X^2$: $$ x\sowie \beta_{L}^*x, \quad{\bf y}\sowie
\beta_{L}^*{\bf y}, \quad{\bf z} \sowie \beta_{L}^*{\bf z},
\qquad\mbox{and}\qquad x'\sowie \beta_{R}^*x', \quad{\bf y}'\sowie
\beta_{R}^*{\bf y}', \quad {\bf z}'\sowie \beta_{R}^*{\bf z}'.$$ A
priori the two coordinate sets on $X$ can be completely
independent. It will usually be assumed though, that $x=x'$ is the
fixed boundary defining function on $X$. Sometimes, depending of
the geometrical situation -- see the proof of Lemma \ref{philift}
below, we also assume that ${\bf y}={\bf y}'$ and ${\bf z}={\bf
z}'$ on $X$. In that case, we will say that the coordinates
 ${\bf y}, {\bf y}'$ or ${\bf z}, {\bf z}'$ on $X^2$ are {\em paired}.
\par
 For an
easier analysis of the different classes of pseudodifferential
operators on $X$ we introduce: The $b$-kernel space is defined as
the blow up (compare \cite{Melaps}) $$ X^2_b:= [X^2, \partial
X\times
\partial X].$$ The corresponding blow down map is denoted
$\beta_b:X^2_b\rightarrow X^2$, and the induced left and right
projections are $\beta_{b,L}$ and $\beta_{b,R}$. We write $$
\bigtriangleup_b:= \overline{\beta_b^{-1}(\bigtriangleup-\partial
      \bigtriangleup)},\quad
 \fabf := \beta_b^{-1}(\partial X\times\partial X).$$\par
\bigskip
\centerline{\epsfig{file=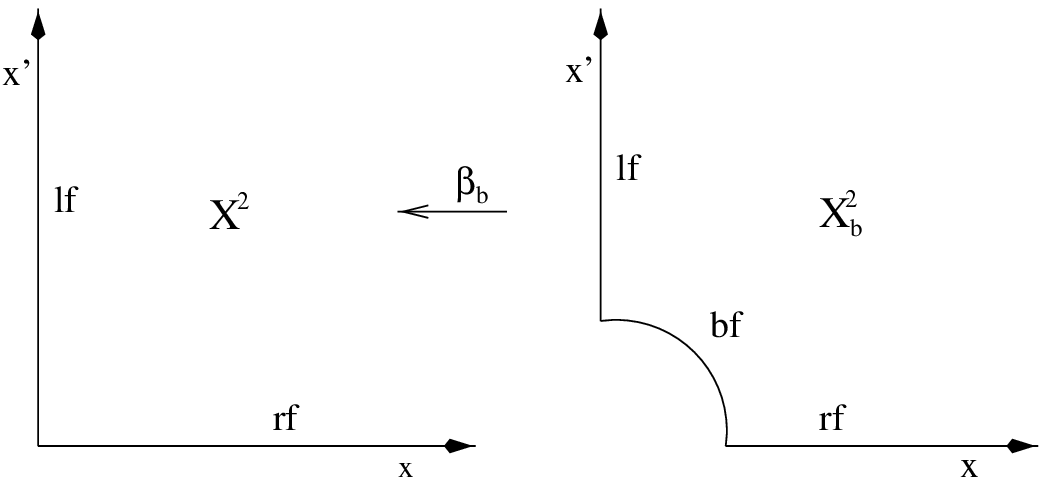}}\par\bigskip
\nopagebreak\centerline{Figure 1}
\par\bigskip
The face $\fabf$ is diffeomorphic to $\partial X\times\partial
X\times [-1,1]_{\sigma}$. Here, the coordinate $\sigma$ is given
by $\sigma=\frac{x-x'}{x+x'}$. The {\em fibre diagonal} $\partial
X\times_Y\partial X$ is the preimage of the diagonal
$\bigtriangleup_Y$ under the map $(\phi,\phi):\partial
X\times\partial X\rightarrow Y\times Y$. Its lift in $\fabf$ is
denoted by $$ F_{\phi}:=\beta_b^{-1}(\partial X\times_Y
\partial X)\times \{\sigma=0\} \cong
\partial X\times_Y \partial X. $$
Note that this definition depends on our (fixed) choice for the
boundary defining function $x$. Blowing up $F_{\phi}$ in $X_b^2$
yields the kernel space for the $\phi$-pseudodifferential
operators (compare \cite{MelMazmwfb})
$$X^2_{\phi}:=[X^2_b,F_{\phi}].$$ Again, we write
$\beta_{\phi}:X^2_{\phi}\rightarrow X^2$ for the blow down map and
$\beta_{\phi,L}$ and $\beta_{\phi,R}$ for the induced left and
right projections onto $X$. Also, the lifted diagonal and the
blown up faces are denoted by $$ \bigtriangleup_{\phi}:=
\overline{ \beta_{\phi}^{-1}(\bigtriangleup-
\partial \bigtriangleup) } ,\quad \faff := \beta_{\phi,b}^{-1}F_{\phi},
\quad \faphibf := \overline{\beta_{\phi,b}^{-1}(\fabf-F_{\phi})}.$$
\par
\bigskip
\centerline{\epsfig{file=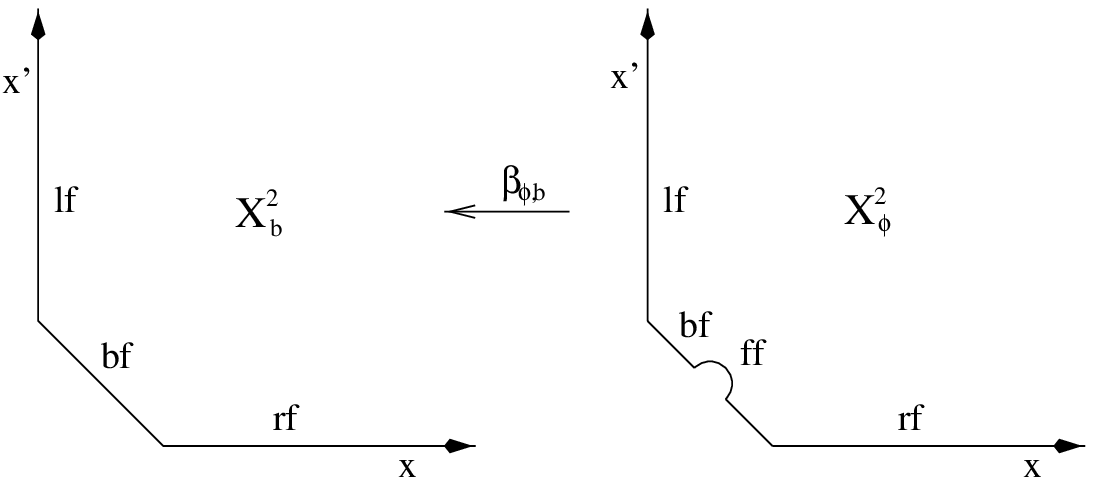}}\par\bigskip \centerline{Figure
2}
\par\bigskip
\noindent As a general rule, given a boundary face $F$, we denote
by $\rho_F$ a boundary defining function for $F$. The following
Lemma describes the fundamental property of the blow up
$X^2_{\phi}$:
\begin{lemma}\label{philift}
\begin{enumerate}
\item Restricted to $\bigtriangleup_{\phi}$ the spaces
$\beta_{\phi, R}^*\phiv$ and $\beta_{\phi, L}^*\phiv$ each span
the normal bundle $N\bigtriangleup_{\phi}$. Especially, recalling
that $\bigtriangleup_{\phi}$ is canonically diffeomorphic to $X$,
we have $TX_{\phi}^2|_{\bigtriangleup_{\phi}} \equiv TX\oplus
\phitx$.
\item  $\beta_{\phi,R}^*W$, for
$W\in \Gamma(\phitx)$,  is a $\phi$-vector field at $\falf$ and
$\faphibf$ and is a $b$-vector field at $\faff$ and $\farf$.
\item $\stackrel{\circ}{\faff}$ is canonically diffeomorphic to
$\phindx\times_Y \partial X$ and $\faff$ to $\RC_{\rm
fib}(\phindx\times_Y \partial X)$.
\end{enumerate}\par\smallskip
\end{lemma}
\begin{proof}
(a) and (b) can be proved by calculating the lift of vector fields
in local coordinates. First, taking local coordinates
{\bf\boldmath near the corner in $X^2$}, we can introduce local
coordinates on $X^2_b$, projective  w.r.t. $x'$:
\begin{equation}\label{bfcoordinates}
s=\frac{x}{ x'},\quad x',\quad{\bf y},\quad{\bf y'},\quad{\bf z},
\quad {\bf z'} \end{equation} These coordinates are valid  near
$\fabf$ but away from $\farf$. There we calculate $$ \dd s
=\frac{\dd x}{x'}-\frac{1}{s x'}\dd x'\quad \mbox{and}\quad
 \beta_{b,R}{}^*x^2 \frac{\partial}{\partial x} =x's^2
 \frac{\partial}{\partial s}, \quad \beta_{b,R}{}^* x\frac{\partial}
{\partial {\bf y}}=x's \frac{\partial}{\partial {\bf y}}.$$
{\bf\boldmath Near the fibre diagonal $F_{\phi}$} we can assume
that the coordinates ${\bf y}$ and ${\bf y}'$ are paired, thus
$F_{\phi}$ is described by $\{x'=0, s=1, {\bf y} ={\bf y'} \}$.
Near $\faff$, but away from $\fabf$, we can introduce coordinates
on $X^2_{\phi}$, projective w.r.t. $x'$:
\begin{equation}\label{ffcoordinates}
x',\quad S=\frac{s-1}{x'}, \quad {\bf Y}=\frac{{\bf y}-{\bf
y'}}{x'}, \quad{\bf y'},\quad{\bf z},\quad {\bf z'}\quad\mbox{and}
\end{equation}
 $$ \dd S=\frac{\dd s}{x'}-\frac{s-1}{(x')^2}\dd x',\quad
    \dd {\bf Y}=\frac{\dd {\bf y}-\dd {\bf y}'}{x'}
    -\frac{{\bf y}-{\bf y}'}{(x')^2}\dd x'. $$
From this we get:
\begin{equation}\label{pairedxy}
\beta_{\phi,R}{}^*x^2\frac{\partial}{\partial x} =
   (1+x'S)^2\frac{\partial}{\partial S},\quad
   \beta_{\phi,R}{}^*x\frac{\partial}{\partial {\bf y}}=
    (1+x'S)\frac{\partial}{\partial {\bf Y}}.
    \end{equation}
From these calculations (b) is obvious. Finally, near the lifted
diagonal $\bigtriangleup_{\phi}$ we can add the condition that the
coordinates ${\bf z}$ and ${\bf z}'$ are also paired. Then
$\bigtriangleup_{\phi}$ is just $\{x'=0, S=0, {\bf Y}=0\}, {\bf
z}={\bf z}'$ and  the assertion of (a) also follows from
(\ref{pairedxy}).
\par
 For (c), we just mention that the map from
$\phindx\times_Y\partial X$ to the interior of $\faff$ is given as
follows. A nonzero element $W\in \Gamma(\phindx)$ lifts via
$\beta_{\phi,L}$ to a nonvanishing vector field
$\beta_{\phi,L}^*W$  on $\faff$, which generates a flow
$\Phi_W(t)$ on that face. This flow is identified with the
linear flow generated by $W$ in the fibres of $\phindx$. \qed
\end{proof}
\par\smallskip
 The {\em small} calculus of $b$-pseudodifferential operators on $X$
 is defined as the space
$$ \Psi^m_{b,{\rm cl}}(X):= \cdotinfty{\fabf}I^m_{\rm cl
}(X^2_{b},\bigtriangleup_{b}; \beta_{b, R}^* {}^b\Omega(X)),$$ of
distributions  conormal  to $\bigtriangleup_b$, but $C^{\infty}$
everywhere else, and vanishing to infinite order at all the faces
in $X^2_b$, except $\fabf$ (compare \cite{Meaomwc} and  Remark
\ref{remdefconormal} for an explanation of these notions). The
{\em small} calculus of $\phi$-pseudodifferential operators on $X$
is defined as $$ \Psi^m_{\phi, {\rm cl}}(X):=
\cdotinfty{\faff}I^m_{\rm cl }(X^2_{\phi},\bigtriangleup_{\phi};
\beta_{\phi, R}^* {}^\phi\Omega(X)).$$ Note that we {\em always}
assume our conormal distributions to be classical at the diagonal.
The ``${\rm cl}$''-notation in the above definitions refers to the
fact that the distributions are also assumed to be $C^{\infty}$ at
the faces $\fabf$ or $\faff$. More generally, given a
$C^{\infty}$-index set ${\cal E}=(E_{\faff}, E_{\faphibf},
E_{\farf},E_{\falf})$, see Appendix \ref{appzeronotation}, we can
define the corresponding space in the big calculus of
$\phi$-pseudodifferential operators by $$ \psiphi^{m,{\cal
E}}(X):= {\cal A}^{\cal E}I_{\rm cl}^m(X^2_{\phi},
\bigtriangleup_{\phi}; \beta_{\phi, R}^* {}^\phi\Omega(X)),$$ and
ditto for b-operators. These definitions also make
 sense for finite index sets, where the index of lowest order is
 interpreted as a conormal bound.
 Note that with this notation
 $$ \Psi_{\phi,{\rm cl}}^m(X)=\psiphi^{m,{\cal E}}(X),\quad
 E_{\faff}=m+\NN,\quad E_{\falf}=E_{\farf}=E_{\faphibf}=\infty,$$
 and the similar statement for $b$-operators is left to the
 reader.
\par
\smallskip
An advantage of choosing the density bundle to be lifted from the
right factor is that for any $\nu\in \Gamma({}^\phi\Omega(X))$ and
$W\in \Gamma(X,\phitx)$ we have $$ L_{\beta_{L}^*W} \beta_{R}^*\nu
=0 \quad \mbox{\rm on}\quad X^2, \qquad L_{\beta_{\phi,L}^*W}
\beta_{\phi,R}^*\nu =0\quad\mbox{\rm on}\quad X^2_{\phi}.$$ Thus,
the lifts to $X^2_{\phi}$ of vector fields in $\phiv$ from the
left factor act on elements in $\psiphi^{m,\cal E}(X)$ (as Lie
derivatives) in a particularly simple manner.
\par
\smallskip
The definitions above can be extended to operators between
sections of vector bundles in the usual manner. In the case of the
Dirac operator on a ($\phi$-)Clifford bundle $E$ we introduce the
``coefficient bundle''
$$\CB_{\phi}=\beta_{\phi,L}^*E\otimes\beta_{\phi,R}^*(E^*\otimes\phiomega(X))$$
on $X^2_{\phi}$ and the analogous definition on $X^2_b$ to
simplify the notation.

\subsection{Mapping Properties}\label{secmapprop}
In order to describe elementary mapping properties of our
pseudodifferential operators we use the pullback and pushforward
results for conormal distributions on blown up spaces as described
in \cite{Meaomwc} (see also \cite{Paul} or Appendix
\ref{appzeronotation}). First, we need a little
\begin{lemma}[Densities and Blow up]\label{liftdensity}
Let $N$ be a manifold with corners, and $F\subset \partial Y$ a
p-submanifold. Let $S_q\subset N^*F$ be a subspace of
$q$-parabolic directions. Then the blow down map
$\gamma:N_F=[N;F,S_q]\longrightarrow N$ has the property
$$\gamma^*\Omega(N)= \rho_F^{n(F)+(q-1)\dim S_q}\Omega(N_F),$$
where $n(F):= \dim N -1-\dim F$ is the codimension of the blow up.
\qed
\end{lemma}
Note that the blow ups in this Chapter are all {\em simple}, i.e. $q=1$.
\par
\noindent Using this, we can describe the action of an element
$P\in \psiphi^{m,{\cal E}}(X)$ on a function $f\in{\cal A}^{\cal
F}(X)$ as follows. Choose a nonvanishing section $u\in
C^{\infty}(X,\Omega(X))$.  Lemma \ref{liftdensity} and the
pullback Theorem \ref{pullback} then tell us that we have a
conormal distribution on $X^2_{\phi}$ of the type
\begin{eqnarray*} \beta_{\phi,L}^*u\cdot P
\cdot \beta_{\phi,R}^*f \in {\cal A}^{\cal
G}I^m(X^2_{\phi},\bigtriangleup_{\phi}, \Omega(X^2_{\phi})) & &
G_{\faff}=E_{\faff}+F,\quad G_{\faphibf}=E_{\faphibf}+F-h-1,\\ &
&G_{\farf}= E_{\farf}+F-h-2,\quad G_{\falf}= E_{\falf}
\end{eqnarray*}
If $G_{\farf}> -1$, we can push forward this expression via
$\beta_{\phi,L}$ to the left factor $X$ (i.e. integrate the
expression over the right factor) using Theorem \ref{pushforward}
and set
 \begin{equation}\label{defphiaction}
 Pf=\beta_{\phi,L \,*}(\beta_{\phi,L}^*u\cdot P\cdot\beta_{\phi,R}^*f)/u,
\quad\mbox{i.e}\quad\langle Pf,u\rangle_X=\langle P,
\beta_{\phi}^*(f\boxtimes u) \rangle_{X^2_{\phi}},
\end{equation}
where $u$ is only used for bookkeeping purposes w.r.t. the
densities. By the pushforward Theorem \ref{pushforward}  this is a
conormal distribution
\begin{equation}\label{mappingphi}
Pf \in {\cal A}^G(X), \quad G=E_{\falf}\overline{\cup}(E_{\faff}+F)
\overline{\cup} (E_{\faphibf}+F-h-1).
\end{equation}
Analogously, if $Q$ is a b-operator in $\Psi_b^{m,{\cal E}}(X)$,
with $E_{\farf}+F-1>-1$ we can set
\begin{equation}\label{mappingb}
 Qf=\beta_{b,L \,*}(\beta_{b,L}^*u\cdot Q\cdot\beta_{b,R}^*f)/u
  \in {\cal A}^G(X), \quad G=E_{\falf} \overline{\cup}(E_{\fabf}+F).
\end{equation}
As a special case note that elements in the small ($b$- or
$\phi$-) calculus always map ${\cal A}^F(X)$ to itself.
\par\medskip
Having defined the action of $\phi$- and $b$-pseudodifferential
operators we can now describe their relationship with the
corresponding classes of differential operators. First, note that
the identity is an operator in the small calculus of order $0$:
\begin{equation}\label{identity}
[\Id]\in \Psi_{b,{\rm cl}}^0(X),\quad\mbox{and}\quad [\Id] \in
\Psi_{\phi,{\rm cl}}^0(X).
\end{equation}
To describe the composition of a $\phi$-differential operator, or
$\phi$-vector field $W$ for simplicity,
 with a $\phi$-pseudodifferential operator $P$,
 recall that for forms $\alpha_1,\alpha_2$ and a $b$-vector field
$W$ on $X^2_{\phi}$ the partial integration rule reads $$\int
L_W\alpha_1\wedge \alpha_2 =-\int \alpha_1\wedge L_W \alpha_2. $$
Hence, from the definition (\ref{defphiaction}) we get $$ \langle
W \circ Pf, u\rangle_X =-\langle Pf,L_W u\rangle= -\langle P,
\beta_{\phi}^*(f\boxtimes L_W u) \rangle_{X^2_{\phi}} = \langle
L_{\beta_{\phi,L}^*W} P, \beta_{\phi}^*(f\boxtimes
u)\rangle_{X^2_{\phi}},$$ from which we can deduce that $ W\circ P
=  L_{\beta_{\phi,L}^*W}  P$.  Here we have used the fact (Lemma
\ref{philift}) that the vector field $\beta_{\phi,L}^*W$ is
tangent to the faces of $X^2_{\phi}$. This argument shows that
composition of a $\phi$-differential operator and a
$\phi$-pseudodifferential operator is well defined:
\begin{equation}\label{diffinpsi}
\Diff_{\phi}^k(X)\circ\psiphi^{m,\cal E}(X) \subset
\psiphi^{m+k,\cal E}(X),
\end{equation}
for any index set ${\cal E}$. Also, from (\ref{identity}) we
immediately infer that $\Diff_{\phi}^k(X)\subset\Psi_{\phi,{\rm cl
}}^{k}(X) $.
\par\bigskip
\noindent{\bf Sobolev Spaces, Compactness and Trace Properties}
\par
\medskip
 We now want to analyze the mapping properties of our
operators w.r.t. $L^2$- and Sobolev spaces. Using the volume form
$\dvol_{\phi}$ which is induced by the exact $\phi$-metric
$g^\phi$ and $\dvol_b=x^{h+1}\dvol_\phi$, we can introduce the
spaces $$ L^2_{\phi}(X,E) = L^2(X,E;\dvol_{\phi}), \qquad
    L^2_b(X,E)= L^2(X,E;\dvol_b),$$
with norms $\|\cdot \|_\phi$ and $\|\cdot\|_b$.
The associated scales of
Sobolev spaces are
\begin{eqnarray*}
H^s_\phi(X,E) &=& \{\xi\in C^{-\infty}(X,E)\quad |\quad \Psi_\phi^s(X,E)\xi
\subset L^2_\phi(X,E) \}  \\
H^s_b(X,E) &=& \{\xi\in C^{-\infty}(X,E)\quad | \quad \Psi_b^s(X,E)\xi
\subset L^2_b(X,E) \},
\end{eqnarray*}
and the Sobolev norms $\|\cdot\|_{\phi,s}$, $\|\cdot\|_{b,s}$ can
be defined as usual using fixed families of elliptic operators
$P_s\in \Psi_{\phi,{\rm cl}}^s(X,E)$, $Q_s\in\Psi_{b,{\rm cl}
}^s(X,E)$ (invertible, or
-- in the b-case -- invertible modulo an element in $\Psi_b^{-\infty}(X,E)$ ).
The main mapping properties of $b$- and $\phi$-operators w.r.t.
these spaces are described in the following
\begin{proposition}
\begin{enumerate}
\item $Q\in\Psi_b^{m,{\cal F}}(X,E)$ with $m\leq 0$
and index sets $F_{\fabf} \geq 0$ and $F_{\falf}, F_{\farf}
> 0$ is a bounded operator on $L^2_b(X,E)$.
\item $P\in\psiphi^{m,{\cal E}}(X,E)$ with $m\leq 0$ and
index sets $E_{\faff} \geq 0$, $E_{\faphibf} \geq h+1$, $E_{\falf}
> \frac{h+1}{2}$ and $E_{\farf}
>\frac{h+1}{2}$ is a bounded operator on $L^2_\phi(X,E)$.
\item $P\in\psiphi^{m,{\cal E}}(X,E)$ with $m\leq 0$ and
index sets $E_{\faff} \geq 0$, $E_{\faphibf} \geq h+1$, $E_{\falf}
> 0$ and $E_{\farf} > h+1$ is a bounded operator on
$L^2_b(X,E)$.
\end{enumerate}
\end{proposition}
\begin{proof}
Reduction to the case $m=-\infty$ is done as usual using
H\"ormander's argument (see \cite{Meaomwc} or Proposition 7 in
\cite{MelMazmwfb}). In (a) it then suffices to verify the mapping
property for an operator $Q$ supported in a coordinate patch.
There, write $Q=a({\bf w},{\bf w}')\dvol_{b}({\bf w}')$. We can
estimate
\begin{eqnarray*}
\|Q f\|^2_b &=&
\int_X |\int_X a({\bf w},{\bf w}')f({\bf w}')\dvol_{b}({\bf w}')
|^2\dvol_{b}({\bf w})\\
 &\leq& \int_X \int_X |a({\bf w},{\bf w}')|\dvol_{b}({\bf w}')
\int_X |a({\bf w},{\bf w}')||f({\bf w}')|^2\dvol_{b}({\bf w}')
\dvol_{b}({\bf w}).
\end{eqnarray*}
Now, it is easy to see that there is an operator $B =b({\bf
w},{\bf w}')\dvol_{b}({\bf w}')$ such that $b({\bf w},{\bf
w}')\geq |a({\bf w},{\bf w}')|$. Then, the above is certainly
smaller than $$ \sup_{{\bf w}} (\int_X b({\bf w},{\bf
w}')\dvol_{b}({\bf w}'))\cdot \sup_{{\bf w}'} (\int_X b({\bf
w},{\bf w}')\dvol_{b}({\bf w}))\cdot \int_X |f({\bf
w}')|^2\dvol_{b}({\bf w}').$$ The $\sup$'s are finite, since $B$
maps the function $1$ to $L^{\infty}(X)$ and the  adjoint kernel
satisfies the same estimates; (b) and (c) are then consequences of
(a) using the fact that
$\beta_{\phi,R}^*\Omega_\phi=(x')^{-h-1}\beta_{\phi,R}^*\Omega_b$
and $L^2_\phi(X)=x^{(h+1)/2}L^2_b(X)$. \qed
\end{proof}
\begin{corollary}
For $P\in \Psi_{\phi,{\rm cl}}^m(X,E)$, $Q\in\Psi_{b,{\rm cl}
}^m(X,E)$ the maps $$ P:H_\phi^{s+m}(X,E)\longrightarrow
H_\phi^{s}(X,E),\qquad
    Q:H_b^{s+m}(X,E)\longrightarrow H_b^{s}(X,E). $$
are continuous.\qed
\end{corollary}\par
The standard  results on compactness- and trace class-properties
of  pseudodifferential operators on a compact manifold generalize
to our context as follows:
\begin{proposition}[Sobolev, Rellich and Lidskii]\quad\hfill
\begin{enumerate}
\item $H_b^{n/2+l+\epsilon}(X,E)\hookrightarrow C^l(\stackrel{\circ}{X},E)$,
continuously.
\item $x^{\delta}H_b^{s+\epsilon}(X,E)\hookrightarrow H_b^{s}(X,E)$ is
      compact.
\item Let $P$ be a bounded operator on $L^2_b(X,E)$ and assume
 $P:H^s_{b}(X,E)\rightarrow x^{\delta}H_b^{s+n+\epsilon}(X,E)
 \quad\mbox{continuously}$.
 Then $P$ is trace class and
$ \tr(P)=\int_{\bigtriangleup}[P].$\qed
\end{enumerate}
\end{proposition}
The above Sobolev spaces and their variants will only play a minor
role in subsequent Chapters of this work, since we will explicitly
construct the integral kernels of the resolvent and the heat
operator of $\Dir^d$ and not rely on abstract existence results.
The main result that we are going to use is the following
Corollary, which gives us a minimum requirement for the smallness
of the error term in the construction.
\begin{corollary}[Trace Class
Operators]\label{traceclass}\quad\par\noindent
 Operators $Q\in
\Psi_b^{-n-\epsilon, (\epsilon,\epsilon,\epsilon)}(X,E)$ are trace
class on $L^2_b(X,E)$. \qed
\end{corollary}
Finally, let us note how to calculate adjoints w.r.t.
$L^2_b(X,E)$: Let $Q$ be an operator in $\Psi_b(X,E)$ or
$\Psi_{\phi}(X,E)$ and write $Q^+$ for the formal adjoint of $Q$
obtained by flipping sides in $X^2$. Then
\begin{lemma}\label{adjoints}
The $L^2_b$-adjoint of $Q$ is given by $Q^*=Q^+\cdot
\dvol_{b,L}^{-1}\otimes \dvol_{b,R}$.\qed
\end{lemma}
\par\bigskip
In order to avoid any confusion due to the variety of different
$L^2$-spaces, and to be able to use the standard constructions of
the $b$-calculus (see \cite{Melaps}) without difficulty, we decide
to  {\bf\boldmath  let all operators act on $L^2_b(X)$}. However,
the goal of this work is the description of the spectral
properties of the Dirac operator $\Dir^d$ on $
L^2(X,E;\dvol_{d})$, where $\dvol_d=x^n\dvol_{\phi}=x^{v}\dvol_b$ is the volume
form corresponding to the metric $g^d=x^2g^{\phi}$. Thus, we have
to consider the diagram (truly valid only on the domains of
definition) $$
\begin{array}{ccc}
   L^2(X,E;x^n\dvol_{\phi})& \stackrel{\Dir^d}{\longrightarrow} & L^2(X,E;
     x^n\dvol_{\phi})\\
   \quad \downarrow x^{v/2}& & \quad \downarrow x^{v/2}\\
   L^2_b(X,E)& \stackrel{\Dir^d-\frac{v}{2}\cl_{d}(\frac{\dd x}{x})}
  {\longrightarrow} & L^2_b(X,E).
    \end{array}. $$
Thus, requiring all our operators to act on $L^2_b(X,E)$ means that
 we are ultimately looking at the operator {\bf\boldmath $$ \Dir
=\Dir^d-\frac{v}{2}\cl_{d}(\frac{\dd x}{x})=
\frac{1}{x}\Dir^{\phi}+\frac{h}{2}\cl_{\phi}(\frac{\dd x}{x^2})
 \quad\mbox{on}\quad L^2_b(X,E).$$}
  This choice will turn out to
be particularly convenient in Sections \ref{normalopsbf} and
\ref{inversionatphibf}.

\subsection{Composition Formula}\label{compositionformula}
In this Section, we will prove that the composition of two
operators in $\Psi_{\phi}(X)$ is again in $\Psi_{\phi}(X)$. To
explain the main idea of the proof, consider  first the ``triple
space'' $X^3$ together with its projections onto $X^2$:
\begin{equation}\label{triplex}
\begin{array}{ccccc} X^2 &\stackrel{\pi_L}{\longleftarrow} &
X^3 &\stackrel{\pi_R}{\longrightarrow}& X^2
\\ & & \quad\downarrow \pi_M
\\ & &  X^2
\end{array}
\end{equation}
The composition of two operators $A, B\in
\cdotinfty{}(X^2,\beta_R^*\Omega(X))$ is then simply given by
$$\pi_{M\,*}(\pi_L^*A\cdot\pi_R^*B).$$ In order to describe the
composition properties of $\phi$-pseudodifferential operators we
need to replace the image spaces in (\ref{triplex}) by
$X^2_{\phi}$ and modify $X^3$ in such a way that (\ref{triplex})
becomes a diagram of $b$-fibrations. This will allow us to use
the pullback- and pushforward results in \cite{Meaomwc} (see also Appendix
\ref{appzeronotation}) to
describe  composition of $\phi$-pseudodifferential operators.
\par
\smallskip
The first step is to blow up the corner of $X^2$ in
(\ref{triplex}), thus replacing it by $X_b^2$. Blowing up the
preimages of the corner under $\pi_L$, $\pi_M$ and $\pi_R$ at
their intersections $$T:=(\partial X)^3, \quad F:=X\times\partial
X\times
\partial X, \quad
   C:=\partial X \times X \times \partial X, \quad S:=\partial X\times
\partial X \times X,$$
leads to the $b$-triple space $X_b^3$ defined as the blow up
$$X^3_b:=\left[ X^3;T,F,C,S\right].$$
\begin{table}[t]
\centerline{\epsfig{file=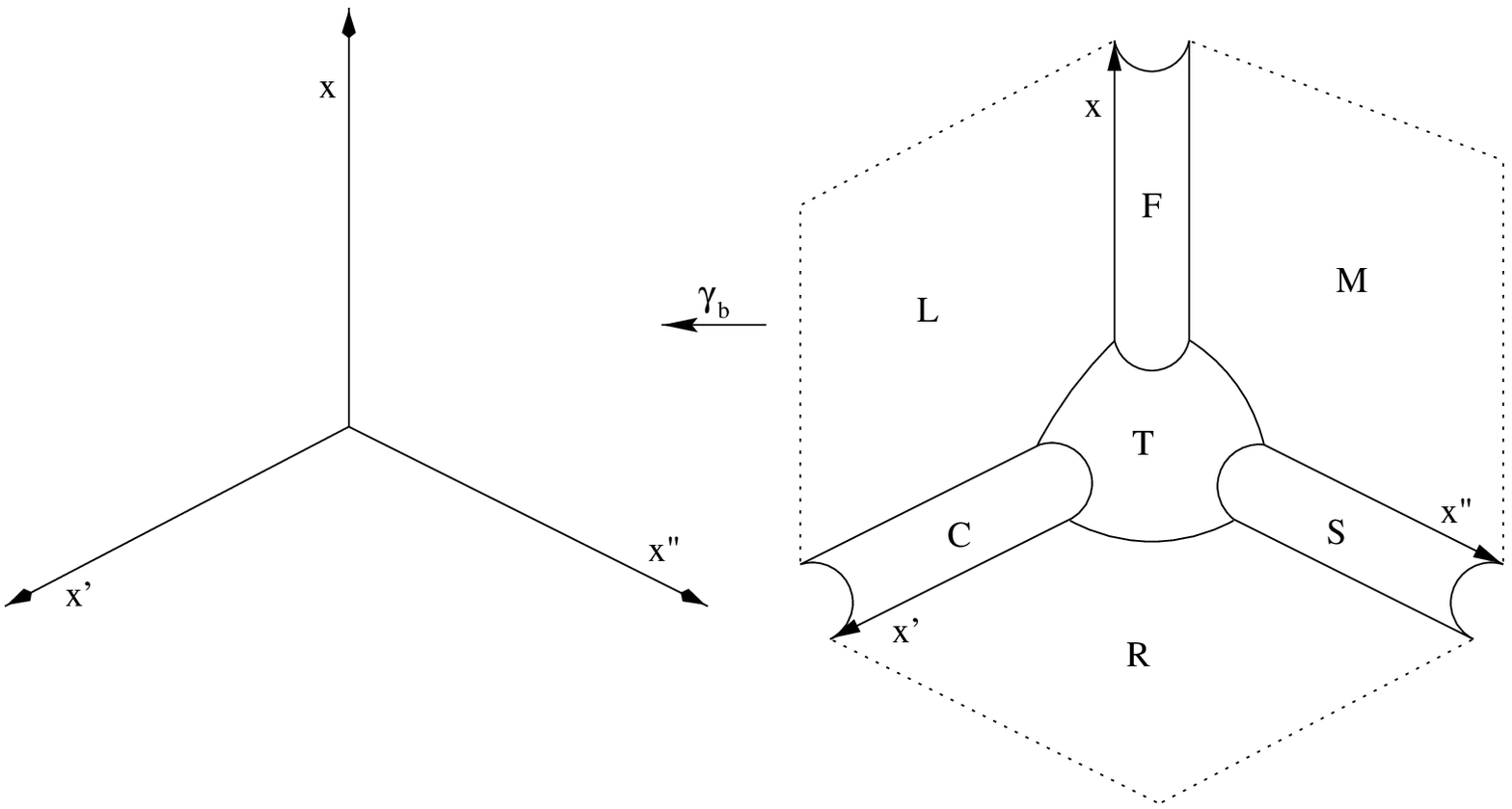}}
 \centerline{Figure 3: $X^3_b$}
\end{table}
The projections in (\ref{triplex}) induce {\em $b$-fibrations}
$\pi_{b,o}:X^3_b\rightarrow X^2_b$ for $o=L,M,R$ such that the
diagram
\begin{eqnarray}
X^3_b&\stackrel{\pi_{b,o}}{\longrightarrow}&X^2_b\nonumber
\\ \downarrow& & \downarrow \label{tripleb}
\\ X^3& \stackrel{\pi_o}{\longrightarrow}& X^2\nonumber
\end{eqnarray}
is commutative. The next step is to replace $X^2_b$ in
(\ref{tripleb}) by $X^2_{\phi}$ by blowing up the fibre diagonal
$F_{\phi}$. To define the new LHS in (\ref{tripleb}),
$X^3_{\phi}$, we have to blow up the lifts of $F_{\phi}$ under
$\pi_L$, $\pi_M$ and $\pi_R$. We denote these faces (as afterwards
the result of their blow up) by $\phi_{F}$, $\phi_{C}$,
$\phi_{S}$, $\phi_{FT}$, $\phi_{CT}$, $\phi_{ST}$, $\phi_{TT}$,
and define $$X^3_{\phi}:=\left[ X^3_b; \phi_{TT} ;\,\phi_{FT}
;\,\phi_{CT} ;\,\phi_{ST};
  \,\phi_{F} ;\,\phi_{C} ;\,\phi_{S}\right].$$
 By $R, M, C, F, S, C, T$ we denote the lifts to $X^3_{\phi}$ of the
non-blown up parts of these faces in $X^3_b$. We will also have to
consider the behavior of our operators at the lifted diagonals
$\bigtriangleup_{\phi, o}:=\pi_{\phi,
o}^{-1}\bigtriangleup_{\phi}$, with $o=L, M, R$.  The following is
proved in \cite{MelMazmwfb}:
\begin{proposition}
The projections in (\ref{triplex},\ref{tripleb}) induce
$b$-fibrations $\pi_{\phi,o}:X^3_{\phi}\rightarrow X^2_{\phi}$ for
$o=L,M,R$ such that the diagram
\begin{eqnarray*}
X^3_{\phi}&\stackrel{\pi_{\phi,o}}{\longrightarrow}&X^2_{\phi}
\\ \downarrow& & \downarrow
\\ X^3& \stackrel{\pi_o}{\longrightarrow}& X^2
\end{eqnarray*}
is commutative.\qed
\end{proposition}
\begin{table}[t]
 \centerline{\epsfig{file=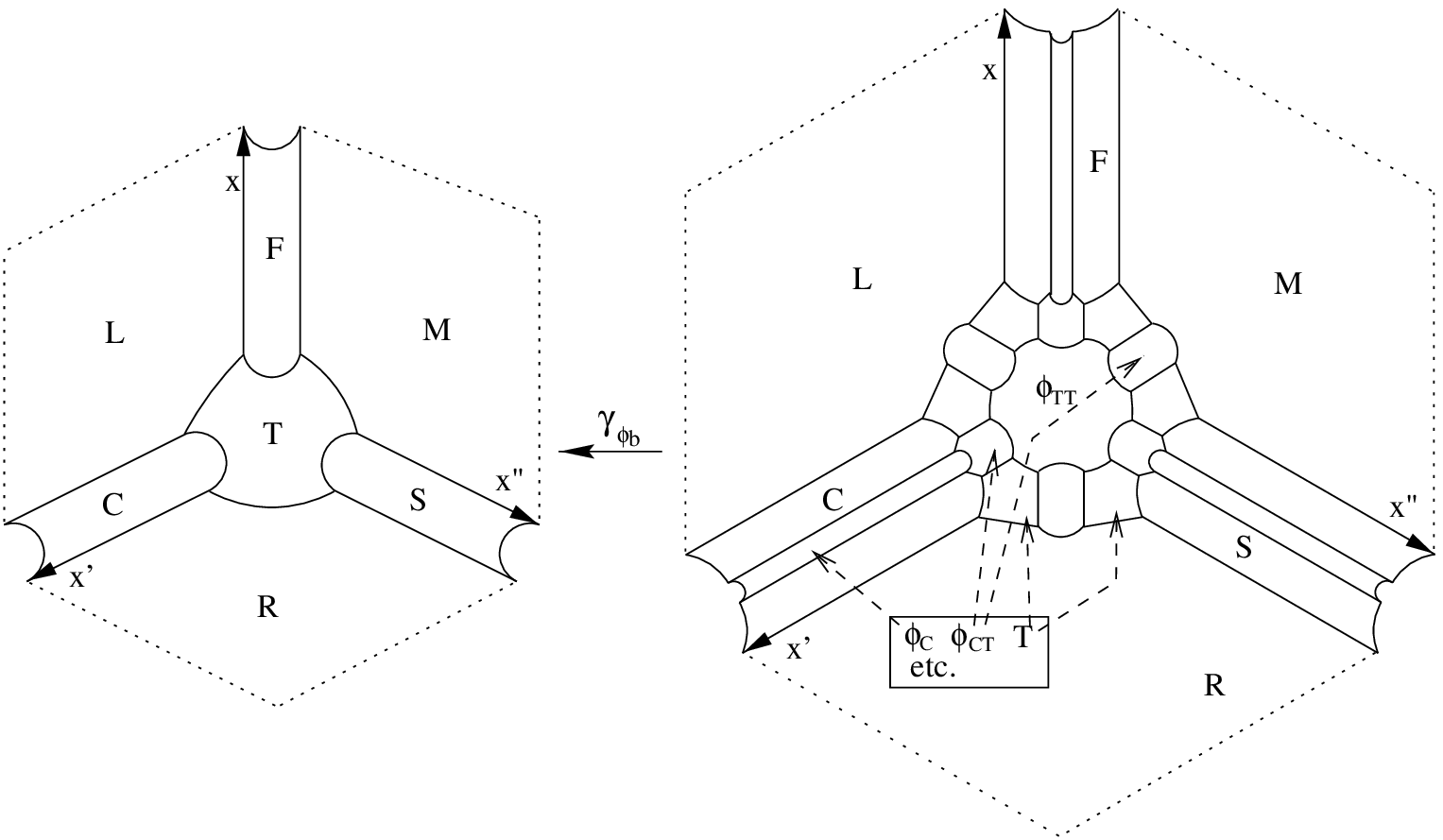}}
 \centerline{Figure 4: $X^3_{\phi}$}
 \end{table}\par
\noindent The coefficient matrices can be read off from Figure 4:
\begin{lemma}\label{coeffmatphi}
\begin{enumerate}
\item $\pi_L$ maps the face $L$ to the interior of
$X^2_{\phi}$. Also
\par $\begin{array}{rclcrcl}
   \pi_L^*\falf &=&\{ R, C, \phi_{C}  \}
& & \pi_L^*\farf&=&\{ M, F, \phi_F \}
\\ \pi_L^*\faphibf &=& \{ S, T, \phi_{FT}, \phi_{CT}\}
& & \pi_L^*\faff &=& \{ \phi_S, \phi_{ST}, \phi_{TT}    \}
\end{array}$
\item $\pi_R$ maps the face $R$ to the interior of
$X^2_{\phi}$. Also
\par $\begin{array}{rclcrcl}
   \pi_L^*\falf &=&\{ M, S, \phi_{S}  \}
& & \pi_L^*\farf&=&\{ L, C, \phi_C \}
\\ \pi_L^*\faphibf &=& \{ F, T, \phi_{ST}, \phi_{CT}\}
& & \pi_L^*\faff &=& \{ \phi_F, \phi_{FT}, \phi_{TT}    \}
\end{array}$
\item $\pi_M$ maps the face $M$ to the interior of
$X^2_{\phi}$. Also
\par $\begin{array}{rclcrcl}
   \pi_M^*\falf &=&\{ R, S, \phi_{S}  \}
& & \pi_M^*\farf&=&\{ L, F, \phi_F \}
\\ \pi_M^*\faphibf &=& \{ C, T, \phi_{ST}, \phi_{FT}\}
& & \pi_M^*\faff &=& \{ \phi_C, \phi_{CT}, \phi_{TT}    \}
\end{array}$\par\qed
\end{enumerate}
\end{lemma}
\par
Starting with two operators $A\in \Psi_{\phi}^{l,{\cal L}}(X)$,
$B\in \Psi_{\phi}^{r,{\cal R}}(X)$, and a placeholder function
$u\in C^{\infty}(X,\Omega(X))$, their lift to $X^3_{\phi}$ is
 \begin{equation}\label{lifttox3phi}
 \gamma_{\phi,L}^*u\cdot\pi_{\phi,L}^*A
\cdot\pi_{\phi,R}^*B\in{\cal A}^{\cal
G}I^{(l-n/4,r-n/4)}(X^3_{\phi}, \bigtriangleup_{\phi,L}\cup
\bigtriangleup_{\phi,R} ;\Omega(X^3_{\phi})) .
\end{equation}
The space on the RHS consists of distributions with conormal
singularity at the union $\bigtriangleup_{\phi,L}\cup
\bigtriangleup_{\phi,R}$ of interior manifolds in $X^3_{\phi}$.
This is described in \cite{Meaomwc}. At the boundary faces its
elements have expansions described by the index set ${\cal G}$ as
in Appendix \ref{appzeronotation}.
 To calculate the index set ${\cal G}$ we need to know the lifting
properties of the density bundle $\phiomega(X)$. The following is
proved using Lemma \ref{liftdensity}:
\begin{lemma}
\begin{enumerate}
\item $\gamma_{b}^*\Omega(X^3)=\rho_{b,F}\, \rho_{b,C}\, \rho_{b,S}\,
      \rho_{b,T}^2 \,\Omega(X^3_b)$
\item $\gamma_{\phi,b}^*\Omega(X^3_b)
    =\rho_{\phi_{T}}^{2h+2}\, (\rho_{\phi_{FT}}\,\rho_{\phi_{CT}}\,
     \rho_{\phi_{ST}}\, \rho_{\phi_{F}}\,\rho_{\phi_{C}}\,
      \rho_{\phi_{S}})^{h+1}\Omega(X^3_{\phi})$
\item $\gamma_{\phi}^*\Omega(X^3)=
    \rho_F\, \rho_C\, \rho_S\, \rho_T^2\,  \rho_{\phi_{F}}^{h+2}\,
    \rho_{\phi_{C}}^{h+2}\, \rho_{\phi_{S}}^{h+2}\,
    \rho_{\phi_{FT}}^{h+3}\, \rho_{\phi_{CT}}^{h+3}\, \rho_{\phi_{ST}}^{h+3}\,
     \rho_{\phi_{TT}}^{2h+4} \Omega(X^3_{\phi})$
\item $ \gamma_{\phi,F}^*\Omega(X)\cdot \pi_{\phi,L}^*\phiomega(X)
\cdot\pi_{\phi,R}^*\phiomega(X)= \rho_L^0\, \rho_M^{-h-2}\,
\rho_R^{-h-2}$\par
 $\qquad\quad\cdot\rho_F^{-h-1}\, \rho_C^{-h-1}\, \rho_S^{-2h-3}\, \rho_T^{-2h-2}\,
   \rho_{\phi_{F}}^{0}\, \rho_{\phi_{C}}^{0}\, \rho_{\phi_{S}}^{-h-2}\,
   \rho_{\phi_{FT}}^{-h-1}\,\rho_{\phi_{CT}}^{-h-1}\,  \rho_{\phi_{ST}}^{-h-1}
   \rho_{\phi_{TT}}^0\Omega(X^3_{\phi})$\qed
\end{enumerate}
\end{lemma}
From this and the pullback Theorem \ref{pullback} (or really its
extension to distributions conormal at the boundary faces {\em
and} to some interior submanifold, see \cite{Meaomwc}) it follows
directly that the index ${\cal G}$ is given by
\par $$\begin{array}{rclcrcl}
G_L&=&R_{\farf} & \quad&G_F&=&L_{\farf}+R_{\fabf}-h-1  \\
G_M&=&L_{\farf}+R_{\falf}-h-2 & & G_C&=&L_{\falf}+R_{\farf}-h-1\\
G_R&=& L_{\falf}-h-2 & &G_S&=&L_{\faphibf}+R_{\falf}-2h-3\\
   & & & &G_T&=&L_{\faphibf}+R_{\faphibf}-2h-2\\
G_{\phi_F}&=&L_{\farf}+R_{\faff}& &
   G_{\phi_{FT}}&=&L_{\faphibf}+R_{\faff}-h-1\\
G_{\phi_C}&=&L_{\falf}+R_{\farf}& &
   G_{\phi_{CT}}&=&L_{\faphibf}+R_{\faphibf}-h-1\\
G_{\phi_S}&=&L_{\faff}+R_{\falf}-h-2& &
   G_{\phi_{ST}}&=&L_{\faff}+R_{\faphibf}-h-1\\
          & & & &G_{\phi_{TT}}&=&L_{\faff}+R_{\faff}.
\end{array} $$
Thus, by the pushforward Theorem \ref{pushforward} (again, it is
really the extension of this result to distributions conormal to
$\bigtriangleup_{\phi,L}\cup \bigtriangleup_{\phi,R}$ in
\cite{Meaomwc} that is needed), the conormal distribution in
(\ref{lifttox3phi}) can be pushed forward via $\pi_{\phi,M}$
whenever $G_M=L_{\farf}+R_{\falf}-h-2>-1$ giving $$
\pi_{M\,*}(\gamma_{\phi,L}^*u\cdot\pi_{\phi,L}^*A
\cdot\pi_{\phi,R}^*B)\in{\cal A}^{\pi_{\phi,M\,\sharp}{\cal
G}}I^{l+r}(X^2_{\phi},\bigtriangleup_{\phi};\Omega(X^2_{\phi})).$$
Here the pushforward of the index set ${\cal G}$ is given by the
rule $$ \pi_{\phi,M\, \sharp\,}{\cal G}(B) = \overline{\bigcup}
_{F\in{\cal M} (X^3_{\phi}), \pi_{\phi,M}(F)\subset B} G_F,  $$
which by  Lemma \ref{coeffmatphi} is $$\begin{array}{rclcrcl}
 \pi_{\phi,M\, \sharp\,}{\cal G}(\falf) &=&G_R\overline{\cup} G_S\overline{\cup}
   G_{\phi_{S}}& &
 \pi_{\phi,M\, \sharp\,}{\cal G}(\farf)&=& G_L\overline{\cup} G_F
\overline{\cup} G_{\phi_F}
\\ \pi_{\phi,M\, \sharp\,}{\cal G}(\faphibf) &=&  G_C\overline{\cup}
  G_T\overline{\cup} G_{\phi_{ST}}\overline{\cup} G_{\phi_{FT}}& &
 \pi_{\phi,M\, \sharp\,}{\cal G}(\faff) &=& G_{\phi_C}
  \overline{\cup} G_{\phi_{CT}}\overline{\cup} G_{\phi_{TT}}
\end{array}$$
Rewriting everything in terms of the $\phi$-density bundle lifted
from the right $$\beta_{\phi,
R}^*\phiomega(X)=\rho_{\falf}^{-h-2}\,
\rho_{\faphibf}^{-h-1}\Omega(X^2_{\phi})/\beta_{\phi,L}^*\Omega(X),$$
we get the final composition formula:
\begin{theorem}[Composition Formula]\label{composition}
Let $A\in \Psi_{\phi}^{l,\cal L}(X)$  and $B\in
\Psi_{\phi}^{r,\cal R}(X) $ and assume
$L_{\farf}+R_{\falf}-h-2>-1$. Then $A\circ B\in \psiphi^{r+l, \cal
M}(X)$ with index set given by
\par\medskip $\begin{array}{lcl}
M_{\falf}&=&L_{\falf}\overline{\cup} (L_{\faphibf}+R_{\falf}-h-1)
\overline{\cup} (L_{\faff}+R_{\falf})\\
M_{\farf}&=&R_{\farf}\overline{\cup} (L_{\farf}+R_{\faphibf}-h-1)
\overline{\cup} (L_{\farf}+R_{\faff})\\
M_{\faphibf}&=&(L_{\falf}+R_{\farf}) \overline{\cup}
(L_{\faphibf}+R_{\faphibf}-h-1) \overline{\cup}
(L_{\faphibf}+R_{\faff})\overline{\cup} (L_{\faff}+R_{\faphibf})\\
M_{\faff}&=&(L_{\falf}+R_{\farf}) \overline{\cup} (L_{\faphibf}+\
R_{\faphibf}-h-1) \overline{\cup} (L_{\faff}+R_{\faff})
\end{array}$\par\qed
\end{theorem}
\par
All the results above also hold for $\phi$-pseudodifferential
operators with coefficients in vector bundles.

\newpage
\section{Normal Operators and the Resolvent Construction}
In this Chapter we will present a construction of the resolvent
$G(\lambda)$ of the Dirac operator $\Dir$.
 We do this by first solving the equation
\begin{equation}\label{resolventeq}
(\Dir-\lambda)Q(\lambda)=I+R(\lambda) \quad \mbox{on}\quad
X^2_{\phi}
\end{equation}
with successively smaller and smaller remainder terms
$R(\lambda)$, then use the composition formula, Theorem
\ref{composition}, to (formally) invert $I+R(\lambda)$ and obtain
a smoothing remainder term. In a final step this remainder term is
then also removed.
\par
The construction will proceed  by the symbolic solution of
(\ref{resolventeq}) and the subsequent removal of error terms
appearing at the faces $\faff$, $\faphibf$, $\falf$. The technical
tools used are the so-called normal operators, which allow to give
a representation of (\ref{resolventeq}) at each of the above
faces. These concepts will be introduced at the corresponding
stages of the construction in Sections \ref{secsymbolmap}, \ref{normalop}
and \ref{normalopsbf}. Prototypes of this kind of
construction can be found in \cite{MelMazmerext} and
\cite{Melaps}.
\subsection{Symbol Map}\label{secsymbolmap}
The symbol map ${}^\phi\sigma$   in the small calculus of
$\phi$-pseudodifferential operators, describes the conormal
singularities at $\bigtriangleup_{\phi}$. It is given by the
composition
\begin{eqnarray}
{}^\phi\sigma_m:\Psi_{\phi,{\rm cl}}^m(X)= \cdotinfty{\faff}I_{\rm
cl }^m(X^2_{\phi},\bigtriangleup_{\phi};\beta_{\phi,R}^*
\phiomega(X))& \longrightarrow& C^{\infty}_{c}I_{\rm cl
}^m(N\bigtriangleup_{\phi},\bigtriangleup_{\phi};
\Omega(N\bigtriangleup_{\phi}/\bigtriangleup_{\phi} ))\nonumber
\\ &\stackrel{F}{\longrightarrow}& {\cal S}_{\rm
cl}^{m}(N^*\bigtriangleup_{\phi})
 \cong  {\cal S}_{\rm cl}^m (\phitsx).\label{symbolphitsx}
\end{eqnarray}
Here, the first map comes from an identification of
$N\bigtriangleup_{\phi}$ with a tubular neighborhood of
$\bigtriangleup_{\phi}$, the second map is Fourier transform and
the identification at the end is just Lemma \ref{philift}. Also,
using radial compactification (for instance by introducing
$\rho:=|\,|^{-1}$ as a defining function for the boundary at
infinity) in the fibres of $\phitsx$ we could have written $${\cal
S}_{\rm cl}^{m}(\phitsx)=\rho^{-m}C^{\infty}(\RC(\phitsx)),$$
which by multiplication with $\rho^m$ and restriction to the
boundary maps to $C^{\infty}({}^{\phi}S^*X)$.
\begin{remark}\rm \label{remdefconormal}
The existence of the map (\ref{symbolphitsx}) can be taken to be
the defining property of the space $I_{\rm cl
}^m(X^2_{\phi},\bigtriangleup_{\phi};\beta_{\phi,R}^*
\phiomega(X))$ of distributions conormal at
$\bigtriangleup_{\phi}$ of order $m$ and $C^{\infty}$ up to
$\faff$.
\end{remark}\par
\medskip
 Putting everything together, the $m$th
order symbol fits into a short exact sequence
\begin{equation}\label{symbolphissx}
 0
\longrightarrow \Psi_{\phi,{\rm cl}}^{m-1}(X)\longrightarrow
\Psi_{\phi,{\rm cl }}^{m}(X)
\stackrel{{}^{\phi}\sigma_m}{\longrightarrow}
C^{\infty}({}^{\phi}S^*X) \longrightarrow 0
\end{equation}
 of {\em filtered
algebra} homomorphisms.
\par
Note that we will use both versions (\ref{symbolphitsx},
\ref{symbolphissx}) of the symbol without distinguishing them
notationally. Also, the generalization of these definitions to the
calculus with coefficients in a vector bundle is left to the
reader. Especially, using the identification of
$N\bigtriangleup_{\phi}$ with $\phitx$:
\begin{proposition} Let $W\in\Gamma(\phitx)$, $Q\in\psiphi^m(X)$,
$Q_E\in\psiphi^m(X,E)$.
\begin{enumerate}
\item ${}^{\phi}\sigma_{1}(W)=iW, \quad$ i.e.
      $\quad{}^{\phi}\sigma_{m+1}(W\circ Q)= iW\otimes {}^{\phi}\sigma_{m}(Q)$
\item ${}^{\phi}\sigma_{1}(\Dir^{\phi})={}^{\phi}\sigma_{1}(x\Dir^d)
={}^{\phi}\sigma_{1}(x\Dir)=i\cl_{\phi},\quad$ i.e.
${}^{\phi}\sigma_{m+1}(\Dir^{\phi}\circ Q_E)=
i\cl_{\phi}\otimes \,{}^{\phi}\sigma_{m}(Q_E),$
\end{enumerate}
 etc., where $W$ (and $\cl_{\phi}$) are interpreted as linear functions
on $\phitsx$.\qed
\end{proposition}
As usual, an element $P$ in $\psiphi^{m>0}(X,E)$ will be called
{\em elliptic}, when its symbol ${}^{\phi}\sigma_m(P)$ is
invertible in $C^{\infty}({}^{\phi}S^*X)$. Thus, for instance, the
Dirac operator $\Dir^{\phi}$ is elliptic.

\subsection{Symbolic Construction}
To perform the symbolic part of the construction of the inverse of
the Dirac operator $\Dir-\lambda$ write
 $$ (\Dir-\lambda)^{-1}=
 (x\Dir-x\lambda)^{-1}\cdot x. $$
Thus, we can as well invert the operator
$P(\lambda)=x\Dir-x\lambda$, which is an elliptic operator in the $\phi$-calculus:
\begin{proposition}\label{symbolicresult}
There are  holomorphic families, $$ \lambda \mapsto
Q_0(\lambda)\in \Psi_{\phi,{\rm cl}}^{-1}(X,E), \qquad \lambda
\mapsto R_0(\lambda)\in\Psi_{\phi,{\rm cl}}^{-\infty}(X,E), $$
such that $P(\lambda)Q_0(\lambda)=1+R_0(\lambda).$
\end{proposition}
\begin{proof}
This is the standard symbolic construction using the symbol
sequence $$  0 \longrightarrow \psiphi^{m-1}(X,E)\longrightarrow
\psiphi^{m}(X,E) \stackrel{{}^{\phi}\sigma_m}{\longrightarrow}
C^{\infty}({}^{\phi}S^*X,\END(E))\longrightarrow 0. $$ The symbol
${}^{\phi}\sigma_1(P(\lambda))$ is independent of $\lambda$ and
invertible. We can therefore find a $g \in
C^{\infty}({}^{\phi}S^*X,\END(E))$ and,
 because of the surjectivity of the symbol
map ${}^{\phi}\sigma_{-1}$, also a $G\in\Psi_{\phi,{\rm cl
}}^{-1}(X,E)$ with $${}^{\phi}\sigma_{-1}(G)= g,\qquad
  P(\lambda)G=1-S(\lambda),\quad
\lambda\mapsto S(\lambda) \in \Psi_{\phi,{\rm cl}}^{-1}(X,E)\quad
\mbox{holomorphic}.$$ We now set as usual
$$Q_0(\lambda)=G\sum_{j}S(\lambda)^j,\qquad
R_0(\lambda)=P(\lambda)Q_0(\lambda)-1.$$ The construction can be
made holomorphic in the parameter $\lambda$.\qed
\end{proof}
\par
In general the remainder $R_0(\lambda)$ will not vanish to any
positive order at the face $\faff$ and therefore (compare Section
\ref{secmapprop}) cannot be expected to be compact. It will
therefore be necessary to ``solve away'' (at least) the leading
term of $R_0(\lambda)$ at $\faff$, i.e. to modify the parametrix
$Q_0(\lambda)$ in such a way that the new remainder also vanishes
at $\faff$. This will be done in the next two Sections.

 \subsection{Normal Operator at ff}\label{normalop}
The normal operator at $\faff$ of an operator $P\in\Psi_{\phi,{\rm
cl}}^*(X)$ is defined as the restriction of the conormal
distribution $P$ to that face. To describe this in more detail,
first recall that the  interior of the front face $\faff$ is
\begin{equation}\label{identifyff}
\stackrel{\circ}{\faff}=\phindx\times_Y\partial X
\end{equation}
and $\faff$ fibres over $\partial X$ via the left and right
projections $$\faff\stackrel{\beta_{\phi,L}}{\longrightarrow}
\partial X \qquad \faff \stackrel{\beta_{\phi,R}}{\longrightarrow}
\partial X.$$
The density bundle for the $\phi$-calculus restricts to $\faff$ as
the density bundle on the fibres of the left projection
\begin{equation}
\beta_{\phi,R}^*\phiomega(X)|_{\faff} \cong
\omfib(\phindx\times_Y\partial
X\stackrel{\beta_{\phi,L}}{\rightarrow}
 \partial X)\cong \beta_{\phi,R}^*\Omega(\phindx/Y)
\end{equation}
This allows us to define the normal map at $\faff$ as follows
\begin{definition}[Suspended Calculus and the Normal Map]\quad
\par\noindent
\begin{enumerate}
\item The space $\psisus^m(\phindx)$ of {\em suspended} pseudodifferential
operators on $\phindx$ is defined as the space $$
\cdotinfty{}I_{\rm cl}^{m+1/4}(\faff,
  \bigtriangleup_{\phi};\omfib(\beta_{\phi,L}))
= {\cal S}_{\rm fib}I_{\rm cl}^{m+1/4}(\phindx\times_{Y}\partial
X, F_{\phi};\beta_{\phi,R}^*\Omega(\phindx/Y))$$
 acting
as ($C^{\infty}(Y)\times \RR^{h+1}$-invariant) convolution
operators on (e.g.) ${\cal S}_{\rm fib}(\phindx)$.
\item The normal map $N_{\faff}:\Psi_{\phi,{\rm cl}}^m(X)\rightarrow \psisus^m(\phindx) $
is given
by $N_{\faff}(A):= A|_{\faff}.$\qed
\end{enumerate}
\end{definition}
 It is again straightforward to extend these  definitions to the
the calculus with coeffients in a vector bundle as well as to any
extended calculus, which allows restriction to $\faff$ to be
defined. Let us note the main properties of this construction:
\begin{lemma}\label{leminducedcalc}
\begin{enumerate}
\item The suspended pseudodifferential operators form a calculus,
i.e. composition gives a map $$\Psi_{\rm
sus}^m(\phindx)\times\Psi_{\rm sus}^l(\phindx)
\stackrel{\circ}{\longrightarrow}\Psi_{\rm sus}^{l+m}(\phindx).$$
\item The normal map fits into a short exact sequence $$
0\longrightarrow \rho_{\faff}\Psi_{\phi,{\rm
cl}}^m(X)\longrightarrow \Psi_{\phi,{\rm cl
}}^m(X)\stackrel{N_{\faff}}{\longrightarrow}
\psisus^m(\phindx)\longrightarrow 0,$$ which  is a sequence of
filtered algebras.
\end{enumerate}
\end{lemma}
\begin{proof}
First, the exactness of the sequence in (b) should be clear. The
algebraic properties follow, since composition in $\Psi_{\rm
sus}(\phindx)$ is induced by composition in $\Psi_{\phi}(X)$:
Taking  two operators $P$, $Q$ in the {\em small} calculus
$\Psi_{\phi,{\rm cl}}^*(X)$,  it follows from the composition
formula, Theorem \ref{composition}, that $(P\circ Q)|_{\faff}$ is
$0$ whenever $P|_{\faff}$ or $Q|_{\faff}$ is $0$. Thus, the
composition of $P|_{\faff}$ and $Q|_{\faff}$  is independent of
the choice of extension to $X^2_{\phi}$. \qed
\end{proof}
\par
\medskip
Again, the description of the normal operators for
$\phi$-differential operators is particularly easy. Since, as we
have seen, the composition of a $\phi$-differential operator, or
$\phi$-vector field $W$, and a $\phi$-pseudodifferential operator
$P$ is given by $$ W\circ P = L_{\beta_{\phi,L}^*W}
P=\beta_{\phi,L}^*W\cdot P, \quad \mbox{and thus}\quad (W\circ
P)|_{\faff} = \beta_{\phi,L}^*W|_{\faff}\cdot P|_{\faff},$$ the
normal operator of $W$ can be viewed as the restriction of
$\beta_{\phi,L}^*W$ to the front face $\faff$, acting on
$\psisus^m(\phindx)$.
\par
This action can be described geometrically using the
identification (\ref{identifyff}). The vector field $\beta_{\phi,L
}^*W$ is tangent to the fibres of $\beta_{\phi,R}:\faff\rightarrow
Y$, i.e. it is tangent to (the left factor  of) $\phindx/Y$ in the
interior of $\faff$.
  Using the identification
\begin{equation}\label{fibrephindx} T\phindx/Y\cong \pi^*\phindx\oplus_{\partial X} V\partial X
\end{equation}
we have
\begin{proposition} Define the map $N_{\phi}:\phitx\rightarrow
T\phindx/Y$ as the composition $$ N_{\phi}: \phitx \longrightarrow
\phitx|_{\partial X} \cong \phindx\oplus V\partial
X\stackrel{(\pi^*,\Id)}{\longrightarrow}T\phindx/Y.$$ Then for
$Q\in\Psi_{\phi,{\rm cl }}^m(X)$, $Q_E\in\Psi_{\phi,{\rm
cl}}^m(X,E)$ and $W\in\Gamma(\phitx)$ we have
\begin{enumerate}
\item $N_{\faff}(W \circ Q)=N_{\phi}(W)\circ N_{\faff}(Q)$.
\item $N_{\faff}((x\Dir-x\lambda)\circ Q_E)=N_{\faff}(\Dir^{\phi}\circ Q_E)=
(N_{\phi}(\cl_{\phi})
\circ\pi^*\nabla^{E,\phi})\circ
 N_{\faff}(Q_E).$
\qed
\end{enumerate}
\end{proposition}
\par
Recall that for two vector fields $A,B\in\Gamma(\partial
X,\phitx)$, the commutator $[{\sf n}A,{\sf v}B]$ lies in
$\phindx$, but is not necessarily $0$. However, the normal
operator
$$N_{\phi}(\Dir^{\phi})=(N_{\phi}(\cl_{\phi})
\circ\pi^*\nabla^{E,\phi})=\Dir^{\phi,V}+{\sf
n}\cl_{\phi}\dd_{\phindx}$$
of $P(\lambda)$ has the following special property
\begin{lemma} \label{commute}
$\Dir^{\phi,V}$ and ${\sf n}\cl_{\phi}\circ \dd$ anticommute over $\partial X$.
\end{lemma}
\begin{proof}
First note that  it follows from \ref{exactconnection}(c) that for
$T$ tangent to the boundary
 $$[\nabla_T^{E,\phi},{\sf
n}]=0\quad\mbox{and}\quad [\nabla_T^{E,\phi},{\sf v}]=0.$$ Using
this, we can write
\begin{eqnarray*} [{\sf v}\cl_{\phi}\pi^*\nabla^{E},{\sf
n}\cl_{\phi}\pi^*\nabla^{E}]&=& ({\sf v}\cl_{\phi}, {\sf
n}\cl_{\phi})[\pi^*\nabla^E,\pi^*\nabla^E]
 +{\sf v}\cl_{\phi}[\pi^*\nabla^E,{\sf
 n}\cl_{\phi}]\pi^*\nabla^{E}\\
  & &\qquad -{\sf v}\cl_{\phi}[\pi^*\nabla^E,{\sf n}\cl_{\phi}]\pi^*\nabla^{E}
+ [{\sf v}\cl_{\phi}, {\sf
n}\cl_{\phi}]\pi^*\nabla^E\pi^*\nabla^E\\ &=& ({\sf v}\cl_{\phi},
{\sf n}\cl_{\phi})\pi^*F^E=0,
\end{eqnarray*}
which proves the claim.\qed
\end{proof}\par
\smallskip\noindent
This Lemma will allow us to calculate solutions for
$N_{\phi}(P(\lambda))$ by separating the horizontal and vertical
variables.

 \subsection{Inversion at ff}
In this Section we use the normal operator defined in the previous
Section to improve the error term obtained in Proposition
\ref{symbolicresult} at the front face $\faff$. This means that we
have to solve the equation
$$N_{\faff}(P(\lambda)Q_1(\lambda))=N_{\phi}(\Dir^{\phi})\circ
N_{\faff}(Q_1(\lambda)) \stackrel{!}{=}-N_{\faff}(R_0(\lambda))\in
\Psi_{{\rm sus}}^{\infty}(\phindx,E),$$ i.e. we have to invert the
operator $N_{\phi}(\Dir^{\phi})=(\Dir^{\phi,V}+{\sf n}\cl_{\phi}
\dd_{\phindx})$ over $\phindx/Y$.
\par
\smallskip
Clearly, sections in the null space of $\Dir^{\phi,V}$ will be of
special importance in this context. At each point $y\in Y$ the
null space of the operator $\Dir^{\phi,V}|_{Z_y}$ is given by
$${\cal K}_y={\rm null}(\Dir^{\phi,V})|_{Z_y}\subset
L^2(Z_y,\dvol_Z).$$ From now on we {\em assume} that
\begin{equation}\label{mainassumption}
\mbox{\bf\boldmath ${\cal K}=({\rm
null}(\Dir^{\phi,V})|_{Z_y})_{y\in Y}\rightarrow Y$ \quad is a
vector bundle over $Y$},
\end{equation}
i.e. ${\rm null}(\Dir^{\phi,V})$ is the space of sections of the
vector bundle ${\cal K}\rightarrow Y$. Denote by $\Pi_{\circ},
\Pi_{\perp}$ the projections onto this null space and its
orthogonal complement. At the boundary this gives us the
decomposition into ``zero modes'' and ``nonzero modes''
$$C^{\infty}(\partial X, E)\cong \Pi_{\circ}C^{\infty}(\partial
X,E) \oplus \Pi_{\perp}C^{\infty}(\partial X,E).$$ Such
decompositons also exist at the faces $\falf$, $\faphibf$, $\faff$
and it follows from Lemma \ref{commute} that $N_{\phi}(\phidir)$
preserves this decomposition over $\faff$.
 As explained in Appendix \ref{appzeronotation} it will be useful
 to consider spaces of polyhomogeneous,
conormal sections whose top (in the sense defined there)
coefficients lie in either part of the decomposition. A more
refined index notation adapted to this situation is also described
in Appendix \ref{appzeronotation}.
\par
\bigskip \noindent{\bf\boldmath The Dirac Operator on $\RR^n$}
\par
\medskip
 We start with the analysis of our model operator restricted to
 $\Pi_{\circ}$, i.e. the (constant coefficient) Dirac operator on $\RR^n$:
 Denote by $g_n$ the euclidean metric on $\RR^n$ and by $\dvol_n$ the euclidean
 volume form. Let ${\cal W}$ be a
Clifford module over $\RR^n$ with  Clifford action
$\cl_n\in \Gamma(T\RR^n,\End({\cal W}))$ and trivial connection $\dd$.
The Dirac operator is then
$\Dir_n=\cl\circ\dd$.
\par
Writing ${\sf R}$ for the tautological (or radial) vector field
in $T\RR^n$ the fundamental solution for $\Dir_n$ is
$$E_n=\Const_n | \cdot |^{-n}\cl_n({\sf R}), \quad \mbox{or}\quad
E_n({\bf y})= \Const_n \sum_{j=1}^nr^{-n}y^j\cl(\dd y^j),$$ when one prefers the use of
the standard euclidean coordinates ${\bf y}\sowie (y^1,\ldots,
y^n)$, $r=|{\bf y}|$, on $\RR^n$. Then
$$\Dir_n E_n = \delta \quad\mbox{in}\quad C^{-\infty}(\RR^n,\End({\cal W})).$$
\par
Our idea is to analyze the behavior of solutions for $\Dir_n$,
by first looking at solutions for the conformally transformed
Dirac operator $\Dir_b$ on the cylinder $\RR_+\times
S^{n-1}$. For this operator we can the use the known methods of the
$b$-calculus from \cite{Melaps}.\par
Thus identify $\RR^n\setminus 0$ with $\RR_+\times S^{n-1}$
and introduce the $b$-metric $g_b=r^{-2}g_n$. The corresponding $b$-volume form
is $\dvol_b=r^{-n}\dvol_n$. It then
follows again as in Appendix \ref{conformal} that ${\cal W}$ is a
Clifford module over $\RR_+\times S^{n-1}$, with Clifford action
$\cl_b= r\cl_n$ and Dirac operator
\begin{equation}\label{euclidcyl}
 \Dir_b=r^{(n-1)/2}r\Dir_nr^{-(n-1)/2}= r\Dir_n -\frac{n-1}{2}\cl_n(\dd r).
 \end{equation}
 This operator is $b$-elliptic and selfadjoint on $L^2_b(\RR_+\times
S^{n-1})$ and
 $r\Dir_n$ is a $b$-elliptic, $b$-differential operator on $\RR_+\times
S^{n-1}$. We have $$ L^2(\RR^n, {\cal W};\dvol_n)\cong r^{-n/2}L^2(\RR_+\times S^{n-1},{\cal W};\dvol_n)
\equiv r^{-n/2}L^2_b(\RR_+\times S^{n-1},{\cal W}).$$
The main result is now
\begin{lemma}\label{basicmap}
The operator $r \Dir_{n,\alpha} : r^{\alpha}H^{k+1}_b(\RR_+\times
S^{n-1},{\cal W})\rightarrow r^{\alpha}H^k_b(\RR_+\times S^{n-1},{\cal W})$ is
invertible if $\alpha\notin \ZZ$. More exactly the extended set of indicial
roots is $\quad i\NN\cup -i((n-1)-\NN)$, especially the indicial roots are all simple.
\end{lemma}
\begin{proof}
First, we know from the general $b$-calculus that the operator $
\Dir_{b,\alpha}$ is invertible except if $i\alpha$ is an indicial root, where it
has a null space in the ``extended $L^2$-sense''.
Since $\Dir_b|_{r=0}$ anticommutes with $\cl_b(\frac{\dd
r}{r})$, one easily finds that the set of indicial roots $\spec_I (\Dir_b)$
is invariant under multiplication with $-1$. Also it follows
immediately from (\ref{euclidcyl}) that
$\spec_I(r\Dir_n)=-i\frac{n-1}{2}+\spec_I(\Dir_b)$.  It therefore
suffices to show that $$\{z\in \spec_I(r\Dir_n)\,|\, \Im(z)\geq -(n-1)\} \stackrel{!}{=}
-i(n-1)\cup i\NN.$$
This is easy. Let $u_{\alpha}\in {\rm null}_-(r\Dir_{n,\alpha})\subset r^{\alpha}L^{\infty}(\RR_+
\times S^{n-1})$. If $\alpha >-n$ then
 the distribution $u_{\alpha}$ becomes
integrable in $0$ w.r.t. $\dvol_n$ and extends to a distribution (also denoted
$u_{\alpha}$) on $\RR^n$. This distribution fulfills $$\Dir_n u_{\alpha}=
\delta w,\quad\mbox{i.e.}\quad \Dir_n(u_{\alpha}-E_n(w))=0$$ in all of $\RR^n$
for some $w\in {\cal W}$. By elliptic regularity, we can conclude that
$u_{\alpha}-E_n(w)$ is in $C^{\infty}(\RR^n,{\cal W})$, from which the
claim follows.\qed
\end{proof}
This result can be used to solve $\Dir_n u=f$ for $f$ with compact
support in $\RR^n$. First, if $f\in C^{\infty}_c(\RR_+\times S^{n-1},{\cal W})$ we also have
$r f\in C^{\infty}_c(\RR_+\times S^{n-1},{\cal W})$. Choosing $\alpha
\notin \ZZ$ set $$u_{\alpha}=(r\Dir_{n,\alpha})^{-1}r f\in
r^{\alpha}H_b^{\infty}(\RR_+\times S^{n-1},{\cal W}).$$ Then the usual
results for elliptic $b$-pseudodifferential operators imply that
$u_{\alpha}$ has expansions at $r=0$ and $1/r=0$ of the form
\begin{eqnarray*}
u_{\alpha}& \sim & \sum_{\ZZ\ni j> \alpha}a_j({\bf \omega})r^j
\quad \mbox{for}\quad r\rightarrow 0\\ u_{\alpha}& \sim &
\sum_{\ZZ\ni j> -\alpha}b_j({\bf \omega})\rho^j \quad \mbox{for}\quad
\rho=r^{-1}\rightarrow 0.
\end{eqnarray*}
If we choose $\alpha>-n$, the distribution $u_{\alpha}$ again extends
to a distribution on $\RR^n $.
Writing $\BB^n:={\RC}(\RR^n)$ for the radial compactification of
$\RR^n$, we have shown the more precise statement
$$u_{\alpha}- E_n(w)\in \rho^{n-1}C^{\infty}(\BB^n,{\cal W}).$$
 Of course, due to the translational invariance of the problem,
 this argument works for any
$f\in C^{\infty}_c(\RR^n,{\cal W})$ by first shifting the support of $f$
away from $0$. We now show
\begin{lemma}\label{dirrn}
\begin{enumerate}
\item $\Dir_n:\rho^{\alpha}H^{\infty}_b(\BB^n,{\cal W})
\rightarrow \rho^{\alpha+1}H^{\infty}_b(\BB^n,{\cal W}).$
\item For any $f\in{\cal A}^F(\BB^n,{\cal W})$ there exists a
$$u\in {\cal A}^U(\BB^n
 ,{\cal W}),\qquad U=(n-1+\NN)\cup ((\ZZ_{\geq F}\overline{\cup} F)-1),$$
which solves $\Dir_n u=f$.
\end{enumerate}
\end{lemma}
\begin{proof}
Choose a cutoff function $\varphi\in C_c^{\infty}(\BB^n)$. Then,
the compactly supported section  $\varphi f$, as treated above, is
responsible for the term $n-1+\NN$ in the index set $U$.\par
Choosing a noninteger $\alpha<F$ we know that $(1-\varphi)f\in
\rho^{\alpha}H^{\infty}_b(\RR_+\times S^{n-1},{\cal W})$, and by Lemma
\ref{basicmap} there is $v_{\alpha}\in
\rho^{\alpha-1}H^{\infty}_b(\RR_+\times S^{n-1},{\cal W})$ satisfying
$\Dir_n v_{\alpha}=(1-\varphi)f$ . By the general theory of
$b$-pseudodifferential operators (compare \cite{Melaps}, Chapter
5) $$ v_{\alpha}\in {\cal A}^G(\RR_+\times S^{n-1},{\cal W}), \quad
G(r=0)=\{1-\alpha<j\in \ZZ\}, \quad
G(\rho=0)=(F\overline{\cup}\{\alpha<j\in\ZZ\})-1.$$ Now, choosing
another cutoff function $\eta\in C^{\infty}_c(\BB^n)$ with
$(1-\eta)(1-\varphi)=1-\varphi$, we can write
$$\Dir_n(1-\eta)v_{\alpha}=(1-\eta)\Dir_n v_{\alpha}-\cl(\dd
\eta)v_{\alpha} =(1-\varphi)f -\cl(\dd \eta)v_{\alpha}.$$ But
$\cl(\dd \eta)$ is compactly supported and can therefore be
treated as in the first part.\qed
\end{proof}
We have solved the model problem in the zero-modes at $\faff$. The
model problem in the nonzero-modes is easy. Noting that
$$\Psi_{{\rm sus}}^{-\infty,F}(\phindx,E)={\cal A}^F
(\RC(\phindx\times_Y\partial X), \CB_{\phi})={\cal
A}^F(\faff,\CB_{\phi})$$ we get
\begin{lemma}[Solution of the Model Problem at {\boldmath $\faff$}]\label{dirff}
\quad\par\noindent
\begin{enumerate}
\item For every $f_{\circ} \in \Pi_{\circ}\Psi_{{\rm sus}}^{-\infty,F}(\phindx,E)$
 there exists  $$u_{\circ} \in \Pi_{\circ}\Psi_{{\rm sus}}^{-\infty,U}(\phindx,E),\qquad
  U=(h+\NN)\cup ((F\overline{\cup}\ZZ_{\geq
F})-1),$$ which solves $N_{\faff}(\phidir) u_{\circ}=f_{\circ}$.
\item For every $f_{\perp} \in \Pi_{\perp}\Psi_{{\rm
sus}}^{-\infty,F}(\phindx,E)$
 there exists
an $u_{\perp} \in \Pi_{\perp}\Psi_{{\rm
sus}}^{-\infty,F-1}(\phindx,E)$ which solves $N_{\faff}(\phidir)
u_{\perp}=f_{\perp}$.
\end{enumerate}
\end{lemma}
\begin{proof}
First, note that $N_{\faff}(\Dir^{\phi})$ restricted to
$\Pi_{\circ}H^{\infty}_b(\faff,\CB_{\phi})$ is just the family of
(constant coefficient) fibre Dirac operators given by ${\sf n
}\cl_{\phi}\circ d$ on $\phindx\times_Y\partial X/F_{\phi}$. Hence
(a) follows by applying Lemma \ref{dirrn} fibre by fibre.\qed
\end{proof}\par\smallskip
Using Lemma \ref{dirff} and the fact that $N_{\faff}(P(\lambda))=
N_{\faff}(\phidir)$ it is now straightforward to prove the
following refinement of Proposition \ref{symbolicresult}
\begin{proposition}\label{invff}
There are holomorphic families $$\lambda \mapsto Q_1^N(\lambda)\in
\psiphi^{-1,\cal E}(X,E),\quad\mbox{and} \quad \lambda \mapsto
R_1^N(\lambda)\in \psiphi^{-\infty,\cal F}(X,E),$$ such that such
that $$P(\lambda)Q_1^N(\lambda)=1+R_1^N(\lambda),$$ and with index
sets ${\cal E}$, ${\cal F}$  given by
 $$F_{\faphibf}=(h+1+\NN)^{\overline{\cup }N},\quad F_{\faff}=N+\NN$$
$$ E_{\faphibf}=(((h,0),\ldots,(h,N-1))^{\circ},
(h+1+\NN)^{\overline{\cup}N})=[\circ](h+\NN)^{\overline{\cup}N},
\quad E_{\faff}=\NN,$$  (all other index sets  are  $\infty$).
\par \medskip \noindent
 Note that the special case $N=1$ will suffice for the construction of the resolvent!
\end{proposition}
\begin{proof}
For purposes of notational simplicity we are going to ignore
coefficients. From the symbolic construction, we have holomorphic
families $$ \lambda \mapsto Q_0(\lambda)\in \Psi_{\phi,{\rm
cl}}^{-1}(X), \qquad \lambda \mapsto
R_0(\lambda)\in\Psi_{\phi,{\rm cl }}^{-\infty}(X),
\quad\mbox{with}\quad P(\lambda)Q_0(\lambda)=1+R_0(\lambda).$$
Setting $S(\lambda)=-R_0(\lambda)$, we want to solve
$$P(\lambda)G(\lambda)=S(\lambda)$$ holomorphically in $\lambda$
to order $N$ at $\faff$. For this, we write $G(\lambda)$ and
$S(\lambda)$ in a (finite) Taylor series (sum): $$
G(\lambda)=G_0(\lambda)+x'G_1(\lambda)+(x')^2G_2(\lambda)+
\ldots+(x')^{N-1} G_{N-1}(\lambda),$$ $$
S(\lambda)=S_0(\lambda)+x'S_1(\lambda)+(x')^2S_2(\lambda)+
\ldots+(x')^{N-1} S_{N-1}(\lambda)+(x')^NS_N(\lambda),$$ and the
maps $$\lambda\mapsto G_j(\lambda)\in \Psi_{\phi}^{-\infty,{\cal
I}_j}(X) \quad\mbox{and} \quad \lambda\mapsto S_j(\lambda)\in
\Psi_{\phi}^{-\infty,{\cal J}_j}(X)$$ are holomorphic. As part of
the proof we will show that $$ I_{j,\faff}=0, \quad
I_{j,\faphibf}= [\circ](h-j+\NN)^{\overline{\cup}j+1}, \quad
J_{j,\faff}=0, \quad J_{j,\faphibf}=(h-j+1+\NN)^{\overline{\cup}j}
.$$ In the first step, we have to solve
$$N_{\faff}(S_0(\lambda))=N_{\faff}(P(\lambda)G_0(\lambda))
 =N_{\faff}(\Dir^{\phi})N_{\faff}(G_0(\lambda)).$$
Write $s_0(\lambda)$ for $N_{\faff}(S_0(\lambda))$. It is a
holomorphic family $\lambda \mapsto s_0(\lambda)\in
\cdotinfty{}(\faff)$. By Lemma \ref{dirff}, we can find a
holomorphic family $$\lambda \mapsto
g_0(\lambda)=g_0(\lambda)^{\perp}+g_0(\lambda)^{\circ} \in
\Pi_{\perp}\cdotinfty{}(\faff)+ \Pi_{\circ}{\cal
A}^{h+\NN}(\faff), \quad
N_{\faff}(\Dir^{\phi})g_0(\lambda)=s_0(\lambda), $$ which can be
extended to a holomorphic family $$\lambda\mapsto G_0(\lambda)\in
\Psi^{-\infty,{\cal I}_0} ,
 \quad I_{0,\fabf}=((h,0)^{\circ},h+1+\NN), \quad
N_{\faff}(\Dir^{\phi}G_0(\lambda)-S_0(\lambda))=0.$$ Thus, we have
shown the case $j=0$.\par Assuming inductively that the equation
has been solved up to order $j$, with index sets as indicated, we
set
 $$ S_{j+1}(\lambda)=(x')^{-1}(P(\lambda)G_j(\lambda)-S_j(\lambda)).$$
Since $P(\lambda)G_j(\lambda)$ has index set
$(h+1-j+\NN)^{\overline{\cup}j+1}$ and $S_j(\lambda)$ has index
set $(h-j+1+\NN)^{\overline{\cup}j}$ at $\fabf$, the index set for
$S_{j+1}(\lambda) $ is just ${\cal J}_j$. Now, the equation
$$N_{\faff}(\Dir^{\phi}G_{j+1}(\lambda)-S_{j+1}(\lambda))=0$$ can
be solved holomorphically in $ \lambda$ by Lemma \ref{dirff} and
$G_{j+1}$ has index set ${\cal I}_{j+1}$. Setting
$Q_1^N(\lambda)=Q_0(\lambda)+G(\lambda)$ and
$R_1^N(\lambda)=(x')^N S_N(\lambda)$ yields the result. \qed
\end{proof}
For reference in the next Sections we note the case $N=1$ as a
Corollary. Set $Q_1(\lambda)= x' Q_1^1(\lambda)$ and
$R_1(\lambda)= \frac{x'}{x}R_1^1(\lambda)$:
\begin{corollary}\label{finalff}
$
D_{\lambda}Q_1(\lambda)=1+R_1(\lambda)$ with holomorphic families
$$\lambda\mapsto R_1(\lambda)\in\Psi_{\phi}^{\infty, {\cal
F}}(X,E),\qquad F_{\faphibf}=(h+1+\NN),\qquad F_{\faff}=1+\NN$$
 $$\lambda\mapsto Q_1(\lambda)\in\Psi_{\phi}^{-1, {\cal
E}}(X,E),\qquad E_{\faphibf}=[\circ](h+1+\NN),\quad
E_{\faff}=1+\NN,$$ the other terms in the index sets are
$\infty$.\qed
\end{corollary}
Thus, we have solved the resolvent equation (\ref{resolventeq}) up
to a  remainder term, which is $C^{\infty}$ in the interior and
which vanishes to first order at the front face. Unfortunately, by
Section \ref{secmapprop}, this new remainder is still not compact,
since it does not vanish to positive order at $\faphibf$ (w.r.t.
to the coefficients $\CB_b$, which are the relevant ones at that
face). Thus the next two Sections will describe how to improve the
remainder at $\faphibf$.

\subsection{Reduced Normal Operators at bf and $\phi$bf}
\label{normalopsbf}
In this Section we define the normal operator at the face
$\faphibf$ and describe the model action of the Dirac operator
$\Dir$ at this face. Recall the relation between the $b$- and the
$\phi$-coefficient bundles: $\CB_b=(x')^{h+1}\CB_{\phi}$. The
normal operators at $\faphibf$ and $\fabf$ are defined by
\begin{definition} Let ${\cal E}$, ${\cal F}$ be index sets for $X^2_{\phi}$
and $X^2_b$ respectively, with $E_{\faphibf}= h+1+\NN$ and
$F_{\fabf}=\NN$.
\begin{enumerate}
\item $N_{\fabf}:\Psi_b^{-\infty, {\cal F}}(X,E)\longrightarrow
  {\cal A}^{(F_{\falf},F_{\farf})}(\fabf,\CB_b)\qquad A\longmapsto A|_{\fabf,\CB_b}$
\item $N_{\faphibf}:\Psi_{\phi}^{-\infty,{\cal E}}(X,E)
\longrightarrow {\cal A}^{(E_{\falf}, E_{\farf},
E_{\faff})}(\faphibf,\CB_b) \qquad A\longmapsto
A|_{\faphibf,\CB_b}$\qed
\end{enumerate}
\end{definition}
The corresponding exact sequences are $$0 \longrightarrow
\rho_{\fabf}\Psi_b^{-\infty, {\cal
F}}(X,E)\longrightarrow\Psi_b^{-\infty, {\cal
F}}(X,E)\stackrel{N_{\fabf}}{\longrightarrow}{\cal
A}^{(F_{\falf},F_{\farf})}(\fabf,\CB_b)\longrightarrow 0$$
 $$0 \longrightarrow
\rho_{\faphibf}\Psi_{\phi}^{-\infty, {\cal
E}}(X,E)\longrightarrow\Psi_{\phi}^{-\infty, {\cal
E}}(X,E)\stackrel{N_{\faphibf}}{\longrightarrow}{\cal
A}^{(E_{\falf},E_{\farf},
E_{\faff})}(\faphibf,\CB_b)\longrightarrow 0.$$ Note also, that
whenever $E_{\faff}=\infty$ we have $$\Psi_{\phi}^{-\infty,{\cal
E}}(X,E)\subset\Psi_{b}^{-\infty,(E_{\farf}-h-1,\NN,E_{\falf})}(X,E),$$
and the normal map $N_{\faphibf}$, restricted to the LHS,  then
equals $N_{\fabf}$. We will not try to describe the general
behavior of the above sequences under composition of operators.
Instead we will restrict ourselves to the analysis of the normal
action of the Dirac operator $\Dir$, which will be sufficient for
our purposes.
\par
\smallskip
Recall that the Dirac operator $\Dir$ does {\em not} map calculi
of the above type to themselves, since its lift to $X^2_b$ or
$X^2_{\phi}$ is not tangent to $\fabf$ or $\faphibf$. From
Appendix \ref{appzeronotation} we know that
 $$ \Dir_{\lambda} :
\Psi_{b}^{-\infty, (E_{\farf} ,[\circ]E_{\falf},
[\circ]\NN_{\fabf}) } (X,E) \longrightarrow \Psi_{b}^{-\infty,
(E_{\farf},E_{\falf}, \NN_{\fabf})}(X,E)$$ $$ \Dir_{\lambda} :
\Psi_{\phi}^{-\infty, (E_{\farf} ,[\circ]E_{\falf},
[\circ]\NN_{\faphibf} ,E_{\faff}) } (X,E) \longrightarrow
\Psi_{\phi}^{-\infty,(E_{\farf},E_{\falf},
\NN_{\faphibf},E_{\faff}-1)}(X,E),$$ i.e. $\Dir_{\lambda}$ maps
nicely, when we restrict to conormal distributions with top
cofficients in the zero modes at $\falf$ and $\faphibf$. This
allows to define the normal action $N_{\fabf}(\Dir_{\lambda})$ on
the sections of ${\cal K}$ over $\fabf$ (resp.
$N_{\faphibf}(\Dir_{\lambda})$ on $\faphibf$), defined as usual by
$$N_{\fabf}(\Dir_{\lambda})N_{\fabf}(A)=\Pi_{\circ}N_{\fabf}(\Dir_\lambda
A ),\quad\mbox{for}\quad A\in \Psi_b^{-\infty}(X,E),\quad
N_{\fabf}(A)\in {\cal K}. $$
\par\medskip
\noindent{\bf The Model Operator at {\boldmath $\fabf$}}
\par\medskip
 Let us describe the action of $N_{\fabf}(\Dir_{\lambda})$ at
$\fabf$. First, identify $\fabf-\farf$ with the fibre bundle $$
\beta_{b,L}:\partial X\times\partial
X\times[0,\infty[_{s}\rightarrow \partial X, \quad s=x/x',\qquad
\beta_{b,R}^*\bomega(X)=\bomega_{\rm fib}(\beta_{b,L}).$$ Thus,
writing $B=Y\times[0,\infty[$, the space of sections in the null
space of $\Dir^{\phi,V}$ over $\fabf$ is
$$\Pi_{\circ}C^{\infty}(\fabf,\CB_b)=C^{\infty}(\partial X,
C^{\infty}(B,{\cal K}\otimes\bomega(B))),$$ and
$N_{\fabf}(\Dir_{\lambda})$ acts as the operator  on
\begin{equation}\label{Ibactson}
C^{\infty}(B,{\cal K}\otimes\bomega(B))=
\Pi_{\circ}C^{\infty}(\partial X\times \RR_+,E_{\partial X }),
\end{equation}
constant in the parameter space $\partial X$, given by fixing
coefficients at the boundary
\begin{equation}
I_b(\Dir_{\lambda}):=\Pi_{\circ}\Dir_{\lambda}|_{\partial
X,{}^bTX}=[\Pi_{\circ}\frac{1}{x}\Dir^{\phi}+\frac{h}{2}\cl_{\phi}(\frac{\dd
x}{x^2})-\lambda]|_{\partial X,{}^bTX}.
\end{equation}
This definition is {\em much less innocent} than it looks at first
sight. The expression actually   may involve hidden endomorphisms originating from
vertical differentiations compressed to ${\cal K}$. We therefore
prefer another way to define this ``indicial operator'':
\begin{lemma}\label{indicialfamily}
\begin{enumerate}
\item Let $\xi \in C^{\infty}(X,E)^{\circ}$, i.e. $\Dir_{\lambda} \xi\in
C^{\infty}(X,E)$.
Then $$I_b(\Dir_{\lambda})(\xi|_{\partial X}):=
      \Pi_{\circ}(\Dir_{\lambda}\xi)|_{\partial X}$$
is well defined and only depends on $\xi|_{\partial X}\in \Pi_{\circ}C^{\infty}
(\partial X,E)$.
\item The family
$$ z\mapsto  I_{b}(\Dir_{\lambda})(z)(\xi|_{\partial X})
:=\Pi_{\circ}(x^{-iz}\Dir_\lambda x^{iz}\xi
|_{\partial X})\in \Diff^1(Y,{\cal K}),$$ is well
defined (i.e. independent of the choice of extension $\xi$ of $\xi|_{\partial X}$)
and holomorphic in $z$.
\item  $I_b(\Dir_\lambda)(z)=
M_{x\rightarrow z} I_b(\Dir_\lambda) M_{z\rightarrow x} $,
where $M$ denotes the Mellin
transform (see (\ref{meldef})).
\item Let $\xi\in C^{\infty}(X,E)$ such that $\Dir_{\lambda}^l\xi\in
C^{\infty}(X,E)$.
Then $$\Pi_{\circ}(\Dir_{\lambda}^l\xi)|_{\partial X}= [I_{b}(\Dir_{\lambda})]^l
(\xi|_{\partial X}).$$
Moreover  any $\xi|_{\partial X}$ can be extended to $\widetilde{\xi}$
such that $(\Dir_{\lambda}^l\widetilde{\xi})|_{\partial X}=
[I_{b}(\Dir_{\lambda})]^l (\xi|_{\partial X}).$
\end{enumerate}
\end{lemma}
\begin{proof}
In (a) and (b), only the independence of the extension of $\xi|_{\partial X}$ might be
worth a little thought. It follows from the fact, that for $\eta
\in C^{\infty}(X,E)$ we have $$\Pi_{\circ}(\Dir_{\lambda} x\eta)|_{\partial X}=
\Pi_{\circ}  (x\Dir_\lambda \eta +x\cl_d(\frac{\dd x}{x})\eta)|_{\partial X}=
\Pi_{\circ}(x\Dir \eta)|_{\partial X}=0.$$
Thus, the expression defining $I_{b}(\Dir_{\lambda})(\xi|_{\partial X})$ vanishes for
the difference of two extensions of $\xi|_{\partial X}$.
\par
Part (d) is proved by induction. We have proved the case $l=1$ in
(a), now assume that
$$\Pi_{\circ}\Dir_{\lambda}^{l-1}\xi|_{\partial
X}=[I_b(\Dir_{\lambda})]^{l-1}(\xi|_{\partial X}).$$
Then starting with $\xi$ such that $\Dir_{\lambda}^l\xi\in
C^{\infty}(X,E)$ means that we have $\Dir_{\lambda}^{l-1}\xi\in
C^{\infty}(X,E)^{\circ}$. Therefore,
$$\Pi_{\circ}\Dir_{\lambda}^{l}\xi|_{\partial
X}=\Pi_{\circ}\Dir (\Dir_{\lambda}^{l-1}\xi)|_{\partial
X}=I_b(\Dir_{\lambda})[I_{b}(\Dir_{\lambda})]^{l-1}(\xi|_{\partial
X}),$$
by (a) and the hypothesis of the induction.
\qed
\end{proof}
\par
\medskip
\noindent {\bf Geometrical Description of the Indicial Operator}
\par
\medskip
Let us take the time to describe the geometry underlying
 our definition of the indicial family $I_b(\Dir)(z)$.
This will not be
strictly necessary for the analysis of the model problem at $\faphibf$
and $\fabf$, but gives some insight in the specialties of
``$\phi$-geometry''.
Recall that bwe write $\nu=\frac{1}{x}{\sf X}_d$. Then, starting with $\xi|_{\partial X}\in \Gamma(Y,{\cal K})$
and an extension$\xi \in C^{\infty}(X,E)^{\circ}$, we can use Lemma
\ref{indicialfamily} to get at  the boundary
\begin{eqnarray}
 I_b(\Dir)(iz)(\xi|_{\partial X})&=&\Pi_{\circ}(x^z
 \Dir x^{-z}\xi|_{\partial X})=
-(z+v/2)\Pi_{\circ}\cl_{d}(\frac{\dd x}{x})\xi
+\Pi_{\circ}\Dir^d\xi\nonumber
\\&=&
-(z+v/2)\Pi_{\circ}\cl_{d}(\frac{\dd x}{x})\xi
+\Pi_{\circ}\nabla_{\nu}^{E,d}x\Dir^{d}\xi\nonumber
\\ &=&
-(z+v/2)\Pi_{\circ}\cl_{d}(\frac{\dd x}{x})\xi
+\Pi_{\circ}\left([\nabla^{E,d}_{\nu},\cl_{d}\circ x \nabla^{E,d}_{\cdot}]
\xi
+x\Dir^d\nabla^{E,d}_{\nu}\xi\right)
\nonumber
\\ &=&
-(z+v/2)\Pi_{\circ}\cl_{d}(\frac{\dd x}{x})\xi
+\Pi_{\circ}\cl_{d}(x F^{E,d}(\nu,\cdot))\xi
-\Pi_{\circ}\cl_{d}\nabla^{E,d}_{x \nabla^{\phi}_{ \cdot}\nu}\xi.
\label{ibintermed}
\end{eqnarray}
In the last step we have used that $\Pi_{\circ}x\Dir^d
\nabla^{E,d}_{\nu}\xi|_{\partial X}=0$, the definition
of the curvature $F^{E,d}$ and the fact that
$A-x\nabla^d_A\nu=-x\nabla^{\phi}_A\nu$. Note that we are only referring
to the ``$d$-structure'' of the bundle $E$ in this
calculation. It is therefore valid for any Dirac operator $\Dir$.
It now follows  from Lemma \ref{halpha} that
$\nabla^{\phi}_W\nu\in{}^bTX$ for $W\in\Gamma(\phitx)$, and
 $$\nabla^{\phi}_A\nu=-{\sf h}\frac{1}{x}A-\frac{\dd x}{x^2}(A){\sf
   X}_d+O(x){}^bTX,\quad\mbox{for}\quad A\in\Gamma([\phindx]).$$ Also,
   using the special extension
described in Lemma \ref{extension}, we find that $ W\mapsto
\nabla^{\phi}_{\overline{W}}\nu$ is a well defined tensor in
$C^{\infty}(\partial X,V^*\partial X\otimes\phitx)$, independent
of the vector field used in the definition of the extension (which
therefore could be different from $\nu$). Writing
\begin{equation}\label{mystery}
2\langle\nabla^{\phi}_{\overline{W}}
\nu,\overline{V}\rangle_{\phi}
=(L_{\nu}g_{\phi})(\overline{W},\overline{V}), \quad
V,W\in\Gamma(V\partial X),
\end{equation} we see that this describes yet another
aspect of the deviation of $g_{\phi}$ from a product metric near
the boundary.\par
 Since the Clifford action $\cl_d$ is parallel w.r.t.
$\nabla^{E,d}$ we can use the special extension introduced in
Section \ref{boundarygeom} and decompose $$\cl_d=\overline{{\sf
n}\cl_d}+\overline{{\sf v}\cl_d},\quad \nabla^{E,d}_{\nu}\overline{{\sf
n}\cl_d}=\nabla^{E,d}_{\nu}\overline{{\sf v}\cl_d}=0.$$
 Plugging this into (\ref{ibintermed}) we get
$$
I_b(\Dir)(iz)\xi=\Pi_{\circ}\left(
 -(z+v/2)\cl_{d}(\frac{\dd
x}{x})\xi +(\overline{{\sf n
}\cl_{d}})\nabla^{E,d}_{\cdot}\xi
+\cl_{d}(x F^{E,d}(\nu,\cdot))\xi
-\overline{{\sf v
}\cl_{d}}\nabla^{E,d}_{x \nabla^{\phi}_{\cdot}\nu}\xi \right).$$
Thus $I_b(\Dir)$ is a $b$-differential
operator on the vector bundle ${\cal K }\rightarrow B$ pulled back
from $Y$ to $B=Y\times\RR_+$.
As explained, the last  summand
is an additional {\em endomorphism} stemming from the fact that we
are looking at a situation where the metric $g_{\phi}$ is not a
product $\phi$-metric.
Let us now analyze the
curvature term  in the formula. Of course,
this term depends on the choice of Clifford bundle $E$, but we
can always look at the general decomposition
$$ F^{E,d}(\nu,A)=\frac{1}{2}\cl_{d}(R^d(\nu,A))+F^{E/S,d}(\nu,A).$$
The contribution coming from the first summand
 can then be further analyzed using Proposition \ref{rdatthedel}:
\begin{eqnarray}
\frac{1}{2}\cl_{d,1}\cl_d(x R^d(\nu,\cdot_1))&=&\frac{1}{4}(x\cl_d,\cl_d,\cl_d)
\langle R^d(\nu, \cdot)\cdot,\cdot \rangle_d\nonumber
\\ &=&
\frac{1}{4}({\sf v}\cl_d,{\sf v}\cl_d,{\sf v}\cl_d)
\langle x R^d(\nu, \cdot)\cdot,\cdot \rangle_d
+\frac{1}{2}({\sf v}\cl_d, {\sf n}\cl_d, {\sf v}\cl_d)
\langle S_d( \cdot)\cdot,\cdot \rangle_{d}\label{rdsdint}
\end{eqnarray}
Using the symmetries of the curvature tensor and the description of
$S_{\phi}$ given in Lemma
\ref{lemnablab},  the reader can check that (\ref{rdsdint})can be
written as
$$
\frac{1}{2} {\sf v}\cl_d({\rm Ric}^{d,V}({\sf X}_d,\cdot))+\frac{v}{2}\cl_d(\frac{\dd x}{x})
+\frac{1}{2}{\sf n}\cl_d(\tr(S_M))
+\frac{1}{4}({\sf n}\cl_d, {\sf v}\cl_d, {\sf v}\cl_d)\langle x B_{\phi}(\cdot,cdot),
\cdot\rangle_d.$$
For $A\in{}^dTX$ we have used the notation
$$ {\rm Ric}^{d,V}({\sf X}_d,A)= ({\sf v}g_d)^{-1}_{aa'}\langle R^d({\sf X}_d,\cdot_a)
\cdot_{a'},A\rangle_d$$
for the ``vertical'' Ricci curvature and $S_M(\cdot){\sf n}A \in \End(V\partial
X)$
is the second fundamental form of the boundary fibration.  Definition
\ref{Sphi} can be used to check that all these terms are  in $
C^{\infty}(\partial X,\End(E))$.
\par
\medskip
To be able to interpret this formula, we need to note some more of
the {\bf geometrical data} of this construction:
\begin{itemize}
\item $g_{B,b}$ is the metric on ${}^b\,TB$ given by $(\frac{\dd x}{x})^2+g_B|_{x=0}$.
\item The Levi Civita connection $\nabla^{B,b}:\Gamma({}^b\,TB)
\rightarrow\Gamma(T^*B\otimes{}^b\,TB)$ is the connection induced
by the connection
$\nabla^d:\Gamma(\phi^*{}^bTB)\rightarrow\Gamma(\phi^*{}^bTB\otimes
T^*Y)$ over the boundary (see also Lemma \ref{lemnablab}). Here,
elements of $TY$ act on $\Gamma(\phi^*{}^bTB)$ via their {\em
horizontal lift} to $TX$
\item $\cl_{B,b}\in\Gamma({}^b\,TB\otimes\End({\cal K}))$ is the
($\RR_+$-invariant) Clifford action induced by $\Pi_{\circ}{\sf h
}\cl_d|_{\partial X}$ on ${\cal K}$ (Here we use
$[\Dir^{\phi,V},{\sf n}\cl_{d}]=0$ at the boundary).
\item The grading operators $\eps=\eps_{\phitsx}$ and
$\eps_B=\eps_{{}^bT^*B}$ induce
 gradings on ${\cal K}$.
\item The hermitian metric $\langle \cdot,\cdot\rangle_{\cal K}$ is the
   $L^2$-metric:
   $\langle\xi,\eta\rangle_{\cal K}=\phi_*(\langle\xi,\eta\rangle_{E}\dvol_Z)$.
\item For $A\in C^{\infty}(Y,\End({\cal K}))$:
$$\tr_{\cal K}(A)=\int_{\partial
X/Y}\tr_E(A)\dvol_Z=\phi_*[\tr_E(A)\dvol_Z] \in C^{\infty}(Y)$$
\item $\widetilde{\nabla}^{\cal K}$ is the pullback to $B$ of the compressed
connection $\Pi_{\circ}\nabla^{E,d}$ over $Y$.
 Here again,
elements of $TY$ act  via their {\em
horizontal lift} to $TX$
\end{itemize}
In this context the usual problem  about the non-unitarity of the
compressed connection $\widetilde{\nabla}^{\cal K}$ arises: For
$T\in \Gamma(TB)$ denote by $\phi^*T\in \Gamma(H(\partial X\times
\RR_+))$ the horizontal lift of $T$. Then
\begin{eqnarray*}
T\cdot\langle \xi,\eta\rangle_{\cal K}&=& T\cdot\phi_*[\langle
\xi,\eta\rangle_E\dvol_Z]\\ &=& \phi_*[\langle
\nabla^{E,d}_{\phi^*T}\xi,\eta\rangle_E\dvol_Z]
    +\phi_*[\langle
    \xi,\nabla^{E,d}_{\phi^*T}\eta\rangle_E\dvol_Z]
    +\phi_*[\langle \xi,\eta\rangle_E\Int(\phi^*T)\dd \dvol_Z].
\end{eqnarray*}
The last summand in this expression (compare Appendix
\ref{geomtens} or \cite{BGV} Chapter 9) is just $$-\phi_*[\tr
(\cdot S_M(\cdot)T)\langle \xi,\eta\rangle_E\dvol_Z]$$ and the
unitary (Clifford-) connection on ${\cal K}$ is given by
$$\nabla^{{\cal K }}_T\xi:=\widetilde{\nabla}^{{\cal
K}}_{T}\xi+\frac{1}{2} \Pi_{\circ}\tr(\cdot S_M(\cdot)T)\xi.$$ We
can summarize our results as follows:
\begin{proposition} \label{localibd}In
terms of the connection $\nabla^{\cal K}$ the indicial operator of
$\Dir$ has the form
\begin{eqnarray*}
I_b(\Dir)&=& \cl_{B,b}\circ\nabla^{\cal
    K}
\\ && +
\Pi_{\circ}\left[ {\sf v}\cl_{d}(F^{E/S,d}({\sf X}_d,\cdot))
+ \frac{1}{2} {\sf v}\cl_d({\rm Ric}^{d,V}({\sf X}_d,\cdot))
+\frac{1}{4}({\sf n}\cl_d, {\sf v}\cl_d, {\sf v}\cl_d)\langle x B_{\phi}(\cdot,cdot),
\cdot\rangle_d
-\overline{{\sf v
}\cl_{d}}\nabla^{E,d}_{x \nabla^{\phi}_{\cdot}\nu}\right].
\end{eqnarray*}
This is an odd operator w.r.t. $\eps$. It is also odd w.r.t.
$\eps_{B}$ whenever $h+1$ is even. Recall that the last
summand is an {\em endomorphism} on ${\cal K}$.\qed
\end{proposition}
Especially, this Proposition shows that in general $I_b(\Dir)$ is
{\em not} the induced Dirac operator on ${\cal K}$! However, the last three
summands on the RHS vanish when the
metric $g_{\phi}$ is a product $\phi$-metric. Also, if $X$ is spin  and
$E=S(X)\otimes E'$, we have
$F^{E/S,d}(\nu,\cdot)=F^{E',d}(\nu,\cdot)$,
which vanishes whenever the twisting bundle $E'$
has product structure near $\partial X$.
\par
\bigskip \noindent{\bf Inverting the Indicial Family}
\nopagebreak[4]
\par\medskip
 We now
want to analyze the structure of the indicial family
$I_b(\Dir_\lambda) (z)$ in more detail. We emphasize that we make
no assumption on the parity of the dimension of the base space
$Y$ at this stage. The operator $I_b(\Dir)$ can be written as
$$I_b(\Dir) :=
\Dir_Y+\cl_{d}(\frac{\dd x}{x})x\frac{\partial }{\partial x}
\quad\mbox{with}\quad
\Dir_Y(\xi|_{\partial X})=I_b(\Dir)(0)(\xi|_{\partial X})
=\Pi_{\circ}(\Dir\xi)|_{\partial X}.$$
Proposition \ref{localibd}, or the calculation
\begin{eqnarray*}
[\Dir^d,\cl_{d}(\frac{\dd
x}{x})]|_{\partial X}&=&
\left(\cl_{d}[\nabla^{E,d},\cl_{d}(\frac{\dd x }{x})]
+[\cl_{d},\cl_{d}(\frac{\dd x }{x})]\nabla^{E,d}\right)|_{\partial X}
\\ &=&(\cl_d,\cl_d){\sf v}g_d(\cdot,\cdot)|_{\partial X}-2\nabla_{{\sf
X}_d}^{E,d}|_{\partial X}=-v,
\end{eqnarray*}
 show that $D_Y$
anticommutes with $\mbox{\boldmath $\sigma$}= \cl_{d}(\frac{\dd
x}{x})$. Also, $\Dir_Y$ is selfadjoint with purely discrete
spectrum. This implies that for any eigenvalue $a$ of $\Dir_Y$ the
negative $-a$ is also an eigenvalue and the corresponding
eigenspaces are intertwined by $\mbox{\boldmath $\sigma$}$. Denote
by $P(\pm a)$ the projector onto the eigenspaces for $\pm a$.
These are invariant subspaces for the indicial family
$I_b(\Dir_{\lambda})(z)$ and for $a\neq 0$ we get the matrix
representation $$ P(\pm a) I_b(\Dir_{\lambda})(z)\sowie
\zwzw{a-\lambda}{-iz}{iz} {-a-\lambda},$$ w.r.t. a basis of
eigenvectors of $\Dir_Y$ of the form $(\mbox{\boldmath
$\psi$},\mbox{\boldmath $\sigma$} \mbox{\boldmath $\psi$})$. The
inverse, whenever it exists, has the form
\begin{equation}\label{indinverse}
 P(\pm a) I_b(\Dir_{\lambda})(z)^{-1}\sowie \frac{1}{\lambda^2-a^2-z^2}
\zwzw{-a-\lambda}{iz}{-iz}{a-\lambda}.
\end{equation}
In the case $a=0$ we have
\begin{equation}\label{indinverseb}
P(0)I_b(\Dir_{\lambda})(z) = iz\mbox{\boldmath $\sigma$}
-\lambda,\qquad
   P(0)I_b(\Dir_{\lambda}(z))^{-1}
   =\frac{1}{z^2-\lambda^2}(iz\mbox{\boldmath $\sigma$}+\lambda).
\end{equation}
Thus, writing $\spec(\Dir_Y)$ for the (necessarily symmetric) set
of eigenvalues of $\Dir_Y$, we find the following description of
the set of indicial roots, i.e. those points $z$, where
$I_b(\Dir_{\lambda})$ is not invertible:
$$\spec_I(I_b(\Dir_{\lambda})(z))=
 \big\{ \pm \sqrt{ \lambda^2-c_j^2 } | \quad c_j\in \spec(\Dir_{Y})
\big\}.$$ Formula (\ref{indinverse},\ref{indinverseb}) shows that
these indicial roots are all simple (i.e. the corresponding poles
are), as long as $\lambda$ is not an eigenvalue of $\Dir_Y$. Also,
for $\lambda=0$ the indicial root in $0$ (which exists if $0$ is
an eigenvalue of $\Dir_Y$) is always simple. If $\lambda$ is a
nonzero eigenvalue of $\Dir_Y$, then $0$ is an indicial root of
order two.
\par
\bigskip \noindent{\bf The Model Operator at {\boldmath
$\faphibf$}}
\par
\medskip
Having analyzed the behavior of the model operator on the fibres
of the face $\fabf$, let us now briefly describe the situation at
the face $\faphibf=\fabf-F_{\phi}$. As usual, this face is a fibre
bundle $$ \faphibf \cong \partial X\times\partial
X\times[-1,1]_{\sigma}-F_{\phi} \stackrel{\beta_{\phi, R}}
{\longrightarrow}
\partial X,\qquad \sigma=\frac{x-x'}{x+x'}.$$
In this compactified picture, the fibre over a point $p\in\partial
X$ is the manifold $\partial X\times [-1,1]_{\sigma}
-(Z_{\phi(p)}\times\{0\})$ with boundaries $$\farf\cap
\faphibf|_p\cong\partial X \times \{1\},\qquad
\falf\cap\faphibf|_p \cong\partial X \times \{-1\},\qquad
 \faff\cap\faphibf|_p\cong S\phindx|_{Z_{\phi(p)}},$$
 and the model operator $N_{\faphibf}(\Dir_{\lambda})$ on the
 fibres is just the restriction
 of $N_{\fabf}(\Dir_{\lambda})$ to this space.
More concretely
\begin{lemma}\label{indZ}
Fix a $p\in \partial X$  and set $$ B_{p}=[Y\times [0,\infty[_s;\,
\{\phi(p) \}\times\{s=1\}],\qquad s=x/x'.$$ By restriction,
 the function
$\rho_{\faff}$ gives a defining function for the face $\faff$ in
$B_p$ obtained by the blow up of $ \{\phi(p) \}\times\{s=1\}$ in
$B_p$. Then
\begin{enumerate}
\item  $N_{\faphibf}(\Dir_{\lambda})$ is the operator
$I_b(\Dir)$ restricted  to  ${\cal K}$ over $B_p$
\item $\rho_{\faff}I_b(\Dir)$ is a $b$-elliptic operator on $B_p$.
\item The extended set of indicial roots at $\faff$ of the operator
$\rho_{\faff}\,I_{b}(\Dir_{\lambda})$ is given by
$$\spec_I(I_{\faff}(\rho_{\faff}\,I_{b}(\Dir_{\lambda}))(z))
=i\NN\cup -i((n-1)-\NN)\subset \ZZ.$$
\end{enumerate}
\end{lemma}
\begin{proof}
This Lemma is a special case of a similar statement for elliptic
operators on a compact manifold minus a point. An instance of this was
presented in Lemma \ref{basicmap} and we use this result for the proof
of (c). Write $x=\rho_{\faff}\cdot\rho_{\faphibf}$ around $\faff$
and consider the little calculation
$$0=\rho_{\faff}\,\cl_{\phi}(\frac{\dd x}{x^2})|_{\faff}=
\cl_{\phi}(\frac{\dd \rho_{\faphibf}}{\rho^2_{\faphibf}})+
\frac{\rho_{\faff}}{\rho_{\faphibf}}\cl_{\phi}(\frac{\dd
\rho_{\faff}} {\rho^2_{\faff}})|_{\faff},\quad\mbox{i.e.}\quad
\frac{1}{x}\cl_{\phi}(\dd
\rho_{\faff})|_{\faff}=-\cl_{\phi}(\frac{\dd
\rho_{\faphibf}}{\rho^2_{\faphibf}})|_{\faff}.$$ This allows us to
calculate the indicial family at $\faff\cap \faphibf$ as
\begin{eqnarray*}
\lefteqn{I_{\faff}(\rho_{\faff}\,I_{b}(\Dir_{\lambda}))(-i z)=
\rho_{\faff}^{-z+1}I_{b}(\Dir_{\lambda})\rho_{\faff}^z
|_{(\faff\cap\faphibf)_{\phi(p)}}
=\rho_{\faff}\,I_{b}(\Dir_{\lambda})|_{(\faff\cap\faphibf)_{\phi(p)}}+
\frac{z}{x}\cl_{\phi}(\dd
\rho_{\faff})|_{(\faff\cap\faphibf)_{\phi(p)}}  }\\ & & \qquad =
\rho_{\faff}\,I_{\faphibf}(\Dir_{\lambda})|_{(\faphibf\cap
\faff)_{\phi(p)}} -z\cdot\cl_{\phi} (\frac{\dd
\rho_{\faphibf}}{\rho_{\faphibf}^2})|_{(\faff\cap\faphibf)_{\phi(p)}}
=\rho_{\faphibf}^{z}\frac{1}{\rho_{\faphibf}}\Dir^{\phi}
\rho_{\faphibf} ^{-z}|_{(\faff\cap\faphibf)_{\phi(p)}},
\end{eqnarray*}
which is just the indicial family whose roots were calculated to
be $i\NN\cup -i((n-1)-\NN)$ in Lemma \ref{basicmap}.\qed
\end{proof}

\subsection{Inversion at $\phi$bf}\label{inversionatphibf}
In Corollary \ref{finalff} we have constructed the resolvent up to
an error term $R_1$, which does not exhibit decay at the face
$\faphibf$. In this Section we want to improve this error at the
face $\faphibf$. To do this, we first remove the Taylor series of
$R_1|_{\faphibf}$ at the face $\faff$ obtaining a remainder lying
in the ``zero modes'' over the blown down face $\fabf$. This term
can then be removed using the inversion of $I_b(\Dir-\lambda)$
described in the last Section.\par
\begin{proposition}[Reduction to {\boldmath $\fabf$}]
\quad\par\smallskip \noindent The equality
$(\Dir-\lambda)Q_2^\phi(\lambda)=1+R_2^\phi(\lambda)+
R_2^{b}(\lambda)$ holds for holomorphic families
\begin{eqnarray*}
\lambda&\mapsto& R^{\phi}_2(\lambda)\in \Psi_{\phi}^{-\infty,{\cal
F}(\lambda)} (X,E),\quad
F(\lambda)_{\faff}=(1+\NN)\overline{\cup}\widehat{(1+\NN)},\quad
   F(\lambda)_{\faphibf}=h+2+\NN
\\ \lambda&\mapsto& R^b_2(\lambda)\in \Psi_{b}^{-\infty,{\cal
G}(\lambda)}(X,E),\quad G(\lambda)_{\fabf}=[\circ]\NN
\\ \lambda&\mapsto& Q_2^\phi(\lambda)\in \Psi_{\phi}^{-1,{\cal
E}(\lambda)}(X,E),\quad E(\lambda)_{\faphibf}=[\circ](h+1+\NN),
\quad E(\lambda)_{\faff}=(1+\NN)
   \overline{\cup}\widehat{(2+\NN)},
\end{eqnarray*}
all other index sets equal $\infty$.
\end{proposition}
\begin{proof}
In this proof, we will uniformly assume that the coefficient
bundle is $\CB_{b}$ and drop it from our notation. Then
 $$s(\lambda)=R_1(\lambda)|_{\faphibf} \in{\cal A}^{\cal
H}(\faphibf), \quad H_{\faff}=-h+\NN,\quad
H_{\falf}=H_{\farf}=\infty.$$ In order to reduce
$s^{\circ}(\lambda)$ to an error term on the face $\fabf$, we have
to solve away its expansion  at $\faff$. Since by Lemma \ref{indZ}
the indicial roots of $\rho_{\faff}I_{\faphibf}(\Dir-\lambda)$ are
just $i\ZZ$, the theory of elliptic $b$-differential operators
(see Lemma 5.44 in \cite{Melaps}, or Lemma \ref{holo2} below for a
similar statement) tells us that the equation
$$I_{\faphibf}(\Dir-\lambda)g^{\circ}(\lambda)-s^{\circ}(\lambda)
=u^{\circ}(\lambda)$$ for the zero-mode part can be solved with
$$g^{\circ}(\lambda)\in\Pi_{\circ}{\cal A}^{\cal I}(\faphibf),
\quad I_{\faff}=(-h+1+\NN)\overline{\cup}
\widehat{(-h+1+\NN)},\quad I_{\falf}= I_{\farf}=\infty,$$
$$\mbox{and}\quad u^{\circ}(\lambda)\in\Pi_{\circ}
\cdotinfty{}(\faphibf) \subset \Pi_{\circ}\cdotinfty{}(\fabf).$$
Extend $g^{\circ}$ and $u^{\circ}$ to kernels $G^{\circ}(\lambda)$
and $U^\circ(\lambda)$ on $X^2_{\phi}$. All this can be done
holomorphically in $\lambda$. Noting that $$G^{\circ}(\lambda)\in
{\cal A}^{\cal I}(X^2_{\phi}), \quad
I_{\faff}=(1-h+\NN)\overline{\cup}\widehat{(1-h+\NN)},\quad
  I_{\faphibf}=[\circ]\NN, \quad I_{\farf}=I_{\falf}=\infty, $$
we get, using Lemma \ref{mapping},  $$
T(\lambda):=(\Dir-\lambda)G^{\circ}(\lambda)-S(\lambda)-U^{\circ}(\lambda)
\,\in{\cal A}^{\cal J}(X),\quad
J_{\faff}=(1+\NN)\overline{\cup}\widehat{(1+\NN)}, J_{\faphibf}=
[\perp](h+1+\NN).$$ \par It remains to solve away $t(\lambda)$,
the restriction of $T(\lambda)$ to $\faphibf$. To do this, just
define $G^{\perp}(\lambda)$ to be $x'$ times the extension of
$(\Dir^{\phi,V})^{-1}t(\lambda)$ to $X^2_\phi$. Then, setting $$
Q_2^{\phi}(\lambda)
:=Q_1(\lambda)+G^{\circ}(\lambda)+G^{\perp}(\lambda),\quad
R_2^b(\lambda)=U^{\circ}(\lambda), \quad
R_2^{\phi}(\lambda)=(\Dir-\lambda) G^{\perp}(\lambda)-T(\lambda)$$
proves the claim.\qed
\end{proof}
\par
\bigskip
\noindent{\bf Meromorphic Solution at bf}
\par
\medskip
 For the formulation of the next Proposition, recall that the
indicial roots of the operator $I_b(\Dir-\lambda)$ were given by
$\pm\sqrt{\lambda^2-c_j^2}$\, , where the $c_j$ denote the
eigenvalues of $\Dir_Y$. We have agreed implicitly  to read this
as $\pm \lambda$ if $c_j=0$. For $\Im(\lambda)>0, a\in \RR-\{0\}$
we have $$
\sqrt{\lambda^2-a^2}=\sqrt{\lambda-a}\sqrt{\lambda+a}.$$ This can
be continued to a meromorphic function on the Riemann surface
$\Sigma_a\rightarrow \CC$, defined as the minimal branched
covering of $\CC$ with that property, with poles only at the
branching points $\pm a$.
\par
Denote by $\Sigma\rightarrow \CC$ the (infinite) branched covering
over $\CC$, onto which all the functions $\sqrt{\lambda^2-c_j^2}$
can be continued to be meromorphic, with poles only at the
branching points. This means that $\Sigma$ has branching points
over $\pm c_j\neq 0$ around which it is uniformized by the
functions $\sqrt{\lambda \mp c_j}$. We emphasize that $0$ is not a
branching point of this covering. We will denote points over
$\lambda\in \CC$ by $\Lambda\in\Sigma$ and write
$\theta_j(\Lambda)$ for the function on $\Sigma$, obtained by
continuation of $\theta_j(\lambda)=-i\sqrt{\lambda^2-c_j^2}$ from
$\Im(\lambda)>0$. Note that on $\Im(\lambda)>0$ the function
$\theta_j(\lambda)$ has positive real part.
\par
More concretely, introduce the {\em physical domain} ${\rm PD}$ as
that part of the preimage of the upper half plane $\{\lambda|
\Im(\lambda)>0\}$ in $\Sigma$ in which
$\theta_j(\Lambda)=\theta_j(\lambda)$ or, equivalently, where all
$\theta_j(\Lambda)$ have positive real part. Define the set
$\Theta(\Lambda)$ as the smallest $C^{\infty}$-index set
containing $(\theta_j(\Lambda)|j\in\NN)$. This is a
generating set for $\Theta(\Lambda)$  and a minimal such {\em away} from accidental
multiplicities, i.e. in the sense of Appendix
\ref{appzeronotation}: $${\rm
lead}(\Theta(\Lambda))=(\theta_j(\Lambda)|j\in\NN),
\qquad\mbox{if}\quad \theta_j(\Lambda)\notin
\theta_k(\Lambda)+\NN\quad\mbox{for}\quad j\neq k. $$ Let us note
some important properties of these index sets as a
\begin{lemma}\label{realparts}
\begin{enumerate}
\item  $0<\Re(\theta_j(\Lambda))<\Re(\theta_k(\Lambda))\quad$ for $\Lambda\in {\rm PD}, j<k$
\item $\theta_j(-r)=\overline{\theta_j(r)} \quad
(\quad=-\theta_j(r) \quad \mbox{if} \quad |r|\geq c_j )$ for $r\in
\RR=\partial{\rm PD}$
\item The set of accidental multiplicities ${\cal M}\subset\Sigma$ is
discrete and ${\cal M}\cap {\rm PD}\subset i]0,\infty[$.
\end{enumerate}
\end{lemma}
\begin{proof}
An accidental multiplicity in ${\rm PD}$ is a point $\lambda$ with
$\Im(\lambda)>0$ such that $$
-i\sqrt{\lambda^2-a^2}=-i\sqrt{\lambda^2-b^2}+m,$$ with $a, b$
positive eigenvalues of $\Dir_Y$ and $m\in \ZZ$. This can be
rewritten as $$(2m\lambda)^2=[(b+m)^2-a^2][a^2-(b-m)^2],$$ from
which the assertion in (c) immediately follows.\qed
\end{proof}
We can now formulate the next step in the construction of the
resolvent
\begin{proposition}\label{holo1}
 The equation
 $(\Dir-\lambda)Q_3^b(\Lambda)=-R_2^b(\Lambda)+R_3^b(\Lambda)$  can
 be solved with
\begin{eqnarray*}
 Q_3^b(\Lambda) &\in& \Psi_{b}^{-\infty,{\cal E}}(X,E), \quad
  E(\Lambda)_{\falf}=[\circ]\Theta(\Lambda),\quad
   E(\Lambda)_{\farf}=\Theta(\Lambda),\quad
   E(\Lambda)_{\fabf}=[\circ]\NN,
\\ R_3^b(\Lambda) &\in&  \Psi_{b}^{-\infty,{\cal F}}(X,E),\quad
F(\Lambda)_{\falf}=1+\Theta(\Lambda),\,
   F(\Lambda)_{\farf}=\Theta(\Lambda),\,
F(\Lambda)_{\fabf}=1+\NN.
\end{eqnarray*}
The maps $$\Sigma\supset B_{R/\sqrt{2}}(0)\ni\Lambda\mapsto
Q_3^b(\Lambda)\in \Psi_{b}^{-\infty,(E_{\fabf},  -R,-R) }(X,E),$$
$$\Sigma\supset B_{R/\sqrt{2}}(0)\ni\Lambda\mapsto R_3^b(\Lambda)
\in \Psi_{b}^{-\infty,(F_{\fabf},-R,-R)}(X,E)$$ into the calculus
 with bounds are meromorphic with poles only at the branching points of $\Sigma$.
\end{proposition}
\begin{proof}
We remind the reader that we will be working with the coefficient
bundle $\CB_b$. Assume for the moment that $\Im(\lambda)>0$, i.e.
$I_b(\Dir-\lambda)$ does not have any real indicial roots. Also
assume for simplicity that $\lambda$ is not an accidental
multiplicity i.e. $\theta_j(\lambda)\notin \theta_k(\lambda)$ for
$j\neq k$.
 Since $$v^{\circ}(\lambda):=-R_2^b(\lambda)|_{\fabf}\in
\Pi_{\circ}\cdotinfty{}(\fabf),$$ we can solve
$I_{\fabf}(\Dir-\lambda)g^{\circ}(\lambda)=-v^{\circ}(\lambda)$
with $$g^{\circ}(\lambda)\in\Pi_{\circ}{\cal A}^{\cal
G}(\fabf),\quad G_{\falf}= \Theta(\lambda), \quad
G_{\farf}=\Theta(\lambda).$$ We can extend $g^{\circ}(\lambda)$ to
$X^2_{b}$ in such a way that all the leading coefficients of ${\rm
lead}(\Theta(\lambda))$ at $\falf$ are in the zero modes and
independent of $x'$ close to $\fabf$.
 This gives an
element $$G^{\circ}(\lambda)\in{\cal A}^{\cal H}(X^2_{b}),\quad
H_{\falf}=[\circ]\Theta(\lambda),\quad
H_{\farf}=\Theta(\lambda),\quad H_{\fabf}=\NN.$$
    This operator can be used to solve away the
    error term
$R_2^b(\lambda)$. By Lemma \ref{mapping} we have $$
-T(\lambda):=(\Dir-\lambda)G^{\circ}(\lambda)-R_2^b(\lambda)
   \in{\cal A}^{\cal I}(X^2_{b}),\quad I_{\fabf}=[\perp]\NN,\,
   I_{\falf}=[\perp]\Theta(\lambda),
\, I_{\farf}=\Theta(\lambda),$$ where the form of $I_{\falf}$ is
due to our special choice of extension.  We are left with the
problem of solving away an error perpendicular to the
zero-modes.\par
 This is done as usual by extending
$(\Dir^{\phi,V})^{-1}T(\lambda)|_{\fabf}$, to $X^2_b$ such that
the coefficients  associated to $[\perp]{\rm
lead}(\Theta(\lambda))$at $\falf$ are just $(\Dir^{\phi,V})^{-1}$
applied to the corresponding coefficients of $T(\lambda)$ at
$\falf$. Define $G^{\perp}(\lambda)$ to be $x$ times this
extension, thus $$G^{\perp}(\lambda)\in {\cal A}^{\cal
J}(X^2_{b}),\qquad
 J_{\falf}=[\perp]\Theta(\lambda)+1,
\quad J_{\farf}=\Theta(\lambda),\quad J_{\fabf}=[\perp](1+\NN).$$
With Lemma \ref{mapping} it is now easy to see that, setting
$Q_3^b(\lambda)=G^{\circ}(\lambda)+ G^{\perp}(\lambda)$,  we have
$$R_3^b(\lambda)=(\Dir-\lambda)Q_2^b(\lambda)-R_2^b(\lambda)
  =(\Dir-\lambda)G^{\perp}(\lambda)-T(\lambda),$$
with index sets as indicated in the Proposition.
\par
\medskip
To prove the claimed meromorphy,
 a more careful analysis of the term $g^{\circ}(\lambda)$ above is needed.
Using $s=x/x'$, recall that since
$v^{\circ}\in\Pi_{\circ}\cdotinfty{}(\fabf)$, the Mellin transform
\begin{equation}\label{meldef}
M(v^{\circ})(z)=\int_0^{\infty}s^{iz}v^{\circ}(s)\frac{\dd s}{s}
\end{equation}
is holomorphic and (still assuming $\Im(\lambda)>0$) we have
\begin{equation}\label{ghoch0}
-g^{\circ}(\lambda)(s)= I_{\fabf}(\Dir-\lambda)^{-1}
v^{\circ}(\lambda)(s)= \frac{1}{2\pi}\int_{\Im(z)=0}s^{-iz}
I_{\fabf}(\Dir-\lambda)(z)^{-1} M(v^{\circ})(z) \dd
z
\end{equation}
The inverse of the indicial family has the form (\ref{indinverse},
\ref{indinverseb}). The integral thus is well defined, as long as
$\lambda$ stays away from the real axis. When $|\lambda| <
R/\sqrt{2}$ approaches the real axis from $\Im(\lambda)>0$, only
the eigenvalues smaller than $ R / \sqrt{2} $ contribute  a
singularity. Denote by $ c_1, \ldots,c_N $ the eigenvalues of
$\Dir_Y$ which are smaller than a given $L>2R$. Writing $P(c_j)$
and $P(\geq L)$  for the projections onto the corresponding
eigenspaces of $\Dir_Y$, we get the decomposition
\begin{eqnarray}
I_b(\Dir-\lambda)(z)^{-1}&=& P(\geq L)I_b(\Dir-\lambda)(z)^{-1}
+P(0)\frac{1}{z^2-\lambda^2}(iz\mbox{\boldmath
$\sigma$}+\lambda)\nonumber\\ & &\qquad +\sum_{j=1}^N P(\pm c_j)
\frac{1}{\lambda^2-c_j^2-z^2}
\zwzw{-c_j-\lambda}{iz}{-iz}{c_j-\lambda}.\label{invdec}
\end{eqnarray}
The first term in this decomposition is uniformly bounded in $z$
and has no poles in the strip $-L/2\leq \Im(z)\leq L/2$. Since the
term $M(v^{\circ})(z)$ is rapidly decreasing at real infinity,
performing the integral in (\ref{ghoch0}) for this part yields a
section in ${\cal A}^{(L/2,L/2)}(\fabf)$, which depends
holomorphically on $|\lambda|<R/\sqrt{2}$.\par Plugged into
(\ref{ghoch0}), the summands of the third term yield integrals of
the form $$\frac{1}{2\pi}\int_{\Im(z)=0}s^{-iz} P(\pm c_j)
  \frac{1}{\lambda^2-c_j^2-z^2} \zwzw{-c_j-\lambda}{iz}{-iz}{c_j-\lambda}
M(v^{\circ})(z) \dd z.$$ For small $s$ this can be evaluated by
shifting contours to $\Im(z)\gg 0$. From the residue theorem
(using the same method for the term stemming from the
0-eigenvalue) we get an expansion for $g^{\circ}(\lambda)(s)$ of
the form
\begin{equation}\label{morethanneeded}
 A_0(\lambda)s^{-i\lambda}+\sum_{j=1}^N A_j(\Lambda)
\frac{s^{\theta_j(\Lambda)}}{\theta_j(\Lambda)} +
O_{\Lambda}(s^{L/2}).
\end{equation}
Here, the coefficients $A_j$ are holomorphic with values in the
$\pm c_j$-eigenspaces of $\Dir_Y$. The expansion for
$1/s\rightarrow 0$ can be similarly obtained by shifting contours
to $\Im(z)\gg 0$.
\par The expression (\ref{morethanneeded}) is meromorphic in the sense
of being an expansion with meomorphic coefficients. Its extension
to $X^2_b$ is still meromorphic in this strong sense, i.e. as an
element in $\Psi_{b}^{-\infty, {\cal E}(\Lambda)}(X)$, as long as
we stay away from accidental multiplicities. Away from these
points, the extensions of the leading coefficients $A_j$ of the
expansion can
 (be clearly distinguished and therefore) still be chosen to lie
 in the $\pm c_j$-eigenspaces of $\Dir_Y$.
\par
Around points of accidental multiplicity, the extension of
(\ref{morethanneeded}) is still meromorphic in the weaker sense
claimed in the Proposition.\qed
\end{proof}
\par\medskip
\noindent{\bf Removal of the Remainder at {\boldmath $\falf$}}
\par
\medskip
 The next step in the parametrix construction is to solve
away the error terms at the left face $\falf$. The approach to
this problem is based on the following Lemma:
\begin{lemma}\label{holo2}
For any family of sections $f(\Lambda)\in {\cal
A}^{1+\Theta(\Lambda)}(X,E)$ there is a family $u(\Lambda)\in
{\cal A}^{1+ \Theta(\Lambda)}(X,E)$, such that
$(\Dir-\lambda)u(\Lambda)-f(\Lambda)\in\cdotinfty{}(X,E)$. The
construction gives a meromorphic map $$\Sigma\supset
B_{R/\sqrt{2}}(0)\ni \Lambda \mapsto u(\Lambda)\in {\cal
A}^{-R}(X,E), $$ whenever $f(\Lambda)$ is meromorphic in that
sense.
\end{lemma}
\begin{proof}
This Lemma is a variation of Lemma 5.44  in \cite{Melaps}, with an
easier proof. Again, we first assume that $\Im(\lambda)>0$. We
have to show how to solve for a section $f$ of the form
$f=x^{\theta_j(\lambda)+m}g(\lambda)$, with $g(\lambda)$
meromorphic in the $C^{\infty}$-sections over $X$ and $m\in
\NN_+$. As usual, decompose $g=g^{\circ}+g^{\perp}$ with
$g^{\circ}\in C^{\infty}(X,E)^{\circ}$ and $g^{\perp}\in
C^{\infty}(X,E)^{\perp}$. Starting with $g^{\circ}$, we can just
set
\begin{equation}\label{ansatz}
u_0^{\circ}(\lambda)=x^{\theta_j(\lambda)+m}
I_b(\Dir-\lambda)(i(\theta_j(\lambda)+m))^{-1} g^{\circ}(\lambda),
\end{equation}
which, according to  (\ref{indinverse}, \ref{indinverseb}),
continues to be meromorphic with extra poles arising only at
points of accidental multiplicity.
\par
Ansatz (\ref{ansatz}) solves our problem to first order in the
zero modes:
$$(\Dir-\lambda)u_0^{\circ}(\lambda)-x^{\theta_j(\lambda)+m}g(\lambda)
=x^{\theta_j(\lambda)+m}v_0(\lambda)^{\perp}.$$
 Again,
solving for $x^{\theta_j(\lambda)+m}v_0(\lambda)^{\perp}$
 is easy. Just
set $$u_0^{\perp}=x^{\theta_j(\lambda)+m+1}(D^{\phi,V})^{-1}
v_0^{\perp}(\lambda)$$ This solves the problem to first order.
Iteration of this procedure finishes the proof. \qed
\end{proof}
\par
The obvious parametrized modification of this Lemma can now be
used to get rid of the expansion of the error $R_3^b(\lambda)$ at
the left boundary of $X^2_b$.
\begin{proposition}\label{finalfaces}
The equation
$D_{\lambda}Q_4^b(\Lambda)=-R_3^b(\Lambda)+R^b_4(\Lambda)$ can be
solved with
\begin{eqnarray*}
 R_4^b(\Lambda)&\in&  \Psi_{b}^{-\infty,{\cal F}}(X,E),\quad F(\Lambda)_{\farf}
 =\Theta(\Lambda),\quad
F(\Lambda)_{\falf}=\infty,\quad F(\Lambda)_{\fabf}=1+\NN,
 \\ Q_4^b(\Lambda) &\in& \Psi_{b}^{-\infty,{\cal E}}(X,E),\quad
 E(\Lambda)_{\falf}=1+\Theta(\Lambda),\,
   E(\Lambda)_{\farf}=\infty,\,
   E(\Lambda)_{\fabf}=1+\NN.
\end{eqnarray*}
The maps $$\Sigma\supset B_{R/\sqrt{2}}(0)\ni\Lambda\mapsto
Q_4^b(\Lambda)\in \Psi_{b}^{-\infty,(E_{\fabf}, -R,-R) }(X,E),$$
$$\Sigma\supset B_{R/\sqrt{2}}(0)\ni\Lambda\mapsto R^b_4(\Lambda)
\in \Psi_{b}^{-\infty,(F_{\fabf},\infty,-R)}(X,E)$$ are
meromorphic.\qed
\end{proposition}
\par\medskip
\begin{remark}\rm 
Note that, eventhough $\lambda=0$ might be a point of accidental
multiplicity, the coefficient of the highest power $s^{-i\lambda}$
in the expansions of $Q_3^b(\Lambda)$  and $Q^b_4(\Lambda)$, will
be meromorphic near $0$. We will use this fact in Section
\ref{secconsequences}, but we abstain from introducing even more
specialized notation for this.
\end{remark}

\subsection{End of the Construction}
So far, we have found a parametrix of the type $$
(\Dir-\lambda)Q(\Lambda)=I+R(\Lambda),\qquad\mbox{where}$$
$$Q(\Lambda):=Q_2^\phi(\lambda)+Q_3^b(\Lambda)+Q_4^b(\Lambda),\quad
R(\Lambda):=R_2^\phi(\lambda)+R_4^b(\Lambda)$$ have the properties
described in the above series of Propositions -- in a loose sense,
we have obtained a remainder term $R$, which vanishes to first
order at the faces $\faff$, $\faphibf$ and to infinite order at
$\falf$. From the  composition formula in Theorem
\ref{composition} its powers are of the form
$$R^j\in\Psi_{\phi}^{-\infty,{\cal E}_j}(X,E),\quad
E_{j,\faff}\geq j-\varepsilon,\quad E_{j,\faphibf}\geq
h+1+j-\varepsilon,\quad E_{j,\falf}=\infty.$$ The $\varepsilon$ is
used to take care of the ``log's'' which appear at the face
$\faff$. The expansion of $R^j$ at the right face does not improve
with increasing $j$, but becomes eventually constant in any fixed
compact range of powers. Thus, the Neumann series for $(I+R)^{-1}$
makes sense as the asymptotic sum $$
(1+R)^{-1}=\sum_{k=0}^\infty(-R)^k=1+\widetilde{R} $$ at the faces
$\falf$, $\faphibf$, and $\faff$, and setting $$Q_{\rm
sm}(\Lambda)= Q(\Lambda)(1+R)^{-1}=
Q(\Lambda)+Q(\Lambda)\widetilde{R}(\Lambda)=:Q(\Lambda)+\widetilde{Q}(\Lambda)
$$ we get
 $$\Dir_\lambda Q_{\rm sm}(\Lambda)=1+R_{\rm sm}(\Lambda),\qquad R_{\rm
sm}(\Lambda)="R^{\infty}(\Lambda)"$$
\par
To easily describe the index sets for these operators, we allow
for a little sloppiness in the notation: Given an index set ${\cal
F}$, denote by \begin{equation} \label{sloppy} T{\cal F} :={\cal
F}|_{ [\inf({\cal F}),\inf({\cal F})+1[}
\end{equation} the finite
index set which coincides with ${\cal F}$ to first order. Also, we
will write down the index set at $\faff$ w.r.t. the
$\phi$-density, and those at the faces $\faphibf$, $\falf$ and
$\farf$ w.r.t. the $b$-density, in order to avoid the notorious
term ``$h+1$''. Then
\begin{eqnarray}
 \widetilde{R}\in \Psi_{\phi}^{-\infty,{\cal
G}}(X,E),& &\!\!\!\! G_{\falf,b}=\infty,\,\,
TG_{\faphibf,b}=\{(1,0)\},\,\,
 TG_{\farf,b}=T\Theta(\Lambda),\nonumber
 \\ & &\!\!\!\! TG_{\faff,\phi}=\{(1,0),(1,1)\}\nonumber
  \\ R_{\rm sm}\in
 \Psi^{-\infty,{\cal H}}(X,E),& &\!\!\!\!
TH_{\farf,b}=T\Theta(\Lambda),\,\, TH_{\falf, b}=\infty\nonumber
\\ \widetilde{Q}\in
\Psi_{\phi}^{-1,{\cal I}}(X,E), & &\!\!\!\!
I_{\falf,b}=[\circ]\Theta(\Lambda), \,\,
T_{<1}I_{\faphibf,b}=[\circ]T_{<1}(\Theta(\Lambda)+\Theta(\Lambda)),\nonumber
\\ & &\!\!\!\!
TI_{\farf,b}=T\Theta(\Lambda),\,\,
T_{<1}I_{\faff,\phi}=[\circ]T_{<1}(\Theta(\Lambda)+\Theta(\Lambda))\nonumber
\\ Q_{\rm sm}\in
\Psi_{\phi}^{-1,{\cal J}}(X,E), & &\!\!\!\!
J_{\falf,b}=[\circ]\Theta(\Lambda), \,\,
T_{<1}J_{\faphibf,b}=[\circ]T_{<1}(\{0\}\cup
\Theta(\Lambda)+\Theta(\Lambda)),\nonumber
\\ & &\!\!\!\!
TJ_{\farf,b}=T\Theta(\Lambda),\,\,
T_{<1}J_{\faff,\phi}=[\circ]T_{<1}(\Theta(\Lambda)+\Theta(\Lambda))\label{calL}
\end{eqnarray}
 Also, to first order, the expansion of $R_{\rm sm}(\Lambda)$ at
the right face is of type (\ref{morethanneeded}) with $s$ replaced
by $s^{-1}$.
\par
The removal of this smoothing remainder is now standard. Since,
for instance, $$R_{\rm sm}(i):L_b^2(X,E)\longrightarrow
\cdotinfty{}(X,E),$$ the kernel of the projection $M$ onto the
null space of $1+R_{\rm sm}(i)$ lies in $\cdotinfty{}(X^2,\CB_b)$.
Setting $Q_{\rm inv}(\Lambda)=Q_{\rm sm}(\Lambda)+M$ yields
$$(\Dir-\lambda)Q_{\rm inv}(\Lambda)=1+R_{\rm
sm}(\Lambda)+(\Dir-\lambda)M =1+R_{\rm inv}(\Lambda),$$ where the
new remainder   $R_{\rm inv}$ has the same expansions as $R_{\rm
sm}$. Therefore the remainder is a compact family (compare
Corollary \ref{traceclass}) of the form: $$R_{\rm
inv}:\Sigma\subset B_{R/\sqrt{2}}(0)\longrightarrow {\cal K
}(x^{R}L^2_b(X,E)).$$
 Since $1+R_{\rm inv}(i)$ is invertible
by construction we infer from analytic Fredholm theory that the
inverse of $1+R_{\rm inv}(\Lambda)$ is a meromorphic family of the
same type:
\begin{proposition}
The inverse $(1+R_{\rm inv}(\Lambda))^{-1}$ is of the form
$1+S(\Lambda)$, where $$
S(\Lambda)\in\Psi^{-\infty,(\infty,H_{\farf})}(X,E), \quad
TH_{\farf}=T\Theta(\Lambda), $$
 is meromorphic as a map
 $$\Sigma\supset B_{R/\sqrt{2}}(0)\ni \Lambda\mapsto
 \Psi^{-\infty,(\infty,-R)}(X,E).$$
\end{proposition}
\begin{proof}
The proof closely follows the proof of the analytic Fredholm
theorem. For a $\Lambda_0$ in $B_{R/\sqrt{2}}(0)$ denote by
$M_0=M(\Lambda_0)\in \cdotinfty{}(X^2)$ the projection (of finite
rank) onto the null space of $1+R_{\rm inv}(\Lambda_0)$. Then
w.r.t. $1-M_0$ and $M_0$ the operator $1+R_{\rm inv}(\Lambda)$ can
be written as the matrix $$ 1+R_{\rm
inv}(\Lambda)=\zwzw{A(\Lambda)}{B(\Lambda)}{C(\Lambda)}{D(\Lambda)}.$$
The operator $A(\Lambda)$ is of the form
$A(\Lambda)=1+T(\Lambda)$, with $$T(\Lambda)= R_{\rm
inv}(\Lambda)-M_0R_{\rm inv}(\Lambda)-R_{\rm inv}(\Lambda)M_0 +
M_0R_{\rm inv}(\Lambda)M_0\in\Psi^{-\infty,(\infty,H_{\farf
})}(X,E).$$ Also $A(\Lambda)$ is invertible for $\Lambda$ close to
$\Lambda_0$ with an inverse of the same type. To see this, we
write $$ (1+T)^{-1}=1-T+T(1+T)^{-1}T,$$ and use the following
``semi-ideal'' property of the space
$\Psi_b^{-\infty,(\infty,a,b)}(X,E)$, which is described in
\cite{Melaps}, and which shows that the last term in the above sum
is of the same type as the operator $T$:
\begin{lemma}
For $Q_1,Q_2 \in \Psi_b^{-\infty,(\infty,a,b)}(X,E)$ and $A\in
{\cal L}(x^{-b}L^2_b(X,E))$ we have
$Q_1AQ_2\in\Psi_b^{-\infty,(\infty,a,b)}(X,E)$. \qed
\end{lemma}
The rest of the proof now follows as usual: Assuming the
invertibility of $A(\Lambda)$, the operator $1+R_{\rm
inv}(\Lambda)$ is invertible exactly if the endomorphism
$D-CA^{-1}B$ is invertible, since $$ \zwzw{A}{B}{C}{D}^{-1}=
\zwzw{A^{-1}+A^{-1}B(D-CA^{-1}B)^{-1}CA^{-1}}{-A^{-1}B(D-CA^{-1}B)^{-1}}
   {-(D-CA^{-1}B)^{-1}CA^{-1}}{(D-CA^{-1}B)^{-1}}.$$
Thus the inverse exists, and is of the claimed form, whenever
$\det(D-CA^{-1}B)\neq 0$. The determinant is not constantly 0,
since
 $1+R_{\rm inv}(i)$  is invertible by construction.\qed
\end{proof}
\par\smallskip
Putting everything together, we obtain one of the central results
in this work
\begin{theorem}[The Resolvent of {$\sf\bf D$}]\label{mainresult1}
The resolvent  of $\Dir$ continued meromorphically to $\Sigma$
from $\Im(\lambda)>0$ is given by
 $$G(\Lambda)\equiv G^-(\lambda):=(Q_{\rm sm}(\Lambda)+M)(1+R_{\rm
inv}(\Lambda))^{-1}\in \Psi_{\phi}^{-1,{\cal J}}(X,E),$$ with
index set ${\cal J}$ as in (\ref{calL}). It is meromorphic as a
map
 $$\Sigma\supset B_{R/\sqrt{2}}(0)\ni\Lambda\mapsto
G^-(\Lambda)\in \Psi_{\phi}^{-1,(-R,-2R,  -2R,-R) }(X,E),$$ \qed
\end{theorem}
\par
\begin{remark}\rm \label{variantsofresolvent}
Let us end this Section with a remark about  some variants of
this construction: First, it is possible to construct the
meromorphic continuation of $(\Dir-\lambda)^{-1}$ from
$\Im(\lambda)<0$. This will be needed in view of Stone's formula
(\ref{stone}). Equivalently, we can use the continuation
$G^+(\Lambda)$ of $(\Dir+\lambda)^{-1}$ from $\Im(\lambda)>0$. The
results in this Chapter hold {\em verbatim} for this continuation
(with $\lambda$ replaced by $-\lambda$, {\em but the terms
$\Lambda$, $\theta_j(\Lambda)$ etc. remain unchanged!}).
\par
\smallskip
Also, a generalized inverse for $\Dir-\lambda_0$ can be
constructed for any fixed $\lambda_0\in \CC$. For $\alpha \in \RR$
not equal to the real part of any $\pm \theta_j(\lambda_0)$ we
define the $C^{\infty}$-index set
$$\Theta_{\alpha}(\lambda_0):=\{ s=\pm \theta_j(\lambda_0)\,|\quad
\Re(s)> \alpha \} +\NN.$$
 Then the construction of the
parametrix as an operator on $x^{\alpha}L^2_b(X,E)$ up to a
smoothing term can be followed through as before. The integral in
(\ref{ghoch0}) is now over $\Im(z)=\alpha$ and the index set
$\Theta(\lambda_0)$ has to be replaced by
$\Theta_{\alpha}(\lambda_0)$ at the left face and by
$\Theta_{-\alpha}(\lambda_0)$ at the right face. The constructions
in this Section and Section \ref{inversionatphibf} then yield
\begin{proposition}\label{variantalpha}
$\Dir_{\lambda_0}Q_{{\rm sm}, \alpha}(\lambda_0)=
1+R_{{\rm sm}, \alpha}(\lambda_0),$ holds
with  $R_{{\rm sm},\alpha}\in
 \Psi^{-\infty,{\cal H}}(X,E)$,  $Q_{{\rm sm},\alpha}\in
\Psi_{\phi}^{-1,{\cal J}}(X,E)$, and index sets
\begin{eqnarray*}
 & &
TH_{\farf,b}=T\Theta_{-\alpha}(\lambda_0),\quad TH_{\falf, b}=\infty
\\ & &
J_{\falf,b}=[\circ]\Theta_{\alpha}(\Lambda), \quad
T_{<1}J_{\faphibf,b}=[\circ]T_{<1}(\{0\}\cup
\Theta_{\alpha}(\lambda_0)+\Theta_{-{\alpha}}(\lambda_0)),
\\ & &
TJ_{\farf,b}=T\Theta_{-\alpha}(\lambda_0),\quad
T_{<1}J_{\faff,\phi}=[\circ]T_{<1}(\Theta_{\alpha}(\lambda_0)
+\Theta_{-\alpha}(\lambda_0)).
\end{eqnarray*}
Especially, the operator $(\Dir-\lambda_0):x^{\alpha}L^2_b(X,E)\rightarrow
x^{\alpha}L^2_b(X,E)$ is Fredholm, if $\alpha\neq
\Re(\pm\theta_j(\lambda_0))$.\qed
\end{proposition}
\end{remark}

\subsection{Consequences}\label{secconsequences}
In this Section we present some  standard applications of our
construction of the resolvent. Our presentation here is an
adaptation of the developments in Chapter 6 of \cite{Melaps} to
the case of our Dirac operator $\Dir$. As a special feature, this
approach allows a description of the continuous part of the
spectral measure for $\Dir$, based solely on Stone's formula,
without using any of the machinery of Hilbert space scattering
theory.\par
The arguments in this Section do not really differ from those in
the $b$-case. They have been included here, since a detailed treatment
of the Dirac operator along these lines does not seem to exist in the
literature. A treatment of the $b$-Laplace operator can be
found in \cite{Tanya}.
\par \bigskip
\noindent{\bf Structure of the Null Spaces }
\par
\medskip
As an immediate application of the resolvent construction, and its
extension described at the end of the last Section, we can
describe the structure of sections in the weighted null spaces of
$\Dir-{\lambda_0}$.
\begin{proposition}\label{basicnull}
Recall the notation introduced in (\ref{sloppy}). For
$\lambda_0\in \CC$, $r_0\in \RR$ and $\alpha\in \RR$ such that
$\alpha\neq \Re(\pm\theta_j(\lambda_0))$ we have
\begin{enumerate}
\item ${\rm null}_{x^{\alpha} L^2_b}(\Dir-\lambda_0)\subset
   {\cal A}^{[\circ]I_{\alpha}(\lambda_0)}(X,E)$,
   \qquad $TI_{\alpha}(\lambda_0)=T\Theta_{\alpha}(\lambda_0)$,
\item ${\rm null}_{L^2_b}(\Dir-r_0)\subset
x^{\varepsilon}H^{\infty}_b(X,E)$, \quad ${\rm null}_{L^2_b}(\Dir^2-r_0)\subset
x^{\varepsilon}H^{\infty}_b(X,E)$
\item ${\rm  null}_-(\Dir-r_0):=\bigcap_{\varepsilon>0}
{\rm null}_{x^{-\varepsilon}L^2_b}(\Dir-r_0)\subset {\cal
A}^{[\circ]I(r_0)}$,\qquad $TI(r_0)=T\Theta(r_0)$.
\item ${\rm  null}_-(\Dir^2-r_0^2):=\bigcap_{\varepsilon>0}
{\rm null}_{x^{-\varepsilon}L^2_b}(\Dir^2-r_0^2)\subset {\cal
A}^{[\circ]J(r_0)}$, with $$TJ(r_0)=T[\Theta_0(-r_0)\overline{\cup}
(\Theta(r_0)\cup(\Theta(r_0)+\Theta_0(-r_0))\overline{\cup}
(\Theta(r_0)+\Theta_0(-r_0)+\Theta(-r_0))].$$
\end{enumerate}
All of these spaces are of finite dimension.
\end{proposition}
\begin{proof}
Let $\xi$ be a section in $x^{\alpha}L^2_b(X,E)$ such that
$(\Dir-\lambda_0)\xi =0$. Then, starting with the parametrix
described in Proposition \ref{variantalpha}, we get the equation $$
(\Dir-\overline{\lambda_0})Q_{{\rm
sm},-\alpha}(\overline{\lambda_0})=1+R_{{\rm sm}, -\alpha}
(\overline{\lambda_0}) $$ on $x^{-\alpha}L^2_b(X,E)$. The adjoint
of this equation can then be applied to $\xi$:
 $$0=Q_{{\rm
sm},{-\alpha}}(\overline{\lambda_0})^*(\Dir-\lambda_0)\xi=
(1+R_{{\rm sm},-\alpha}(\overline{\lambda_0})^*)\xi.$$ Since by
Lemma \ref{adjoints}  the remainder $R_{{\rm
sm},-\alpha}(\overline{\lambda_0})^*$ has an expansion with index
set of type $T\Theta_{\alpha}(\lambda_0)$ at the left face, we
find $\xi \in {\cal A}^{[\circ]I_{\alpha}(\lambda_0)}(X,E)$, with
$I_{\alpha}(\lambda_0)$ as above. Part (b) is an immediate
consequence of this, since all indeces in $\Theta_0(r_0)$ have
positive real part.
 \par
  The proof of (c) then follows from this, once one recalls that $\Theta(r_0)$ is
obtained from $\Theta(\lambda)$ by continuation from
$\Im(\lambda)>0$. Alternatively
$$\Theta(r_0)=\bigcap_{\alpha>0}
\Theta_{-\alpha}(r_0).$$
All of these null spaces
are of finite dimension, since the operator
$(\Dir-\lambda_0)$ is Fredholm on $x^{\alpha}L^2_b(X,E)$.
\par
To analyze the structure of the null space of $\Dir^2-r_0^2$, with
$r_0\in \RR$, we temporarily write $G(s):=Q_{{\rm sm},\alpha}(s)$
and $R(s):=R_{{\rm sm},\alpha}(s)$. Then $G(r_0)G(-r_0)$ is our
first guess at a parametrix:
\begin{eqnarray*}
(\Dir^2-r_0^2)G(r_0)G(-r_0)&=&(\Dir+r_0)(\Dir-r_0)G(r_0)G(-r_0)
=(\Dir+r_0)(1+R(r_0))G(-r_0)
\\ &=&1+R(-r_0)+(\Dir+r_0)R(r_0)G(-r_0)
\end{eqnarray*}
with $R(-r_0)\in\Psi^{-\infty, {\cal E}}(X,E)$,
 $(\Dir+r_0)R(r_0)G(-r_0)\in \Psi^{-\infty, {\cal F}}(X,E)$
and index sets given by
$$
 E_{\falf}=\infty,\quad TE_{\farf}=T\Theta_{-\alpha}(-r_0) \quad \mbox{and}
 \quad F_{\falf}=\infty,$$
$$ TF_{\farf}=T[\Theta_{\alpha}(-r_0)\overline{\cup}
 (\Theta_{-\alpha}(r_0)\cup(\Theta_{-\alpha}(r_0)+\Theta_{\alpha}(-r_0)))\overline{\cup}
 (\Theta_{-\alpha}(r_0)+\Theta_{\alpha}(-r_0)+\Theta_{-\alpha}(-r_0))].$$
 Taking the limit $\alpha\searrow 0$ as before proves the claim.
 \qed
\end{proof}
\par
\smallskip
Note that this Proposition only makes a statement about the {\em
form} of null sections, not about their existence. For instance,
since $\Dir$ is selfadjoint on $L^2_b$, the operator
$\Dir-\lambda_0$ does not have an $L^2_b$-null space when
$\lambda_0$ is not purely real. Also, part (d) should only be regarded as an
``a priori''-result. We will show below that the top
coefficients in ${\rm null}_-(\Dir^2-r_0^2)$ are much simpler.
\par
\smallskip
Proposition \ref{basicnull} implies that the pairing
\begin{equation}\label{epspairing}
x^{-\varepsilon}L^2_b(X,E)\times x^{\varepsilon}L^2_b(X,E)\longrightarrow \CC
\end{equation}
induces a pairing between the $L^2_b$-nullspaces and the extended
$L^2_b$-nullspaces. The relationship between the extended
nullspaces ${\rm null}_-(\Dir\pm s)$ and ${\rm null}_-(\Dir^2-r^2)$
can therefore be described as follows
\begin{lemma}\label{exactdd}
The sequence
$$0\longrightarrow {\rm null}_-(\Dir-r)\longrightarrow {\rm null}_-(\Dir^2-r^2)
\stackrel{\Dir-r}{\longrightarrow}[{\rm
null}_{L^2_b}(\Dir-r)^{\perp}\cap
{\rm null}_-(\Dir+r)]\longrightarrow 0$$
is exact.
\end{lemma}
\begin{proof}
Injectivity of the left map is clear, so is exactness in the
middle.
We only have to show that the map on the right
is (well defined and) surjective. Also, it should be clear that for $\xi\in {\rm
null}_-(\Dir^2-r^2)$,
the expression $(\Dir-r) \xi$ lies in ${\rm null}_-(\Dir +r)$ and
 pairs to $0$ with ${\rm null}_{L^2_b}(\Dir-r)$ under
(\ref{epspairing}).
Surjectivity then follows from the fact that
$${\rm null}_{L^2_b}(\Dir-r)^{\perp}={\rm null}_{x^{\varepsilon}L^2_b}(\Dir-r)^{\perp}
=\overline{{\rm im}_{x^{-\varepsilon}L^2_b}(\Dir-r)}={\rm
im}_{x^{-\varepsilon}L^2_b}(\Dir-r),
$$
where we have chosen $\varepsilon$
 such that $(\Dir-r)$ is Fredholm on
 $x^{-\varepsilon}L^2_b(X,E)$. Since $\varepsilon >0$ with this
 property can be chosen arbitrarily small, the claim follows.
 \qed
\end{proof}
\par
For a more detailed description of the null spaces  and the
generalized eigenspaces we need to analyze the generalized
eigenspaces of the model operator $I_b(\Dir)$, first. Let $P\in
\Diff_b(B,{\cal K})$ be an $\RR_+$-invariant operator (think of
$P=I_b(\Dir)$ or $P=I_b(\Dir)^2$).
The {\em formal null space} of $P$ at
$z\in \spec_I(P)$ is defined as the extended $x^{-iz}L^2_b(B,{\cal K})$-null
space of $P$, i.e. $$ F(P,z) =\{ u=\sum_{l=0}^{{\rm ord}(z)-1} a_l
x^{-iz}\log(x)^l\quad|\quad a_l\in\Gamma(Y,{\cal K}), \quad
Pu=0\}.$$ For $z\notin \spec_b(P)$ we have $F(P,z)=\{0\}$. This
can be seen, noting that any $u$ with an expansion as above
fulfills
\begin{eqnarray*}
0&=& Px^{iz}\sum_{l\leq p}a_l\log(x)^l=x^{iz}I_b(P)(z)\sum_{l\leq
p}a_l\log(x)^l +x^{iz}\sum_{l\leq p-1}b_l \log(x)^l.
\end{eqnarray*}
Since $I_b(P)(z)$ is invertible, it follows that the coefficient
of the highest $\log$-power has to be $0$. Thus, the whole
expansion must be $0$.
\par
In the case of the Dirac operator $P=I_b(\Dir\mp\lambda)$
the structure of the formal eigenspaces can
be inferred from (\ref{indinverse}) and (\ref{indinverseb}). Let
us start with a
\begin{definition}
Let $\Im(\lambda)>0$. For $\xi_j$ in the $c_j$-eigenspace of
$\Dir_Y$ set
\begin{eqnarray*}
\pi^{\pm}(\xi_j,\lambda)&:=&\sqrt{\lambda+c_j}\xi_j\mp
        i\sqrt{\lambda-c_j}\mbox{\boldmath $\sigma$}\xi_j\\
\pi^{\pm}(\xi_0,\lambda)&:=&\frac{1}{2}(\xi_0\pm i\mbox{\boldmath
$\sigma$}\xi_0)\quad\mbox{(the projection onto the $\mp
i$-eigenspace of $\mbox{\boldmath $\sigma$}$)}.
\end{eqnarray*}
\end{definition}
The sections $\pi^{\pm}$ can be continued to $\Sigma$. Recall
however that $\theta_j(-r)=\overline{\theta_j(r)}$ for
$\lambda=r\in \partial{\rm PD}=\RR$ and
\begin{equation}\label{picontinued}
\pi^{\pm}(\xi_j,-r)=\left\{
      \begin{array}{ll}
       \pm \sgn(r)\mbox{\boldmath $\sigma$}\pi^{\pm}(\xi_j,r)& |r|> c_j\neq 0
     \\ \pm \mbox{\boldmath $\sigma$}\pi^{\mp}(\xi_j,r)& |r|\leq c_j \neq 0
     \end{array}\right.
\end{equation}
Now the formal eigenspaces of $\Dir\pm\lambda$ have the following
form
\begin{proposition}\label{formalnull}
Let $\Im(\lambda)\geq 0$.
\begin{enumerate}
\item For any $c_j$ and $\lambda\neq \pm c_j$:
 \begin{eqnarray*}
 F(I_b(\Dir-\lambda),\mp
i\theta_j(\lambda) )
        &=& \left\langle\left. x^{\mp\theta_j(\lambda)}
  \pi^{\mp}(\xi_j,\lambda)\,\right|\,\Dir_Y\xi_j=c_j\xi_j\right\rangle,
\\ F(I_b(\Dir+\lambda),\mp i\theta_j(\lambda) )
        &=& \left\langle\left. x^{\mp\theta_j(\lambda)}
  \mbox{\boldmath $\sigma$}
  \pi^{\pm}(\xi_j,\lambda)\,\right|\,\Dir_Y\xi_j=c_j\xi_j\right\rangle,
\end{eqnarray*}
Especially, recalling that $\theta_0(\lambda)=-i\lambda$:
\begin{eqnarray*}
F(I_b(\Dir-\lambda),\pm\lambda\neq 0)
      &=& \langle \pi^{\pm}(\xi_0)x^{\mp
      i\lambda}\,|\,\Dir_Y\xi_0=0\rangle,
\\F(I_b(\Dir+\lambda),\pm\lambda\neq 0)
      &=& \langle \pi^{\pm}(\xi_0)x^{\pm
      i\lambda}\,|\,\Dir_Y\xi_0=0\rangle.
\end{eqnarray*}
\item For $\lambda\in\{\pm c_j\}$ we get:
\begin{eqnarray*}
F(I_b(\Dir-c_j),0)&=& \langle \xi_j x^0 \log(x)
+c_j/2\,\mbox{\boldmath $\sigma$} \xi_j x^0,\xi_jx^0\, |\, \Dir_Y\xi_j=c_j\xi_j\rangle, \quad
j\neq 0
\\F(I_b(\Dir+c_j),0)&=& \langle \mbox{\boldmath $\sigma$}\xi_j x^0 \log(x) +c_j/2\, \xi_j
x^0, \mbox{\boldmath $\sigma$}\xi_j x^0\, |\, \Dir_Y\xi_j=c_j\xi_j\rangle,\quad j\neq 0
\\F(I_b(\Dir),0)&=& \langle \xi_0 x^0\,|\,\Dir_Y\xi_0=0\rangle.
\end{eqnarray*}
\item For any $c_j^2$ and $\lambda^2\neq c_j^2$:
$$
 F(I_b(\Dir)^2-\lambda^2),\mp
i\theta_j(\lambda) )
        = \left\langle\left. x^{\mp\theta_j(\lambda)}
  \eta_j\,\right|\,\Dir_Y^2\eta_j=c_j^2\eta_j\right\rangle,$$
\item For $\lambda^2=c_j^2$
$$F(I_b(\Dir)^2-c_j^2,0)=\left\langle x^0\log(x)\eta_j^1+x^0\eta^0_j\,|\, \Dir_Y^2\eta_j^i
=c_j^2\eta^i_j\right\rangle $$
\end{enumerate}\qed
\end{proposition}
These spaces encode the ``top order behavior'' of generalized
eigensections of $\Dir$ or $\Dir^2$.
To describe this, introduce the notation
$G(P,r):=\sum_{\Im(z)=-r}F(P,z)$. By Lemma \ref{realparts} we know
that for $\lambda\in \CC$
$$G(\Dir-\lambda,0)=\bigoplus
_{|c_j|\leq |\lambda|}F(\Dir-\lambda,\pm\theta_j(\lambda)).$$
It then follows from Lemma \ref{indicialfamily} that the maps
\begin{equation}
{\rm null}_-{(\Dir-r_0)}\rightarrow G(I_b(\Dir)-r_0,0), \quad
{\rm null}_-{(\Dir^2-r_0^2)}\rightarrow G(I_b(\Dir)^2-r_0^2,0)
\end{equation} given
by restriction to the top coefficients, are well defined and have
kernel contained in
$x^{\varepsilon}H^{\infty}_b$.  The images of these maps are  denoted by
$G'(I_b(\Dir)-r_0,0)$ and $G'(I_b(\Dir)^2-r_0^2,0)$.
They can be identified with the quotient between the corresponding
extended $L^2_b$-nullspace and the true $L^2_b$-nullspace
\begin{eqnarray*}
G'(I_b(\Dir)-r_0,0)&\cong& [{\rm null}_{L^2_b}(\Dir-r_0)^{\perp}\subset{\rm null}_-(\Dir-r_0)]
\\ G'(I_b(\Dir)^2-r_0^2,0)&\cong& [{\rm null}_{L^2_b}(\Dir^2-r_0^2)^{\perp}\subset
{\rm null}_-(\Dir^2-r_0^2)],
\end{eqnarray*}
where $\perp$ again refers to the pairing (\ref{epspairing})
above.
\par
\smallskip
Using part (a) of
Proposition \ref{formalnull}, we can now abstractly ``construct'' generalized
eigensections for $\Dir$, representing elements in  $G'(I_b(\Dir)-r_0,0)$,
 as follows: Fixing a cutoff function
$\chi$ near the boundary we have seen that for $\Im(\lambda)>0$:
\begin{eqnarray*}
 (\Dir-\lambda) x^{-\theta_j(\lambda)}\pi^-(\xi_j,\lambda)\chi =
O(x^{1-\theta_j(\lambda)})
\\ (\Dir+\lambda) x^{-\theta_j(\lambda)}\mbox{\boldmath $\sigma$}\pi^+(\xi_j,\lambda)\chi =
O(x^{1-\theta_j(\lambda)})
\end{eqnarray*}
 For small $\Im(\lambda)>0$ this is
certainly in $L^2_b$ and we can apply the resolvent to it. Thus
define
\begin{eqnarray} U^-(\xi_j,\lambda)&:=&
x^{-\theta_j(\lambda)}\pi^-(\xi_j,\lambda)\chi -
(\Dir-\lambda)^{-1}(\Dir-\lambda)
x^{-\theta_j(\lambda)}\pi^-(\xi_j,\lambda)\chi\label{Uminus}
 \\ U^+(\xi_j,\lambda)&:=&
x^{-\theta_j(\lambda)}\mbox{\boldmath
$\sigma$}\pi^+(\xi_j,\lambda)\chi -
(\Dir-\lambda)^{-1}(\Dir-\lambda)
x^{-\theta_j(\lambda)}\mbox{\boldmath
$\sigma$}\pi^+(\xi_j,\lambda)\chi.\label{Uplus}
\end{eqnarray}
 These sections can be continued meromorphically over
$\Sigma$.
Before noting more of their properties we have to take another
look at the resolvent of $\Dir$.
\par\medskip
Let us give a
rough description of the real poles of the resolvent:
\begin{proposition}[{\bf Poles of the Resolvent}]\label{poles}
\begin{enumerate}
\item Let $\Lambda_0\in\Sigma$ be located above $\lambda_0$ on
   the real axis,
    away from  any $\pm c_j$ but possibly near $0$.
    Then $G^{\mp}(\Lambda)$ has a pole of order at most $1$ in
    $\Lambda_0$:
    $$G^{\mp}(\Lambda)=\frac{A^{\mp}_{\Lambda_0}}{\lambda-\lambda_0}
    +B^{\mp}_{\Lambda_0}+
    (\lambda-\lambda_0)C^{\mp}_{\Lambda_0}(\Lambda),$$
    and $ A^{\mp}_{\Lambda_0}$ is the projection onto the $L^2_b$-null space of
    $\Dir\mp\lambda_0$.
\item If $\lambda_0$ equals some $\pm c_j\neq 0$ then $G^{\mp}(\Lambda)$
       has a pole of order
      at most two. Use $\theta_j$ as a uniformizing function
      of $\Sigma$ around this branching point by writing
      $\Lambda=\Lambda_0+\theta_j(\Lambda)^2$. Then
            $$G^{\mp}(\Lambda)=\frac{A^{\mp}_{\Lambda_0}}{\theta^2}
      +\frac{B^{\mp}_{\Lambda_0}}{\theta}
      +C^{\mp}_{\Lambda_0}(\theta).$$
 Again, $ A^{\mp}_{\Lambda_0}$ is the projection onto the
$L^2_b$-null space of $\Dir\mp\lambda_0$.
\end{enumerate}
\end{proposition}
\begin{proof}
 To prove (a), note that $\lambda_0$ is not a
branching point and $G^-(\lambda_0+i\varepsilon)=(\Dir-(\lambda_0
+i\varepsilon))^{-1}$. It follows from the selfadjointness of
$\Dir$ and the spectral radius formula that
$\|G^-(\lambda_0+i\varepsilon)\|\leq \varepsilon^{-1}$. Thus
$G^-(\lambda)$ can have a pole of order at most $1$ with a Laurent
expansion as in (a). The operators $A^-_{\Lambda_0},
B^-_{\Lambda_0}:x^{\varepsilon}L^2_b(X,E)\rightarrow
x^{-\varepsilon}L^2_b(X,E)$, fulfill the equalities
$$(\Dir-\lambda_0)A^-_{\Lambda_0}=0, \quad
(\Dir-\lambda_0)B^-_{\Lambda_0}=1-A^-_{\Lambda_0},\quad
(\Dir-\lambda_0)C^-_{\Lambda_0}=-B^-_{\Lambda_0}$$ from the first
two of which the claim follows.
\par
The proof of (b) is analogous, when one uses the fact that around
a branching point $\pm c_j$ the surface $\Sigma$ is uniformized by
$\theta_j(\Lambda)$.\qed
\end{proof}
\par\smallskip
The main properties of our generalized eigensetions (\ref{Uminus}), (\ref{Uplus})
are now listed in the following
\begin{proposition}[Generalized Eigenfunctions]\label{geneig}
\hfill
\begin{enumerate}
\item $(\Dir\mp\lambda)U^{\mp}(\xi_j,\Lambda)=0$
\item $U^{\mp}(\xi_j,\Lambda)$ is regular near $\Lambda=0$
\item $U^{\mp}(\xi_j,0)\perp {\rm null}_{L^2_b}(\Dir)$.
\item $U^-(\xi_j,\Lambda)= x^{-\theta_j(\Lambda)}\pi^-(\xi_j,\Lambda)
+\sum_{c_k\leq |\lambda|}
x^{\theta_k(\Lambda)}A_{kj}(\Lambda)\pi^-(\xi_j,\Lambda)+O(x^{\varepsilon})$
\item $U^+(\xi_j,\Lambda)=x^{-\theta_j(\lambda)}\mbox{\boldmath $\sigma$}\pi^+(\xi_j,\lambda)
+\sum_{c_k\leq |\lambda|}
x^{\theta_k(\Lambda)}B_{kj}(\Lambda)\mbox{\boldmath
$\sigma$}\pi^+(\xi_j,\Lambda)+O(x^{\varepsilon})$
\par
\noindent
 close to $\partial{\rm PD}$ and where the families of
endomorphisms $$A_{kj}(\Lambda):\langle
\pi^-(\xi_j,\Lambda)\rangle\rightarrow\langle\pi^+(\xi_k,\Lambda)\rangle,\quad
B_{kj}(\Lambda):\langle \mbox{\boldmath
$\sigma$}\pi^+(\xi_j,\Lambda)\rangle\rightarrow\langle\mbox{\boldmath
$\sigma$}\pi^-(\xi_k,\Lambda)\rangle$$ are meromorphic in that
area.
\end{enumerate}
\end{proposition}
\begin{proof}
(a) should be clear, (b) follows from Proposition \ref{poles},
since
$(\Dir-\lambda)x^{-\theta_j(\lambda)}\pi^-(\xi_j,\lambda)\chi$ is
orthogonal to the $L^2_b$-null space of $\Dir-\lambda$.
\par
Also (d) and (e) follow from Proposition \ref{formalnull}, since
$U^{\pm}$ certainly belong to the formal eigenspaces to first
order and the sums on the RHS are just the most general linear
combination of formal eigensections, which are $L^2_b$ for
$\Lambda$ in the physical domain. \qed
\end{proof}
\par
\smallskip
Restricting to  $r\in\partial{\rm PD}$ we can rewrite (d) and (e)
in the above Proposition in the following way. First, $$
\Pi_r:=E_{]-|r|,|r|[}(\Dir_Y)=\Pi_r^+\oplus\Pi^-_r,
\quad\mbox{with}\quad \Pi^{\pm}_r:= \langle \pi^{\pm}(\xi_l,r)\, |
\, c_l<|r|\rangle.$$ Then we can write the endomorphisms in (d)
and (e) as: $$ A(r) := \bigoplus_{c_j,c_k<|r|}A_{kj}(r):\Pi_r^-
\longrightarrow \Pi^+_r,\qquad B(-r) :=
\bigoplus_{c_j,c_k<|r|}B_{kj}(-r): \Pi^{+}_r\longrightarrow
\Pi_r^{-},$$ and, introducing
 the definition
 $\mathbb{X}(r)(\pi^{\pm}(\xi_j,r)):=\pi^{\pm}(\xi_j,r)x^{\pm\theta_j(r)},$
 and using (\ref{picontinued}), the generalized eigensections of
 $\Dir$ in (d) and (e) can now be
 written
\begin{equation}\label{dande}
U^{-}(\xi_j,r)=\mathbb{X}(r)(1+A(r))\pi^-(\xi_j,r), \quad
  U^{+}(\xi_j,-r)=\mathbb{X}(r)(1+B(-r))\pi^+(\xi_j,r),
  \end{equation}
  modulo $O(x^{\varepsilon})$.
More properties of the ``scattering matrices'' $A$ and $B$ will be
described in Proposition \ref{abunitary} below.
\par
\medskip
Choosing an orthonormal basis of eigenvectors $\xi_j^a$ of the
$c_j$-eigenspace of $\Dir_Y$ (with the additional condition that
$\mbox{\boldmath $\sigma$}\xi_0^a=i\xi_0^a$) we get an independent
set of sections
\begin{equation}\label{genbas}
U^-(\xi_0^1,\Lambda),\ldots,U^-(\xi_0^{a_0},\Lambda),\ldots,
U^-(\xi_j^1,\Lambda),\ldots,U^-(\xi_j^{a_j},\Lambda)\quad c_j\leq
|\lambda|
\end{equation}
in ${\rm null}_-(\Dir-\lambda)/{\rm null}_{L^2_b}(\Dir-\lambda)$.
The set of generalized eigensections (\ref{genbas}) spans a
subspace of dimension $\dim(\Pi_{|\lambda|})/2$. To see that this
span is the whole space, we need a second way to obtain
information about the dimension.
\par
\medskip
What is needed here is a scalar product on the spaces of formal
eigensections of $\Dir$. To define this, let $P$ be in
$\Diff_b(B,{\cal K})$ as before. The {\em boundary pairing} for
$P$ is defined as a the sesquilinear pairing $$ b_P:F(P,z)\times
F(P^*,\overline{w})
   \longrightarrow \CC, \quad
    b_P(\xi,\eta):=i\int_B\langle P(\chi\xi), \chi \eta\rangle_E
            -\langle \chi\xi, P^*(\chi\eta)\rangle_E \dvol_b, $$
where $\chi$ is a cutoff function supported near $x=0$ and
$\langle\cdot,\cdot\rangle_E$ is antilinear in the {\em second}
variable; $b_P$ is independent of the choice of $\chi$. It is
nondegenerate iff $z=w$ and $0$ otherwise (see \cite{Melaps},
Chapter 6). The boundary pairing for the Dirac operator has the
following fundamental properties:
\begin{lemma}\label{boundarypairing}
Let $u\in  F(I_b(\Dir-\lambda),a)$, $v\in F(I_b(\Dir-\lambda),b)$
be elements of the general form $u=x^{-ia}(u_0 + u_1 \log(x))$,
$v=x^{-ib}(v_0 + v_1 \log(x))$ with extensions to $X$ denoted by
$U,V$. Then for real $\lambda$ $$ b_{(\Dir-\lambda)}(u,v)=i\int_X
[\langle \Dir U,V\rangle_E -\langle U,\Dir V\rangle_E]\dvol_{X,b}
= \RES_{z=0}\,i\int_X x^z\langle U,\mbox{\boldmath $\sigma$}
V\rangle \dvol_{X,b},$$ which equals
 $i\langle u_0,\mbox{\boldmath $\sigma$} v_0\rangle_{\partial X}$ when $a=\overline{b}$.
\end{lemma}
\begin{proof}
This is a simple consequence of ``Green's formula'', which looks
like $$[\langle\Dir U,V\rangle_E-\langle U,\Dir
V\rangle_E]\dvol_b= [\dd^{*,d}x^{-v}\langle
U,\cl_d(\cdot)V\rangle_E]\dvol_d$$ in our context.\qed
\end{proof}
\par
The boundary pairing can be used to obtain relations between the
formal and true null spaces of $\Dir-\lambda$.  One can show (as in
Chapter 6 in \cite{Melaps}):
\begin{proposition}\label{lagrange}
$G'(\Dir-\lambda,0)^{\perp,b_{\Dir}}=G'(\Dir-\overline{\lambda},0)$,
especially $G'(\Dir-r,0)$, with $r$ real,  is a Lagrangian
subspace for the symplectic form $b_{\Dir}$.\qed
\end{proposition}
\begin{corollary}
The set (\ref{genbas}) is a basis in ${\rm
null}_-(\Dir-r)/{\rm null}_{L^2_b}(\Dir-r)$.\qed
\end{corollary}
\par
Proposition \ref{lagrange} can also be used to obtain more
information about the scattering operators $A(r)$, $B(r)$ in
(\ref{dande}).
 We
start with a simple
\begin{lemma}\label{piscalarproduct}
Let $\xi_j, \eta_k$ be eigensections of $\Dir_{Y}$ with
eigenvalues $c_j$, $c_k$ and $0\neq c_j,c_k<|r|$. Then w.r.t. the
$L^2$-scalar product over $Y$
\begin{enumerate}
\item $\langle \pi^{\pm}(\xi_j,r),\pi^{\pm}(\eta_k,r)\rangle_{L^2(Y)}
=\langle\xi_j,\eta_k\rangle_{L^2(Y)}( |r+c_j|+|r-c_j|)$
\item $\langle
\pi^{\pm}(\xi_j,r),\mbox{\boldmath
$\sigma$}\pi^{\pm}(\eta_k,r)\rangle_{L^2(Y)} =\mp
2i\langle\xi_j,\eta_k\rangle_{L^2(Y)}  \sqrt{r^2-c_j^2} $\qed
\end{enumerate}
\end{lemma}
Thus introduce $$ w_0^{\pm}(r):=\mp 1,\qquad w_j^{\pm}(r):=\pm 2
\frac{ \sqrt{r^2-c_j^2}}{|r+c_j|+|r-c_j|},$$ and define the scalar
product $\llangle\, ,\,\rrangle_r$ on $\Pi_{|r|}$ such that
$\Pi_r^+\perp\Pi_r^-$ and $$\llangle \varphi^{\pm}_j,\eta^{\pm}_k
\rrangle_r:=i\langle \varphi^{\pm}_j,\mbox{\boldmath $\sigma$}
\eta^{\pm}_k
\rangle_{L^2(Y)}=w_j^{\pm}(r)\langle\varphi^{\pm}_j,\eta^{\pm}_k\rangle_{L^2(Y)}\quad
\mbox{for}\quad\varphi^{\pm}_i,\eta^{\pm}_i\in \langle
\pi^{\pm}(\xi_i,r)\rangle.$$ This is a well defined scalar product
away from $r=\pm c_j$ and
\begin{proposition}\label{abunitary}
$A(r)$, $B(-r)$ are unitary w.r.t. $\llangle \, ,\,\rrangle_r$ and
$A(r)^*=B(-r)$.
\end{proposition}
\begin{proof}
First, the above Corollary implies that the $U^-(\xi_j,r)$ and the
$U^+(\xi_j,-r)$  span the generalized eigenspace of $\Dir-r$. The
boundary values of any element in the generalized eigenspace have
a unique decomposition into their $\Pi^{\pm}_r$-parts, which can
be written in two ways according to (\ref{dande}): $$
 (1+A(r))B(-r)\pi^+(\xi_j,r)= (1+B(-r))\pi^+(\xi_j,r).$$
This shows that $A(r)$ and $B(-r)$ are inverses of each other!
\par
To prove unitarity, we use Proposition \ref{lagrange} and Lemma
\ref{boundarypairing}. For any $\eta\in \Pi^-_r$, for instance
$\eta\in \langle \pi^{-}(\xi_j,r)\rangle$ for a fixed eigenvalue
$c_j$, we get
\begin{eqnarray*}
0&=&
b_{\Dir-r}(\mathbb{X}(r)(1+A(r))\eta,\mathbb{X}(r)(1+A(r))\eta)
\\ &=&i\langle \eta,\mbox{\boldmath $\sigma$}\eta\rangle_{L^2(Y)}
+i\langle A(r)\eta,\mbox{\boldmath $\sigma$}
A(r)\eta\rangle_{L^2(Y)} = \llangle \eta,\eta\rrangle_r+\llangle
A(r)\eta, A(r)\eta\rrangle_r
\end{eqnarray*}
and so on. \qed
\end{proof}
\par
\bigskip \noindent{\bf Spectral Measure around 0}
\par\medskip
\noindent Let $f\in C(\RR)$ have support in around $0$. Then,
according to Stone's formula
\begin{eqnarray}
f(\Dir)\varphi&=&\lim_{\varepsilon\rightarrow 0}
          \frac{1}{2\pi i}\int_{\RR}f(r)\left[ (\Dir-r-i\varepsilon)^{-1}-
          (\Dir-r+i\varepsilon)^{-1}\right]\varphi\dd r\label{stone}
\\&=& \lim_{\varepsilon\rightarrow 0}
          \frac{1}{2\pi i}\int_{\RR}f(r)[ G^-(r+i\varepsilon)-
          G^+(-r+i\varepsilon)]\varphi\dd r\nonumber
\\&=&    \frac{1}{2\pi i}\int_{\RR}f(r)[ G^-(r)-
          G^+(-r)]\varphi\dd r\nonumber
\end{eqnarray}
The last integrand can have singularities. However, fixing
$\varphi\in x^{\varepsilon}L^2_b$ {\em orthogonal} to the
$L^2_b$-null space of $(\Dir-r)$, we get, according to Proposition
\ref{poles}
 $$ [ G^-(r)-
          G^+(-r)]\varphi=[B_r^--B^+_{-r}]\varphi
          =(\Dir-r)[C_r^--C^+_{-r}]\varphi$$
which shows that the LHS exists and lies in ${\rm null}_-(\Dir-r)$
and is orthogonal to the $L^2_b$-null space of $(\Dir-r)$. This
allows us to make a general {\em Ansatz} of the form
\begin{equation}\label{measureansatz}
[ G^-(r)-G^+(-r)]\varphi=\sum_{|c_j|,|c_k|<r;a,b}a_{j,k}^{a,b}(r)
          \langle \varphi, U^-(\xi_j^a,r)\rangle U^-(\xi_k^b,r),
\end{equation}
where we have used the orthonormal basis introduced before
(\ref{genbas}). This is a linear combination of finite rank linear
operators, continuous on $x^{\varepsilon}L^2_b$, and has null
space in $x^{\varepsilon}L^2_b$ of codimension at least
$\dim(\Pi_{|r|})$. Thus it will be enough to calculate
(\ref{measureansatz}) on a (special ) finite set of sections
$\varphi$ as follows:
\begin{lemma}
Let $\varphi$ be the test function of the form $\varphi=
(\Dir-r)x^{-\theta_j(r)}\pi^-(\xi_j,r)$. Then
\begin{enumerate}
\item $[ G^-(r)- G^+(-r)]\varphi=U^-(\xi_j,r)$
\item $\langle \varphi , U^-(\eta_k,r)\rangle=
\langle \pi^-(\xi_j,r),\mbox{\boldmath $\sigma$}
\pi^-(\eta_k,r)\rangle_{L^2(Y)}$, ($=2i\langle
\xi_j,\eta_i\rangle_{L^2(Y)} \sqrt{r^2-c_j^2}$ for $k,j\neq 0$)
\end{enumerate}
\end{lemma}
\begin{proof}
This proof is a good exercise to test the reader's  versality with
the concepts introduced in this Section. We give the details for
reference.
\par
First, using Lemma \ref{realparts}(b) and
(\ref{picontinued}), we know that
\begin{equation}\label{aboveandbelow}
-\theta_j(r)=\theta_j(-r),\qquad\pi^-(\xi_j,r)=\sgn(r)\mbox{\boldmath
$\sigma$} \pi^-(\xi_j,-r),\qquad\mbox{for $|r|>c_j$}.
\end{equation}
Therefore $\varphi$ can be written as the limiting function of a
function on ${\rm PD}$ in two ways: \begin{eqnarray} \varphi&=&
\lim\limits_{{\rm PD}\ni\lambda\rightarrow
r}(\Dir-\lambda)x^{-\theta_j(\lambda)}\pi^-(\xi_j,\lambda)\label{above}
\\ \varphi&=&
\lim\limits_{{\rm PD}\ni\lambda\rightarrow
-r}(\Dir+\lambda)x^{\theta_j(\lambda)}\sgn(r)\mbox{\boldmath
$\sigma$ }\pi^-(\xi_j,\lambda).\label{below}
\end{eqnarray}
Note that in (\ref{below}) the Dirac operator is applied to a
function which is already in $L^2_b$! Thus
\begin{eqnarray*}
G^+(-r)\varphi &=& \lim\limits_{{\rm PD}\ni\lambda\rightarrow
-r}(\Dir+\lambda)(\Dir+\lambda)x^{\theta_j(\lambda)}\sgn(r)\mbox{\boldmath
$\sigma$ }\pi^-(\xi_j,\lambda)
\\ &=& x^{\theta_j(-r)}\sgn(r)\mbox{\boldmath
$\sigma$ }\pi^-(\xi_j,-r)\stackrel{(\ref{aboveandbelow})}{=}
x^{-\theta_j(r)}\pi^-(\xi_j,r).
\end{eqnarray*}
Recalling the definition (\ref{Uminus}) of $U^-(\xi_j,r)$, part
(a) is proved. To prove (b), simply recall that $U^-(\eta_k,r)$ is
in the extended $L^2_b$-null space of $(\Dir-r)$. Thus, from the
definition of the boundary pairing $$ \langle \varphi ,
U^-(\eta_k,r)\rangle=-ib_{(\Dir-r)}(x^{-\theta_j(r)}\pi^-(\xi_j,r),
x^{-\theta_k(r)}\pi^-(\eta_k,r))=\langle \pi^-(\xi_j,r),
\mbox{\boldmath $\sigma$}\pi^-(\eta_k,r)\rangle_{L^2(Y)},$$
 as
claimed. The last equation then follows from Lemma
\ref{piscalarproduct}(b).
 \qed
\end{proof}
Plugging this into our Ansatz (\ref{measureansatz}) we get
\begin{theorem}\label{spectmeasure} The continuous part
of the spectral measure for $\Dir$ is of the form
\begin{eqnarray*}
\lefteqn{\dd E_{\Dir, c}(r) \varphi = \frac{1}{2\pi i}[ G^-(r)-
          G^+(-r)]\varphi \dd r}
\\ &=& \frac{1}{2\pi}\sum_{a} \langle \varphi , U^-(\xi_0^a,r)\rangle
U^-(\xi_0^a,r)\dd r - \frac{1}{4\pi}\sum_{0\neq |c_j|<r;
a}(r^2-c_j^2)^{-1/2}_+
 \langle \varphi , U^-(\xi_j^a,r)\rangle U^-(\xi_j^a,r)\dd r
\end{eqnarray*}
\qed
\end{theorem}
\par
\medskip
\noindent {\bf Regularized Trace of the Heat Kernel for Large
Times and the Index}
\par\medskip
\noindent Using the spectral measure $\dd E_{\Dir}$, the heat
kernel for $\Dir^2$ can be written $$
e^{-t\Dir^2}=\int_{\RR}e^{-tr^2}\dd E_{\Dir}(r). $$ Choose a cut
off function $\chi$, which is $1$ near $0$ and whose support lies
in $]-c_1,c_1[$. Then, using Theorem \ref{spectmeasure}, we can
calculate the regularized graded trace:
\begin{eqnarray*}
\lefteqn{\Str_{d}(e^{-t\Dir^2}\chi(\Dir))
  =\REG_{z=0}\Str(x^z e^{t\Dir^2}\chi(\Dir))}\\
&=&\REG_{z=0}\Str(x^z[P_{\Dir}(0)+
\sum_{a}\frac{1}{2\pi}\int_{\RR}e^{-tr^2}\chi(r)U^-(\xi_0^a,r)\langle\cdot,
U^-(\xi_0^a,r)\rangle_E \dd r])
\\&=&{\rm ind}_{L^2_b}(\Dir)+\REG_{z=0}
\sum_{a}\frac{1}{2\pi}\int_Xx^z\int_{\RR}e^{-tr^2}\chi(r)\str_E(U^-(\xi_0^a,r)\langle\cdot,
U^-(\xi_0^a,r)\rangle_E )\dd r\dvol_b
\end{eqnarray*}
Now, using the fact that the grading operator $\epsilon\,$
intertwines the $\pm i$-eigenspaces of $\mbox{\boldmath
$\sigma$}$, the supertrace can be written as
\begin{eqnarray*}
\str_E(U^-(\xi_0^a,r)\langle\cdot, U^-(\xi_0^a,r)\rangle_E )&=&
\langle \epsilon\,
(x^{ir}\xi_0^a+x^{-ir}A(r)\xi_0^a),x^{ir}\xi_0^a+x^{-ir}A(r)\xi_0^a
\rangle_E +O_r(x^{\delta})
\\ &=&x^{2ir}\langle \epsilon\,\xi_0^a,A(r)\xi_0^a\rangle_E+
      x^{-2ir}\langle \epsilon\, A(r)\xi_0^a, \xi_0^a\rangle_E+O_r(x^{\delta})
\end{eqnarray*}
 The term $O_r(x^{\delta})$ is
regular in $r=0$. Thus the corresponding part of the above
integral vanishes for $t\rightarrow \infty$. The remaining part
gives integrals of the type $$ \REG_{z=0}
\frac{1}{2\pi}\int_{\RR_+}x^z\int_{\RR}e^{-tr^2}\chi(r)
x^{2ir}\langle \epsilon\,\xi_0^a,A(r)\xi_0^a\rangle_{\partial X}
\chi(x)\dd r \frac{\dd x}{x}=\frac{1}{2}\langle
\epsilon\,\xi_0^a,A(0)\xi_0^a\rangle_{\partial X},$$ since this is
just the Mellin transform back and forward and evaluated in $z=0$.
Thus, recalling that the $U(\xi_0^k,0)/\sqrt{2}$ are orthonormal
in ${\rm null}(\Dir_Y)$, we finally obtain
\begin{equation}\label{extindex}
\lim_{t\rightarrow \infty}\Str_{d}(e^{-t\Dir^2})
   ={\rm ind}_{L^2_b}(\Dir)+\frac{1}{2}\sum_a\langle \epsilon\, U(\xi_0^a,0),
   U(\xi_0^a,0)\rangle_{\partial X}={\rm sdim}({\rm
   null}_-(\Dir)),
\end{equation}
the extended $L^2_b$-index of $\Dir$, which we also denote by
${\rm ind}_-(\Dir)$. To calculate this index we will have to
consider the limit $t\rightarrow 0$ of the heat supertrace. This
problem will be attacked in the next two Chapters.

\newpage
\section{The Heat Kernel at Finite Time}\label{heatchapter}
\par
In this Chapter we describe the construction of the heat kernel of
$\Dir^2$ for finite times. We introduce a heat space $X^2_H$ and
heat calculus adapted to our problem and show by direct
construction that the heat kernel lies in that calculus. As in the
case of the resolvent, this construction will be performed by
first solving the heat equation for $\Dir^2$ to first order (at
the relevant faces in $X^2_H$) and then using an iteration
argument, based on a composition formula for our calculus. The
presentation given here leans on \cite{Daimel}.
\par
\subsection{$d$-Heat Space and $d$-Heat Calculus} \label{dheatbasics}
In this Section we describe the heat space and the corresponding
calculus
 adapted to our Dirac operator.
 We define the heat space as the parabolic blow up
$$ X_H^2:= \left[ X^2_b\times [0,\infty[_t;(F_{\phi}\times
\{0\}),\dd t;
 \bigtriangleup_X\times\{0\},\dd t \right] .$$
 The invariant notion of parabolic blow up is explained in \cite{Meaomwc},
 \cite{Melaps} and \cite{Daimel}. For our purposes here, let us
 just mention that  local descriptions of the space $X^2_H$ are
 given by the sets of coordinates (\ref{localhbf}, \ref{localhbftff}, \ref{localtff},
 \ref{localtfftf}) below. See also Remark \ref{globalsquare}.\par
Now, the corresponding blow down map is denoted by
$\beta_{H}:X_H^2\rightarrow X^2\times [0,\infty[$ and the left and
right projections by $\beta_{H,L}:X_H^2\rightarrow X$ and
$\beta_{H,R}:X_H^2\rightarrow X$. The faces of $X_H^2$ are denoted
by $\farf, \falf$ and $$\begin{array}{rclcrcl} \fahbf &:=&
\overline{\beta_{H}^{-1}(\partial X ^2\times ]0,\infty[)} &\quad&
  \fatff &:=& \beta_{H}^{-1}(F_{\phi} \times \{0\})\\  \\
 \fatf &:=& \beta_{H}^{-1}( \bigtriangleup_X\times\{0\})&\quad&
    \fatb&:=&\overline{\beta_H^{-1}((X^2-\bigtriangleup_X)\times\{0\})}.
\end{array}$$
\par
\bigskip
\centerline{\epsfig{file=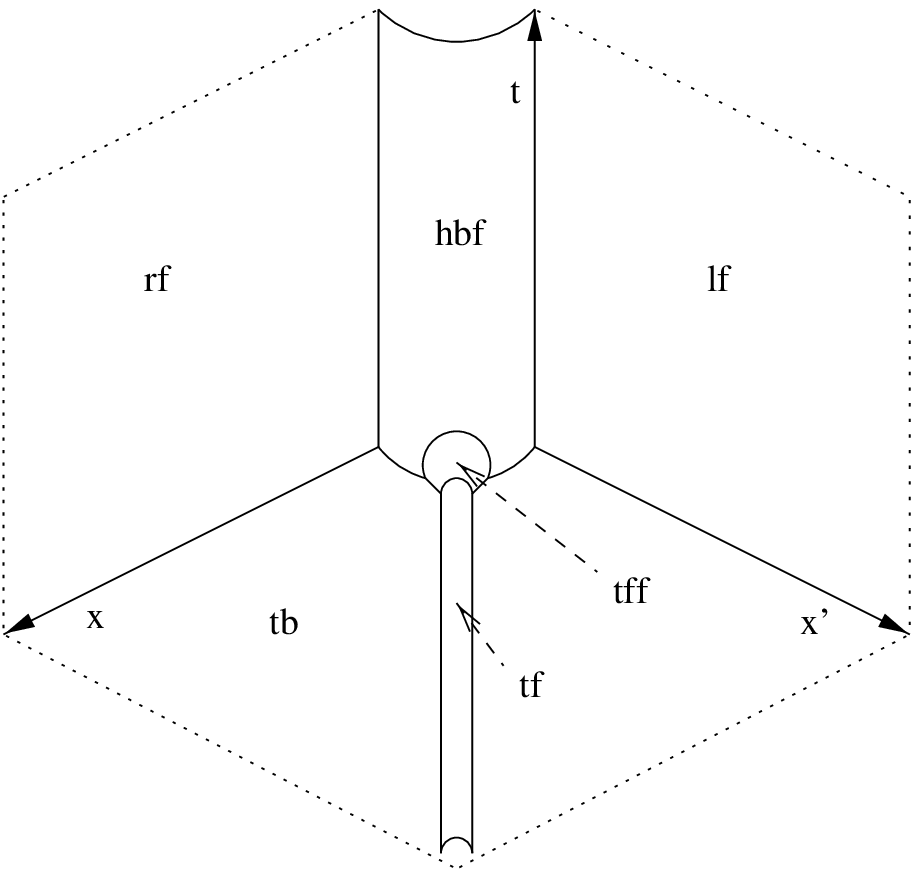}} \centerline{Figure 5:
$X_H^2$}\par\bigskip Using the coordinates $x, x',{\bf y},{\bf
y'},{\bf z},{\bf z'},t $ on $X^2\times [0,\infty[$  we can
describe the properties of lifted vector fields under the blow
down map $\beta_{H}$ and near the above faces as follows: In the
interior of {\boldmath $\fahbf$} we can use projective coordinates
w.r.t. $x'$, where $x$ is the boundary defining function in $X$,
as usual: \begin{equation}\label{localhbf} s=\frac{x}{ x'},\quad
x',\quad{\bf y},\quad{\bf y'},\quad{\bf z}, \quad {\bf z'},\quad
t.
 \end{equation}
 Then $$ \dd s =\frac{\dd x}{x'}-\frac{1}{s
x'}\dd x'\qquad \mbox{and}\quad
 \beta_H^*x \frac{\partial}{\partial x} =s \frac{\partial}{\partial s}. $$
Near {\boldmath $\fahbf$ and $\fatff$} we have to consider the
$\dd t$-parabolic blow up of $F_{\phi}\times\{0\}$. We do this by
introducing coordinates \begin{equation}\label{localhbftff}
S_t=\frac{s-1}{\sqrt{t}},\quad X_t=\frac{x'}{\sqrt{t}},\quad {\bf
Y}_t=\frac{{\bf y}-{\bf y}'}{\sqrt{t}}, \quad \sqrt{t},\quad{\bf
y}', \quad {\bf z}, \quad {\bf z}',
\end{equation} where now ${\bf
y}$, ${\bf y}'$ are assumed to be {\bf paired}. Then $$ \dd
S_t=\frac{\dd s}{\sqrt{t}}-\frac{s-1}{2t^{3/2}}\dd t,\quad
   \dd X_t=\frac{\dd x'}{\sqrt{t}}-\frac{x'}{2t^{3/2}}\dd t,\quad
   \dd {\bf Y}_t=\frac{\dd {\bf y}-\dd {\bf y}'}{\sqrt{t}}
    -\frac{{\bf y}-{\bf y}'}{2t^{3/2}}\dd t  , $$
and therefore
$$\beta_H^*t \frac{\partial}{\partial t}=\frac{1}{2}
 \left( \sqrt{t}\frac{\partial}{\partial \sqrt{t}}-
  S_t\frac{\partial}{\partial S_t}-X_t\frac{\partial}{\partial X_t}
  -{\bf Y}_t \frac{\partial}{\partial {\bf Y}_t}\right),$$
$$ \beta_H^*x \frac{\partial}{\partial x}=
  (S_t+\frac{1}{\sqrt{t}})\frac{\partial}{\partial S_t},\quad
  \beta_H^*\frac{\partial}{\partial {\bf y}}
   = \frac{1}{\sqrt{t}}\frac{\partial}{\partial {\bf Y}_t},\quad
   \frac{1}{x}\frac{\partial}{\partial {\bf z}}=
      \frac{1}{sX_t\sqrt{t}}\frac{\partial}{\partial {\bf z}}.$$
In the interior of {\boldmath $\fatff$} we can use projective
coordinates w.r.t. $x'$:
\begin{equation}\label{localtff}
S_x = \frac{s-1}{x'}, \quad {\bf Y}_x=\frac{{\bf y}-{\bf
y}'}{x'},\quad x', \quad \tau=\frac{t}{(x')^2}, \quad{\bf y}',
\quad{\bf z}, \quad{\bf z}'.
\end{equation}
We then have $$ \dd
S_x=\frac{\dd s}{x'}-\frac{s-1}{(x')^2}\dd x',\quad
    \dd {\bf Y}_x=\frac{\dd {\bf y}-\dd {\bf y}'}{x'}
   -\frac{{\bf y}-{\bf y}'}{(x')^2}\dd x',\quad
    \dd \tau=\frac{\dd t}{(x')^2}-\frac{t}{2(x')^3}\dd x',$$
and $$\beta_{H}^*t\frac{\partial}{\partial
t}=\tau\frac{\partial}{\partial \tau}, \quad \beta_{H}^*\delyy =
\frac{1}{x'}\frac{\partial}{\partial {\bf Y}_x},
\quad\beta_{H}^*\frac{1}{x}\delzz = \frac{1}{x}\delzz,\quad
\beta_{H}^*x\delx =(S_x+\frac{1}{x'})\frac{\partial}{\partial
S_x}$$ At {\boldmath $\fatf$ and near $\fatff$} (i.e. after
parabolically blowing up $\bigtriangleup_{X}$) we have projective
coordinates w.r.t. $\tau^{1/2}$
\begin{equation}\label{localtfftf}
{\sf S}=\frac{S_x}{\tau^{1/2}}=S_t,\quad
 {\bf\sf Y}=\frac{{\bf Y}_x}{\tau^{1/2}}={\bf Y}_t,\quad
{\bf\sf Z}=\frac{{\bf z}-{\bf z'}}{\tau^{1/2}}
    =\frac{x'({\bf z}-{\bf z'})}{t^{1/2}}, \quad
\sqrt{\tau}, \quad x', \quad {\bf y}', \quad {\bf
z}'.
\end{equation}
 Here, we have added the requirement that ${\bf
z}$, ${\bf z}'$ also be {\bf paired}. Now $$\dd {\sf S}=\frac{\dd
s}{\sqrt{t}}-\frac{1}{2t^{3/2}}(s-1)\dd t,\quad
   \dd {\bf\sf Y}_t=\frac{\dd {\bf y}-\dd {\bf y}'}{\sqrt{t}}
    -\frac{1}{2t^{3/2}}({\bf y}-{\bf y}')\dd t,$$
$$   \dd {\bf \sf Z}= \frac{{\bf z}-{\bf z}'}{\sqrt{t}}\dd x'
  -\frac{x'({\bf z}-{\bf z}')}{2t^{3/2}} \dd t +\frac{x'}{\sqrt{t}}
(\dd {\bf z}-\dd{\bf z}'),\quad \dd
\sqrt{\tau}=\frac{1}{2x'\sqrt{t}}\dd t+\frac{\sqrt{t}}{2(x')^3}\dd
x', $$ and the standard vector fields lift as $$
\beta_H^*x\frac{\partial}{\partial x}= ({\sf
S}+\frac{1}{\sqrt{t}})\frac{\partial}{\partial{\sf S}}, \quad
\beta_{H}^*\delyy = \frac{1}{\sqrt{t}}\frac{\partial}{\partial
{\bf\sf Y}},\quad \beta_{H}^*\frac{1}{x}\delzz =
\frac{1}{s\sqrt{t}} \frac{\partial} {\partial{\bf\sf Z}}, $$
$$\beta_{H}^*t\frac{\partial}{\partial t}=
 \frac{1}{2}\left( \tau^{1/2}\frac{\partial}{\partial \tau^{1/2}}
 -{\bf\sf Z} \frac{\partial}{\partial \bf\sf Z}
-{\sf S}\frac{\partial}{\partial \sf S} - {\bf\sf
Y}\frac{\partial} {\partial \bf\sf Y}\right).$$
\par
These calculations  show that for $V\in\Gamma(\phitx)$ the vector
field $\frac{t^{1/2}}{x}\beta_{H,L}^*V$ is tangent  to (the fibres
under $\beta_{H,R}$ of) $\fatff$ and $\fatf$. We thus get the
following analogue of Lemma \ref{philift}:
\begin{lemma}
The faces of $X^2_H$ are isomorphic to standard
spaces:
\begin{enumerate}
\item $ \stackrel{\circ}{\fatf}\, \cong  \phitx $\\
      $\fatf \cong\RC^2(\phitx)$
\item $ \stackrel{\circ}{\fatff} \,\cong  \phindx\times_Y\partial X
\times ]0, \infty[_{\tau}$
 \\  $\fatff\cong \left[ \RC^2(\phindx\times_Y\partial X) \times
[0,\infty[_{\tau}; \bigtriangleup_{\partial X}\times \{0\},\dd
\tau \right]$
\item  $\stackrel{\circ}{\fahbf}\,\cong (\partial
X)^2\times ]0,\infty[_{s} \times ]0,\infty[_{t}$
\\  $\fahbf\cong\left[ (\partial X)^2\times [-1,1]_{\sigma} \times
[0,\infty[_{t} ; F_{\phi}\times \{0\},\dd t\right]$\qed
\end{enumerate}
\end{lemma}
Here, given a euclidean vectorspace $(W,g)$, we denote by  ${\rm
RC}^2(W)$ the ``quadratic'' radial compactification of $W$,
obtained by introducing $\rho=|\cdot |^{-2}$ as a boundary
defining function at infinity. The following gives a description
of the restriction of $\phi$-vector fields to these faces:
\begin{lemma} \label{liftface}
The lifts of vector fields $V\in \Gamma(\phitx)$,
$t\frac{\partial} {\partial t}$, restricted to the different faces
are as follows:
\begin{enumerate}
\item At $\fatf$ :
$ \beta_{H,L}^*\, t\frac{\partial}{\partial t}|_{\fatf} =
-\frac{1}{2}R^{\phi}\in T\phitx/X, \qquad
\beta_{H,L}^*\,\frac{t^{1/2}}{x}V|_{\fatf} = V \in T\phitx/X.$
\item  At $\stackrel{\circ}{\fatff}$ :
$ \beta_{H,L}^*\, t\frac{\partial}{\partial t}|_{\fatff} =
\tau\frac{\partial}{\partial \tau}, \quad
\beta_{H,L}^*\,\frac{t^{1/2}}{x} V|_{\fatff}=
\tau^{1/2}N_{\phi}(V) \in T(\phindx\times_{Y}\partial X/Y).$\qed
\end{enumerate}
\end{lemma}\par\noindent
Now define the density and the coefficient bundle for the heat
kernel as $$ \KD_H:= \frac{(x')^n}{t^{n/2}}\frac{\dd t}{t}
\beta_{H,R}^*\phiomega(X),\quad
   \CB_H:= \KD_H \otimes\beta_{H,L}^*E\otimes \beta_{H,R}^*E^*=: \KD_H \otimes
   \END(E) .$$
Since we have defined the density $\KD_H$ to be lifted from the
right, its Lie derivative w.r.t. a vector field lifted from the
left is trivial: $$L_{\beta_{H,L}^*V}\KD_H =0, \qquad
L_{\beta_{H}^* t\frac{\partial}{\partial t}}
\KD_H=-\frac{n}{2}\KD_H.$$ This will again prove useful in Section
\ref{normalops}, when we calculate the action of such vector
fields on our calculus.
\par
The $d$-heat calculus is defined to be $$\Psi^{l,m,p}_{H,{\rm
cl}}(X,E):=
 \rho_{\fatf}^l\,\rho_{\fatff}^m\,\rho_{\fahbf}^{p-v}
\cdotinfty{\fatf,\fatff,\fahbf}(X_H^2,\CB_H)$$ for $l>0,m\geq
0,p\geq 0$. For $l=0$, we set $$ \Psi^{0,m,p}_{H,{\rm cl}}(X,E):=
\rho_{\fatff}^m\,\rho_{\fahbf}^{p-v}
\cdotinfty{\fatf,\fatff,\fahbf}(X_H^2,\CB_H)^{\fatf}\oplus \CC,$$
where, using (\ref{KDtf=}),
 we require for $A\in \cdotinfty{\fatf,\fatff,\fahbf}(X_H^2,\CB_H)^{\fatf}$
that the following mean value condition holds
\begin{equation}\label{meanvaluecond}
\int_{\phitx/X}\Int(t\frac{\partial}{\partial t}) A|_{\fatf}=0.
\end{equation}
\par
Let us show, how we can define the action of an element $A$ in the
calculus on some function $g \in  \cdotinfty{c}
(X\times[0,\infty[,E) $. First, for any function  $f \in
\cdotinfty{c}(X\times[0,\infty[, \Omega(X)\dd t\otimes E^*)$ we
consider the convolution product
$$f*g(t)=\int_0^{\infty}f(t+s)g(s) \in \cdotinfty{c} (X\times
X\times [0,\infty[,\beta_{H,L}^*\Omega(X)\otimes
\beta_{H,L}^*E^*\otimes \beta_{H,R}^*E  ).$$ Note that for
functions vanishing rapidly at the faces $\farf, \falf, \fatb$ we
have $$\beta_{H,L}^*\Omega(X)\otimes\KD_H \sim
\rho_{\fatff}^0\,\rho_{\fatf}^{-1}\,
\rho_{\fahbf}^{v}\Omega(X_H^2),$$ which implies that for an
operator $A\in\Psi^{l,m,p}_{H,{\rm cl}}(X,E)$ the pairing $$
(A,\beta_{H}^*(f*g))\in \rho_{\fatff}^m\,\rho_{\fatf}^{l-1}\,
\rho_{\fahbf}^{p}
\cdotinfty{\fatf,\fatff,\fahbf}(X_H^2,\Omega(X_H^2)),$$ is
integrable for $l>0$. In the case $l=0$ we have to use the mean
value condition (\ref{meanvaluecond}) to make sense of the
integral. This is done by writing, in the neighborhood
$U(\fatf):=\{\rho_{\fatf}<\varepsilon\}\cong [0,\varepsilon[\times
\phitx$: $$ \int_{U(\fatf)}(A, \beta_{H}^*(f*g))
   = \lim_{\delta\rightarrow 0} \int_{\delta}^{\varepsilon}\int_X
\int_{\phitx/X} a \dvol_{\phitx/X} \dvol_{X} \frac{\dd
\rho_{\fatf}} {\rho_{\fatf}}. $$ Since the innermost integral
converges to $0$ for $\rho_{\fatf}\rightarrow 0$, this limit
exists.\par Now we can define the action of $A$ on $g$ by setting
$$ \langle A g,f\rangle_{X\times[0,\infty[} := \langle A,
\beta_{H}^*(f*g) \rangle_{X_H^2}.$$ As usual, the task of writing
down the definitions for the extended heat calculus is left to the
reader. It is important to note, however, that we can only allow
conormal functions for which the above mean value condition makes
sense.

\subsection{Composition Formula}
To explain our constructions in this Section consider simplified
convolution operators $A, B\in
\cdotinfty{c}(X^2\times[0,\infty[,\Omega_R(X)\dd t )$. Their
composition is given by
\begin{equation}\label{convolution}
A\circ B({\bf w}, {\bf w}'', t)= \int_{{\bf w}',t'} A({\bf w},
{\bf w}',t-t')B({\bf w}', {\bf
w}'',t')\in\cdotinfty{c}(X^2\times[0,\infty[,\Omega_R(X)\dd t ).
\end{equation}
Using  the space
$X^3\times[0,\infty[_{t'}[\times[0,\infty[_{t-t'}$ and the three
projections $$\begin{array}{ccccc}
X^2\times[0,\infty[_{t-t'}&\stackrel{\pi_L}{\longleftarrow}&
X^3\times[0,\infty[_{t'}[\times[0,\infty[_{t-t'}&\stackrel{\pi_R}{\longrightarrow}&
X^2\times[0,\infty[_{t'}\\
 & & \quad\downarrow \pi_M \\
& & X^2\times[0,\infty[_{t}
\end{array}$$
this can also be written as: $$A\circ
B=\pi_{M,*}(\pi_L^*A\cdot\pi_R^*B).$$ Now, in our case, the image
spaces of the projections will really have to be $X^2_H$. To be
able to use the pullback- and pushforward-results in
\cite{Meaomwc} (see also Appendix \ref{appzeronotation}) we will
have to modify the ``triple space'' such that the corresponding
lifted projections become $b$-fibrations.
\par
\begin{remark}\rm \label{globalsquare}
To simplify things a bit, we agree to define the $d$-heat space
$X^2_H$ using  $\sqrt{t}$ as a {\em global} variable, i.e. as the
simple blow up $$ X^2_H=\left[ X^2_b\times
[0,\infty[_{\sqrt{t}};(F_{\phi}\times \{0\});
 \bigtriangleup_X\times\{0\} \right]\quad\mbox{away from tb}.$$
Since elements in the heat calculus vanish at $\fatb$ by
definition, nothing is lost by considering the RHS above, instead
of the original space $X^2_H$.
\end{remark}
\par
\medskip
We now define $${\sf T}^2:=
[0,\infty[_{\sqrt{t'}}\times[0,\infty[_{\sqrt{t-t'}},$$ and denote
its faces by ${\rm lt}$ and ${\rm rt}$. The first thing to note
now, is that the map $$\pi_M:{\sf T}^2  \rightarrow
[0,\infty[_{\sqrt{t}} \quad (a,b)\mapsto \sqrt{a^2+b^2}$$ is {\em
not} a $b$-fibration. Thus we will have to take care to also blow
up the corner ${\rm bt}={\rm lt}\cap{\rm rt}$ whenever this
situation arises. So our starting point is the diagram
\begin{equation}\label{startdiagram}
\begin{array}{ccccc}
X^2_b\times[0,\infty[_{\sqrt{t-t'}}&\stackrel{\pi_L}{\longleftarrow}&
X^3_b\times{\sf T}^2&\stackrel{\pi_R}{\longrightarrow}&
X^2_b\times[0,\infty[_{\sqrt{t'}}\\
 & & \quad\downarrow \pi_M \\
& & X^2_b\times[0,\infty[_{\sqrt{t}}
\end{array}
\end{equation}
where we keep in mind that the horizontal maps are still
$b$-fibrations but $\pi_M$ is not.
\par
\smallskip
The definition of the heat triple space now proceeds in two steps.
The first step consists in replacing the image spaces in
(\ref{startdiagram}) by the partial blow up
$\widetilde{X}^2_H=[X^2_b\times[0,\infty[;F_{\phi}\times\{0\},\dd
t ]$, which we interpret as simple blow up w.r.t. the global
variable $\sqrt{t}$.
\begin{table}[h]
\centerline{\epsfig{file=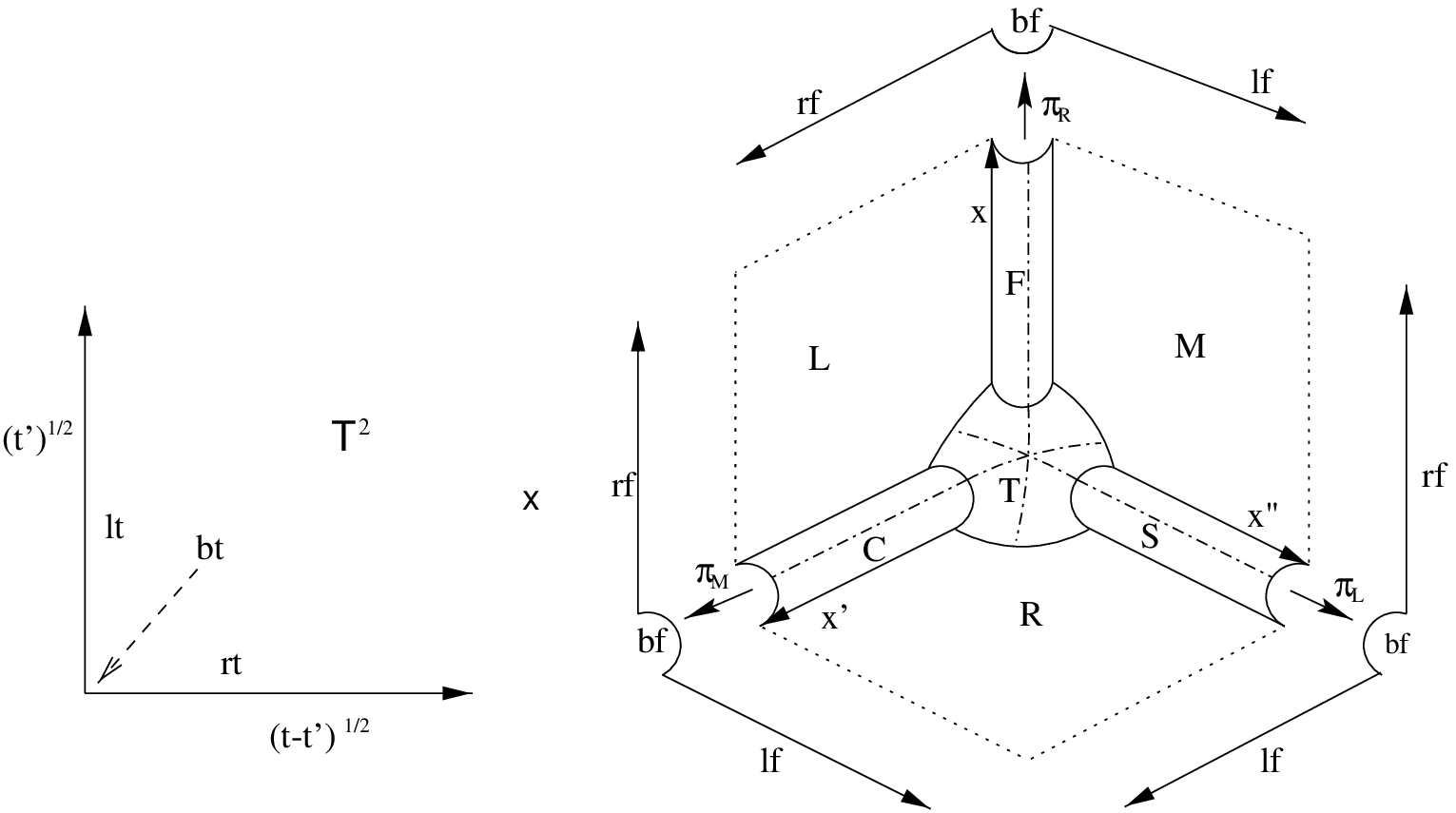}} \centerline{Figure 7: The
${\sf T}^2$- and $X^3_b$-factor}
\end{table}
\par
 As in the case of the
$\phi$-triple space  we denote the faces to be blown up (as well
as the resulting faces) by $${\rm bt}\times\phi_C, \,{\rm
bt}\times\phi_{CT}, \,{\rm bt}\times\phi_{TT},\,{\rm rt}\times
\phi_F,\,{\rm lt}\times \phi_S,
 \,{\rm rt}\times\phi_{FT},  \,{\rm lt}\times\phi_{ST},$$
  and define the (intermediate) heat triple space to be
$$ \widetilde{X}^3_H := \left[
   X^3_{b} \times {\sf T}^2; {\rm bt}\phi_{TT};
 {\rm bt}\phi_{CT};  {\rm rt}\phi_{FT}; {\rm lt}\phi_{ST};
 {\rm bt} \phi_C; {\rm rt}\phi_F; {\rm lt}\phi_S  \right] .
$$ To avoid any misunderstandings, recall that this is the space
$X^3_{b} \times {\sf T}^2$, with faces blown up as listed from
left to right. Each face in that list thus really represents the
{\em preimage} of the corresponding face in $X^3_{b} \times {\sf
T}^2$ under the previous blow ups.
\par
For the second step, denote the different lifts of the diagonal
$\bigtriangleup_{\phi}\subset X^2_{\phi}$ to $X^3_{\phi}$ by
$\bigtriangleup_L$, $\bigtriangleup_M$, $\bigtriangleup_R$ and
analogously their intersection by $\bigtriangleup_{LMR}$. Then the
true ``triple space'' is defined as
$$X^3_H:=[\widetilde{X}^3_H;{\rm bt}\bigtriangleup_{LRM}; {\rm
bt}\bigtriangleup_{M}; {\rm lt}\bigtriangleup_L,{\rm
rt}\bigtriangleup_R].$$ Besides those faces obtained by the above
blow ups, this space has the faces $$M\sowie {\sf T}^2M, L, R, S,
C, T,\quad \mbox{and}\quad {\rm bt}\sowie {\rm bt}X^3_b, {\rm rt},
{\rm lt}.$$ Again, these are understood to be the lifts of the
corresponding faces in $X^3_b\times{\sf T}^2$ minus those faces,
which were obtained in the previous blow up. The important result
now is
\begin{lemma}
There are induced $b$-maps $
X^3_H\stackrel{\pi_{H,o}}{\longrightarrow}
  X^2_H$, $o=L,M,R$, (where we agreed to define the RHS to be
  defined with global time variable $\sqrt{t}$) such that the diagrams
 $$ \begin{array}{ccc}
X^3_H&\stackrel{\pi_{H,o}}{\longrightarrow}&X^2_H\nonumber
\\ \downarrow& & \downarrow
\\ X^3_b\times{\sf T}^2& \stackrel{\pi_o}{\longrightarrow}& X^2_b\times
{\sf T}\nonumber
\end{array} $$
commute.
\end{lemma}
\begin{proof}
We show how the map $\pi_{H,L}$ is defined. Starting with $X^3_H$
we first blow down all those faces which do not intersect any of
the faces lifted from the left factor:
 $$ {\rm
rt}\bigtriangleup_R,\, {\rm bt}\bigtriangleup_M, \,{\rm
rt}\phi_F,\, {\rm rt}\phi_{FT},\, {\rm bt}\phi_C,\,{\rm
bt}\phi_{CT},\, F,\,C, $$ thus getting
 $$[X^2_b\times
X\times{\sf T}^2; \fabf\times\partial X\times{\sf T}^2; {\rm
bt}\phi_{TT}; {\rm lt}\phi_{ST};{\rm lt}\phi_{S};{\rm
bt}\bigtriangleup_{LRM}; {\rm lt}\bigtriangleup_{L}].$$ Now,
obviously, ${\rm bt}\bigtriangleup_{LRM}$ is a p-submanifold of
${\rm lt}\bigtriangleup_{L}$, i.e. their blow ups can be exchanged
and the triple diagonal can be blown down. The same argument let
us us exchange the blow ups of   ${\rm bt}\phi_{TT}$ and  ${\rm
lt}\phi_{ST}$, yielding
 $$[X^2_b\times
X\times{\sf T}^2; \fabf\times\partial X\times{\sf T}^2; {\rm
lt}\phi_{ST};{\rm lt}\phi_{S}; {\rm bt}\phi_{TT}; {\rm
lt}\bigtriangleup_{L}].$$ Now the face $ {\rm bt}\phi_{TT}$ in
question here is really the preimage of this face under the blow
up of $ {\rm lt}\phi_{ST}$, and one easily checks that it is
transversal to $ {\rm lt}\bigtriangleup_{L}$. Being disjoint from
the other blow ups, we can switch the corresponding blow up to the
front and then blow down ${\rm bt}\phi_{TT}$ altogether. Following
this, it is easy to see that the blow up of $\fabf\times\partial
X\times{\sf T}^2$ can be ``pulled to the front'', if one switches
the order of the blow ups $ {\rm lt}\phi_{ST}$ and ${\rm
lt}\phi_{S}$. Thus this face can be blown down, too, leaving us
with
 $$[X^2_b\times
X\times{\sf T}; F_{\phi}\times X \times \{0\}; F_{\phi}\times
\partial X \times \{0\};\bigtriangleup_{b}\times X \times\{0\}]\times
{\sf T},$$ where we have written the faces $\phi_{S}$, $\phi_{ST}$
in their blown down form in $X^2_b\times X$. It is important to
note that, again, $\phi_{ST}\times\{0\}$ is really the preimage of
that face under the blow up of $\phi_S\times \{0\}$, making it
transversal to the lifted left diagonal. Thus $\phi_{ST}\{0\}$ can
be blown down and we are left with $$[X^2_b\times{\sf T};
F_{\phi}\times\{0\} ; \bigtriangleup_{b}\times\{0\}]\times
X\times{\sf T},$$ which projects nicely to $X^2_H$.
 \qed
\end{proof}
Since all their entries are either $0$ or $1$, we can  write down
the ``coefficient matrices'' for $\pi_L$, $\pi_M$, $\pi_R$ (see
Appendix \ref{appzeronotation}) by just listing the preimages of
the different faces of $X^2_H$:
\begin{lemma}
\begin{enumerate}
\item $\pi_L$ maps the faces $\{L, {\rm rt}, {\rm
rt}\bigtriangleup_R \}$ to the interior of $X^2_H$. Also
\par
 $\begin{array}{rcl}
   \pi_L^*\fahbf&=&\{ D, T, {\rm bt}\phi_{CT}, {\rm rt}\phi_{FT} \}
\\ \pi_L^*\fatff&=&\{  {\rm lt}\phi_{S},
 {\rm lt}\phi_{ST}, {\rm bt}\phi_{TT}  \}
\\ \pi_L^*\fatf&=&\{    {\rm lt}\bigtriangleup_L,
{\rm bt}\bigtriangleup_{LMR} \}
\\ \pi_L^*\fahlf&=&\{ R, C, {\rm bt}\phi_C   \}
\\ \pi_L^*\fahrf&=&\{ M, F, {\rm rt}\phi_F   \}
\\ \pi_L^*\fatb&=&\{ {\rm bt}, {\rm lt},  {\rm bt}\phi_C, {\rm
bt}\phi_{CT}, {\rm bt}\bigtriangleup_M \}
\end{array}$
\item $\pi_R$ maps the faces $\{R, {\rm lt}, {\rm
lt}\bigtriangleup_L \}$ to the interior of $X^2_H$. Also
\par
 $\begin{array}{rcl}
   \pi_R^*\fahbf&=&\{ F, T, {\rm bt}\phi_{CT}, {\rm lt}\phi_{ST} \}
\\ \pi_R^*\fatff&=&\{  {\rm rt}\phi_{F},
 {\rm rt}\phi_{FT}, {\rm bt}\phi_{TT}  \}
\\ \pi_R^*\fatf&=&\{ {\rm rt}\bigtriangleup_R,
{\rm bt}\bigtriangleup_{LMR} \}
\\ \pi_R^*\fahlf&=&\{ M, S,  {\rm lt}\phi_S   \}
\\ \pi_R^*\fahrf&=&\{ L, C, {\rm bt}\phi_C   \}
\\ \pi_R^*\fatb&=&\{ {\rm bt}, {\rm rt}, {\rm bt}\phi_C, {\rm
bt}\phi_{CT}, {\rm bt}\bigtriangleup_M \}
\end{array}$
\item $\pi_M$ maps the faces $\{S, {\rm rt}, {\rm lt}, {\rm
rt}\bigtriangleup_R, {\rm lt}\bigtriangleup_L \}$ to the interior
of $X^2_H$. Also
\par
$\begin{array}{rcl}
   \pi_M^*\fahbf&=&\{ C, T, {\rm rt}\phi_{FT}, {\rm lt}\phi_{ST} \}
\\ \pi_M^*\fatff&=&\{ {\rm bt}\phi_{C}, {\rm bt}\phi_{CT}, {\rm bt}\phi_{TT} \}
\\ \pi_M^*\fatf&=&\{ {\rm bt}\bigtriangleup_M, {\rm bt}\bigtriangleup_{LMR} \}
\\ \pi_M^*\fahlf&=&\{ R, M,  {\rm lt}\phi_S   \}
\\ \pi_M^*\fahrf&=&\{ L, F, {\rm rt}\phi_F   \}
\\ \pi_M^*\fatb&=&\{ {\rm bt} \}
\end{array}$
\par \qed
\end{enumerate}
\end{lemma}
This Lemma shows that $\pi_L$, $\pi_R$, are not $b$-normal, and
thus cannot be $b$-fibrations, since the faces ${\rm bt}\phi_{C}$,
${\rm bt}\phi_{CT}$ appear in the preimage of several faces of
$X^2_H$. However, for two heat operators $$A\in \Psi_H^{\cal
L}(X)={\cal A}^{\cal L}(X^2_H,\rho_{\fahbf}^{-v}\KD_H),\quad B\in
\Psi_H^{\cal R}(X)={\cal A}^{\cal
R}(X^2_H,\rho_{\fahbf}^{-v}\KD_H)$$ their lift from the left and
the right to $X^3_H$ is $$\pi^*_LA\cdot \pi^*_RB \in {\cal
A}^{\cal H}(X^3_H,\Omega(X^3_H)),$$ and the above Lemma tells us
that
\begin{eqnarray}
\lefteqn{H_L= H_M=H_R=H_F=H_S=H_C=H_{\rm bt}=H_{\rm rt}=H_{\rm
lt}} & & \nonumber\\ & &\qquad\qquad\qquad\qquad =H_{{\rm
bt}\phi_C}=H_{{\rm
    lt}\phi_S}
 =H_{{\rm rt}\phi_F}
    =H_{{\rm bt}\phi_{CT}} =H_{{\rm bt}\bigtriangleup_M}
    =\infty.\label{Htrivial}
\end{eqnarray}
To describe the expansions at the remaining faces, we need to
caculate the lifts of our densities under the blow down map
$\gamma_H:X^3_H\rightarrow X^3\times [0,\infty[^2$. Writing
$s=t-t'$ for the (right) variable in ${\sf T}^2$ note first the
coefficients which appear when lifting to $X^3_b\times{\sf T}^2$:
\begin{equation}\label{lifttob}
\gamma_b^*\Omega(X^3)\otimes\,{}^b\Omega({\sf T}^2)=
   \rho_{T{\sf T}^2_b}^2\,\rho_{F{\sf T}^2_b}\,\rho_{S{\sf T}^2_b}
   \,\rho_{C{\sf T}^2_b}\,\rho_{X^3_b{\rm rt}}^{-1}\,
\rho_{X^3_b{\rm lt}}^{-1}\,\gamma_{H,b}^*\Omega(X^3_b\times{\sf
T}^2).
\end{equation}
 To determine the relevant coefficients of the lift of
$$\Omega_L(X)\pi_L^*(\rho_{\fahbf}^{-v}\KD_H)\pi^*_R(\rho_{\fahbf}^{-v}\KD_H)
=(s^{1/2}t^{\prime 1/2})^{-n}(x'x'')^{v-1}\,
\pi_L^*\rho_{\fahbf}^{-v}\,\pi_R^*\rho_{\fahbf}^{-v}
\Omega(X^3)\bomega({\sf T}^2),$$ we use Lemma \ref{liftdensity},
i.e we need to list the contributions of the factors in
(\ref{lifttob}) at the faces of $X^3_H$ and to keep track of the
codimensions of the faces blown up in $\partial(X^3_b\times{\sf
T}^2)$:
 $$\begin{array}{|l|l|l|l|l|l|l|}
\hline & & & & & &
\\ \mbox{face} &\mbox{codim} &
(\ref{lifttob})&
 (s^{1/2} t^{\prime 1/2})^{-n} &(x' x'')^{v-1}
 & \pi_L^*\rho_{\fahbf}^{-v}\pi_R^*\rho_{\fahbf}^{-v}
 &\mbox{Total}
 \\ \hline & & & & & &
 \\ T& 0 & 2 & 0  & 2v-2&-2v  &0
\\ {\rm bt}\phi_{TT}& 2h+4 & 0 & -2n  & 2v-2&  &0
\\{\rm rt}\phi_{FT} & h+2 & 1 & -n &2v-2  &  -v  &0
\\{\rm lt}\phi_{ST} & h+2 & 1 & -n &2v-2 &  -v  & 0
\\{\rm bt}\bigtriangleup_{LMR} &2n+1 & -2 & -2n &   &   & -1
\\{\rm rt}\bigtriangleup_{R} & n & -1 & -n &  &   &-1
\\{\rm lt}\bigtriangleup_{L} & n & -1 & -n &  &   &-1
\\ \hline
 \end{array} $$
Thus, at the relevant faces, the density is
\begin{eqnarray}
\lefteqn{\gamma_b^*\Omega_L(X)(\pi_L^*\rho_{\fahbf}^{-v}\KD_H)(\pi^*_R
\rho_{\fahbf}^{-v}\KD_H) }\nonumber
\\  & &\qquad\qquad=\rho_T^{0}\,\rho_{{\rm rt}\phi_{FT}}^{0}\,\rho_{{\rm
lt}\phi_{ST}}^{0} \,\rho_{{\rm bt}\phi_{TT}}^{0}\,\rho_{{\rm
lt}\bigtriangleup_{L}}^{-1}\, \rho_{{\rm rt}
\bigtriangleup_{R}}^{-1}\,  \rho_{{\rm bt}
\bigtriangleup_{LMR}}^{-1} \,\Omega(X^3_H),\label{x3hdensity}
\end{eqnarray}
and the nontrivial terms of the index set ${\cal H}$ are given by
$$\begin{array}{lclclcl}
 H_{T} &=&L_{\fahbf}+R_{\fahbf}  \\
H_{{\rm rt}\phi_{FT}}&=& L_{\fahbf}+R_{\fatff} & &  H_{{\rm
lt}\phi_{ST}}&=& L_{\fatff}+R_{\fahbf}
\\ H_{{\rm bt}\phi_{TT}} &=&L_{\fatff}+R_{\fatff}
\\ H_{{\rm lt}\bigtriangleup_L}&=& L_{\fatf}-1 & & H_{{\rm
rt}\bigtriangleup_R}&=& R_{\fatf}-1
\\  H_{{\rm bt}\bigtriangleup_{LMR}}&=&L_{\fatf}+R_{\fatf}-1
\end{array}$$
Though $\pi_M$ is not a $b$-fibration it can be made into one in
this special case, because (\ref{Htrivial}) allows the faces
causing the nuisance to be blown down. Thus the pushforward
$$\pi_{M,*}(\gamma_{H,L}^*f\pi^*_LA\cdot \pi^*_RB)\in {\cal
A}^{\pi_{M \sharp}{\cal H}}(X^2_H,\Omega(X^2_H))$$ is well defined
when it is integrable at those faces mapped to the interior by
$\pi_M$, i.e. whenever $$ H_{{\rm rt}\bigtriangleup_R}, H_{{\rm
lt}\bigtriangleup_L} > -1 \quad\mbox{or equivalently}\quad
R_{\fatf},L_{\fatf}>0 $$ (the condition at the other faces again
is trivial due to (\ref{Htrivial})), and the pushforward of the
index set is given by: $$ \pi_{M \sharp}{\cal H}(\fahbf)=H_T
\overline{\cup} H_{{\rm lt}\phi_{ST}} \overline{\cup} H_{{\rm
rt}\phi_{TF}},\quad \pi_{M \sharp}{\cal H}(\fatff)=  H_{{\rm
bt}\phi_{TT}}, \quad \pi_{M \sharp}{\cal H}(\fatf)= H_{{\rm
bt}\bigtriangleup_{LMR}}.$$ Also, in the case
$R_{\fatf},L_{\fatf}\geq 0$, the pushforward can still be defined
using the mean value condition (\ref{meanvaluecond}).
 Rewriting
everything in terms of
$\beta_{H,L}^*\Omega(X)\rho_{\fahbf}^{-v}\KD_H=\rho_{\fatf}^{-1}\Omega(X^2_H)$
we finally get
\begin{theorem}[Composition Formula]
  \label{heatcomposition} Let $R_{\fatf},L_{\fatf}\geq 0$. Then
composition gives a map $$ \Psi^{\cal L}_{H}(X) \times\Psi^{\cal
R}_{H}(X) \stackrel{\circ}{\longrightarrow}\Psi^{\cal M}_{H}(X)$$
with index sets given by $$
M_{\fahbf}=(L_{\fahbf}+R_{\fahbf})\overline{\cup}
(L_{\fatff}+R_{\fahbf})\overline{\cup}(L_{\fahbf}+R_{\fatff}),$$
$$ M_{\fatff}=L_{\fatff}+R_{\fatff},\quad
M_{\fatf}=L_{\fatf}+R_{\fatf}.$$ \qed
\end{theorem}

\subsection{Normal Operators and Construction of the Heat Kernel}  \label{normalops}
Using the local coordinate representations described in Section
\ref{dheatbasics} it is easy to see that the restriction of the
density $\KD_H$ to the different faces has the following form.
\begin{eqnarray}
 \KD_H|_{\stackrel{\circ}{\fatf}}&\cong &\frac{\dd t}{t}\Omega(\phitx/X)\\
   \KD_H|_{\stackrel{\circ}{\fatff}}&\cong& \frac{\dd \tau}{\tau}
\tau^{-n/2}\omfib( \phindx\times_Y\partial X
\stackrel{\beta_{H,L}}{\longrightarrow}{\partial X})
\label{KDtf=}\\ \KD_H|_{\stackrel{\circ}{\fahbf}}&\cong& \frac{\dd
t}{t}X_t^v\,{}^b\omfib(\fahbf\stackrel{\beta_{H,L}}{\rightarrow}\partial
X )
\end{eqnarray}
For  $A\in\Psi_{H,{\rm cl}}^{l,m,p}(X)$ we now define the normal
operators at the faces $\fatf$ and $\fatff$ and $\fahbf$ by
\begin{eqnarray*}
N_{H,\fatf,l>0}(A)&:=&
  t^{1-l/2}A|_{\fatf}\in \rho_{\fatff}^{m-l}
   \cdotinfty{\fatff}(\fatf,t \KD_H)
\cong x^{m-l}{\cal S}_{\rm fib} (\phitx,\Omega(\phitx/X)\dd t)
\\ N_{H,\fatf,0}(A\oplus c)&:=&
  t A|_{\fatf}\oplus c \in \rho_{\fatff}^{m}
 \cdotinfty{\fatff}(\fatf,t \KD_H)^{\fatf}\oplus \CC
 \\ N_{H,\fatff,m}(A)&:=&(x')^{-m}A|_{\fatff}\in
\rho_{\fatf}^{l}\,\rho^{p-v-m}_{\fahbf}\,\cdotinfty{\fatf,\fahbf}
(\fatff,\KD_H)
\\ N_{H,\fahbf,p}(A)&:=& (x')^{-p}A|_{\fahbf}
                 \in \rho_{\fatff}^{m-p} \cdotinfty{\fatff}(\fahbf,
                 {}^b\omfib(\fahbf)\dd t)
\end{eqnarray*}\par
For $l\geq 1$ these normal operators fit into the usual short
exact sequences:
\[ 0  \longrightarrow  \Psi^{l-1,m,p}_{H,{\rm cl}}(X)  \longrightarrow
  \Psi^{l,m,p}_{H,{\rm cl}}(X)  \stackrel{ N_{H,\fatf,l}}{\longrightarrow}
   x^{m-l}{\cal S}_{\rm fib}(\phitx,\Omega(\phitx/X)\dd t)
   \longrightarrow 0 \]
\begin{equation}
 0  \longrightarrow  \Psi^{l,m-1,p}_{H,{\rm cl}}(X)  \longrightarrow
 \Psi^{l,m,p}_{H,{\rm cl}}(X)  \stackrel{ N_{H,\fatff,m}}{\longrightarrow}
\rho_{\fatf}^{l}\,\rho^{p-v-m}_{\fahbf}\,\cdotinfty{\fatf,\fahbf}(\fatff,
\KD_H) \longrightarrow 0\label{RHSsusp}
\end{equation}
\[ 0  \longrightarrow  \Psi^{l,m,p-1}_{H,\rm cl}(X)  \longrightarrow
  \Psi^{l,m,p}_{H,\rm cl}(X)  \stackrel{ N_{H,\fahbf,p}}{\longrightarrow}
  \rho_{\fatff}^{m-p} \cdotinfty{\fatff}(\fahbf,{}^b\omfib(\fahbf)\dd t
  ) \longrightarrow 0, \]
and we leave the formulation of the case $l=0$ to the reader, as
well as the usual extension to the calculus with coefficients.
\par
\smallskip
As usual, the spaces on the RHS in the above sequences represent
some calculi of operators whose composition rules are induced by
the composition rules for operators in the ``big'' calculus
$\Psi_{H,{\rm cl}}(X)$ (compare Lemma \ref{leminducedcalc}). For
instance, in the case (\ref{RHSsusp}) we can define $$
\Psi_{H,{\rm
sus}}^{l,p}(\phindx):=\rho_{\fatf}^{l}\,\rho^{p-v}_{\fahbf}
\,\cdotinfty{\fatf,\fahbf}(\fatff,\KD_H).$$ This is the  heat
space associated with the suspended calculus introduced in Section
\ref{normalop}. Some related calculi are described in more detail
in Appendices \ref{Mehler} and \ref{vertfam}.
\par\medskip
At this point we are not interested in the general algebraic
properties of the normal maps. All we really need is to describe
the normal action of the Dirac operator $\Dir$ at the different
faces, for which it will suffice  -- as usual -- to  analyze
compositions of vector fields with elements in the heat calculus.
To do this, we remind the reader of the following simple fact
about Lie derivatives of forms and densities: Let
$j_F:F\hookrightarrow M$ be a face of a manifold with corners with
defining function $f$, W a vector field on $M$ tangent to $F$, and
$\omega$ a form on $M$. Then $$
j_F^*L_{fW}\omega=j_F^*(fL_W\omega+\ext(\dd f)\Int(W)\omega)=0,$$
which means that $$j_F^*L_W\omega =L_{W|_F}j_F^*\omega,$$ i.e. the
Lie derivative restricted to the face can be calculated {\em on}
the face. The normal action of vector fields on the elements in
the calculus is now described in
\begin{proposition}\label{normcalc}
Let $V$ be in $\Gamma(\phitx)$.  Composition gives maps $$
\Psi_{H,{\rm cl}
}^{l,m,p}(X)\stackrel{\frac{1}{x}V\circ}{\longrightarrow}
   \Psi_{H,{\rm cl}}^{l-1,m-1,p-1}(X),\qquad \Psi_{H,{\rm cl}}^{l,m,p}(X)
   \stackrel{\frac{\partial}{\partial t}\circ}{\longrightarrow}
   \Psi_{H,\rm cl}^{l-2,m-2,p}(X)$$
Applied to $A\in\Psi^{l,m,p}_{H,{\rm cl}}(X)$, we get at the
different faces
\begin{enumerate}
\item $N_{H,\fatf,l-2}({\frac{\partial}{\partial t} }\circ A)=
\left\{
 \begin{array}{ll}
     -\frac{1}{2}( L_{R^{\phi}}+2-l )N_{H,\fatf,l}(A) & \mbox{for}\quad l>2 \\
     -\frac{1}{2}( L_{R^{\phi}}+2-l )N_{H,\fatf,2}(A)\oplus
    -\int_{\phitx/X} N_{H,\fatf,2}(A) & \mbox{for}\quad l=2\\
 \end{array}   \right.$
\item $   N_{H,\fatf,l-1}({\frac{1}{x}V}\circ A)=
L_V  N_{H,\fatf,l}(A)$
\item  $N_{H,\fatff,m-2}(\frac{\partial}{\partial t}\circ A)=
 L_{\frac{\partial}{\partial \tau}}N_{H,\fatff,m}(A) $
\item $   N_{H,\fatff,m-1}({\frac{1}{x}V}\circ A)=
L_{N_{\phi}(V)}  N_{H,\fatff,m}(A)$
\item  $N_{H,\fahbf,p}({\frac{\partial}{\partial t}}\circ A)=
 L_{\frac{\partial}{\partial t}}N_{H,\fahbf,p}(A) $
\item $   N_{H,\fahbf,p}({\frac{1}{x}V}\circ A)=
L_{\beta_H^*\frac{1}{x}V}  N_{H,\fahbf,p}(A)$, for
$V\in\Gamma([\phindx])$.
\end{enumerate}
\end{proposition}
\begin{proof}
First, recall the partial integration formula for the Lie
derivative w.r.t. a vector field $V$ on a manifold with boundary
$N$: $$ \int_N \alpha \wedge L_V\beta = -\int_N(L_V\alpha)\wedge
\beta +\int_{\partial N}\Int(V)\alpha\wedge\beta .$$ Now let  $f
\in  \cdotinfty{c}(X\times[0,\infty[, \Omega(X)\dd t)$ and $g \in
\cdotinfty{c} (X\times[0,\infty[) $. We then have
\begin{eqnarray*}
\langle L_{\frac{\partial}{\partial t}} (A*g),f\rangle_{X\times[0,\infty[} &=&
-\langle A*g, L_{\frac{\partial}{\partial t}} f\rangle_{X\times[0,\infty[}
+\int_{X\times\{0\}}\Int(\frac{\partial}{\partial t})A*g\cdot f\\
&=& -\langle A, \beta_{H}^*(L_{\frac{\partial}{\partial t}} f)*g
\rangle_{X_H^2}=-\langle A,
L_{\beta_{H}^*\frac{\partial}{\partial t}} \beta_{H}^*(f*g)
\rangle_{X_H^2}\\
&=&\langle  L_{\beta_{H}^*\frac{\partial}{\partial t}}A ,
\beta_H^*(f*g)\rangle_{X_H^2}
-\int_{\fatf\cup\fatff}\Int(\frac{\partial}{\partial t})A\cdot \beta_H^*
( f*g).
\end{eqnarray*}
Since near $\fatf$ the density factor in $A\cdot \beta_H^*( f*g)$
is of the form $\rho_{\fatf}^l\, \rho_{\fatff}^m \,\frac{\dd
t}{t}\beta_{H,L}^*\Omega(X)\Omega\phitx/X$ we see that the last
integral in the above expression vanishes for $l>2$, but
contributes the constant $$\int_X
\left(\int_{\phitx/X}N_{H,\fatf,2}(A)\right) \langle
f,g\rangle_{[0,\infty[}$$ for $l=2$. Thus, recalling that
$$tL_{\frac{\partial}{\partial t}}
        = L_{t\frac{\partial}{\partial t}}-1,\qquad L_{\frac{\partial}{\partial t}}
\dd t = \dd t, \qquad \mbox{and}\quad L_{R_{\phi}}\dd t =0,$$ and
using Lemma \ref{liftface}, we find for $l>2$:
\begin{eqnarray*}
N_{H,\fatf,l-2}(L_{\frac{\partial}{\partial t}}\circ A) &=&
t^{1-l/2} t L_{\beta_H^*\frac{\partial}{\partial t}}A|_{\fatf}
=t^{1-l/2} ( L_{\beta_H^*t\frac{\partial}{\partial t}}-1)A|_{\fatf}\\
& =& ( L_{\beta_H^*t\frac{\partial}{\partial t}}-2+l/2)t^{1-l/2}A|_{\fatf}
= -\frac{1}{2}( L_{R^{\phi}}+2-l)N_{H,\fatf,l}(A).
\end{eqnarray*}
Note also that in the case $l=2$ we have $$
N_{H,\fatf,0}(\frac{\partial}{\partial t}\circ A)=
 -\frac{1}{2} L_{R^{\phi}}N_{H,\fatf,2}(A)=
 -\frac{1}{2}\dd t L_{R^{\phi}}(N_{H,\fatf,2}(A)/\dd t),$$
therefore $$ \int_{\phitx/X}
N_{H,\fatf,0}(\frac{\partial}{\partial t}\circ A)=0 ,$$ i.e. the
mean value condition is fulfilled.\qed
\end{proof}\par
As usual, the above results can all be formulated for operators
with coefficients in a vector bundle. We leave this to the reader.
Let us just make the remark that the construction is independent
of the choice of connection made. Also, for $W\in\Gamma(\phitx)$,
the normal operators of $\tau^{1/2}\nabla_W^d$ and
$\tau^{1/2}\nabla_W^{\phi}$ at $\fatf$ and $\fatff$ are the same.
\begin{corollary}
\begin{enumerate}
\item $N_{H,\fatf}(t^{1/2}\Dir)=
  N_{H,\fatf}(t^{1/2}\Dir^d)=
  N_{H,\fatf}(t^{1/2}\frac{1}{x}\Dir^{\phi})=\cl_{\phitx}\circ \dd_{\phitx}$
\item $N_{H,\fatf}(t\Dir^2)=\Delta_{\phitx}$,
   the fibrewise Laplacian for the fibre metric $g_{\phitx/X}$.
\item $N_{H,\fatff}(t^{1/2}\Dir)=t^{1/2}\frac{1}{x}N_{H,\fatff}(\Dir^{\phi})
=\tau^{1/2} N_{\phi}(\Dir^{\phi})$
\item $N_{H,\fatff}(t\Dir^2) =\tau N_{\phi}(\Dir^{\phi})^2=\tau\Delta_{\phindx}
+\tau(\Dir^{\phi,V})^2$.\qed
\end{enumerate}
\end{corollary}
\par
\smallskip
The description of the action of $\Dir^2$  can be described in the same way
using Lemma \ref{indicialfamily} (d). We present this as part of the proof of
\begin{theorem}[Parametrix for the Heat Equation]\label{heatnormal}
There exists an element $K\in \Psi^{2,2,0}_{H,{\rm cl}}(X,E)$ such
that
\begin{equation}\label{heat1}
\left( \frac{\partial}{\partial t}+\Dir^2\right)*  K= I+R ,\qquad
   R\in \Psi^{\infty,1,1}_{H,{\rm cl}}(X,E).
\end{equation}
\end{theorem}
\begin{proof}
We start by solving the above equation to first order at $\fatf$.
Applying $N_{H,\fatf,0}$ to both sides, we get first
$$\left(-\frac{1}{2}L_{R^{\phi}}+\Delta_{\phitx} \right)
N_{H,\fatf,2}(K_2)=0\quad\mbox{and}\quad
\int_{\phitx/X}(N_{H,\fatf,2}(K_2)/\dd t)=1.$$ This is to be
understood as an equation on the fibres of $\phitx$. Writing
$N_{H,\fatf,2}(K_2)=k_2^{\fatf} \dd t\dvol_{\phitx/X}$ (and
forgetting about the ${}^{\fatf}$-index until it is really needed)
we get
$$(-\frac{1}{2}(R_{\phi}+n)+\Delta_{\phitx})k_2=0,\qquad
     \int_{\phitx/X}k_2 \dvol_{\phitx/X}=1.$$
Taking the Fourier transform in $\phitx$ yields
$$(\frac{1}{2}({\bf \xi}\frac{\partial}{\partial \bf\xi})
+\Betrag{\xi}^2)\widehat{k_2}=0, \qquad \widehat{k_2}(0)=1
$$ the solution of which is easily seen to be $$
k_2=\frac{1}{\sqrt{4\pi}^n}\exp(-\Norm{.}_{\phi}^2/4)\in x^0{\cal
S }_{\rm fib}(\phitx),$$ which also shows that the solution is of
order $2$ at $\fatff$ in the calculus. Thus the compatibility
condition at the corner reads
$$N_{H,\fatff,2}(K_2)_{\fatf\cap\fatff}=(x')^{-2}
N_{H,\fatf,2}(K_2)_{\fatf\cap\fatff}=k_2^{\fatf} \dd
\tau\dvol_{\phitx/X}$$ and the heat equation on the interior of
$\fatff$ is
\begin{equation}\label{heat1tff}
 \tau( L_{\tau\frac{\partial}{\partial \tau}}+
 \Delta_{\phindx}+(\Dir^{\phi,V})^2)
N_{H,\fatff,2}(K_2)=0
\end{equation}
Writing $N_{H,\fatff,2}(K_2)=k_2^{\fatff} \dvol_{\phindx/Y}\dd
\tau$ we get the equation $$ ( \frac{\partial }{\partial \tau}+
\Delta_{\phindx}+(\Dir^{\phi,V})^2 )k_2^{\fatff}=0,$$ which is
solved by $$
 k_2^{\fatff} =
\sqrt{4\pi\tau}^{-h-1}e^{-\Norm{.}^2_{\phindx}/{4\tau}}
([\exp(-\tau(\Dir^{\phi,V})^2)]/\dd \tau\beta_{H,R}^*\dvol_{M/Y}
).$$ The compatibility condition at $\fahbf\cap\fatff$ is then
\begin{equation}
  N_{H,\fahbf,0}(K_2)|_{\fatff}=(x')^2N_{H,\fatff,2}(K_2)|_{\fahbf}
=\frac{e^{-\Norm{.}^2_{\phindx}/4}}{\sqrt{4\pi}^{h+1}}
[\Pi_{\circ}]\frac{\dd s}{s}\dvol_{\phindx/\partial X}\dd
t,\label{initialhbf}
\end{equation}
where $[\Pi_{\circ}]\in C^{\infty}(M\times_Y
M,\END\otimes\beta_{R}^*\Omega(M/Y))$ is the kernel of the
projection onto the null space of $\Dir^{\phi,V}$. Then
(\ref{initialhbf}) is the initial value for the problem at
$\fahbf$.
\par\smallskip
In order to solve the heat equation at $\fahbf$, recall that
$$x\Dir: C^{\infty}(X,E)\rightarrow C^{\infty}(X,E),
\quad x\Dir|_{\partial X}=\Dir^{\phi,V}
,\quad\mbox{and}\quad \Dir: C^{\infty}(X,E)^{\circ}\rightarrow
C^{\infty}(X,E).$$
Thus for $G\in \Psi_{H,{\rm cl}}^{2,2,0}(X)$  the expression $\Dir^2 G$
a priori only lies in
$\Psi_{H,{\rm cl}}^{0,0,-2}(X)$. However, it follows from Lemma
\ref{indicialfamily} (d) that $G|_{\fahbf}$ can be extended to the
interior in such a way that $\Dir^2G \in \Psi_{H,\rm
cl}^{\infty,1,0}(X)$ and
\begin{equation}\label{heatclaim}
\Dir^2 G|_{\fahbf}
=(\Pi_0\Dir\Pi_0)^2G_0|_{\fahbf}=I_b(\Dir)^2G|_{\fahbf}
\end{equation}
 It is then clear that the heat
equation at $\fahbf$ can be solved by choosing $G|_{\fahbf}$ more
concretely, namely $$G|_{\fahbf}=[\exp(-tI_b(\Dir)^2)]
=[\exp(-t\Dir_Y^2)]\frac{1}{\sqrt{4\pi
t}}e^{-(\ln(s))^2/4t}\frac{\dd s}{s},$$
and then use the extension which fulfills (\ref{heatclaim}).
\par\noindent
Since the solutions at the faces $\fatf$, $\fatff$, $\fahbf$ are
compatible, we can choose an extension $K_2\in \Psi_{H,{\rm cl
}}^{2,2,0}(X)$ with $$N_{H,\fatf,2}(K_2)=k_2^{\fatf}\dd t
\dvol_{\phi},\quad N_{H,\fatff,2}(K_2)=k_2^{\fatff}\dd \tau
\dvol_{\phindx/Y},\quad N_{H,\fahbf,-v}(K_2)=g_0.$$ This gives a
parametrix, which solves the heat equation to {\em first} order at
all the faces
 $$ \left( \frac{\partial}{\partial t}+\Dir^2\right)*
K_2= I+R_2 ,\qquad
   R_2\in \Psi^{1,1,1}_{H,{\rm cl}}(X).$$
\par\smallskip
In order to further improve our parametrix at $\fatf$, we have to
successively find $K_{l}\in \Psi^{l,3,\infty}_{H,{\rm cl}}(X)$,
$l=3,4,\ldots$ with $$N_{\fatf,l-1}\left(
(\frac{\partial}{\partial t}+\Dir^2)* K_{l+1}
\right)=-N_{\fatf,l-1}(R_{l}) \in (x')^{2-l}{\cal S}(\phitx,\dd
t\dvol_{\phi}) ,$$ and $$  (\frac{\partial}{\partial
t}+\Dir^2)(K_2+\ldots+K_{l}) = I+R_{l},\qquad R_{l}\in
\Psi_{H,{\rm cl}}^{l-1,1,\infty}(X),$$ Again writing $N_{H,\fatf,
l}(K_l)=k_l\dd t\dvol_{\phitx/X}$ and so on, we have $$
(\frac{1}{2}({\bf \xi}\frac{\partial}{\partial \bf\xi})
+\Betrag{\xi}^2+l/2-1/2+\Lambda)\widehat{k_{l+1}}=\widehat{r_{l}}\in
(x')^{2-l} {\cal S}(\phitx).$$ That this can be uniquely solved is
the content of the following
\begin{lemma}
Let $f\in{\cal S}(\RR)$ and $m>0$. Then there is a unique solution
$u\in {\cal S}(\RR)$ with
$$((s\frac{\partial}{\partial s})+2s^2+m)u=f.$$
\end{lemma}
\begin{proof}
First, the homogeneous equation
$$ ((s\frac{\partial}{\partial s})+2s^2+m)u_0=0$$
is uniquely solved by $u_0=e^{-s^2}s^{-m}$. Writing $u=u_0u_1$ we see that
$$ f= ((s\frac{\partial}{\partial s})+2s^2+m)u_0u_1=
   u_0 s\frac{\partial}{\partial s} u_1.$$
Thus $u_1$ solves the equation
$$\frac{\partial}{\partial s}u_1 =(s u_0)^{-1}f, $$
i.e. $u$ is of the form
$$ u = e^{-s^2}s^{-m}\int_0^s f(r) e^{r^2}r^{m-1}\dd r, $$
which is obviously in ${\cal S}(\RR)$.\qed
\end{proof}
This also finishes the proof of the Theorem.\qed
\end{proof}
\par
\smallskip
 This result can be further improved  by formally inverting the
term $I+R$. From the composition formula, Theorem
\ref{heatcomposition}, we get: $$ R^j \in \Psi_H^{{\cal
G}_j}(X),\quad G_{\fahbf,j}\geq j, \quad G_{\fatff}= j+\NN, \quad
G_{\fatf,j}=\infty .$$
 Thus the powers $R^j$ can be asymptotically
summed, and the Neumann series $$(1+R)^{-1}=\sum_{k=0}^{\infty}
(-R)^k=1+\widetilde{R}$$ is valid in that sense. We can then
define
\begin{eqnarray}
K_{\rm sm}&:=&K(1+R)^{-1}=K+K\widetilde{R}\in \Psi_H^{{\cal H
}}(X),\quad \mbox{with}\nonumber
\\ & &\quad H_{\fatf}=2+\NN,\quad
H_{\fatff}=2+\NN,\quad H_{\fahbf} \subset \NN\times\NN,\quad
T(H_{\fahbf})=\{0\}. \label{heatindexset}
\end{eqnarray} This
gives a parametrix of the heat equation up to a smoothing
remainder $$(\frac{\partial}{\partial t}+\Dir^2)K_{\rm sm}=
1+R_{\rm sm}, \qquad R_{\rm sm} \in \cdotinfty{}(X^2_H,\KD_H). $$
Now, as in the final step of the resolvent construction, the
operator $1+R_{\rm sm}$ is a {\em Volterra} operator, which can be
inverted and whose inverse $(1+R_{\rm sm})^{-1}$ is again of that
type. We have therefore shown \nopagebreak[4]
\begin{theorem}[Heat Kernel]\label{mainresult2}
 The heat kernel for $\Dir^2$ is  given by
$$ {\sf K} = K_{\rm sm}(1+R_{\rm sm})^{-1} \in\Psi_H^{{\cal H
}}(X,E),\qquad {\cal H} \quad\mbox{as in}\quad
(\ref{heatindexset}).$$ Its normal operators at the different
faces are (well defined and) given by
\begin{enumerate}
\item $N_{H,\fatf,2}({\sf K})=
\frac{1}{\sqrt{4\pi}^n}\exp(-\Norm{.}_{\phi}^2/4)\dvol_{\phitx/X}\dd
t$
\item $N_{H,\fatff,2}({\sf K})=
\frac{1}{\sqrt{4\pi\tau}^{h+1}}\exp(-\Norm{.}^2_{\phindx}/{4\tau})\dvol_{\phindx/\partial
X } ([\exp(-\tau(\Dir^{\phi,V})^2)]$
\item $N_{H,\fahbf,0}({\sf K})=[\Pi_{\circ}][\exp(-t\Dir_Y^2)]\frac{1}{\sqrt{4\pi
t}}e^{-(\ln(s))^2/4t}\frac{\dd s}{s}$
\end{enumerate}
\qed
\end{theorem}
\par
\begin{remark}\rm 
There are several possible extensions of this result which are not
presented here, since they are not crucial to our argument. First,
one can eliminate half of the coefficients in the expansion of the
heat kernel at $\fatf$, by analyzing their behavior under the
``flip map'' in $X^2_H$. This procedure is described in Chapter 7
of \cite{Melaps}. Also, direct analysis of the model problem at
$\fahbf$ shows that the heat kernel expansion does not exhibit
log-terms at that face. Thus really, $H_{\fahbf}=\NN$
\end{remark}
\par
Still, Theorem \ref{mainresult2} gives a rather crude description
of the expansions of the heat kernel at $\fatf$ and $\fatff$.
Especially, it does not yet allow us to calculate the regularized
heat supertrace at small times. For this, we have to introduce a
refinement of the heat calculus, which keeps track of the Clifford
degree of the coefficients in the expansions. This will be done in
the next Chapter.

\newpage
\section{Rescaled d-Heat Calculus}\label{rescaledchapter}
In this Chapter we present a refinement of the $d$-heat calculus
presented in Chapter \ref{heatchapter}. This ``rescaled'' calculus
encodes the behavior of the Taylor coefficients of the heat kernel
at $\fatf$ and $\fatff$ w.r.t. the Clifford filtration of the
endomorphism bundle at these faces. This generalization of
Getzler's technique (compare \cite{BGV}) was introduced in
\cite{Melaps}.  A good description of this philosophy can also be found in
\cite{Daimel}.
\par
While the results in earlier Sections are independent of the
parity of the dimension $n$ of $X$ we now assume that $n=v+h+1$ is
{\em even}.
\par
\subsection{The Double Rescaling}\label{basicresc}
In this Section we show how the natural filtration by Clifford
degree of the bundle $\END(E)$ over the faces $\fatf$ and $\fatff$
can be used to define a new bundle ${\rm Gr}(\END)$ {\em rescaled}
at these faces.  First, we have to describe the Clifford
action on $\END(E)$ in more detail.
\par
\smallskip
 Elements $\alpha \in
\Cl({}^dT^*X)$ act on $\beta_L^*E$ via Clifford left
multiplication and on $\beta_R^*E^*$ via Clifford right
multiplication $$ \cl_{d,L}(\alpha)\beta_L^*\xi=
\beta_L^*\cl_d(\alpha)\xi,\qquad
\cl_{d,R}(\alpha)\beta_R^*\xi^*=\beta_R^*\xi^*\circ
\cl_d(\alpha).$$ The right action can be turned into a left
action, using the transpose map $(\cdot)^{\dag}$ on
$\Cl({}^dT^*X)$: $$\overline{\cl}_{d,R}(\alpha)\beta_R^*\xi^*
=(-1)^{|\alpha||\xi^*|}\beta_R^*\xi^*\circ\cl_d(\alpha)^{\dag}.$$
\par
 Since the fibres of ${}^dT^*X$ are even dimensional, the
Clifford modules $E$ and $E^*$, as well as $\END(E)$ are naturally
$\ZZ_2$-graded by the action of the volume element in
$\Cl({}^dT^*X)$. The idea of introducing a rescaling into our heat
calculus is based on the observation that the bundle $\END(E)$
over the heat space $X^2_H$ carries a filtration when restricted
to the different faces.
\begin{lemma}\label{filtration} At $\fatf$, $\fatff$ and $\fahbf$ the bundle $\END$
has the following form
\begin{enumerate}
\item $\END(E)|_{\fahbf}=\Cl( \frac{\dd x}{x}
)\otimes \END_{\Cl(\frac{\dd x}{x})}(E)$
\item $\END(E)|_{\fatff}= \phi^*\Cl({}^{\rm b}T^*B|_Y)\otimes
       \END_{\Cl({}^{\rm b}T^*B|_Y)}(E)$.
\item $\END(E)|_{\fatf}=\Cl({}^{\rm d}T^*X)\otimes \End_{\Cl({}^{\rm d
}T^*X)}(E)$
\end{enumerate}
The restrictions of the connection $\nabla^{\END,d}$ split
accordingly. Note that (b) is true independent of the parity of
the dimension $h$ of the base $Y$. Also, part (a) will be of
no importance to us. We state it because its proof is exemplary
for the proofs of (b) and (c).
\end{lemma}
\begin{proof}
  First, since the bundle $\END$ is just pulled back
from $X^2\times [0,\infty[$, it suffices to prove the above
statements in the blown down picture.
\par To prove (a), restrict to the face $\fahbf$. There,
we can assume the coordinates $x$, $x'$ to be {\em paired}
(compare Section \ref{bandphisection}). For a section
$Q\in\Gamma(X^2_H,\END)$, homogeneous w.r.t. to the grading of
$\END$, we can therefore define
 $$[\cl_d(\frac{\dd x}{x}),Q]:=\cl_{d}(\frac{\dd x}{x})\circ Q
- (-1)^{|Q|}Q\circ \cl_{d}(\frac{\dd x'}{x'}).$$ Then
$$\left[\cl_d(\frac{\dd x}{x}),[\cl_d(\frac{\dd
x}{x}),Q]\right]=0,\qquad \left[\cl_d(\frac{\dd x}{x}),
Q+\frac{1}{2}\cl_d(\frac{\dd x}{x})\circ [\cl_d(\frac{\dd
x}{x}),Q]\right]=0,$$ which shows that the map
 $$Q\longmapsto \left(-\frac{1}{2}\cl_d(\frac{\dd x}{x})\circ [\cl_d(\frac{\dd x}{x}),Q]
 \right)+\left(Q+ \frac{1}{2}\cl_d(\frac{\dd x}{x})\circ [\cl_d(\frac{\dd
 x}{x}),Q]\right)$$
has image in the RHS of (a), and is an isomorphism of algebras.
This proves (a).
\par
Parts (b) and (c) can be obtained by iteration of this scheme
using an orthonormal basis in $\Cl({}^bT^*B)$ and $\Cl({}^dT^*X)$
respectively. The important observation is that an element
$\alpha$ in $\Gamma(Y,{}^bT^*B)$ can be lifted from the left or
the right to $\fatff$, and we can define as above
 $$[\cl_d(\alpha),Q]:=\cl_{d}(\beta_{H,L}^*\alpha)\circ Q
- (-1)^{|Q|}Q\circ \cl_{d}(\beta_{H,R}^*\alpha),$$ and the same
can be done for {\em any} $\alpha\in\Gamma({}^dT^*X)$ at
$\fatf$.\qed
\end{proof}
\par\smallskip
To rescale the bundle $\End(E)$ we want to continue the filtrations at
$\fatf$ and $\fatff$ to the interior using parallel transport
along a normal vector field, and define the space of
sections of our rescaled bundle as the set of sections, whose
$l$th normal derivative has degree at most $l$ at the
corresponding face. This encodes the (so far: supposed) behavior
of the heat kernel at these faces.
\par
To explain the construction, assume for a moment that the bundle
and the connection are trivial. A local basis for the rescaled
bundle near $\fatf\cap\fatff$ could then look like
\begin{equation}\label{localbasis}
 \rho_{\fatf}\,\rho_{\fatff}\,\frac{\dd x}{x},\quad
  \rho_{\fatf}\,\rho_{\fatff}\,\dd {\bf y},\quad
  \rho_{\fatf}\, x\dd {\bf z} \,(!),
\end{equation}
  and (wedge-) products thereof. For simplicity of notation,
   we have also dropped the ``coefficients''
   $\END_{\Cl({}^dT^*X)}(E)$, i.e. we have taken $E$ to be the
   spinor bundle.
\par
In reality however, neither our bundle nor the connection are
trivial. Also, we would like to avoid the choice of a normal
vector field, at least for the mere definition of the rescaled
bundle.
 The relevant connection on $\END(E)$ in
our context is the lifted connection
$$\nabla=\beta_{H,L}^*\nabla^{E,d}\otimes\beta_{H,R}^*\nabla^{E^*,d}\otimes
\beta_H^*\dd t\frac{\partial}{\partial t}: \Gamma(\END(E))
\longrightarrow \Gamma(T^*X^2_H\otimes\END(E)),$$ and the
associated curvature will be denoted by $K$. Also, recall that the
identification of the faces $\fatf$ and $\fatff$ was based on the
following identifications for $A\in \Gamma(\phitx)$:
\begin{eqnarray}
\tau^{1/2}\beta^*_LA|_{\fatf}\in T\fatf/X &\mbox{with}& A\in
\phitx \subset T\phitx/X\label{tfidrep}
\\ \beta^*_LA|_{\fatff}\in
T\fatff/\partial X&\mbox{with}& N_{\phi}(A)\in
T\phindx/Y.\label{tffidrep}
\end{eqnarray}
\par
\bigskip
\noindent{\bf Rescaling at tff}
\par
\medskip
We agree to write
\begin{eqnarray*}
 \Cl^k_{\fatff}&\sowie& \phi^*\Cl^k({}^{\rm b}T^*B|_Y)\otimes
       \END_{\Cl({}^{\rm b}T^*B|_Y)}(E)\quad\mbox{at}\quad\fatff
 \end{eqnarray*}
for the filtration space in Lemma \ref{filtration}. Our
construction then starts from the following simple observations
(not all of which are going to be used right away):
\begin{lemma} [Basic Properties of the Curvature and the Filtration at {\boldmath $\fatff$}]
\label{curvprop}{\quad}
\begin{enumerate}
\item Let $A$ be a vector tangent to $\fatff$. Then $\nabla_A$
preserves the filtration of $\END$ at $\fatff$.
\item $K(T_{\fatff}, \cdot)|_{\fatff}\in
N^*\fatff\otimes\Cl^1_{\fatff}$ \quad whenever $T_{\fatff}$ is
tangent to the fibres $\fatff/ Y$.
\item $(\nabla_{N}K)(T_{\fatff},U_{\fatff})\in
\Cl^1_{\fatff}$, for $T_{\fatff}, U_{\fatff}$ tangent to the
fibres  $\fatff/Y$.
\item The covariant derivatives of the curvature $\nabla^l K$ are
 at most of order 2 at $\fatf$ and $\fatff$.
\end{enumerate}
\end{lemma}
\begin{proof}
Since the connection $\nabla$ is the connection on $\END$ pulled
back from $X^2$, it suffices to prove (a) for vectors $A$, tangent
to the fibre diagonal $F_{\phi}$. In view of the proof of Lemma
\ref{filtration} and the fact, shown in Lemma
\ref{exactconnection}(c) and (\ref{gnice}), that $\nabla^d_A$
preserves the decomposition ${}^dT^*X={}^dN^*\partial X\oplus
{}^dV^*\partial X$ over the boundary, the assertion is then clear.
\par
For part (b) to (d), look at the decomposition of the curvature
(and see the remark at the end of Appendix \ref{cliffordconv})
\begin{eqnarray}
K&=&\beta_{L}^*(R^{S,d}+F^{E/S,d})+\beta_{R}^*(R^{S^*,d}+F^{E/S,d})\nonumber\\
  & & \qquad\qquad\qquad =\frac{1}{2}\beta_{L}^*\cl_{d}(R^d)
    -\frac{1}{2}\beta_{R}^*\cl_d(R^d)
   +\beta_L^*F^{E/S,d} +\beta_R^*F^{E/S,d}.\label{curvature}
\end{eqnarray}
Parts (b)-(d) then follow immediately from Proposition
\ref{rdatthedel}.\qed
\end{proof}
\par
As a first
step in the rescaling procedure, we can now define the space of
sections, {\em rescaled} at $\fatff$ as
 $$
{\cal D}_{\fatff}:=\left\{ u\in \Gamma(X^2_H,\END) |\,
  \nabla^k u|_{\fatff}\in (T^*X^2_H)^k\otimes \Cl^k_{\fatff} \right\}.$$
First note that, although it looks like this definition uses a
connection
  on $TX^2_H$, a  second look convinces us that no such choice is
  necessary. In fact, by choosing any vector field $N_{\fatff}$
  (strictly) normal to $\fatff$, and using the compatibility
  of the connection and the filtration at $\fatff$ as stated in
  Lemma \ref{curvprop}(a) and (d), we could also write
\begin{equation}\label{Dtffwithnormal}
{\cal D}_{\fatff}=\left\{ u\in \Gamma(X^2_H,\END) |\quad
  (\nabla_{N_{\fatff}})^k u|_{\fatff}\in \Cl^k_{\fatff} \right\}.
 \end{equation}
As in \cite{Melaps} this is the space of sections of a vector
bundle ${\rm Gr}_{\fatff}(\END)$. We note this as
\begin{proposition}\label{firststep}
\begin{enumerate}
\item ${\cal D}_{\fatff}=\Gamma(X^2_H,{\rm Gr}_{\fatff}(\END))$.
\item ${\rm Gr}_{\fatff}(\END)|_{\fatf}=\Cl([{}^dV^*\partial X])
\otimes \End_{\Cl({}^{\rm
d}T^*X)}(E)= \Cl(x{}^{\rm b}T^*X)\otimes \End_{\Cl({}^{\rm
d}T^*X)}(E)$
\end{enumerate}
\end{proposition}
\begin{proof}
(a) is proved as in Chapter 7 of \cite{Melaps} by writing down a
local basis for ${\rm Gr}_{\fatff}(\END)$. To prove (b), we first
check that sections in $\Cl([{}^dV^*\partial X])$ satisfy the
condition in the definition (\ref{Dtffwithnormal}) of ${\cal
D}_{\fatff}$. Note that for sections $\alpha_i\in
\Gamma([{}^dV^*\partial X])$ the derivative $\nabla_N^d\alpha_i$
still lies in ${}^dT^*X$. Thus $$
(\nabla_N)^k(\alpha_1\wedge\ldots\wedge\alpha_l)|_{\partial
X}\in \Lambda^{k\wedge l}\,{}^dT^*X\otimes\Lambda^{0\vee
l-k}[{}^dV^*\partial X],$$ and $(\nabla_N)^k(\cl_d(\alpha_1)\ldots\cl_d(\alpha_l)$
 is of order at most $k$ at $\fatff$. The
proof of the reverse inclusion proceeds by induction and is left
to the reader.\qed
\end{proof}
Let us describe the structure of the bundle ${\rm
Gr}_{\fatff}(\END)$ in more detail. To see what the bundle looks
like at $\fatff$, note that the space of sections  over $\fatff$
is $${\cal D}_{\fatff}/\rho_{\fatff}{\cal D}_{\fatff},$$ and the
``denominator'' space can be described as follows
\begin{lemma}
$\rho_{\fatff}{\cal D}_{\fatff}=\{u\in{\cal D}_{\fatff}\,|\,
\nabla^ku|_{\fatff}\in (T^*X^2_H|_{\fatff})^k\otimes
\Cl_{\fatff}^{k-1}\}$.\qed
\end{lemma}
Thus {\em choosing} a projection $P_{\fatff}$ of
$TX^2_H|_{\fatff}$ onto a subspace (or, if one prefers, a vector
field $\nu_{\fatff}=P_{\fatff}/\dd x$) normal to $\fatff$, the map
$$u\longmapsto u|_{\fatff,{\rm Gr}}=(u|_{\fatff},P_{\fatff}\nabla
u|_{\fatff},\ldots,\frac{1}{(h+1)!}P_{\fatff}(\nabla^{h+1}
u)|_{\fatff}) $$ yields an exact sequence $$0 \longrightarrow
{\cal D}_{\fatff}/\rho_{\fatff}\,{\cal D}_{\fatff}\longrightarrow
{\cal D}_{\fatff}\stackrel{|_{\fatff,{\rm Gr}}}{\longrightarrow}
C^{\infty}(\fatff,\mbox{$\bigoplus$}_k
(N^*\fatff)^k\otimes\Cl^k_{\fatff}/\Cl^{k-1}_{\fatff})\longrightarrow
0.$$ The vector bundle appearing on the right represents the
restriction of ${\rm Gr}_{\fatff}(\END)$ to $\fatff$. It will be
denoted by $${\cal N} \Lambda(\fatff)=\mbox{$\bigoplus$}_k{\cal N
}\Lambda^k(\fatff).$$ However, the reader should keep in mind that
the restriction {\em map} depends on the choice of projection
$P_{\fatff}$.
\par
We have chosen to keep track of the factors $N^*\fatff$ for a
while (and to remind the reader of this fact by introducing
``${\cal N}$'' into the notation as above) because this allows for
an easy description of the induced connections: Denote by
${}^bT_{\fatff}X^2_H$ the vector bundle whose sections are vector
fields over $X^2_H$ which are tangent to $\fatff$. Also
${}^bT_{\fatff/Y}X^2_H$ is defined as the vector bundle of vector
fields tangent to the {\em fibres} of $\fatff\rightarrow Y$. On
$N\fatff$ we will use the (b-) connection
\begin{equation}\label{definennabla}
\nabla^{N\fatff}:\Gamma(\fatff,N\fatff)\longrightarrow
\Gamma(\fatff,{}^bT_{\fatff}^*X^2_H|_{\fatff}\otimes N\fatff)\quad
\nabla^{N\fatff}_TN:=[T,N].
\end{equation}
The reader should note that this is only a b-connection, since for
instance, $[xN',N]=-(Nx)N'$. However, we can  restrict it to an
action of $T\fatff/Y$ on $N\fatff$ in a reasonably natural way, by
using the inclusion map $$[]:T\fatff/Y\longrightarrow
{}^bT_{\fatff}^*X^2_H|_{\fatff}, \quad T\longmapsto [T],\quad
\mbox{with}\quad \frac{\dd x}{x}([T])=O(x).$$
\par
We can now describe the connection induced on the rescaled bundle
${\rm Gr}_{\fatff}(\END)$ and its representation under restriction
to the boundary.
\begin{remark}\rm \label{remconventions}
It has been above this author's capabilities to devise a good
notation for the following developments.  A connection $\nabla$
will appear in two distinct functions, first as the connection
appearing in the definition of the space of sections (e.g. ${\cal
D}_{\fatff}$) or of the restriction map (e.g. $|_{\fatff,{\rm
Gr}}$), second as the connection {\em acting} on sections of the
rescaled bundle in question (e.g. ${\rm Gr}_{\fatff}(\END)$). In
the latter case, the connection will usually have to be restricted
to a specific class of vectors (e.g. ${}^bT_{\fatff/Y}X^2_H$) and
will then be written in {\em boldface} $\nnabla$.  As an
additional piece of notation,  we will write ${\mathbb
K}:=[\nnabla,\nabla]$ and ${\mathbf K}:=[\nnabla,\nnabla]$. The
complication with this notation arises from the fact that we have
to rescale {\em twice}, thus the connections change their role in
the process.
\end{remark}
\par
\smallskip\noindent
Using these conventions we get
\begin{proposition}\label{inducedtff}
Let $T, U\in\Gamma({}^bT_{\fatff/Y}X^2_H)$. Then
\begin{enumerate}
\item $\nnabla_T$ maps ${\cal D}_{\fatff}$ to itself and thus induces
a connection $\nnabla_T\sowie\nnabla^{{\rm Gr}(\fatff)}_T$ on
${\rm Gr}_{\fatff}(\END)$ with $$(\nnabla_Tu)|_{\fatff,{\rm Gr}}=
(\nnabla_T-P_{\fatff}{\mathbb
K}(T,\cdot)-\frac{1}{2}P_{\fatff}\nabla_{\cdot}{\mathbb K
}(T,\cdot))(u|_{\fatff,{\rm Gr}})=:\nnabla^{{\cal N
}\Lambda(\fatff)}_T(u|_{\fatff,{\rm Gr}}),$$ where
$P_{\fatff}{\mathbb K}(T,\cdot)\in {\cal N }\Lambda^1(\fatff)$ and
$P_{\fatff}\nabla_{\cdot}{\mathbb K }(T,\cdot)\in {\cal
N}\Lambda^2(\fatff).$
\item ${\bf K}(T,U)\in{\cal D}_{\fatff}$. Its (rescaled)
restriction to $\fatff$ is just the curvature of $\nnabla^{{\cal
N}\Lambda(\fatff)}$: $$ ({\bf K}(T,U)u)|_{\fatff,{\rm Gr}}={\bf
K}(T,U)|_{\fatff,{\rm Gr}}\cdot u|_{\fatff,{\rm Gr}}={\bf
K}^{{\cal N}\Lambda(\fatff)}(T,U)\cdot u|_{\fatff,{\rm Gr}}.$$
\end{enumerate}
\end{proposition}
\begin{proof}
We make abundant use of the following commutator identity:
\begin{lemma}\label{calc1}
Let $T\in\Gamma({}^bT_{\fatff/Y}X^2_H)$, $N\in TX^2_H$. Then over
$X^2_H$
\begin{enumerate}
\item $\nabla_T\nnabla(N)=\nnabla_T\nabla_N-\nabla_{[T,N]}
-{\mathbb K}(T,N)$ ( $=\nnabla\nabla(T,N)-{\mathbb K}(T,N)$
restricted to $\fatff$)
\item $P_{\fatff}\nnabla_T\nabla-\nnabla_TP_{\fatff}\nabla$ is of
order $0$ at $\fatff$
\end{enumerate}
\end{lemma}
\begin{proof}
(a) is just a consequence of our definition of $\nnabla, {\mathbb
K }$ and $\nnabla^{N^*\fatff}$. For the proof of (b) consider $$
 \nnabla_T\nabla(P_{\fatff}N)-\nnabla_TP_{\fatff}\nabla(N)=
 \nabla_{P_{\fatff}[T,N]}-\nabla_{[T,P_{\fatff}N]}$$
This is of order $0$ at $\fatff$, since it is covariant
differentiation w.r.t. a vector tangent to $\fatff$ as can be seen
from $$ P_{\fatff}([T,P_{\fatff}N]-P_{\fatff}[T,N])
   = P_{\fatff}[T,P_{\fatff}N-N]=0.$$ This
finishes the proof of the Lemma.\qed
\end{proof}
\par
\smallskip
\noindent
The proof of part (a) of the Proposition now proceeds as follows:
Let $u\in{\cal D}$. Then using the identies from Lemma \ref{calc1}
 we get modulo lower order at  $\fatff$:
$$\frac{1}{k!}P(\nabla^k)\nnabla_T u=\frac{1}{k!}\nnabla_T P
(\nabla^k)u -P{\mathbb K}(T,\cdot)\frac{1}{(k-1)!}P(\nabla^{k-1})u
-\frac{1}{2}(P\nabla{\mathbb
K}(T,\cdot))\frac{1}{(k-2)!}P(\nabla^{k-2})u.$$ Higher derivatives
of the curvature do not contribute because of Lemma
\ref{curvprop}(d).
\par
To show (b), note that is clear from Lemma \ref{curvprop} that
$${\bf K}(T,U)|_{\fatff}\in\Cl^0_{\fatff},\qquad (\nabla_N {\bf
K})(T,U)\in \Cl^1_{\fatff},$$ all the higher derivatives being of
order at most $2$. Thus ${\bf K}(T,U)$ is in fact an element in
${\cal D}_{\fatff}$. The first equality in (b) is now an immediate
consequence of Leibniz' rule while the second equality can be
obtained by playing around with the definition in part (a):
\begin{eqnarray*}
{\bf K}^{{\cal N}\Lambda(\fatff)}(T,U)\cdot u|_{\fatff,{\rm
Gr}}&=& (\nnabla^{{\cal N}\Lambda(\fatff)}_T\nnabla^{{\cal
N}\Lambda(\fatff)}_U-\nnabla^{{\cal
N}\Lambda(\fatff)}_U\nnabla^{{\cal
N}\Lambda(\fatff)}_T-\nnabla^{{\cal
N}\Lambda(\fatff)}_{[T,U]})\cdot u|_{\fatff,{\rm Gr}}
\\ &=&
[(\nnabla_T\nnabla_U-\nnabla_U\nnabla_T-\nnabla_{[T,U]})u]|_{\fatff,{\rm
Gr}}
\end{eqnarray*}
which is just $({\bf K}(T,U)u)|_{\fatff,{\rm Gr}}$.
 \qed
\end{proof}
\par
\medskip
\noindent{\bf Rescaling at tf}
\par \medskip
 In a second step, we
now want to rescale the bundle ${\rm Gr}_{\fatff}(\END)$ at
$\fatf$ using its induced filtration as described in Proposition
 \ref{firststep}(b)
\begin{eqnarray*}
 \Cl^k_{\fatf} &\sowie& \Cl^k(x{}^bT^*X)\otimes \End_{\Cl({}^{\rm
 d}T^*X)}(E),
\end{eqnarray*}
over $\fatf$ and the induced (partial) connection $\nnabla$
described in Proposition \ref{inducedtff}. To do this, we need the
following addendum to Lemma \ref{curvprop}, which shows that the
filtration and the connection are compatible:
\begin{lemma} [Basic Properties of the Curvature and the Filtration at {\boldmath $\fatf$}]
\label{curvprop2}{\quad}\par\smallskip \noindent Let $A,T_{\fatf},N\in
\Gamma({}^bT_{\fatff/Y}X^2_H)$. Then
\begin{enumerate}
\item Let $B$ be tangent to $\fatf$. Then $\nnabla_B\sowie
\nnabla^{{\rm Gr}(\fatff)}_B$
preserves the filtration of $\,{\rm Gr}_{\fatff}(\END)$ at
$\fatf$.
\item ${\bf K}(T_{\fatf},N)|_{\fatf}=0$ \quad whenever $T_{\fatf}$ is tangent to the
fibres of $\fatf\rightarrow X$.
\end{enumerate}
\end{lemma}
\begin{proof}
Again, since $\nabla^{{\rm Gr}(\fatff)}$ is induced from
$\nabla^{\END}$, away from $\fatff$ it suffices to prove (a) for
vectors $B$, tangent to the diagonal $\bigtriangleup\times{0}$. In view
of the proof of Lemma \ref{filtration} and Lemma
\ref{exactconnection}(c) and (\ref{gnice})) the assertion is then
clear. For (b), just note that $T_{\fatf}$ can be thought of as
being of the form $\tau^{1/2}\beta_L^*V$ for $V\in\Gamma(\phitx)$.
Away from $\fatff$  (\ref{curvature}) then implies $$ {\bf K
}(T_{\fatf},N)={\bf K}(\beta_{L}^*\tau^{1/2}V,N)=
\tau^{1/2}K(\beta_L^*V,N),$$ which is obviously 0 when restricted
to $\fatf$. The assertions extend to $\fatff$ by continuity.\qed
\end{proof}
\par
 Now the space of sections of the bundle rescaled at $\fatff$ and
$\fatf$ can be defined, using a vector field $N_{\fatf}$ normal to
$\fatf$ but tangent to the fibres of $\fatff$:
 $$
{\cal D}\equiv{\cal D}_{\fatff,\fatf}:=\left\{ u\in {\cal
D}_{\fatff} |\,
  (\nabla_{N_{\fatf}}^{{\rm Gr}(\fatff)})^k u|_{\fatf}\in \Cl^k_{\fatf} \right\}.$$
Again, this definition gives the space of sections of a vector
bundle ${\rm Gr}(\END)$ and is independent of the choices made:
\begin{proposition}
\begin{enumerate}
\item ${\cal D}=\Gamma(X^2_H,{\rm Gr}(\END))$.
\item ${\cal D}=\left\{ u\in \Gamma(\END) |\,
  \nabla^k u|_{\fatff}\in (T^*X^2_H)^k\otimes\Cl^k_{\fatff},\quad
  \nabla^l u|_{\fatf}\in (T^*X^2_H)^l\otimes\Cl^l_{\fatf}
   \right\}.$
\end{enumerate}
\end{proposition}
\begin{proof}
(a) follows from iterating the argument in the case of a simple
rescaling, (b) follows from Lemma \ref{curvprop2} and the fact
that $\nabla^{{\rm Gr}(\fatff)}$ on ${\rm Gr }_{\fatff}(\END)$ is
just induced from the connection on $\END$.\qed
\end{proof}
From now on, part (b) in this Proposition will be our preferred
definition of ${\cal D}$.
\par
\smallskip
The rescaling ${\rm Gr}_{\fatf}({\cal N}\Lambda(\fatff))$ of the
bundle ${\cal N}\Lambda_{\fatff}$ at $\fatff\cap\fatf$ can be
defined in the same way, using the connection $\nabla^{{\cal
N}\Lambda(\fatff)}$ and the filtration at $\fatf\cap\fatff$
induced from Lemma \ref{firststep}(b). The description of the
bundle ${\rm Gr}(\END)$ over the faces $\fatf$ and $\fatff$ then
follows the lines set out above: The space of sections over the
faces  $F=\fatf$, or $F=\fatff$ is given by
 $${\cal D}/\rho_F {\cal D},$$
which can be described using
\begin{lemma}
$ \rho_F{\cal D}=\{ u\in{\cal D}|\quad  \nabla^k u|_F\in
 (T^*X^2_H|_F)^k\otimes\Cl^{k-1}_F \}.$ \qed
\end{lemma}
\par
\smallskip
Then {\em choosing}  projections $P=P_F$ of $TX^2_H|_F$ onto a
normal subspace  to $F$ -- of course $P_{\fatf}$ needs to be
tangent to the fibres of $\fatff$ -- the maps $$u\longmapsto
\left\{
\begin{array}{l}
u|_{\fatf,{\rm Gr}}=(u|_{\fatf},P_{\fatf}\nabla
u|_{\fatf},\ldots,\frac{1}{n!}P_{\fatf}(\nabla^nu) |_{\fatf})\\
u|_{\fatff,{\rm Gr}}=(u|_{\fatff},P_{\fatff}\nabla
u|_{\fatff},\ldots,\frac{1}{(h+1)!}P_{\fatff}(\nabla^{h+1}
u)|_{\fatff})
\end{array}\right.$$
give exact sequences $$0 \longrightarrow {\cal
D}/\rho_{\fatff}\,{\cal D}\longrightarrow {\cal
D}\stackrel{|_{\fatff,{\rm Gr}}}{\longrightarrow}
C^{\infty}(\fatff,{\rm Gr}_{\fatf}({\cal N
}\Lambda({\fatff})))\longrightarrow 0,$$ $$0 \longrightarrow {\cal
D}/\rho_{\fatf}\,{\cal D}\longrightarrow {\cal
D}\stackrel{|_{\fatf,{\rm Gr}}}{\longrightarrow}
C^{\infty}(\fatf,{\cal N}\Lambda({\fatf}))\longrightarrow 0,$$
where we have written ${\cal
N}\Lambda(\fatf):=\mbox{$\bigoplus$}_k
(N^*\fatf)^k\otimes\Cl^k_{\fatf}/\Cl^{k-1}_{\fatf}$. Of course,
there are also induced restriction maps
 $$|_{\fatff,{\rm Gr
}}:C^{\infty}(\fatf,{\cal N}\Lambda({\fatf}))\rightarrow
C^{\infty}(\fatf\cap\fatff,{\cal N}{\cal N}\Lambda(\fatf,\fatff)),
$$ $$|_{\fatf,{\rm Gr }}:C^{\infty}(\fatff,{\rm Gr}_{\fatf}({\cal
N }\Lambda({\fatff}))) \rightarrow
C^{\infty}(\fatf\cap\fatff,{\cal N}{\cal N
}\Lambda(\fatf,\fatff)),$$ defined by $$ \quad u|_{\fatff,{\rm
Gr}}=(\frac{1}{j!}P_{\fatff}(\nabla^ju)|_{\fatff})_{j=1,\ldots,h+1}
\quad\mbox{and} \quad u|_{\fatf,{\rm
Gr}}=(\frac{1}{j!}P_{\fatf}(\nabla^{{\cal N }\Lambda(\fatff)})^ju
|_{\fatf})_{j=1\ldots,n}.$$
 The vector bundle on the RHS is ${\cal
N}{\cal N}\Lambda(\fatf,\fatff)=\mbox{$\bigoplus_{m=1}^n$}{\cal N
}{\cal N}\Lambda^m(\fatf,\fatff)$. Its part of grade $m$ is
defined as
\begin{equation}\label{cornerbundle}
(N^*\fatf)^m\otimes\left[\mbox{$\bigoplus_{j+l=m}$}
\Lambda^j{}^dV^*\partial X \otimes
(N^*{\fatff})^l\otimes\Lambda^l\phi^*{}^bT^*B\right]
\otimes\END_{\Cl({}^dT^*X)}(E).
\end{equation}
The restriction maps have been defined in such a way that
restriction to the corner $\fatf\cap\fatff$ is well defined
independent of the order
\begin{equation}\label{restcorner}
u|_{\fatf\cap\fatf,{\rm Gr}}:=u|_{\fatff,{\rm Gr}}|_{\fatf,{\rm
Gr}}=u|_{\fatf,{\rm Gr}}|_{\fatff,{\rm Gr}}.
\end{equation}
Defining the (b-) connection on $N\fatf$  in analogy with
(\ref{definennabla}) the analogue of Proposition \ref{inducedtff}
reads
\begin{proposition}\label{inducedtfftf}
Let $T\in\Gamma({}^bT_{\fatff/Y,\fatf/X}X^2_H)$, i.e. $T$ is
tangent to the fibres  $\fatf/X$ and $\fatff/Y$. Then $\nnabla_T$
maps ${\cal D}$ to itself and thus induces a connection
$\nnabla_T\sowie\nnabla^{{\rm Gr}}_T$ on ${\rm Gr}(\END)$.
Restricted to the faces $\fatf$ and $\fatff$ this has the form
\begin{enumerate}
\item $(\nnabla_Tu)|_{\fatf,{\rm Gr}}= (\nnabla_T-\frac{1}{2}P_{\fatf}\nabla_{\cdot}
{\mathbb K }(T,\cdot)) (u|_{\fatf,{\rm Gr}})=:\nabla_T^{{\cal N
}\Lambda(\fatf)}(u|_{\fatf,{\rm Gr}})$
\item $(\nnabla_Tu)|_{\fatff,{\rm Gr}}= (\nnabla_T-P_{\fatff}{\mathbb
K}(T,\cdot)-\frac{1}{2}P_{\fatff}\nabla_{\cdot}{\mathbb K
}(T,\cdot))(u|_{\fatff,{\rm Gr}})=\nabla_T^{{\cal N
\Lambda}(\fatff)}(u|_{\fatff,{\rm Gr}})$
\item $(\nnabla_Tu)|_{\fatff\cap\fatf{\rm Gr}}
= (\nnabla^{{\cal N}\Lambda(\fatff)}_T
-\frac{1}{2}P_{\fatf}\nabla^{{\cal N
}\Lambda(\fatff)}_{\cdot}{\mathbb K}^{{\cal
N}\Lambda(\fatff)}(T,\cdot))(u|_{\fatff\cap \fatf,{\rm Gr}})$
\end{enumerate}
\end{proposition}
\begin{proof}
As in Proposition \ref{inducedtff} this follows from the
commutator identity in Lemma \ref{calc1} and the properties of the
curvature described in Lemma \ref{curvprop2}. In addition, one has
to verify that the term $\nabla_{\cdot}{\mathbb K }(T,\cdot)$ (is
of order 2 and) preserves $\Cl(x{}^bT^*X)$ at $\fatf$! This will
follow from the more concrete calculation of this term in
Proposition \ref{curvterms} below. Note also that ${\mathbb K
}(T,\cdot)$ {\em vanishes} at $\fatf$. \qed
\end{proof}

\subsection{The Rescaled Calculus and Normal Operators}\label{secrescnormal}
Using the rescaled heat bundle ${\rm Gr}(\END)$ we can introduce
the corresponding ``coefficient bundle'' $$\CB_G:={\rm
Gr}(\END)\otimes\KD_H.$$ It is then straightforward to define the
rescaled heat calculus as $$\Psi^{l,m,p}_{G,{\rm cl}}(X,E):=
 \rho_{\fatf}^l\,\rho_{\fatff}^m\,\rho_{\fahbf}^{p-v}
\cdotinfty{\fatf,\fatff,\fahbf}(X_H^2,\CB_G),$$ for $l>0, m\geq 0,
p\geq -v$. For $l=0$, we set $$ \Psi^{0,m,p}_{G,{\rm cl}}(X,E):=
\rho_{\fatff}^m\,\rho_{\fahbf}^{p-v}
\cdotinfty{\fatf,\fatff,\fahbf}(X_H^2,\CB_G)$$
 and we
remind the reader  that
 $A\in\cdotinfty{\fatf,\fatff,\fahbf}(X_H^2,\CB_G)^{\fatf} $
 satisfies the mean value condition $$\int_{\phitx/X}\Int(t\frac{\partial}{\partial t}) [A|_{\fatf}]^0=0. $$
The notation $[]^0$ is meant to indicate the order $0$-part of the
restriction. As usual, there are definitions for the extended
(rescaled) calculi.
\par
The rescaled heat calculus  is a refinement of the usual heat
calulus in the sense that we have an inclusion of the spaces of
sections $$ \Psi^{l,m,p}_{G,{\rm cl}}(X,E)\subset
\Psi^{l,m,p}_{H,{\rm cl}}(X,E).$$ Therefore the action on sections
and the composition of rescaled operators is well defined. It is
important to know, however, that the composition of rescaled
operators is again a rescaled operator. It is an advantage of our
invariant definition of ${\cal D }$ that it allows for a
relatively easy proof of the following refinement of Theorem
\ref{heatcomposition} for rescaled operators:
\begin{theorem}[Rescaled Composition Formula] \label{rescheatcomp}
Let $R_{\fatf},L_{\fatf}\geq 0$. Then composition gives a map $$
\Psi^{\cal L}_{G}(X) \times\Psi^{\cal R}_{G}(X)
\stackrel{\circ}{\longrightarrow}\Psi^{\cal M}_{G}(X)$$ with index
set given by $$ M_{\fahbf}=(L_{\fahbf}+R_{\fahbf})\overline{\cup}
(L_{\fatff}+R_{\fahbf})\overline{\cup}(L_{\fahbf}+R_{\fatff})$$
$$M_{\fatff}=L_{\fatff}+R_{\fatff},\qquad
M_{\fatf}=L_{\fatf}+R_{\fatf}$$
\end{theorem}
\begin{proof}
We do not give the details of this proof here. Basically, the bundle $\END(E)$
on $X^2_H$ can be lifted to $X^3_H$ from the left and the right, using
$\pi_{H,L}$, $\pi_{H,R}$. Composition of endomorphisms gives a
well defined map
$$\pi_{H,L}^*\END(E) \otimes \pi_{H,R}^*\END(E)\longrightarrow
 \pi_{H,M}^*\END(E).$$
 The connections
used in the definition of the rescaling can be lifted to $X^3_H$
in the same way.
 Using
compatibility of these connections with this product, as well as with
the pullback and pushforward
operations, the proof then boils down to Leibniz' rule. \qed
\end{proof}
\par
\smallskip
Normal operators for our rescaled heat calculus can now be defined
in complete analogy with Section \ref{normalops}. For
$A\in\Psi^{l,m,p}_{G,{\rm cl}}(X)$ define the normal operators of
$A$ at the faces $\fatf$ and $\fatff$  by
\begin{eqnarray*}
N_{G,\fatf,l>0}(A)&:=&
  t^{1-l/2}A|_{\fatf, \rm Gr}\in \rho_{\fatff}^{m-l}
   \cdotinfty{\fatff}(\fatf,t \KD_H\otimes \Lambda_{\fatf})\\
    N_{G,\fatf,0}(A\oplus c)&:=&
  t A|_{\fatf,\rm Gr}\oplus c \in \rho_{\fatff}^{m}
 \cdotinfty{\fatff}(\fatf,t \KD_H\otimes
 \Lambda_{\fatf})^{\fatf}\oplus \CC \\
 N_{G,\fatff,m}(A)&:=&(x')^{-m}A|_{\fatff,\rm Gr}\in
\rho_{\fatf}^{l}\,\rho^{p-v-m}_{\fahbf}\,\cdotinfty{\fatf,\fahbf}
(\fatff,\KD_H\otimes{\rm Gr}_{\fatf}(\Lambda_{\fatff})).
\end{eqnarray*}
\par\smallskip
Our aim is to construct the heat kernel for $\Dir^2$ in the
rescaled calculus. This means that we have to describe
the action of $t\Dir^2$ and $t\frac{\partial}{\partial t}$ on ${\cal D}$ at
the faces $\fatf$ and $\fatff$. To see what this requires,
consider the Lichnerowicz formula for $\Dir^d$
\begin{equation}\label{Lichnerowicz}
 t(\Dir^d)^2
\sowie t\beta_L^*g_d^{-1}(\nnabla,
\nnabla)+t\beta_L^*\frac{\kappa^d_X}{4}
+t\cl_{d,L}(\beta_L^*F^{E/S,d}),
\end{equation} where the notation is
meant to emphasize that we are looking at the operator lifted to
$X^2_H$ from the left. Thus we will have to calculate the normal
operators for \begin{equation}\label{tasks}
 t\frac{\partial}{\partial
t},\quad \nnabla,\quad t\beta_L^*g_d^{-1},\quad
t^{1/2}\cl_{d,L},\quad t\kappa^d_X,\quad
t\cl_{d,L}(\beta_L^*F^{E/S,d}).
\end{equation}
As in Proposition \ref{inducedtfftf} we will show that all of
these operators map ${\cal D}$ to ${\cal D}$ in an appropriate
sense by calculating their commutators with $\nabla^l$.
\par \smallskip
Let us start with the analysis of the action of $\nnabla$.
The main point is to calculate the curvature terms in Proposition
\ref{inducedtfftf}. However, instead of
considering the action of general sections in
${}^bT_{\fatff/Y,\fatf/X}X^2_H$ on ${\cal D}$ we now restrict to
the special class of sections of the type
\begin{equation}\label{calW}
{\cal W} := \langle t^{1/2}\beta_{H,L}^*{}^{\rm d}TX,
t\frac{\partial}{\partial t}\rangle=\langle
\tau^{1/2}\beta_{H,L}^*\phitx, t\frac{\partial}{\partial
t}\rangle.
 \end{equation}
 The essential property of sections in ${\cal W}$ is that they
are tangent to the fibres in $\fatf/X$ and $\fatff/Y$. At
$\fatff$ this is also true for sections in $\beta_{H,L}^*\phitx$,
which is why we will occasionaly  drop the factor $\tau^{1/2}$ on
the RHS in (\ref{calW}). Restricting to the class ${\cal W}$
allows us to make sense of the formulas in Proposition
\ref{inducedtfftf} on all of $X^2_H$ and to follow the conventions
described in Remark \ref{remconventions} more consistently. Thus,
allowing from now on the shorthand $\beta_L\sim\beta_{H,L}$, we
introduce the following notation for our connections:
\begin{itemize}
\item $\nabla$ describes the action of $TX^2_H$ on ${\cal W}$ and
$\END(E)$ and represents the connection used in the definition of
${\cal D}$.
 $$\nabla: \Gamma({\cal W})\rightarrow \Gamma(T^*X^2_H\otimes
       {\cal W}) \mbox{ with }
   \nabla_N  t^{1/2}\beta_{L}^*V=t^{1/2}\beta_{L}^*\nabla^d
      \beta_{L}^*V \quad \mbox{and}\quad
    \nabla_N t\frac{\partial}{\partial t}=0.$$
\item $\nnabla$ describes the action of ${\cal W}$
on $T X^2_H$ and $\END(E)$. It represents the connection used in
the definition of our operators acting on ${\cal D}$.
\begin{eqnarray}
\lefteqn{\nnabla:\Gamma(TX^2_H)\rightarrow \Gamma({\cal W}^*
\otimes TX^2_H)\!\mbox{ with}}\hfill  \nonumber
\\ & & \nnabla_{t\frac{\partial }{\partial
t}}N=[t\frac{\partial}{\partial t},N] \,\mbox{ and
}\,\nnabla_{t^{1/2}\beta_{L}^*V}N=t^{1/2}\,\beta_{L}^*\nabla^d_N\,
\beta_{L}^*V
   +[t^{1/2}\beta_{L}^*V,N]\label{bconnectiona}
\end{eqnarray}
\end{itemize}
It is easy to see that these definitions are compatible with the
definition (\ref{definennabla}) of $\nabla^{N\fatff}$ in the last
Section.
 Note that the ($b$-) connections on $N\fatf$ and
$N\fatff$ are curvature-free by the Jacobi identity. Also in some
cases they are rather easy to calculate:
\begin{lemma}\label{calcntf}
Let $\alpha \in \Gamma(T^*X)$,
$[\phi^*\omega]\in\Gamma([H^*\partial X] )$, let
$W\in\Gamma(\phitx)$ and $A\in\Gamma([\phindx])$. Then
$P_{\fatf}\beta_{L}^*\alpha\in N^*\fatf$,
$P_{\fatff}\beta_{L}^*[\phi^*\omega]\in N^*\fatff$ and
\begin{enumerate}
\item
$\nnabla^{N^*\fatf}_{\tau^{1/2}\beta_{L}^*W}(P_{\fatf}\beta_{L}^*\alpha)
=\alpha(W)\dd \tau^{1/2}$, and
$\nnabla^{N^*\fatf}_{t\frac{\partial}{\partial
t}}(P_{\fatf}\beta_{L}^*\alpha)=0$
 \item
$\nnabla^{N^*\fatff}_{\beta_{L}^*A}(P_{\fatff}\beta_{L}^*[\phi^*\omega])
=\phi^*\omega(A/x)\dd x$, and
$\nnabla^{N^*\fatff}_{t\frac{\partial}{\partial
t}}(P_{\fatff}\beta_L^*\alpha)=0$
\item $\nnabla^{N^*\fatf}_{\tau^{1/2}\beta_{L}^*W}\dd \tau^{1/2}=0$
and $\nnabla^{N^*\fatf}_{t\frac{\partial}{\partial t}}\dd
\tau^{1/2}=\frac{1}{2}\dd \tau^{1/2} $
\item $\nnabla^{N^*\fatff}_{\beta_{L}^*W}\dd x=0$
and $\nnabla^{N^*\fatff}_{t\frac{\partial}{\partial t}}\dd x=0 $.
\end{enumerate}
\end{lemma}
\begin{proof}
As an example, let us do the first part in (a) (leaving out the
ubiquitous $\beta_L^*$):
\begin{eqnarray*}
\nnabla_{\tau^{1/2}W}\alpha
\,(N)&\stackrel{(\ref{bconnectiona})}{=}&
\tau^{1/2}W\,\alpha(N)-\alpha\left(\tau^{1/2}x
\nabla^d_N\frac{1}{x}W+ [\tau^{1/2}W,N]\right)
\\ &\stackrel{(*)}{=}&
\dd \alpha(\tau^{1/2}W, N)+
N\,\alpha(\tau^{1/2}W)\stackrel{(**)}{=}\dd
\tau^{1/2}(N)\alpha(W),
\end{eqnarray*} where, in $(*)$ and $(**)$, we have used the fact
that $\tau^{1/2}$ vanishes at the face $\fatf$. \qed
\end{proof}
\par\medskip
In order to make the calculation of normal
operators look as invariant as possible, we will consider rescaled
``calculi'' with coefficients in ${\cal W}$ and ${\cal W}^*$, the
definition proceeding just as before.
Also
write  ${\mathbb K}=[\nnabla,\nabla]$, which  is the restriction
of the curvature $K$ to $\Gamma({\cal W}^*\otimes
T^*X^2_H\otimes\END(E))$, and similarly  ${\bf
K}=[\nnabla,\nnabla]$.
We can now restate
Proposition \ref{inducedtfftf} as a refinement of Proposition
\ref{normcalc}:
\begin{proposition}\label{normcalcresc}
The maps $$ \Psi_{G,{\rm cl}}^{l,m,p}(X)\stackrel{
\nnabla}{\longrightarrow}
   \Psi_{G,{\rm cl}}^{l,m,p}(X,{\cal W}^*),\qquad \Psi_{G,{\rm cl}}^{l,m,p}(X)
   \stackrel{\frac{\partial}{\partial t}}{\longrightarrow}
   \Psi_{G,{\rm cl}}^{l-2,m-2,p}(X)$$
   have the following general form when restricted to the different faces
\begin{enumerate}
\item $N_{G,\fatf,l}(\nnabla A)=
(\nnabla-\frac{1}{2}P_{\fatf}\nabla{\mathbb K})
 N_{G,\fatf,l}(A)= \nnabla^{{\cal N}\Lambda(\fatf)} N_{G,\fatf,l}(A)$
\item  $N_{G,\fatff,m}(\nnabla A)=
(\nnabla-P_{\fatff}{\mathbb K}-\frac{1}{2}P_{\fatff}\nabla{\mathbb
K})
 N_{G,\fatff,m}(A)= \nnabla^{{\cal N}\Lambda(\fatff)} N_{G,\fatff,m}(A)$
 \item  $N_{G,\fatf, l-2}(\frac{\partial}{\partial t}A)
 =\left\{
 \begin{array}{ll}
     -\frac{1}{2}( L_{R^{\phi}}+2-l-{\sf N} )N_{G,\fatf,l}(A) & \mbox{for}\quad l>2 \\
     -\frac{1}{2}( L_{R^{\phi}}+2-l -{\sf N})N_{G,\fatf,2}(A)\oplus
    -\int_{\phitx/X }N_{H,\fatf,2}(A) & \mbox{for}\quad l=2\\
 \end{array}   \right.$
 \item $N_{G,\fatff, m-2}(\frac{\partial}{\partial t} A)=
 L_{\frac{\partial}{\partial \tau}}N_{G,\fatff,m}(A)$
\end{enumerate}
Here ${\sf N}$ denotes the number operator in $\Lambda_{\fatf}$.
Note also that $\nnabla^{N^*\fatf}_{R^{\phi}}\dd \tau^{1/2}=0$.
\end{proposition}
\begin{proof}
This is really just a reformulation of Proposition
\ref{inducedtfftf}.
 For (c) and (d) just recall that we now
also have  to consider the action of
$\nnabla_{t\frac{\partial}{\partial t}}$ on $N^*\fatf$ and
$N^*\fatff$, which is $1/2$ on $N^*\fatf$ and $0$ on $N^*\fatff$
according to Lemma \ref{calcntf}.\qed
\end{proof}
\par
\smallskip
 For the calculation of these curvature terms it will be useful to
 list some further properties of the curvature ${\bf K}$ at the
 faces $\fatf$ and $\fatff$. It has been shown in  Lemmas
 \ref{curvprop}(a) and \ref{curvprop2}(a) that
 \begin{equation}\label{KTN}
 K(T_{\fatf,}\,\cdot\,)|_{\fatf}\in
 N^*\fatf\otimes\Cl^2_{\fatf}/\Cl^1_{\fatf}\,\mbox{ and }\,
 K(T_{\fatff},\,\cdot\,)|_{\fatff}\in
 N^*\fatff\otimes\Cl^2_{\fatff}/\Cl^1_{\fatff}
 \end{equation}
 are well defined for vectors $T_{\fatf}$, $T_{\fatff}$ tangent
 to $\fatf$, $\fatff$. Using the identifications
 (\ref{tfidrep},\ref{tffidrep}) we now show a more precise result:
\begin{lemma}[More Properties of the Curvature] \label{morecurvprop}
For the faces $F=\fatf, \fatff$, let $T_F$, $N_F$ be vector fields
on $X^2_H$ tangent respectively normal to $F$. Also, we assume
that $T_F$ are lifted from vector fields on
$X^2_b\times[0,\infty[$ tangent to the blown down faces $F$.
\begin{enumerate}
\item $K(T_{\fatf},N_{\fatf})|_{A\in\phitx}= \dd \tau^{1/2}(N_{\fatf})
K(T_{\fatf}, \beta_L^*A)$ \quad modulo $\Cl_{\fatf}^1$
\item $K(T_{\fatff},N_{\fatff})|_{B\in\phindx}= \dd x(N_{\fatff})
K(T_{\fatff}, \frac{1}{x}\beta_L^*[B])$ \quad modulo
$\Cl^1_{\fatff}$.
\end{enumerate}
Here $[B]$ denotes a choice of extension of $B$ to $[\phindx]$.
\end{lemma}
\begin{proof}
The proofs of (a) and (b) follow the same pattern. To prove (b),
note that
\begin{equation}\label{zerosection}
K(T_{\fatff},N_{\fatff})|_{0\in\phindx}\in\Cl^1_{\fatff}.
\end{equation}
This can be seen by looking at the map $\pi$ flipping sides in
$X^2_H$. Since, by assumption, the vector $T_{\fatff}$ is lifted
from a vector on $X^2_b\times[0,\infty[$, tangent to
$F_{\phi}\times\{0\}$, we can restrict ourselves to the following
situation: $$\pi(N_{\fatff})=N_{\fatff}\,\mbox{ and
}\,\pi(T_{\fatff})= T_{\fatff}+ VF_{\phi}\quad\mbox{ at the
zero-section of }\, \fatff \rightarrow F_{\phi}.$$  Using (b)
above, we get $$K(T_{\fatff},N_{\fatff})=-K(\pi(T_{\fatff}), \pi(
N_{\fatff}))\in -K(T_{\fatff},N_{\fatff})+\Cl^1_{\fatff},$$ which
proves (\ref{zerosection}). It now suffices to calculate the
derivative of $K(T_{\fatff},\,\cdot\,)$ in $\Lambda^2_{\fatff}$
along the fibres of $\fatff$. But from (\ref{curvature}), the
assumption on $T_{\fatff}$ and Corollary \ref{curvboundary}, this
is just a linear combination of sections of the type $$
\beta_L^*[\phi^*\omega]\otimes\beta_L^*\xi,\quad
\beta_R^*[\phi^*\omega]\otimes\beta_R^*\xi,\quad\mbox{with}\quad
[\phi^*\omega]\in\Gamma(X,[H^*\partial X]),\quad \xi \in
\Gamma(X,\End(E)) $$ Thus for $B\in\phindx$ representing a point
in the fibre, we can  use Lemma \ref{calcntf}(b) (and its analogue
for sections lifted from the right) to get $$(\nnabla_{\beta_L^*B}
K(T_{\fatff},\,\cdot\,))(N_{\fatff})= \dd x(N_{\fatff})
K(T_{\fatff},\frac{1}{x}[\beta_L^*B]),$$
 which proves the claim.\qed
\end{proof}
\par
\smallskip
The curvature terms in $N_{G,\fatf}(\nnabla)$,
$N_{G,\fatff}(\nnabla)$ can now be calculated. Readers, who feel
uneasy about the signs or the different factors of $x$ in the
formulas should (hope that the author has made no mistake and)
consult Definition \ref{Sphi} and Appendix \ref{cliffordconv}.
\begin{proposition}\label{curvterms}
\begin{enumerate}
\item Using the identification (\ref{tfidrep}) we have for $A\in\phitx$
 \begin{eqnarray*}
 &&(P_{\fatf}\nabla {\mathbb K})
 (\cdot,\tau^{1/2}\beta_{L}^*A,\cdot)|_{B\in\phitx}=
  (\dd\tau^{1/2})^2\frac{1}{2}\ext_d[\langle x^2R^d(\cdot,\cdot)
  (B+\nu_{\fatf}|_0) ,A\rangle_{\phi}]
 \\&&\qquad\qquad =:(\dd
  \tau^{1/2})^2\left(\frac{1}{2}\ext_d[\langle x^2R^d(\cdot,\cdot) B,A
  \rangle_{\phi}]
  +Q_{\fatf}(A)\right) \in \Lambda^2(\fatf).
  \end{eqnarray*}
   Here
  $\nu_{\fatf}|_0$ is the vector field $P_{\fatf}/\dd{\tau^{1/2}}$ at the zero section
  in $\phitx$.
\item Using the identification of (\ref{tffidrep}) we find for $W\in
\phitx$
\begin{eqnarray*}
&& (P_{\fatff}{\mathbb K})(\beta_{L}^*W,
\cdot)|_{B\in\phindx}=-\frac{1}{2}\dd x
(\ext_d,\cl_d)(\langle x S_{d}(W)\cdot,\cdot\rangle_{\phi})
\\&&\qquad\qquad =-\frac{1}{2}\dd x
\left[(\ext_d,\cl_d)\langle S_{\phi}(W)\cdot,\cdot\rangle_d+\ext(\frac{\dd
x}{x})\cl_d(W/x)\right]\in\Lambda^1(\fatff)
\end{eqnarray*}
covariant constant along the fibres  $\fatff/F_{\phi}$.
\item Let $A\in \Gamma(\phindx)$. Then
\begin{eqnarray*}
&& (P_{\fatff}\nabla {\mathbb
K})(\cdot,\beta_{L}^*A,\cdot)|_{B\in\phindx}
  =(\dd x)^2 \frac{1}{2}\ext_d[\langle \phi^*R^{B,b}(\cdot,\cdot)
  ([B]+x\nu_{\fatff}|_0),[A]\rangle_{\phi}]
\\&&\qquad\qquad =:(\dd x)^2
\left(\frac{1}{2}\ext_d[\langle\phi^*R^{B,b}(\cdot,
\cdot)[B],[A]\rangle_{\phi}]+Q_{\fatff}(A)\right)
  \in \Lambda^2(\fatff).
  \end{eqnarray*}
Here,
  $\nu_{\fatff}|_0$,
  is the vector field $P_{\fatff}/\dd x$ evaluated at the zero section
  in $\fatff\rightarrow F_{\phi}$ and $[A], [B]$ denote choices of extension of
   $A, B$ into $[\phindx]$. The result is clearly independent of that
   choice.
\item
Now let $V\in\Gamma( V\partial X)$. Then
\begin{eqnarray*}
(P_{\fatff}\nabla
{\mathbb K})(\cdot,\beta_{L}^*V,\cdot)_{B\in\phindx}
   &=& \frac{1}{4}(\dd x)^2(\ext_{d},\ext_d)
   \langle x^2\Omega_{\phi}(\cdot,\cdot),V\rangle_{\phi}
 \\   &=&\frac{1}{4}(\dd x)^2(\ext_d,\ext_d)\langle \Omega_{M}(\cdot,\cdot)
 ,V\rangle_{M,b}]
    \in \Lambda^2(\fatff).
 \end{eqnarray*}
   Again, this term is covariant constant along the fibres $\fatff/F_{\phi}$.
\end{enumerate}\par
\smallskip\noindent
Note that, except for the term in (b), none of these terms contain
a $\ext(\frac{\dd x}{x})$-contribution.
\end{proposition}
\begin{proof}
The proof of (a) at $\fatf$ is easy. First, recall again that $$
{\mathbb K}(\tau^{1/2}\beta_L^*A,N)=\tau^{1/2}K(\beta_L^*A,N)$$
vanishes at $\fatf$ since $\tau^{1/2}$ vanishes there. Thus $$
(P_{\fatf}\nabla {\mathbb K})(N,\tau^{1/2}\beta_{L}^*A,N) = \dd
\tau^{1/2}(N) K(\beta_{L}^*A,P_{\fatf}N),$$ and we can calculate
the derivative along the fibre as in the proof of Lemmas
\ref{curvprop2}(a) and \ref{morecurvprop}. In this case, however,
the constant term along the fibres does not vanish, accounting for
the term $Q_{\fatff}(A)$ in the formula. Now using
(\ref{curvature}) and writing $$ g_{d}(R^d(A,B)\cdot, \cdot)
=g_{\phi}(x^2R^d(\cdot,\cdot)B,A),$$ gives the result.
\par
 It remains to show that the RHS in (a) really is
an element in $\Cl_{\fatf}^2$. This should be clear for the linear
term in the fibres by Proposition \ref{inducedtff}(c). It is also
true for $Q_{\fatf}(A)$, since $\nu_{\fatf}$ is tangent to the
fibres  $\fatff/Y$, and can therefore be written in the form
$\nu_{\fatf}=\beta_L^*C+xN$ for $C\in \Gamma(\phitx), N\in
\Gamma(TX)$. Then
$$Q_{\fatf}(A)=K(\beta_L^*A,(\beta_L^*C+xN)|_{0\in\phitx}) =
xK(\beta_L^*A,N)+\Cl_{\fatf}^2,$$ but the order-2 part of
$K(\beta_L^*A,N)$ is seen to be in $x\Lambda^2{}^bT^*X$ in part
(b)!
\par\
Part (b) follows from (\ref{curvature}) and Proposition
\ref{rdatthedel} (a).
The proof of (c) is completely analogous to (a) once we restrict
ourselves to $A,B\in \phindx$ with extensions $[A], [B]$ in
$[\phindx]$ as stated. Then again
\begin{equation}\label{firstN}
(P_{\fatff}\nabla {\mathbb
 K})(N,\beta_{L}^*A,N)|_{\fatff} =
  \dd x(N) K(\frac{1}{x}\beta_{L}^*[A],P_{\fatff}N),
\end{equation}
and the result follows by calculating the derivative along the
fibres as in  Lemma \ref{curvprop2}(d) and Lemma \ref{calcntf}.
Finally (d): For $V\in V\partial X$, $ A_1, A_2 \in \phindx$ the
term at $\fatff$,
 $$ \nabla _N{\mathbb K}(N,\beta_{L}^*V)=
\frac{1}{2}\beta_L^*\cl_d((\nabla^d_NR^d)(V,N)),$$
was calculated in Proposition \ref{rdatthedel} (b).
\qed
\end{proof}
\par
\medskip
Having described the role of the conormal bundles $N^*\fatf$ and
$N^*\fatff$, and used them in our calculations, we allow ourselves
to assume from now on that $$ \mbox{{\bf\boldmath
$N^*\fatf$ is trivialized by $\dd \tau^{1/2}$ and $N^*\fatff$ is
trivialized by $\dd x$ ,}} $$
to simplify our formulas a bit.
 Then, recalling Lemma
\ref{liftface}, the rescaled bundle $\fatf$
looks like $$\Lambda(\fatf)=\Lambda(x{}^bT^*X)\otimes
\END_{\Cl({}^dT^*X)}(E)\longrightarrow \phitx.$$ Also, $\nnabla$
restricts to the fibres of $\fatf\rightarrow X$ as the trivial
connection $(\dd_{\END},\dd_{\phitx/X})$. The rescaled bundle over
the {\em interior} of $\fatff$ (i.e.  without the second
rescaling) looks like $$
\Lambda(\fatff)=\Lambda(\phi^*{}^bT^*B)\otimes
\END_{\Cl(\phi^*{}^bT^*X)}(E)\longrightarrow
\phindx\times_Y\partial X\times]0,\infty[_{\tau^{1/2}}.$$ The
induced connection on the fibres  $\fatff/ Y$ is trivial along the
fibres of $p:\phindx\times]0,\infty[_{\tau^{1/2}}\rightarrow
\partial X$. One could write it  as $(p^*\nabla^{\END,d}_{{\sf v}\cdot},
p^*\nabla^{\phindx}_{{\sf v}\cdot})$.
\begin{corollary}[Explicit Form of the Induced Rescaled
Connection]\label{explicitconnection}{\quad}\par
\smallskip\noindent
 Denote by ${\bf x}$, ${\bf y}$ the linear
variables in the fibres of $\phitx$ and $\phindx$ respectively.
Then for $A\in \phitx$, $B\in\phindx$ and  $W\in V\partial X$
\begin{enumerate}
\item $\nnabla^{\Lambda(\fatf)}_A=\Int(A)\dd_{\END} -
\frac{1}{4}\ext[\langle x^2R^d(\cdot,\cdot) {\bf
x},A\rangle_{\phi}]-\frac{1}{2}Q_{\fatf}(A))$
\item  $\nnabla^{\Lambda(\fatff)}_{W}= \nnabla^{\END,d}_{{\sf v}W}
-\frac{1}{2} (\ext_d, {\sf v}\cl_{d})\langle S_{\phi}(W)\cdot,\cdot\rangle_d-
 \frac{1}{4}\ext_d(x^2\langle\Omega_{\phi}(\cdot,\cdot),W\rangle_{\phi})+\frac{1}{2}\ext_d(\frac{\dd
x}{x})\cl_d(\frac{1}{x}{\sf v}W)$
\par\qquad\qquad\qquad\qquad\qquad\qquad$=:\widehat{\nnabla}^{\Lambda(\fatff)}_{W}
  +\frac{1}{2}\ext_d(\frac{\dd x}{x})\cl_d(\frac{1}{x}{\sf v}W)$
 \item $\nnabla^{\Lambda(\fatff)}_{B}=
 \Int(B) \dd_{\END} - \frac{1}{4}\ext_d[\langle \phi^*R^{B,b}(\cdot,\cdot){\bf
y},B\rangle_{\phi}]-\frac{1}{2}Q_{\fatff}(B))$.\qed
\end{enumerate}
\end{corollary}
\par
\smallskip
Let us now describe the action of the Clifford action and of the
metric in (\ref{tasks}).
\begin{proposition}
The metric gives a map $$ \Psi_{G,{\rm cl}}^{l,m,p}(X)\stackrel{
t\beta_L^*g_d^{-1}\,}
   {\longrightarrow}\Psi_{G,{\rm cl}}^{l,m,p}(X,{\cal W}^*\otimes{\cal W}^*)$$
such that, in view of (\ref{tfidrep}, \ref{tffidrep}),
$N_{G,\fatf}(tg_d^{-1})=g^{-1}_{\phitx/X}$ and
 $N_{G,\fatff}(tg_d^{-1})=\tau g^{-1}_{\phindx/Y}.$
\qed
 \end{proposition}
To describe the Clifford action at $\fatf$ and $\fatff$, we use
the partial Clifford action ${\sf m}_d$ on $\Lambda(\fatff)$,
given for $\alpha\in {}^dT^*X$ by $${\sf m}_d(\alpha):=\cl_d({\sf
v}\alpha)+\ext_d({\sf n}\alpha).$$ Also, we recall that the map
$\delta_h^*:x{}^bT^*X|_{\partial X}\rightarrow{}^dT^*X|_{\partial
X}$ was described in Lemma \ref{deltah}.
\begin{proposition}\label{cliffordact}
 The (left) Clifford action of an element
 $x\alpha\in\Gamma(x{}^bT^*X)=\Gamma([{}^dV^*\partial X])$
  gives a map
  $$ \Psi_{G,{\rm cl}}^{l,m,p}(X)\stackrel{ \tau^{1/2}\beta_L^*\cl_d(x\alpha)\,}
   {\longrightarrow}\Psi_{G,{\rm cl}}^{l,m,p}(X). $$
It restricts to the faces as
\begin{enumerate}
 \item $N_{G,\fatf,l}(\tau^{1/2}\beta_L^*\cl_d(x\alpha)A)=\dd
 \tau^{1/2}\ext(x\alpha)\circ N_{G,\fatf,l}(A)\sowie \ext(x\alpha)\circ
 N_{G,\fatf,l}(A)$
 \item $N_{G,\fatff,m}(\tau^{1/2}\beta_L^*\cl_d(x\alpha)A)=
 \tau^{1/2}(\cl_d(x{\sf v}\alpha)+\dd x \ext({\sf h}\alpha))\circ
 N_{G,\fatff,m}(A)$
 \\ \rule{0pt}{0pt} \qquad\qquad\hfill $\sowie \tau^{1/2} {\sf
 m}_d(\delta_{h}(x\alpha))\circ N_{G,\fatff,m}(A).$
\end{enumerate}
\end{proposition}
\begin{proof}
Part (a) follows from the fact that
$\tau^{1/2}\beta_L^*\cl_d(x\alpha)$ vanishes at $\fatf$ and $$
\nabla^{\End,d}_N\tau^{1/2}\beta_L^*\cl_d(x\alpha)|_{\fatf} =\dd
\tau^{1/2}(N)\beta_L^*\cl_d(x\alpha)|_{\fatf}+\tau^{1/2}\nabla^{\End,d}
\beta_L^*\cl_d(x\alpha)|_{\fatf},$$ where the last summand again
vanishes at $\fatf$. To prove (b), first note that
$\tau^{1/2}\cl_d(x\alpha)$ restricts as $\tau^{1/2}\cl_d(x{\sf v
}\alpha)$ at $\fatff$. The first order contribution of
$\tau^{1/2}\cl_d(x\alpha)$ at $\fatff$ can be calculated using
Lemma \ref{halpha}
$$\nabla_N^{\End,d}\tau^{1/2}\cl_{d,L}(x\alpha)=
\tau^{1/2}\cl_{d,L}(\nabla_N^dx\alpha)=\tau^{1/2}\dd
x(N)\cl_{d,L}({\sf h}\alpha),$$
 as claimed.\qed
\end{proof}
\par
A slightly different formulation of the last result will also be
of some use:
\begin{corollary}\label{easyrestrict}
At $\fatf$, the rescaled restriction to the corner
is given by
\begin{enumerate} \item $|_{\fatff,{\rm Gr}}={\sf m}_d\circ
\delta_h:\Cl(x{}^bT^*X)\longrightarrow
 \Cl({}^dV^*\partial X)\otimes \Lambda \phi^*{}^bT^*B|_{\partial X}$
 in  ${\rm Gr}_{\fatff}(\END)$,
\item  $|_{\fatff,{\rm
Gr}}=\delta_h:\Lambda(x{}^bT^*X)\longrightarrow
 \Lambda{}^dT^*X|_{\partial X}$
in ${\rm Gr}(\END)$.\qed
\end{enumerate}
\end{corollary}
\smallskip
We are left with the analysis of the behavior of the scalar
curvature term $t\kappa_X$ in (\ref{tasks}). Now near $\partial X$
we can write, using Corollary \ref{curvboundary}
\begin{eqnarray}
t\kappa_X&=& t\tr_{aa'}(\langle
(g_d^{-1})_{bb'},{}_{a}R^d(\cdot_{a'},\cdot_b)_{b'}\rangle)\nonumber
\\ &=& \tau\tr_{aa'}(\langle (g_{\phi}^{-1})_{bb'},{}_{a}{\sf v
}R^{\phi}(\cdot_{a'},\cdot_b){\sf v}_{b'}\rangle)+ \tau
\tr_{aa'}(\langle
(g_{\phi}^{-1})_{bb'},{}_{a}\phi^*R^{B,b}(\cdot_{a'},\cdot_b)_{b'}\rangle).
\label{minorrole2}
\end{eqnarray}
This is a $C^{\infty}$-function on $X^2_H$, which vanishes at
$\fatf$. At $\fatff$ the second summand on the RHS vanishes, since
the curvature $\phi^*R^{B,b}$ vanishes when evaluated on a vertical
vector. Since the $B_{\phi}$-term in (\ref{minorrole}) also
vanishes on vertical vectors, we are left with $$
t\kappa_X|_{\fatf}=0, \qquad t\kappa_X|_{\fatff}=\tau
\kappa_{M/Y}.$$
Putting everything together, we have found the following
description of $\Dir^2$ at the different faces:
\begin{corollary}\label{diracsquare}
The square of the Dirac operator gives a map $$\Psi_{G,{\rm cl}
}^{l,m,p}(X) \stackrel{t\Dir^2}{\longrightarrow}\Psi_{G,{\rm cl}
}^{l,m,p}(X)$$ with
\begin{enumerate}
\item $N_{G,\fatf}(t\Dir^2)=N_{G,\fatf}(t(\Dir^d)^2)\sowie g_{\phitx/X}^{-1}
(\nnabla^{\Lambda(\fatf)},\nnabla^{\Lambda(\fatf)})
+\ext_d(x^2F^{E/S,d}(\cdot,\cdot))\quad=:{\cal H}_X$
\item $
N_{G,\fatff}(t\Dir^2)=N_{G,\fatff}(t(\Dir^d)^2)  $
\par\noindent\qquad $=\tau g_{\phindx/\partial X
}^{-1}(\nnabla^{\Lambda(\fatff)},\nnabla^{\Lambda(\fatff)})\qquad\qquad
\qquad\qquad\qquad\qquad\qquad\qquad\quad=:\tau{\cal H}_B$
\par\noindent\qquad $\qquad + \tau\left( g_{M/Y}^{-1}
(\nnabla^{\Lambda(\fatff)},\nnabla^{\Lambda(\fatff)})+
\frac{\kappa_{M/Y}}{4}+ {\sf
m}_d(\delta_h(x^2F^{E/S,d}(\cdot,\cdot)))\right)\qquad =:\tau{\bf B}^2$
\end{enumerate}
\end{corollary}
\begin{remark}\rm \label{oldfibredetanew}
Using the definition of  $\widehat{\nnabla}^{\Lambda(\fatff)}$
given in Corollary \ref{explicitconnection} we set
 $$\widehat{\bf
B}^2:=g_{M/Y}^{-1}(\widehat{\nnabla}^{\Lambda(\fatff)},\widehat{\nnabla}^{\Lambda(\fatff)})+
\frac{\kappa_{M/Y}}{4}+ {\sf m}_d(\delta_h(x^2F^{E/S,d}(\cdot,\cdot))).$$ Then
$\widehat{\bf B}^2$ (truly) commutes with $\ext(\frac{\dd x }{x})$
and we can write ${\bf B}^2$ as
\begin{equation}\label{B2hat}
\widehat{\bf
B}^2+\ext(\frac{\dd
x}{x})\left[\Dir^{\phi,V}-\frac{1}{2}(\ext_{d},\ext_d,
\cl_d) \langle xB_{\phi}(\cdot,\cdot),\cdot \rangle_d
-\frac{1}{8}(\ext_d,\ext_d,\cl_d)\langle
x\Omega_{\phi}(\cdot,\cdot),\cdot\rangle_d
\right].
\end{equation}
 Recall that our model space
$\phindx\rightarrow Y$ can be identified with ${}^bHM\rightarrow
Y$ using the map $\delta_h$. Also the Clifford algebra
$\Cl({}^dT^*X|_{\partial X})$ will be identified with
$\Cl({}^bT^*X|_M)$ in that way. Using Definition \ref{Sphi} and
Lemma \ref{lemnablab} we find $$ \Omega_M=x^2\Omega_{\phi}, \quad
"xS_M =S_{\phi}+\frac{1}{2}B_{\phi}",\quad R^Y=\phi^*R^{B,b},\quad
{\sf m}={\sf m}_d\circ\delta_h\circ x.$$ In the case $B_{\phi}=0$,
for instance if $g_{\phi}$ is a product near the boundary, we thus
have $S_{\phi}=xS_M$ and $\widehat{\bf B}^2={\bf A}^2_M$, i.e.
(\ref{B2hat}) is of the classical form.
\end{remark}
\begin{proof}
These are direct consequences of
 the Lichnerowicz formula (\ref{Lichnerowicz}) for
the Dirac operator $\Dir^d$ and the calculations made in the above
series of Propositions, especially \ref{normcalcresc} and
\ref{curvterms}. To see that we can use the Dirac operator
$\Dir^d$ instead of $\Dir$, note that by the Lichnerowicz formula
and Lemma \ref{cab}(c) the difference of their squares:
\begin{eqnarray*}t(\Dir^2-(\Dir^d)^2)&=&tx^{v/2}(\Dir^d)^2x^{-v/2}-t(\Dir^d)^2=
tg_d^{-1}(x^{v/2}\nnabla x^{-v/2},x^{v/2}\nnabla
x^{-v/2})-tg_d^{-1}(\nnabla,\nnabla)\\ &=&
 -vt\nnabla_{X_d}+t\frac{v^2}{4}-t\frac{v}{2}\langle
 g_d^{-1},(\nnabla\frac{\dd x}{x})\rangle =
 -vt\nnabla_{X_d}-\frac{v^2}{4}t+O(x,t),
\end{eqnarray*}
vanishes at $\fatff$ and $\fatf$. \qed
\end{proof}

\subsection{The Rescaled Heat Kernel and the Index}\label{rescmodprob}
The aim of this Section is to prove an extension of Theorem
\ref{heatnormal} in the rescaled heat calculus. We start writing
down the solutions of the model problems at $\fatf$, $\fatff$ by
adapting the results in Appendices \ref{Mehler} and \ref{vertfam}
to our setting.
\par
\begin{proposition}
\begin{enumerate}
\item The heat kernel for ${\cal H}_X$
$${\sf K}_{{\cal H}_X}({\bf x},t):=p(x^2R^d,{\bf
x},t)\exp(\frac{1}{2}\langle x^2 R^d\cdot\nu_{\fatf}|_0,{\bf x
}\rangle_{\phi})\exp(-x^2F^{E/S,d})\dd t\dvol_{\phitx/X}$$ is
 in $ \Psi_{G,{\rm s},{\rm cl}}^{2,0}(\phitx/ X, \Lambda(\fatf))$
with normal operator at $\fatf(\phitx)=\fatf(X^2_H)$ given by
$$N_{G,\fatf,2}({\sf K }_{{\cal H}_X})({\bf x})=p(x^2R^d,{\bf
x},1)\exp(\frac{1}{2}Q_{\fatf}({\bf x})) \exp(-x^2F^{E/S,d})\dd
t\dvol_{\phitx/X}.$$ Recall that $Q_{\fatf}({\bf x})=\langle x^2
R^d\cdot\nu_{\fatf}|_0,{\bf x }\rangle_{\phi}$.
\item The heat kernel for ${\cal H}_B$
$${\sf K}_{{\cal H}_B}({\bf y},\tau):=p(\phi^*R^{B,b},{\bf
y},\tau)\exp(\frac{1}{2}Q_{\fatff}({\bf y}))\dd
\tau\dvol_{\phindx/\partial X}$$ is in $\Psi_{G,{\rm s},{\rm
cl}}^{2,0}(\phindx/
\partial X, \Lambda(\fatff))$
with normal operators at $\fatf=\phindx$ and ${\rm ti}$ given by
 $$N'_{G,\fatf,2}({\sf K}_{{\cal H}_B})({\bf y})
 =p(\phi^*R^{B,b},{\bf y},1)
\exp(\frac{1}{2}Q_{\fatff}({\bf y}))
 \dd \tau\dvol_{\phindx/\partial X},$$ $$N'_{H,{\rm ti},0}({\sf
K}_{{\cal H}_B})({\bf y})=\sqrt{4\pi}^{-h-1}\exp(-\|{\bf y}\|^2/4)
\dd \tau\dvol_{\phindx/\partial X}.$$
\end{enumerate}\par
\medskip\noindent  Recall that
$Q_{\fatff}({\bf y})=\langle x
\phi^*R^{B,b}\cdot\nu_{\fatff}|_0,{\bf y }\rangle_{\phi}$. Also,
the notation $N'_{G,\fatf}$ is used to remind the reader that
these are the normal operators defined using $\tau$ as the
temporal variable as  in Appendix \ref{Mehler}. Also, we emphasize
that all the normal operators above  commute with $\ext(\frac{\dd
x}{x})$.
\end{proposition}
\begin{proof}
We prove (a) using (the generalization to families of) Proposition
\ref{mehlerresult}. The connection $\nnabla^{\Lambda(\fatf)}_A$,
for $A\in\phitx$, appearing in the definition of ${\cal H }_X$ was
calculated in Corollary \ref{explicitconnection}(a):
\begin{eqnarray*}
\nnabla^{\Lambda(\fatf)}_A &=&\Int(A)\dd_{\END} -
\frac{1}{4}g_{\phitx/X}(x^2R^d(\cdot,\cdot){\bf
x},A)-\frac{1}{2}Q_{\fatf}(A).
\end{eqnarray*}
This is of the form required in Theorem \ref{mehlerformula} with
$R=x^2R^d$ and $q=Q_{\fatf}$. The underlying vector space $W$ is
any fibre of the bundle $\phitx\rightarrow X$, the metric $g_W$ is
given by the fibre metric $g_{\phitx/X}$ and the coefficient
algebra ${\cal A}$ is $\Lambda^{\rm ev}(x{\,}^bT^*X)\otimes
\END_{\Cl({}^dT^*X)}(E)$.\qed
\end{proof}
\begin{theorem}[Rescaled Heat Kernel]\label{rescaledhk}
The heat kernel for $\Dir^2$ is an element $${\sf K}\in
\Psi_{G}^{\cal H}(X,E), \qquad {\cal
H}=(2+\NN,2+\NN,0+H_{\fahbf})\quad\mbox{as in}\quad
(\ref{heatindexset}).$$ Its rescaled normal operators at the
different faces are (well defined and) given by
\begin{enumerate}
\item $N_{G,\fatf,2}({\sf K})({\bf x})
   =p(x^2R^d,{\bf
x},1)\exp(\frac{1}{2}Q_{\fatff}({\bf x})) \exp(-x^2F^{E/S,d})\dd
t\dvol_{\phitx/X}$
\item $N_{G,\fatff,2}({\sf K})({\bf y},\tau)=
{\sf K}_{{\cal H}_B}({\bf y},\tau)/\dd \tau
\cdot {\sf K}_{{\bf B}^2}(\tau)$
\par \qquad\qquad\qquad\qquad $=p(\phi^*R^{B,b},{\bf y},\tau)
\exp(\frac{1}{2}Q_{\fatff}({\bf y}))
 \dvol_{\phindx/\partial X}{\sf K}_{{\bf B}^2}(\tau), $
\item $N_{H,\fahbf,0}({\sf K})=[\Pi_{\circ}][\exp(-t\Dir_Y^2)]\frac{1}{\sqrt{4\pi
t}}e^{-(\ln(s))^2/4t}\frac{\dd s}{s}$
\end{enumerate}
\end{theorem}
\begin{proof}
The solutions at the faces patch together because of Corollaries
\ref{curvboundary} and \ref{easyrestrict}. We thus get a
parametrix which solves the heat equation to first order in the
rescaled heat calculus. This is transformed into a   parametrix
with a smooting remainder using the Composition Formula
\ref{rescheatcomp}. Then a final inversion step as in Section
\ref{normalops} is performed.\qed
\end{proof}
\begin{remark}\rm
As in Theorem \ref{mainresult2}, we can show that
$H_{\fahbf}=\NN$. In the following paragraphs we assume that we
have shown the Theorem in this form, if only for reasons of
notational simplicity.
\end{remark}
\par
\bigskip
\noindent{\bf The Regularized Supertrace}
\par\medskip
Let us note what we have achieved w.r.t. the local supertrace of
the heat kernel by introducing the rescaling. Denote again by
$\bigtriangleup_H$ the lift of the diagonal
$\bigtriangleup_{X}\times[0,\infty[$ to $X^2_H$. First, using our
notation for blowups we get $$
\bigtriangleup_H=[X\times[0,\infty[_{t^{1/2}},\partial
X\times\{0\}].$$ The behavior of general elements in the rescaled
heat calculus w.r.t.  the local supertrace is now described in
parts (c) and (d) of the following
\begin{lemma}\label{traceproperties}
\begin{enumerate}
\item $t\KD_H|_{\bigtriangleup_H}=
\tau^{-n/2}\beta_{\bigtriangleup_H,R}^*\phiomega(X)\dd t
=\rho_{\fahbf}^{v-1}\rho_{\fatff}^{-h-2}\rho_{\fatf}^{-n}
\beta_{\bigtriangleup_H,R}^*\Omega(X)\dd t$
\item $\str:\Psi_H^{2,2,0}(X,E)\longrightarrow \rho_{\fahbf}^{-1}
\rho_{\fatff}^{-h-2}\rho_{\fatf}^{-n}
C^{\infty}(\bigtriangleup_H,\beta_{\bigtriangleup_H,R}^*\Omega(X)\dd
t )$
\item $\str:\Psi_G^{2,2,0}(X)\longrightarrow
C^{\infty}(\bigtriangleup_H,\beta_{\bigtriangleup_H,R}^*{}^b\Omega(X)\dd
t) \quad \left( = tC^{\infty}(\bigtriangleup_H,
\bomega(\bigtriangleup_H)) \right)$
\item $\beta_{\bigtriangleup_H,R}^*{}^b\Omega(X)|_{\fatf\cap\bigtriangleup_H}
\sowie{}^b\Omega(X),\qquad
\beta_{\bigtriangleup_H,R}^*{}^b\Omega(X)|_{\fatff\cap\bigtriangleup_H}\sowie
\frac{\dd \tau^{1/2}}{\tau^{1/2}}\Omega(\partial X)=
{}^b\Omega(\fatff) .$\qed
\end{enumerate}
\end{lemma}
\par\medskip
\noindent
 Of course, an element $L(t)\in
C^{\infty}(\bigtriangleup_H,\beta_{\bigtriangleup_H,R}^*{}^b\Omega(X))$
will in general not be integrable for fixed $t=T$. We therefore
have to consider the {\em regularized} integral
$$\REG_{z=0}\int_{X,T} x^z L(T),$$ as in Section
\ref{secconsequences}. There we showed that the limit of the
regularized supertrace $$ \Str_d({\sf
K}(T))=\REG_{z=0}\Str(x^z{\sf K}(T))=\REG_{z=0}\int_{X,T}x^z
\str({\sf K}(T))/\dd t$$ of the heat kernel at {\em large} times
equals the extended $L^2_b$-index of $\Dir$, denoted by ${\rm
ind}_-(\Dir)$. To analyze the behavior of this supertrace for {\em
small} times we have to improve Lemma \ref{traceproperties}
further:
\begin{lemma}[Small Time Limit of the Heat Supertrace]
\label{smalltimelem}\quad\par
\begin{enumerate}
\item The local supertrace $\str({\sf K}(t))/\dd t\in
C^{\infty}(\bigtriangleup_H,\beta_{\bigtriangleup_H,R}^*
 {}^b\Omega(X))$ vanishes in $\fatf\cap\fatff$ and
 $\fatff\cap\fahbf$.
\item $\lim\limits_{T\rightarrow 0}\Str_d({\sf K}(T))=
\lim\limits_{T\rightarrow 0}\REG\limits_{z=0}\int_{X} x^z
\str({\sf K}(T))/\dd t =
\int_{\fatff\cap\bigtriangleup_H}\str({\sf
K})|_{\fatff}+\int_{\fatf\cap\bigtriangleup_H}\str({\sf
K})|_{\fatf}$,\par the RHS being well defined because of (a).
\end{enumerate}
\end{lemma}
\begin{proof}
(a) follows from the fact, proved in Theorem \ref{rescaledhk}(a)
and (c) (see also the remark in Proposition \ref{curvterms}), that
the restriction of ${\sf K}$ to these corners contains no
contribution of type $\ext_d(\frac{\dd x}{x})$. Thus the
supertrace has to vanish there.
\par
 To prove (b) we localize
around the corner $\fatf\cap\fatff$ (the argument at
$\fatff\cap\fahbf$ is even easier), where we can use the
coordinates $x$ and $\tau^{1/2}$. Then, using the vanishing at the
corner shown in (a), we can calculate the $T\rightarrow 0$-limit
as follows: $$\REG_{z=0}\int x^z a(x,\frac{T^{1/2}}{x})\frac{\dd
x}{x}=\int a(x,\frac{T^{1/2}}{x})\frac{\dd x}{x}=-\int
a(\frac{T^{1/2}}{\tau^{1/2}},\tau^{1/2})\frac{\dd
\tau^{1/2}}{\tau^{1/2}},$$ for $a\in
C_c^{\infty}([0,\infty[\times[0,\infty[)$ with compact support and
$a(0,0)=0$. We can decompose $ a=b +a_1+a_2$ such that $b$
vanishes at $\{x=0\}\cup\{\tau^{1/2}=0\}$, and
$a_1=a(0,\tau^{1/2})$ and $a_2=a(x,0)$. Thus
\begin{eqnarray*}
\REG_{z=0}\int x^z a(x,\frac{T^{1/2}}{x})\frac{\dd x}{x}&=& \int
b(x,\frac{T^{1/2}}{x})\frac{\dd x}{x}+\int a(x,0)\frac{\dd
x}{x}-\int a(0,\tau^{1/2})\frac{\dd \tau^{1/2}}{\tau^{1/2}}\\ &
&\qquad \stackrel{T\rightarrow 0}{\longrightarrow}\int
a(x,0)\frac{\dd x}{x} -\int a(0,\tau^{1/2}) \frac{\dd
\tau^{1/2}}{\tau^{1/2}},
\end{eqnarray*}
which proves the claim.\qed
\end{proof}
\par
\smallskip
We can now calculate the (extended $L^2_b$-) index of $\Dir$
(resp. the extended $L^2_d$-index of $\Dir^d$), using an extension
of the McKean-Singer formula. Starting from (\ref{extindex}) we
have
\begin{eqnarray}
{\rm ind}_-(\Dir)&=&\lim_{T\rightarrow \infty} \REG_{z=0}
\Str(x^z{\sf K}(T)/\dd t)\nonumber
\\ &=&\lim_{T\rightarrow 0} \REG_{z=0}\Str(x^z{\sf K}(T)/\dd t) +\int_{0}^{\infty}
 \REG_{z=0}\frac{\partial}{\partial T}\Str(x^z{\sf
K}(T))\label{smallandfinite}
\end{eqnarray}
For the second summand, we apply the usual trick:
\begin{eqnarray*}
\frac{\partial}{\partial T}\Str(x^z{\sf
K}(T))&=&\frac{1}{2}\Str(x^z[\Dir,\Dir]{\sf K}(T))=
\frac{1}{2}\Str(x^z[\Dir,\Dir{\sf K}(T)])
\\ &=&\frac{1}{2}\Str([\Dir,x^z\Dir{\sf K}(T)]-[\Dir,x^z]\Dir{\sf K}(T))
=-\frac{1}{2}\Str(zx^z\cl_d(\frac{\dd x}{x})\Dir{\sf K}(T))
\end{eqnarray*}
The integrand in the second  summand in (\ref{smallandfinite})
therefore becomes
\begin{eqnarray*}
&&-\frac{1}{2}\RES_{z=0}\Str(x^z\cl_{d}(\frac{\dd x}{x})\Dir {\sf
K}(T)) =-\frac{1}{2}\int_{\partial
X}\str(N_{H,\fahbf,0}(\cl_{d}(\frac{\dd x}{x})\Dir {\sf K}(T)))
\\&&\qquad =-\frac{1}{2}\int_{\partial X}
\tr\left(\eps\Pi_{\circ}\cl_{d}(\frac{\dd
x}{x})[\Dir_Y+\cl_d(\frac{\dd x}{x})s\frac{\partial}{\partial
s}][\exp(-T\Dir_Y^2)]\frac{1}{\sqrt{4\pi
T}}e^{-(\ln(s)^2)/4T}\right)|_{s=1}
\\ &&\qquad =-\frac{1}{4\sqrt{\pi T}}\int_{\partial X}
\tr(\eps\Pi_{\circ}\cl_{d}(\frac{\dd
x}{x})\Dir_Y[\exp(-T\Dir_Y^2)])
\\ &&\qquad =-\frac{1}{4\sqrt{\pi T}}\int_{Y}
\tr_{\cal K}(\eps\cl_{d}(\frac{\dd x}{x})\Dir_Y[\exp(-T\Dir_Y^2)])
=:\frac{1}{2}\eta(\Dir_Y,{\cal K},T).
\end{eqnarray*}
This is the usual eta-integrand for the operator $\Dir_Y$, see for
instance \cite{Vaillant}. Note however that the underlying vector
bundle ${\cal K}\rightarrow Y$ is graded by $\eps$ and not just
$\eps_B\sowie\eps_{{}^bT^*B}$. Also recall from Section
\ref{normalopsbf} that $\Dir_Y$ might have a rather complicated
structure in general. In any case the eta-term vanishes when $h+1$
is odd, since
\begin{eqnarray*}
&&\tr_{\cal K}(\eps\cl_{d}(\frac{\dd
x}{x})\Dir_Y[\exp(-T\Dir_Y^2)])= \tr_{\cal
K}(\eps\eps_{B}^2\cl_{d}(\frac{\dd
x}{x})\Dir_Y[\exp(-T\Dir_Y^2)])\\ &&\qquad\qquad=\tr_{\cal
K}(\eps_{B}\eps\eps_{B}\cl_{d}(\frac{\dd
x}{x})\Dir_Y[\exp(-T\Dir_Y^2)]) =-\tr_{\cal
K}(\eps\cl_{d}(\frac{\dd x}{x})\Dir_Y[\exp(-T\Dir_Y^2)]).
\end{eqnarray*}
We have used that $\cl_{d}(\frac{\dd
x}{x})\Dir_Y[\exp(-T\Dir_Y^2)]$ is even w.r.t. $\eps_{B}$ and that
$\eps\eps_B=-\eps_B\eps$, when $h+1$ is odd.
\par\smallskip
 The calculation of
the other summand in (\ref{smallandfinite}) is straightforward,
using Lemma \ref{smalltimelem}. Recall that $Q_{\fatf}(0)=0$ as
well as $Q_{\fatff}(0)=0$, so the corresponding terms in Theorem
\ref{rescaledhk} do not contribute in the supertrace:
\begin{lemma}
The contributions at $\fatf$ and $\fatff$ are
\begin{enumerate}
\item $\str({\sf K})/\dd t|_{\fatf\cap\bigtriangleup_H}
=\sqrt{2\pi i}^{-n}{\rm top}_{\Lambda {}^d
T^*X}[\widehat{A}(R^d)\Ch(F^{E/S,d})]$
\item $\str({\sf K})/\dd t|_{\fatff\cap\bigtriangleup_H}=
\sqrt{2\pi i}^{-h-1}{\rm top}_{\Lambda {}^bT^*B} [\widehat{A}(
R^Y)\tau^{-{\sf N}_B/2}\str_{E/B}(\eps_{V^*\partial X}{\sf
K}_{{\bf B}^2}(\tau))]$
\end{enumerate}
\end{lemma}
\begin{proof}
First, at $\stackrel{\circ}{\fatf}\cong\phitx$ the density factor
looks like
$\tau^{n/2}t\KD_H|_{\fatf\cap\bigtriangleup_H}=\phiomega(X)\dd t$.
Therefore
\begin{eqnarray*}
\str({\sf K})/\dd t|_{\fatf\cap\bigtriangleup_H}&=&
\sqrt{2/i}^n\frac{{\rm top}_{\Lambda {}^dT^*X}}{
\dvol_{{}^dT^*X}}(\tau^{n/2}N_{G,\fatf,2}({\sf K})/\dd t)
\\ &=&\sqrt{2/i}^n {\rm
top}_{\Lambda {}^dT^*X} [p(R^d,0,1)\str_{E/S}(\exp(F^{E/S, d}))]
\\&=&\sqrt{2\pi i}^{-n}{\rm top}_{\Lambda {}^d
T^*X}[\widehat{A}(R^d)\Ch(F^{E/S,d})].
\end{eqnarray*}
In (b), using the counting function ${\sf N}_B$ on $\Lambda
{}^bT^*B$ and Lemma \ref{tracedecomplemma} we get:
\begin{eqnarray*}
\str({\sf K})/\dd t|_{\fatff\cap\bigtriangleup_H}&=&
\sqrt{2/i}^{h+1}\frac{{\rm top}_{\Lambda
{}^bT^*B}}{\dvol_{{}^bT^*B}}(\frac{x^{h+1}}{\dd
\tau}\str_{E/B}(N_{G,\fatff,2}({\sf K})) |_{\bigtriangleup_H})
\\ &=&\sqrt{2/i}^{h+1}{\rm top}_{\Lambda {}^bT^*B}\left[
p(R^Y,{\bf y}=0,\tau)\str_{E/B}({\sf K}_{{\bf
B}^2}(\tau))\right]|_{\bigtriangleup_H}
\\&=& \sqrt{2\pi i\tau}^{-h-1}{\rm top}_{\Lambda {}^bT^*B}
[\widehat{A}(\tau R^Y)\str_{E/B}({\sf K}_{{\bf B}^2}(\tau))]
\\ &=&\sqrt{2\pi i}^{-h-1}{\rm top}_{\Lambda {}^bT^*B}
[\widehat{A}( R^Y)\tau^{-{\sf N}_B/2}\str_{E/B}({\sf K}_{{\bf
B}^2}(\tau))].
\end{eqnarray*}
\qed
\end{proof}
\par Using the definition of the fibred eta invariant in
(\ref{deffibredeta}) and putting everything together we get the
following index theorem:
\begin{theorem}[Index Theorem]\label{indextheorem}
For the Dirac operator ${\Dir}$ fulfilling (\ref{mainassumption})
the index formula $${\rm ind}_-(\Dir)= \frac{1}{(2\pi
i)^{n/2}}\int_X\widehat{A}(R^d)\Ch(F^{E/S,d})+\frac{1}{(2\pi
i)^{(h+1)/2}}\int_Y\widehat{A}(R^Y)\widehat{\eta}(\Dir^{\phi,V},E)
+\frac{1}{2}\eta(\Dir_Y)$$ holds. \qed
\end{theorem}
\par
 Note that since
$R^d\in\Gamma(\Lambda^2T^*X\otimes\End({}^dT^*X))$, the first
integrand on the RHS is a {\em true} form on the compact manifold
$X$ and integrability is automatic. The second and third term on
the RHS should be compared with the adiabatic limit formula (see
\cite{BC}, \cite{Daithesis}) for the eta invariant of a geometric
Dirac operator on $\partial X\rightarrow Y$ under collapsing of
the metric in the fibres. We describe this in the case of the
signature operator in the next Section.
\par
Theorem \ref{indextheorem} contains the usual APS (or b-) index
theorem (see \cite{APS1} or \cite{Melaps}) for Dirac operators as
the special case $Z=\{{\rm pt}\}$. Also, note that $\Dir$  has
pure point spectrum when $\Dir^{\phi,V}$ is invertible. Then on
the LHS we have the true index of $\Dir$ and the eta-term on the
RHS vanishes. An extreme instance of this is  the case $B=\{{\rm
pt}\}$, when only the local term appears on the RHS.

\subsection{The Signature Operator ${\sf S}$}\label{signatureop}
 The signature operator on $X$ for the metric $g_d$ is the operator
$${\sf S}^d=\Dir^d=\dd_d+\dd^*_d\quad \mbox{on} \quad E=\Lambda\,{}^dT^*X.$$
It has the property that its square, $({\sf S}^d)^2$, is equal to
the Laplacian on forms
$\bigtriangleup_{\Lambda\,{}^dT^*X}$. Even if $X$ is not spin, there is
a {\em local} isomorphism $E\cong
\End(S({}^dT^*X))$.
The grading operator $\eps$ on $E$ is given by the Clifford action of the volume
element on the left spinor factor. A short calculation shows that
$$\eps = i^{n/2}\cl_d(\dvol_d) = i^{n/2} \star_d (-1)^{{\sf N}({\sf N}-1)/2}.
$$
For $n=\dim(X)=4k$ we have $\star_d^2=1$ and the $L^2_d$-signature
 ${\rm sgn}_{L^2}(X,g_d)$ of $X$ is defined as the signature of the quadratic form
$$\alpha\longmapsto
(\alpha,\star_d\alpha) \quad\mbox{on}\quad H^{2k}_{\rm Ho}(X, g_d)
:={\rm null}_{L^2_d}(\bigtriangleup_{\Lambda^{2k}\,{}^dT^*X}).$$
As before, the operator given by
$$ {\sf S}= x^{v/2} {\sf S}^d x^{-v/2}=x^{v/2} \dd_d x^{-v/2}+ x^{v/2} \dd_d^* x^{-v/2}
   =:\widehat{\dd}+\widehat{\dd}^*$$
is selfadjoint on $L^2_b(X,\Lambda\,{}^dT^*X)$ and its $L^2_b$-index equals
${\rm sgn}_{L^2}(X, g_d)$.
The $\phi$-differential operators
$$ x\widehat{\dd},
x\widehat{\dd}^*:C^{\infty}(X,\Lambda\,{}^dT^*X)\rightarrow
C^{\infty}(X,\Lambda\,{}^dT^*X),$$
have well defined  restrictions
to vertical differential operators over  the boundary, which
define the vertical cohomology of the boundary fibration:
\begin{lemma}
\begin{enumerate}
\item $\dd_Z=x^{-{\sf N}_Z}x\widehat{\dd}x^{{\sf N}_Z}:C^{\infty}(\partial X, \Lambda
V^*\partial X)\rightarrow C^{\infty}(\partial X, \Lambda
V^*\partial X)$
\item $\dd_Z^*=x^{-{\sf N}_Z}x\widehat{\dd}^*x^{{\sf N}_Z}:C^{\infty}(\partial X, \Lambda
V^*\partial X)\rightarrow C^{\infty}(\partial X, \Lambda
V^*\partial X)$
\item $H(Z)=H(\dd_Z)\cong H_{\rm Ho}(\dd_Z)$ is a vector bundle
over $Y$.
\item $H(x\widehat{\dd})\cong H_{\rm Ho}(x\widehat{\dd})
={\rm null}(x{\sf S}|_{\partial X})
\cong C^{\infty}(Y,x^{{\sf N}_Z}H_{\rm Ho}(\dd_Z)\otimes \Lambda\,{}^bT^*B)$
etc.
\end{enumerate}
\end{lemma}
\begin{proof}
Let $\alpha\in [{}^dN^*\partial X]={}^bT^*B$. Then $\dd \alpha \in
\Lambda^2\,{}^bT^*X$ and
$$[x\widehat{\dd},\ext(\alpha)]|_{\partial X}=[x\widehat{\dd},\Int(\alpha)]|_{\partial X}
=0 \quad\mbox{in}\quad C^{\infty}(\partial X,\Lambda\,{}^dT^*X),$$
from which (a) and (b) follow.
Parts (c) and (d) follow from the usual Hodge theoretic arguments.\qed
\end{proof}
\par\smallskip
\noindent
This Lemma shows that our main assumption (\ref{mainassumption})
automatically holds for the signature operator.
 Recall
the notation
$$C^{\infty}(X,\Lambda\,{}^dT^*X)^{\circ}=\{\xi\in
C^{\infty}(X,\Lambda\,{}^dT^*X)\,|\quad \xi|_{\partial X}\in H_{\rm
Ho}(x\widehat{\dd})\}.$$
The indicial operators $I_b(\widehat{\dd})$,
$I_b(\widehat{\dd}^*)$, as well as the corresponding indicial
families can be defined as in Section \ref{normalopsbf}.  For
instance
$$\widehat{\dd}_Y(\xi|_{\partial X})\equiv I_b(\widehat{\dd})(0)(\xi|_{\partial X}):=
(\Pi_{\circ}\widehat{\dd}\xi)|_{\partial
X},$$
with $(\xi|_{\partial X})\in H_{\rm Ho}(x\widehat{\dd})$ extended to
$\xi\in C^{\infty}(X,\Lambda\,{}^dT^*X)^{\circ}$
is well defined on $H_{\rm Ho}(x\widehat{\dd})$.
The operators $\widehat{\dd}_Y$, $(\widehat{\dd}^*)_Y$ have the
following properties
\begin{lemma}
\begin{enumerate}
\item $(\widehat{\dd}_Y)^*=(\widehat{\dd}^*)_Y$ and we just
write this as $\widehat{\dd}_Y^*$
\item $\widehat{\dd}_Y$ is of total degree 1 on $H_{\rm Ho}(x\widehat{\dd}) =
C^{\infty}(Y,x^{{\sf N}_Z}H_{\rm Ho}(\dd_Z)\otimes \Lambda \,{}^bT^*B)$
\item $(H_{\rm Ho}(x\widehat{\dd}), \widehat{\dd}_Y)$ is an
elliptical complex and $H(\widehat{\dd}_Y)\cong H_{\rm Ho}(\widehat{\dd}_Y)$ is finite
dimensional.
\item $C^{\infty}(\partial X,\Lambda\,{}^dT^*X)\supset H_{\rm Ho}(x\widehat{\dd}) =
C^{\infty}(Y, x^{{\sf N}_Z}H_{\rm Ho}(\dd_Z)\otimes \Lambda \,{}^bT^*B)\supset
H_{\rm Ho}(\widehat{\dd}_Y)$\qed
\end{enumerate}
\end{lemma}
\par\smallskip
 We have
 the following version of Proposition
\ref{localibd}, describing the operators $\widehat{\dd}_Y$,
$\widehat{\dd}_Y^*$:
\begin{lemma}\label{useless}
The operators $\widehat{\dd}_Y$,  $\widehat{\dd}_Y^*$ act on
$H_{\rm Ho}(x\widehat{\dd}) =
C^{\infty}(Y, x^{{\sf N}_Z}H_{\rm Ho}(\dd_Z)\otimes \Lambda \,{}^bT^*B)$
as operators of total order $1$, resp. $-1$. We have
\begin{eqnarray*}
 \widehat{\dd}_Y \quad=\quad\Pi_{\circ}\left[ ({\sf h}\ext_d)\nabla_{\cdot}^{\Lambda\,{}^dT^*X}
-{\sf v}\ext_d\lambda(x{\sf n}S_{\phi}^*(\cdot){\sf v})+\ext_d(\frac{\dd x}{x})
({\sf N}_{Z}-\frac{v}{2})\quad\right. & &(1,0)
\\ \qquad -{\sf v}\ext_d\nabla_{x \nabla_{\cdot}^{\phi}\nu}^{\Lambda\,{}^dT^*X}
   -({\sf v}\ext_d, {\sf v}\ext_d, {\sf v}\Int_d)
   \langle R^d({\sf X}_d,\cdot)\cdot,\cdot\rangle_d\quad& &(0,1)
\\ \qquad- \left .\frac{1}{2}({\sf v}\ext_d, {\sf v}\ext_d, {\sf n}\Int_d)
\langle x B_{\phi}(\cdot,\cdot)
,\cdot\rangle_{\dd}\right]
   & &(-1,2)
   \end{eqnarray*}
   \begin{eqnarray*}
 \widehat{\dd}^*_Y \quad=\quad\Pi_{\circ}\left[ -({\sf h}\Int_d)\nabla_{\cdot}^{\Lambda\,{}^dT^*X}
-{\sf v}\Int_d\lambda(x{\sf v}S_{\phi}(\cdot){\sf n})+\Int_d(\frac{\dd x}{x})
({\sf N}_{Z}-\frac{v}{2})\quad\right. & &(-1,0)
\\ \qquad +{\sf v}\Int_d\nabla_{x \nabla_{\cdot}^{\phi}\nu}^{\Lambda\,{}^dT^*X}
   +({\sf v}\Int_d, {\sf v}\ext_d, {\sf v}\Int_d)
   \langle R^d({\sf X}_d,\cdot)\cdot,\cdot\rangle_d\quad& &(0,-1)
\\ \qquad+ \left .\frac{1}{2}({\sf v}\Int_d, {\sf v}\Int_d, {\sf n}\ext_d)
\langle x B_{\phi}(\cdot,\cdot)
,\cdot\rangle_{\dd}\right]
   & &(1,-2)
   \end{eqnarray*}
where the column on the right indicates the $({\sf n},{\sf v})$-
degrees of the different parts. Note once again that the
contributions of type $(-1,2)$, $(0,1)$, $(1,-2)$ and $(0,-1)$
simply vanish when $g_{\phi}$
is a product $\phi$-metric.
\end{lemma}
\begin{proof}
Recall that
$F^{E}=R^{\Lambda\,{}^dT^*X}=-\lambda(R^{{}^dTX})$. We leave it to the
more enthusiastic readers
to
use this with Proposition \ref{rdatthedel} and the method of proof
of Proposition \ref{localibd} to prove the claim.
The appearance of the last summands in the $(1,0)$ and $(-1,0)$-terms
above is due to the fact that
$$[\widehat{\dd},\Int_d({\sf X}_d)]|_{\partial X}=[\widehat{\dd}^*,
\ext_d(\frac{\dd x}{x})]|_{\partial X}={\sf N}_Z-\frac{v}{2}\quad\mbox{on}\quad
\Lambda\,{}^dT^*X.$$
\qed
\end{proof}
\par
\smallskip
As in \cite{BiLo}, the fact that $\widehat{\dd}_Y^2=0$ can be used
to show that
 $$\nabla^{{\cal H}(Z)}_T:=\nabla_{T}^{\Lambda\,{}^dT^*X}
+\lambda(\cdot{\sf v}S_{\phi}({\sf v}\cdot)xT)$$
is a flat connection on ${\cal H}(Z)= x^{{\sf N}_Z}H_{\rm Ho}(\dd_Z)
\otimes\langle\frac{\dd x}{x}\rangle\rightarrow Y$. In the case of a $\phi$-product metric the model
of the signature operator then is
$${\sf S}_Y=\overline{\dd}_Y+\overline{\dd}_Y^*=
\ext_Y \nabla_{\cdot}^{{\cal H}(Z)}-\Int_Y \nabla_{\cdot}^{{\cal H}(Z)^*}
+(\ext_d(\frac{\dd x}{x})+\Int_d(\frac{\dd x}{x}))({\sf
N}_Z-\frac{v}{2})$$
on ${\cal H}(Z)$.
The
operator ${\sf S}_Y$ thus differs by an endomorphism
 from the twisted signature operator ${\sf
S}(Y,{H}(\dd_Z))$. However, we will only be interested in calculating
the associated eta invariants. We know that the eta invariant vanishes in the
case of odd dimensional fibres. In the case of even dimensional
fibres, an easy symmetry argument shows
 that only the
restriction of the operator ${\sf S}_Y$ to ${\cal H}^{v/2}(Z)$
contributes to its eta invariant. Therefore $$\eta({\sf
S}_Y)=\eta({\sf S}(Y,{H}(\dd_Z))),$$
\par
\smallskip
Also, the signature operator ${\sf S}$ has the usual property that its extended
$L^2_b$-index equals its true $L^2_b$-index:
\begin{proposition} \label{extnonext} $
   \ind_{-}({\sf S})=\ind_{L^2_b}({\sf S})={\rm sgn}_{L^2}(X, g_d)$.
\end{proposition}
\begin{proof}
This is shown just as in the $b$-case. Recall that we know from
Proposition \ref{basicnull} that any element $u\in {\rm null}_-({\sf S})$
is of the form $u=u_0+O(x^{\varepsilon})$ with $u_{0}\in
C^{\infty}(X,\Lambda\,{}^dT^*X)^{\circ}$. Defining $J$ to be the
corresponding restriction map $J(u):=u_0|_{\partial X}$ we get
the sequence
$$ 0\longrightarrow {\rm null}_{L^2_b}({\sf S})\longrightarrow {\rm
null}_-({\sf S})\stackrel{J}{\longrightarrow} H_{\rm Ho}(\widehat{\dd}_Y)
\subset H(x\widehat{\dd}),$$
which is exact in the left and in the middle and compatible with the
grading $\epsilon$. Thus
$$ \ind_{-}({\sf S})-\ind_{L^2_b}({\sf S})= {\rm sdim}({\rm null}_-(
{\sf S}))
-{\rm sdim}({\rm null}_{L^2_b}({\sf S}))={\rm sdim}({\rm im}(J)),$$
and we  proceed to show that the RHS vanishes.
The idea is to use the decomposition
\begin{equation}\label{dxornodx}
H(x\widehat{\dd})
=C^{\infty}(Y,x^{{\sf N}_Z}H(\dd_Z)\otimes \Lambda T^*Y)\oplus
C^{\infty}(Y,x^{{\sf N}_Z}H(\dd_Z)\otimes\frac{\dd x}{x}\wedge\Lambda T^*Y),
\end{equation}
and to show that ${\rm im}(J)$ decomposes accordingly, i.e. the
projections on the left and right summand are maps in ${\rm
im}(J)$. Once this is shown, the fact that the grading $\eps$
switches the two summands in (\ref{dxornodx}) implies that
the  $\pm 1$ eigenspaces of $\eps$ in ${\rm im}(J)$ have equal
dimension.
\par\smallskip
Start with $\xi\in {\rm null}_-({\sf S})$ such that
$J(\xi)\neq 0$. We may then assume that $\xi$ pairs to $0$ with
${\rm null}_{L^2_b}({\sf S})$ under the usual
$L^2_b$-pairing (\ref{epspairing}). It then follows from Lemma
\ref{exactdd} that $\xi = {\sf S} \eta$ for some
$\eta \in {\rm null}_-({\sf S}^2)$. According to Lemma
\ref{basicnull} this has the form
$$\eta=  x^0\eta_0+  x^0\log(x) \eta_1+O(x^{\varepsilon}),\quad
\mbox{with}\quad \eta_0, \eta_1\in
C^{\infty}(X,\Lambda\,{}^dT^*X)^{\circ}. $$
Then $J(\xi)=J(\widehat{\dd}\eta)+J(\widehat{\dd}^*\eta)$ is a
decomposition of the desired type, since
$$
\widehat{\dd}\eta = \widehat{\dd} \eta_0 +\ext(\frac{\dd x}{x})\eta_1
+ \log(x) \widehat{\dd}\eta_1 +O(x^{\varepsilon})= \ext(\frac{\dd
x}{x})\eta_1+O(x^{\varepsilon})$$
and analogously $\widehat{\dd}^*\eta=-\Int(\frac{\dd
x}{x})\eta_1+O(x^{\varepsilon})$.\qed
\end{proof}
\par
\smallskip
It then follows immediately that Theorem \ref{indextheorem} has the following form
\begin{theorem}[\mbox{\boldmath $L^2$}-Signature Theorem]\label{signaturetheorem}
$${\rm sgn}_{L^2}(X, g_d)= \frac{1}{(\pi
i)^{n/2}}\int_X L(R^d)+\frac{1}{(\pi
i)^{(h+1)/2}}\int_Y L(R^Y)\widehat{\eta}({\sf S}_Z,\Lambda V^*M)
+\frac{1}{2}\eta({\sf S}(Y,H(\dd_Z)).$$
\end{theorem}
\begin{proof}
The LHS in Theorem \ref{signaturetheorem} equals ${\rm sgn}_{L^2}(X, g_d)$
as shown in Proposition \ref{extnonext}.
For the terms on the RHS use the fact that for $E=\Lambda\, {}^dT^*X$
 the relative curvature $F^{E/S}$
is  right Clifford multiplication  $\cl_{d,R}=\ext_d+\Int_d$ with
$\langle - R^d
\cdot,\cdot\rangle_d/2$, i.e. (see \cite{BGV}, Chapter 4)
$$ F^{E/S}=-\frac{1}{2}\cl_{d,R}\langle R^d
\cdot,\cdot\rangle_d
= -\frac{1}{4}({\sf n}\cl_{d,R},{\sf n}\cl_{d,R})
\langle R^{B,b}
\cdot,\cdot\rangle_d-\frac{1}{4}({\sf v}\cl_{d,R}, {\sf v}\cl_{d,R})
\langle R^d \cdot,\cdot\rangle_d.$$
Recalling (see again \cite{BGV}) that
$\Ch(F^{E/S,d})=2^{-n/2}\str(\eps F^{E/S})$
one shows that
\begin{eqnarray*}
L(R^d)&=&2^{-n/2}\widehat{A}(R^d)\Ch(F^{E/S,d})
\\ L(R^{Y})&=&L(R^{B,b})=2^{-(h+1)/2} \widehat{A}(R^{B,b})\tr_{S^*}
(\exp(-\frac{1}{4}{\sf n}\cl_{d,R}
\langle R^{B,b}
\cdot,\cdot\rangle_d))
\end{eqnarray*}
Now rewrite  (\ref{bsquare}) in Appendix \ref{vertfam} as
\begin{eqnarray*}
{\bf B}^2&=&
g_{VM/M}^{-1}(\nnabla^{\Lambda(\fatff)},\nnabla^{\Lambda(\fatff)})
+\frac{\kappa_{M/Y}}{4}
-\frac{1}{8}({\sf m}_d,{\sf m}_d,\cl_{d,R},\cl_{d,R})
\langle \delta_h(x^2R^d(\cdot,\cdot))\cdot,\cdot\rangle_d
\\ &=:&{\bf B}^2_{Z}-\frac{1}{8}({\sf h}\ext,{\sf h}\ext,\cl_{d,R},\cl_{d,R})
\langle \delta_h(x^2R^{B,b}(\cdot,\cdot))\cdot,\cdot\rangle_d,
\end{eqnarray*}
and define  $\widehat{\eta}({\sf S}_Z,\Lambda V^*M)$ as in Appendix
\ref{vertfam}, but using ${\bf B}^2_Z$ instead of ${\bf B}^2$.\qed
\end{proof}
\par
\bigskip
\noindent{\bf Collapsing Cylinders}
\par
\medskip
Fix a product structure,
$$C(M)=M\times [0,\varepsilon[_x \stackrel{\varphi}{\longrightarrow} Y\times
[0,\varepsilon[_x,$$
in a collar neighborhood of the boundary and consider the family
of exact $b$-metrics given by
$$g_{b,\varepsilon}=\left(\frac{\dd
x}{x}\right)^2+\varphi^*g_Y+(x^2+\varepsilon^2)g_Z,$$
on $C(M)$ and fixed on the interior of the manifold.
Denote the associated curvature by $R^{b,\varepsilon}\in
\Lambda^2\,{}^bT^*X\otimes \End({}^bTX)$ and the associated signature
operator on $\Lambda {}^bT^*X$ by ${\sf S}_{b,\varepsilon}$.
Then the $b$-signature theorem (or our Theorem \ref{signaturetheorem}
in the case with no fibre) tells us that
\begin{equation}\label{signatureeps}
{\rm sgn}_{L^2}(X, g_{b,\varepsilon})
=\frac{1}{(\pi
i)^{n/2}}\int_X L(R^{b,\varepsilon})+
\frac{1}{2}\eta({\sf S}_{b,\varepsilon}|_{\partial X}).
\end{equation}
It is known that
$${\rm sgn}_{L^2}(X, g_{b,\varepsilon})={\rm sgn}(X,\partial X),$$
is independent of $\varepsilon$. As $\varepsilon$
approaches $0$, the metric $g_{b,\varepsilon}$ approaches the $d$-metric $g_d$
and $L(R^{b,\varepsilon})$ approaches $L(R^d)$.
Also, it follows from Dai's adiabatic limit formula for the
signature operator in \cite{Daithesis} that
\begin{equation}\label{ausrede}
\lim_{\varepsilon\searrow 0}\frac{1}{2}
\eta({\sf S}_{b,\varepsilon}|_{\partial X})=
\frac{1}{(\pi
i)^{(h+1)/2}}\int_Y L(R^Y)\widehat{\eta}({\sf S}_Z,\Lambda V^*M)
+\frac{1}{2}\eta({\sf S}(Y,{H}(\dd_Z)))+\mbox{\boldmath $\tau$}
\end{equation}
The invariant $\mbox{\boldmath $\tau$}$ is defined in
\cite{Daithesis}, section 4.
It is a topological invariant associated to the
Leray spectral sequence of the fibration $\phi:M\rightarrow Y$.
Comparing the RHS in (\ref{ausrede}) with the RHS in our
signature formula, we find the following relation for the
$L^2_d$-signature, generalizing another formula in \cite{Daithesis}:
$${\rm sgn}_{L^2}(X,g_d)={\rm sgn}(X,\partial X)-\mbox{\boldmath
$\tau$}.$$
It is very likely that this result can be proved directly by
analyzing the $L^2_d$-cohomology of $X$ in more detail and
establishing its relationship with the spectral sequence.
\par
\smallskip
We leave considerations of this kind to the future.
\par\bigskip
\begin{remark}\rm \label{shorterproof}
The proof of the index formula, Theorem \ref{indextheorem},
presented
in this text is far from ``minimal''. We therefore want finish this work
with a brief discussion of possible shortcuts to the proof.
\par
Recall that  via the McKean-Singer
argument, the LHS of the index formula corresponds to the large time limit
of the regularized heat supertrace, whereas the RHS corresponds to the
small time limit, plus the difference of the two in the form of the eta
invariant.  It is obvious from Section \ref{rescmodprob}
 that a  full construction of
the heat kernel is not necessary for the calculation of the RHS.
 All that is needed is a
{\em parametrix}, solving the heat equation to first order at the faces
$\fatf$, $\fatff$ and $\fahbf$ within the {\em rescaled} heat
calculus. This means  that we only need to solve the rescaled model
problems at these faces. Especially, for the index result itself, no composition
formulas for the heat calculus are needed.
\par
On the LHS, it is easily seen that the parametrix for the
resolvent constructed in Proposition \ref{holo1} is compact for $\Lambda$
close to $0$, when $\Dir_Y$ is invertible. Thus $\Dir$ is Fredholm
in that case, i.e. its spectral measure has a gap at $0$
and $\Dir$ has a true index. The calculation of this index as the large
time limit of the heat supertrace is then straightforward.
Note that in the case of the signature operator ${\sf S}$ this
assumption means that the flat cohomology bundle $H(Z)\rightarrow Y$
has to be {\em acyclic}, which seems somewhat too restrictive.
However, given the index formula for Fredholm operators, we can
easily deduce the general index formula by using a weight shifting
argument, see \cite{Melaps} or \cite{Vaillant}.
\end{remark}

\newpage
\appendix
\section{Appendix}
\subsection{Clifford Algebra Conventions}\label{cliffordconv}
Let $(W, g_W)$ be a euclidean vectorspace of dimension $n$ (usualy
$W$ will be the fibre of  some tangent vector bundle) and
$(W^*,g_{W^*}=g_W^{-1})$ its dual . We denote by $\Cl(W^*)$ the
associated Clifford algebra over $W^*$. This is defined as the
quotient of the tensor algebra over $W^*$ by the relation
$$[\cl(\alpha),\cl(\beta)]=-2g_{W^*}(\alpha,\beta), \qquad \alpha,
\beta \in W\subset \Cl(W^*).$$ Here, $\cl$ is the inclusion of
$W^*$ into $\Cl(W^*)$ which can be regarded as an element $\cl\in
W\otimes \Cl(W^*)$. Elements $\alpha\in W^*$ act on the exterior
algebra by
$$\lambda(\alpha)\omega:=\ext(\alpha)\omega-\Int(g_{W^*}(\alpha))\omega.$$
and it is straightforward to see that this induces an action
$\lambda$ of $\Cl(W^*)$ on $\Lambda W^*$. We will also use the
action $\cl \in W^*\otimes \Cl(W^*)$ by writing $$
\cl(A):=\cl(g_W(A,\cdot)), \qquad\mbox{for}\quad A\in W.$$ We will
include the background metric into the notation, whenever this
seems necessary. The {\em symbol map} is the the isomorphism of
vector bundles given by $$\sigma: \Cl(W^*)\longrightarrow \Lambda
W^*,\qquad \sigma(\xi)=\lambda(\xi)1.$$ The inverse of this map
(sometimes called ``quantization'' map) is a map $$\cl:\Lambda
W^*\rightarrow \Cl(W^*),$$ which extends the above inclusion of
$W^*$ into $\Cl(W^*)$. On elements $\omega\in \Lambda^2 W^*$ this
looks like
$$\cl(\omega)=\sum_{i<j}\cl(e_i^*)\cl(e_j^*)\omega(e_i,e_j)
=\frac{1}{2}\sum_{i,j}\cl(e_i^*)\cl(e_j^*)\omega(e_i,e_j)=\frac{1}{2}(\cl,\cl)(\omega).$$
The notation on the RHS will also be used for general tensors
$\eta\in W^*\otimes W^*$:
$$(\cl,\cl)(\eta):=\sum_{i,j}\cl(e_i^*)\cl(e_j^*)\eta(e_i,e_j).$$
Beware of the following common source for sign errors: For an {\em
endomorphism} $\Phi\in\End(W)$ we write $$(\cl,\cl)(\Phi):=
(\cl,\cl)(g_W( \Phi \cdot,\cdot)) \quad (=2\cl(\Phi) \quad\mbox{if
$\Phi$ is skew}),$$ i.e. the ``intuitive'' order of the Clifford
operations is reversed! Whenever the endomorphism $\Phi$ is skew,
we have $$[\frac{1}{2}\cl(\Phi),\cl(A)]=\cl(\Phi(A)).$$
\par
\medskip
\pagebreak
\noindent{\bf The Supertrace}
\par
\medskip
Fix an orientation on $W$. Then the involution
\begin{equation}\label{gradingoperator}
\eps_W:=i^{[\frac{n+1}{2}]}\cl(\dvol_W),\qquad
\eps_W^2=1,
\end{equation} defines a (possibly trivial) grading on
every representation space of $\Cl(W)$. In the following, let
$n=\dim(W)$ be even and denote by $S(W)$ the spinor space for
$\Cl(W)$ (a possible construction of $S(W)$ is described briefly
in the next Appendix). This is the unique irreducible
representation of $\Cl(W)$. The action of $\Cl(W)$ on $S(W)$
defines an isomorphism
\begin{equation}\label{cliffordendo}\Cl(W)\cong \End(S(W)).
\end{equation}
 This representation is truly graded by the operator $\eps_W$
 and
the identification (\ref{cliffordendo}) allows the corresponding
supertrace of an element $Q$ in $\End(S(W))$ to be calculated as
follows
\begin{equation}\label{tracedecomp}\str_{S(W)}(Q)=\tr(\eps_WQ)=(\frac{2}{i})^{n/2}{\rm
top}_{\Lambda W^*}(\sigma(Q))/\dvol_W.
\end{equation}
We will use this formula for the first factor in the decomposition
$$\End(E)=\phi^*\Cl({}^{\rm b}T^*B)\otimes
       \END_{\Cl({}^{\rm b}T^*B)}(E),$$
       described in Lemma \ref{filtration}.
 Write $\eps=\eps_{{}^dT^*X}$, $\eps_B=\eps_{{}^bT^*B}$ and
       $\eps_Z=\eps_{{}^dV^*\partial X}$.
Then, depending on the parity of the dimensions
$n=\dim({}^dT^*X)$, $h+1=\dim({}^{b}T^*B)$ and
$v=\dim(V^*\partial X)$, with $n=h+1+v$ even, we have
$$ \eps=\eps_B\eps_Z \quad\mbox{if $h+1$ is even},\quad
   \eps=-i\eps_B\eps_Z \quad\mbox{if $h+1$ is odd}.$$
For an element $Q_Z\in\END_{\Cl({}^{\rm b}T^*B)}(E)$ we
define the {\em $B$-relative} supertrace as
\begin{equation}\label{tracetrick}
\str_{E/B}(Q_Z):=\Const(h)2^{-(h+1)/2}\str_E(\eps_B Q_Z),
\quad\mbox{with}\quad \Const(h)=\left\{\begin{array}{ll}
\sqrt{i}^3 & h\quad\mbox{even}
\\ 1 & h \quad\mbox{odd}
\end{array}\right.
\end{equation}
We can then show the following formula for the supertrace on $E$.
\begin{lemma}\label{tracedecomplemma}
For $Q_B\otimes Q_Z \in \Cl(\phi^*{}^{b}T^*B)\otimes
\END_{\Cl({}^{\rm b}T^*B)}(E)$ we have
$$\str_{E}(Q_B\otimes Q_Z)= (\frac{2}{i})^{(h+1)/2}{\rm
top}_{\Lambda {}^{b}T^*B} (\sigma(Q_B))/\dvol_{{}^{b}T^*B}
\cdot\str_{E/B}( Q_Z).$$
\end{lemma}
\begin{proof}
It suffices to prove this for $E=S({}^dT^*X)$ and $Q_Z\in \Cl({}^dV^*\partial
X)$. If both $v$ and $h+1$ are {\em even}, the Lemma is a direct
consequence of (\ref{tracedecomp}):
 \begin{eqnarray*}&&\str_{S({}^dT^*X)}(Q_B\otimes Q_Z)=
 \str_{S({}^{b}T^*B)}(Q_B)
 \str_{S({}^dV^*\partial X)}(Q_Z)\\
 &&\qquad =(\frac{2}{i})^{(h+1)/2}{\rm top}_{\Lambda {}^{b}T^*B}
(\sigma(Q_B))/\dvol_{{}^{b}T^*B}\cdot 2^{-(h+1)2)}\tr_{S({}^dT^*X)}
(\eps_{{}^dV^*\partial X} Q_Z).
\end{eqnarray*}
When  $v$ and $h+1$ are {\em odd}, we write
$Q_B=\widehat{Q}_B\cl_{d}(\frac{\dd x}{x})$ to see that
\begin{eqnarray*}
&&\str_{S({}^dT^*X)}(Q_B\otimes Q_Z)\\
 &&\qquad =(\frac{2}{i})^{h/2}{\rm
top}_{\Lambda {}^{T^*Y}}
(\sigma(\widehat{Q}_B))/\dvol_{T^*Y}\cdot\tr_{S(\langle\frac{\dd x}{x}\rangle
\oplus
{}^dV^*\partial X)}(\eps_{\langle\frac{\dd x}{x}\rangle\oplus
{}^dV^*\partial X} \cl_{d}(\frac{\dd x}{x}) Q_Z)
\\ &&\qquad =(\frac{2}{i})^{h/2}{\rm
top}_{\Lambda {}^{b}T^*B}(\sigma(Q_B))/\dvol_{ {}^{b}T^*B}
2^{-h/2}\tr_{S({}^dT^*X)}(\eps_{Z}Q_Z)
\\ &&\qquad =(\frac{2}{i})^{(h+1)/2}{\rm
top}_{\Lambda {}^{b}T^*B}(\sigma(Q_B))/\dvol_{ {}^{b}T^*B}
2^{-(h+1)/2}\tr_{S({}^dT^*X)}(i^{3/2}\eps\eps_B Q_Z)
\end{eqnarray*}
where we have used $\eps_{\langle\frac{\dd x}{x}\rangle\oplus
{}^dV^*\partial X}=\cl_d(\frac{\dd x}{x})\eps_{{}^dV^*\partial X}$.\qed
\end{proof}

\subsection{Dirac Operator and Conformal Transformations}\label{conformal}
In this Appendix we consider the behavior of Dirac operators under
a conformal transformation of the metric.
\par
Let $M$ be an evendimensional manifold with metric $g$ and
$E\rightarrow M$ a Clifford bundle over $M$. Since our
considerations will be purely local, we may assume
right away that $M$ is an open ball and $E$ is trivial over $M$.
Then $E$ is of the form $$ E\cong S\otimes E',\quad \nabla^E\cong
\nabla^S\otimes\nabla^{E'}$$ The spinor bundle $S$ can be
described as follows: There exists a polarization $ T^*M_{\CC} =
W\oplus \overline{W}$ of the complexified cotangent bundle, i.e.
$g(W,W)=0$. Then $S=\Lambda W$ with metric $g^S$ induced by $g$
and  Clifford action of an element $\alpha \in T^*M$ given by
$$\cl(\alpha)=\cl(\alpha_1\oplus\overline{\alpha}_2)
       = \frac{1}{\sqrt{2}}(\ext(\alpha_1)-\Int_g(\overline{\alpha}_2)).$$
Choosing a local orthonormal frame $e^1,\ldots, e^n$ for $T^*M$,
the Levi Civita connection on $T^*M$ is of the form $\dd + \omega$
with $\omega\in \Gamma(T^*M\otimes{\bf so}(T^*M))$. The spinor
connection $\nabla^S$ is then given by
$\dd+\frac{1}{2}\cl(\omega)$.
\par
Now, for a positive function $f\in C^{\infty}(M)$ consider the conformal
 metric $g_f=f^2g$. The corresponding Levi Civita connection on $T^*M$
is given for $T\in \Gamma(TM)$ and $\alpha \in \Gamma(T^*M)$ by
\begin{equation}\label{levif}
 \nabla^f_T\alpha=\nabla_T\alpha
     -\frac{T f}{f}\alpha+g^{-1}(\alpha, \frac{\dd
        f}{f})g(T)-\alpha(T)\frac{\dd f}{f}=: \nabla_T\alpha+G_f(T)\alpha.
 \end{equation}
As before, the spinor bundle is $S_f=S=\Lambda W$ but the Clifford
action is now given by $$\cl_f (\alpha)=
\frac{1}{\sqrt{2}}(\ext(\alpha_1)
   -\Int_f(\overline{\alpha}_2))=\frac{1}{\sqrt{2}}
   (\ext(\alpha_1)-f^{-2}\Int_g(\overline{\alpha}_2)).$$
The map $$f:(T^*M,g)\ni \alpha\mapsto f\cdot \alpha \in
(T^*M,g_f)$$ is an isometry and thus induces an isometry of the
spinor bundles $f^{\sf N}:S\rightarrow S_f$, where ${\sf N}$ is
the number operator in $\Lambda W$. We now have $$f^{-{\sf
N}}\cl_f(\alpha) f^{\sf N} =f^{-{\sf N}}(\ext(\alpha_1)
   -f^{-2}\Int_g(\overline{\alpha}_2))f^{\sf N}
   = f^{-1}\cl(\alpha),$$
and from the description of $\nabla^S, \nabla^{S_f}$ we deduce $$
f^{-{\sf N}}\nabla^{S_f}_Tf^{\sf
N}=\nabla^{S}_T+\frac{1}{2}\cl(\ext(\frac{\dd f}{f})
\Int(T)-\ext_g(T)\Int_g(\frac{\dd f}{f})) = \nabla^S_T+\widehat{G}(T),$$
where we have defined $G(T)$ to be
$\frac{1}{2}(\cl_g(T)\cl(\frac{\dd f}{f})+\frac{\dd f}{f}(T))$.
This just means that the following diagrams are commutative $$
\begin{array}{ccccccc} S\otimes
T^*M&\stackrel{f^{\sf N}\otimes {\rm id}} {\longrightarrow}&
S_f\otimes T^*M&\quad&\Gamma(S)&\stackrel{f^{\sf
N}}{\longrightarrow}&\Gamma(S_f)\\ \qquad\downarrow f^{-1}\cl&
&\quad\downarrow \cl_f &\quad &\quad \qquad \downarrow \nabla^S+\widehat{G}&
&\quad\downarrow \nabla^{S_f}\\
 S &\stackrel{f^{\sf N}}{\longrightarrow}& S_f &\quad& \Gamma(S\otimes
T^*M)&\stackrel{f^{\sf N}}{\longrightarrow}& \Gamma(S_f\otimes
T^*M),
\end{array} $$
and with these isomorphisms the corresponding Dirac operator is of
the form $$\Dir^f:=\cl_f\circ \nabla^{{S_f}} \sowie
f^{-1}\Dir+\frac{n-1}{2}\cl(\frac{\dd f}{f^2})
=f^{-\frac{n-1}{2}}f^{-1}\Dir f^{\frac{n-1}{2}}.$$
\par
Note that for some Clifford bundles like $E=\Lambda T^*M$ the
bundle $E'$ in the local decomposition $E=S\otimes E'$ might also
contain copies of $S$ or $S^*$. The above formulas remain true in
that case, at the cost of a nonuniform treatment of the different
copies of $S$ in $E$.\par
\smallskip
Since it will be of some relevance in applications,
see Section \ref{signatureop}, we describe
the case of the signature operator ${\sf S}$ on the
 bundle $E=\Lambda T^*M$ in detail.
Clifford action is given
$\cl(\alpha)=(\ext(\alpha)-\Int_g(\alpha))$ and our
identifications above would give $\cl_f(\alpha)=f^{-1}\cl(\alpha)$.
However, the natural isometry,
$$ f^{\sf N}:(\Lambda T^*M, g)\longrightarrow (\Lambda T^*M,
g_f)$$
yields an identification
$$f^{\sf N}\ext(\alpha)f^{-{\sf N}}=\ext(f\alpha) \quad\mbox{and}\quad
f^{\sf N}\Int(\alpha)f^{-{\sf N}}=\Int(f^{-1}\alpha)=f\Int_f(\alpha)).$$
Using the generalization to differential forms of (\ref{levif})
one can also show that
\begin{equation}\label{ddstarf}
f^{\sf N}(\dd +\dd^*)f^{-{\sf N}}=f (\dd_f+\dd^*_f)-f\ext(\frac{\dd
f}{f}){\sf N}-f\Int_f(\frac{\dd
f}{f})({\sf N}-n).
\end{equation}

\subsection{Geometry of Fibre Bundles}\label{geomtens}
In this Appendix we give a brief overview of the main notions in
the geometry of riemannian fibre bundles. Our exposition follows
\cite{BGV} closely.
\par\medskip
Let $\phi:M\rightarrow Y$ be a fibre bundle with fibre $Z$, with
$m:=\dim M, v:=\dim Z.$ The vertical bundle is denoted by $ V M:=
\ker (\phi_*)$. By choosing a horizontal subbundle $ H M\subset
TM$ we get a decomposition $$ TM = VM \oplus HM,$$ and associated
projections ${\sf v}:TM\rightarrow VM $, ${\sf h}:TM \rightarrow
HM$. This decomposition is defined to be orthogonal with respect
to a riemannian metric on $M$, of the form
\begin{equation}
\langle .,.\rangle_{M}= \phi^{*}\langle ., .\rangle_Y + \langle
.,.\rangle_{M/Y},\label{specialmetric}
\end{equation} where
$\langle .,. \rangle_Y$ and $\langle ., .\rangle_{M/Y}$
 are given metrics on the base and on the fibres respectively. \par
 Recall that the map $\phi_*:HM\rightarrow TY$ is a pointwise isomorphism.
This allows us to freely identify $HM$ with the lifted bundle
$\phi^*TY$. Denoting by $\nabla^M$ the Levi-Civita connection on
$M$ we then define the  induced connections
\[ \nabla^V\equiv\nabla^{M/Y}:= {\sf v}\circ \nabla^M\circ {\sf v}
 :\Gamma(M,VM)\rightarrow \Gamma(M,T^*M\otimes VM)\]
\[ \nabla^{H} := \phi^* \nabla^{Y}\circ {\sf h}:\Gamma(M,HM)\rightarrow \Gamma(
M,T^*M \otimes HM).\] These connections are invariant under
multiplication of $\langle.,.\rangle_Y$ and
$\langle.,.\rangle_{M/Y}$ with positive constants. The connection
$\nabla^{M, \oplus} := \nabla^V\oplus\nabla^H$ is thus also
invariant under this operation and compatible with the metric,
 but it can have torsion.\par\medskip
We now define the geometric tensors associated to the fibration:
\begin{definition} Let $X,Y$ be vector fields in $TM$
  \begin{enumerate}
  \item The {\em difference tensor} $\omega_M\in
  \Gamma(M,T^*M\otimes \so(TM))$
  is just
$$\omega_{M}(X)Y:= \nabla^{M}_XY-\nabla^{M, \oplus}_XY.$$
  \item The {\em second fundamental form} $S_{M}
  \in \Gamma(M, TM\otimes T^*M\otimes H^*M)$ of the fibres $Z$ is given by
 $$ S_{M}(X)Y:={\sf v}\nabla^{M}_{{\sf v}X} {\sf h}Y.$$
  \item The {\em curvature} $\Omega_M\in
\Gamma(\Lambda^2T^*M,VM)$ of the
  horizontal space is defined as $$ \Omega_M(X,Y)
             :={\sf v}[{\sf h}X,{\sf h}Y].$$
  \end{enumerate}
\end{definition}
The connection $\nabla^{M, \oplus}$ commutes with the projections
${\sf h}, {\sf v}$. Therefore $$ S_M(X)Y={\sf v}\omega_M({\sf v
}X){\sf h}Y, \qquad \Omega_M(X,Y)={\sf v }(\omega_M({\sf h }X){\sf
h }Y-\omega_M({\sf h}Y){\sf h}X).$$
\par
By definition the term ${\sf v}\omega_M{\sf v}$ vanishes. Using
the Koszul formula, we find furthermore that $\langle
\omega_M(X){\sf h}Y,{\sf h}Z \rangle = -\frac{1}{2}\langle
\Omega({\sf h }Y,{\sf h}Z), {\sf v}X \rangle$. With these formulas
it is now easy to show:
\begin{lemma}\label{omegaform}
Let $X,Y,Z$ be vector fields in $TM$. Then the following formula
for $\omega_{M}$ holds:
\begin{eqnarray*}
\lefteqn{  \langle \omega_{M}(X)Y,Z\rangle_M = \langle S_M({\sf
v}X){\sf h}Y,{\sf v}Z\rangle_M -\langle S_M({\sf v}X)
  {\sf h}Z,
{\sf v}Y\rangle_M} \\ & & - \frac{1}{2}\left(\langle \Omega_M
({\sf h}X,{\sf h}Z),
 {\sf v}Y \rangle_M -
 \langle \Omega_M ({\sf h}X,{\sf h}Y),{\sf v}Z \rangle_M +
 \langle \Omega_M ({\sf h}Y,{\sf h}Z),{\sf v}X \rangle_M\right).\qquad \qed
\end{eqnarray*}
\end{lemma}
\par
Let now ${\cal E}\rightarrow M$ be a Clifford bundle over $M$. The
connection on ${\cal E}$ is denoted by $\nabla^{{\cal E}}$, the
hermitian metric by $\langle.,. \rangle_{{\cal E}}$ and Clifford
multiplication by $\cl$.  We assume that there is another
connection $\nabla^{{\cal E}, \oplus}$ on ${{\cal E}}$ which is a
Clifford connection with respect to $\nabla^{M, \oplus}$, i.e.
$$[\nabla^{{\cal E}, \oplus}_X,\cl(\alpha)]=\cl(\nabla^{M,
\oplus}_X \alpha).$$ We then require that $\nabla^{{{\cal E}}}$ be
of the form $$\nabla^{{\cal E}}_X=\nabla^{{{\cal
E}},\oplus}_X+\frac{1}{2}\cl(\omega_M(X)).$$ This is clearly a
Clifford connection with respect to $\nabla^M$, as $$
[\nabla^{{\cal E}}_X,\cl(Y)]= [\nabla^{{{\cal
E}},\oplus}_X,\cl(Y)] +[\frac{1}{2}\cl(\omega_M(X)),\cl(Y)]=
\cl(\nabla^{M, \oplus}_XY+\omega_M(X)Y) ,$$ and the last term is
equal to $\cl(\nabla^M_XY)$.
\par\smallskip
As an example, these requirements are satisfied for the bundle
${{\cal E}}=\phi^*\Lambda T^*Y\otimes {{\cal E}}'$, where ${{\cal
E}}'$ is a vertical Clifford bundle, the Clifford action of
$\alpha \in T^*Y$ is given by
$\cl(\alpha)=\ext(\alpha)-\Int(\alpha)$,
 and $\nabla^{{\cal E}, \oplus}=\phi^*\nabla^{\Lambda T^*Y}\otimes
 \nabla^{{{\cal E}}'}$.
For these bundles we can also examine the degenerate Clifford
action ${\sf m}$ given by $${\sf m}({\sf h}\alpha)=\ext({\sf
h}\alpha), \qquad {\sf m}({\sf v}\alpha)=\cl({\sf v}\alpha),\qquad
[{\sf m}(\alpha),{\sf m}(\alpha)]=-2\langle{\sf v}\alpha,{\sf
v}\alpha\rangle_M.$$ This degenerate Clifford action can be
extended to more general tensors as in Appendix
\ref{cliffordconv}. The associated `degenerate' connection is then
given by $$ {\widehat{\nabla}}^{{\cal E}}_X :=\nabla_X^{{{\cal
E}},\oplus}+\frac{1}{2}{\sf m} (\omega_M(X)).$$ The identity $
[\frac{1}{2}{\sf m}(\omega(X)),{\sf m}(\alpha)]={\sf
m}(\alpha({\sf v}\omega(X)))$ implies that this connection is a
${\sf m}$-Clifford connection with respect to the connection
$${\widehat{\nabla}}^{T^*M}:= \nabla^{M, \oplus}+({\sf v
}\omega_M)^*\quad\mbox{ or} \quad{\widehat{\nabla}}^{TM}:=
\nabla^{M, \oplus}+{\sf v}\omega_M.$$ This connection can be
described more explicitly:
\begin{lemma} \label{horver}Let $W\in VM$, $A\in HM$
\begin{enumerate}
\item ${\widehat{\nabla}}^{{\cal E}}_W =\nabla_W^{{{\cal
E}},\oplus}+\frac{1}{2}(\ext,{\sf v}\cl)(\langle S_M(W){\sf
h}\cdot,{\sf v}\cdot \rangle_M) -\frac{1}{4} \ext(\langle
\Omega_M({\sf h}\cdot,{\sf h}\cdot),W\rangle_M)$
\item ${\widehat{\nabla}}^{{\cal E}}_A = \phi^*\nabla^{\Lambda
T^*Y}_A +\frac{1}{4}(\ext,{\sf v }\cl)(\langle
\Omega_M(A,\cdot),\cdot\rangle_M)$.
\end{enumerate}
\end{lemma}
\begin{proof}
For (a) note that
\begin{eqnarray*}
 {\sf m}(\omega(W))&=& \ext({\sf
h}\omega(W){\sf h })+(\ext, {\sf v}\cl)({\sf v}\omega(W){\sf h})\\
&=&-\frac{1}{2} \ext(\langle \Omega({\sf h}\cdot,{\sf
h}\cdot),W\rangle_M)+(\ext,{\sf v}\cl)(\langle S_M(W){\sf
h}\cdot,{\sf v}\cdot \rangle_M)
\end{eqnarray*}
\par Analogously in (b), using
$$\langle\omega_M(A)Y,Z\rangle_M= - \frac{1}{2}(\langle \Omega_M
(A,{\sf h}Z), {\sf v}Y \rangle_M -
 \langle \Omega_M (A,{\sf h}Y),{\sf v}Z \rangle_M),$$
 we get
$$ {\sf m}(\omega(A))=\frac{1}{2}({\sf m},{\sf
m})(\langle\omega(A)\cdot,\cdot\rangle_M)\\ =\frac{1}{2}(\ext,{\sf
v }\cl)(\langle \Omega_M(A,\cdot),\cdot\rangle_M) ,$$ and dividing
everything by 2 gives the result.\qed
\end{proof}
\par
\smallskip
Analogous to the definition of the Dirac operator $\Dir =\cl \circ
\nabla^{{\cal E}}$  on ${{\cal E}}$ we can define the
superconnection corresponding to ${\widehat{\nabla}}^{{\cal E}}$
by ${\bf A}_M:= {\sf m} \circ {\widehat{\nabla}}^{{\cal E}}$. In
the case of a fibre bundle with base space $Y$ a point this is
again just the Dirac operator. Using Lemmas \ref{omegaform} and
\ref{horver} we get
\begin{lemma}
Write $\tr(S_M)(A)=\tr(\cdot S_M({\sf v}\cdot)A)\in HM$. Then
 $${\bf A}_M:= {\sf m} \circ {\widehat{\nabla}}^{{\cal E}}=
 \dd_{H}+\Dir^{{\cal E},V}+\frac{1}{2}\ext(\tr(S_M))
 +\frac{1}{8}(\ext,\ext,\cl)\langle\Omega_M(\cdot,\cdot),\cdot\rangle_M.$$
\end{lemma}
 Finally, using the fact that ${\sf h}\omega(A)B$ is symmetric
in $A$ and $B$ it is easy to see that $\widehat{\nabla}^{TM}$ is
torsionfree! Thus, the following generalization of the
Lichnerowicz formula holds (cf. \cite{BGV}):
\begin{proposition} \label{lichn}
Set  $F^{{\cal E}/S}(X,Y):= F^{{\cal E}}(X,Y)-\frac{1}{2} {\sf
m}(R^{M}(X,Y))$. Then the following analogue of the Lichnerowicz
formula holds: $${\bf A }_M^2=\widehat{\bigtriangleup}^{{{\cal
E}},v}+{\sf m}( \widehat{F}^{{\cal E}})=
\widehat{\bigtriangleup}^{{{\cal E}},v}+\frac{\kappa_{M/Y}}{4}
+{\sf m}(F^{{\cal E}/S}(\cdot,\cdot)).$$
\end{proposition}
\begin{proof}
The proof proceeds as usual. $$ 2{\bf A}_M^2=[{\sf
m}\widehat{\nabla},{\sf m}\widehat{\nabla}]\\
 ={\sf m}{\sf m}[\widehat{\nabla},\widehat{\nabla}]
    +{\sf m}[\widehat{\nabla},{\sf m}]\widehat{\nabla}
    -{\sf m}[{\sf m},\widehat{\nabla}]\widehat{\nabla}
    +[{\sf m},{\sf m}]\widehat{\nabla}\widehat{\nabla}$$
The first summand on the RHS gives $2{\sf m}(\widehat{F}^{{\cal
E}})$ since $\widehat{\nabla}^{TM}$ is torsionfree. The second and
third summand vanish since ($\widehat{\nabla}^{\cal E}$,
$\widehat{\nabla}^{TM}$) is a Clifford connection. The last
summand gives $-2 g_M^{-1}({\sf v}\widehat{\nabla}, {\sf
v}\widehat{\nabla})$, which is just two times the vertical
Laplacian in the above formula. This proves the first equation in
the Proposition, for the proof of the second equation we refer to
\cite{BGV} Chapter 10. \qed
\end{proof}

\subsection{Mehler's Formula and the Simple Heat Space}\label{Mehler}
 In this Appendix we adapt  the well known Mehler formula to our setting.
  For this, let $(W,g_W)$ be a euclidean vector space and ${\cal A}$ a commutative
algebra (usually this will be the space of even forms on some
vector space). Let $R\in \End(W)\otimes{\cal A}$ such that
$g_W(R)\in \Lambda^2W^*\otimes {\cal A}$. Also let $q\in
W^*\otimes {\cal A} $. We consider $g_W(R {\bf y},\cdot)$ as a
linear map $W\rightarrow W^*\otimes{\cal A}$ and $q$ as  a
constant map of the same type. Then
\begin{theorem}\label{mehlerformula} The operator
\begin{equation}\label{operatorh}
H =g_W^{-1} (\dd
-\frac{1}{4}g_W(R)-\frac{1}{2}q,\dd-\frac{1}{4}g_W(R)-\frac{1}{2}q)
\end{equation}
on $C^{\infty}(W,{\cal A})$ has the heat kernel $[e^{-tH}](0,{\bf
y})=p(R,{\bf y},t)\exp(q({\bf y})/2)\dd t\dvol_W$ given by the
formula \begin{eqnarray*} p(R,{\bf y},t)&=& \sqrt{4\pi
t}^{-n/2}\det\left( \frac{tR/2}{e^{tR/2}-e^{-tR/2}}
\right)\exp\left(-\frac{1}{4t}\left\langle {\bf y}\left|
\frac{tR}{2} \coth\left(\frac{tR}{2}\right)\right|{\bf
y}\right\rangle_W\right), \\ p(R,0,t)&=&
 \sqrt{4\pi
t}^{-n/2}\det\left( \frac{tR/2}{e^{tR/2}-e^{-tR/2}} \right)=:
 \sqrt{4\pi
t}^{-n/2}\widehat{A}(tR).
\end{eqnarray*}
\end{theorem}
\begin{proof}
The Mehler formula for the operator $H_0$ with $q=0$ is well known
and proved for example in \cite{BGV} Chapter 4. For a general
operator $H$ just note that $$H=e^{q({\bf y})/2}H_0 e^{-q({\bf
y})/2},$$ from which the result follows immediately.\qed
\end{proof}
\par
\smallskip To be able to transfer this result to the model
problems at a specific ``temporal front face'' we define the heat
space, on which the above solution lives.
\par
\begin{table}[h]
 \centerline{\epsfig{file=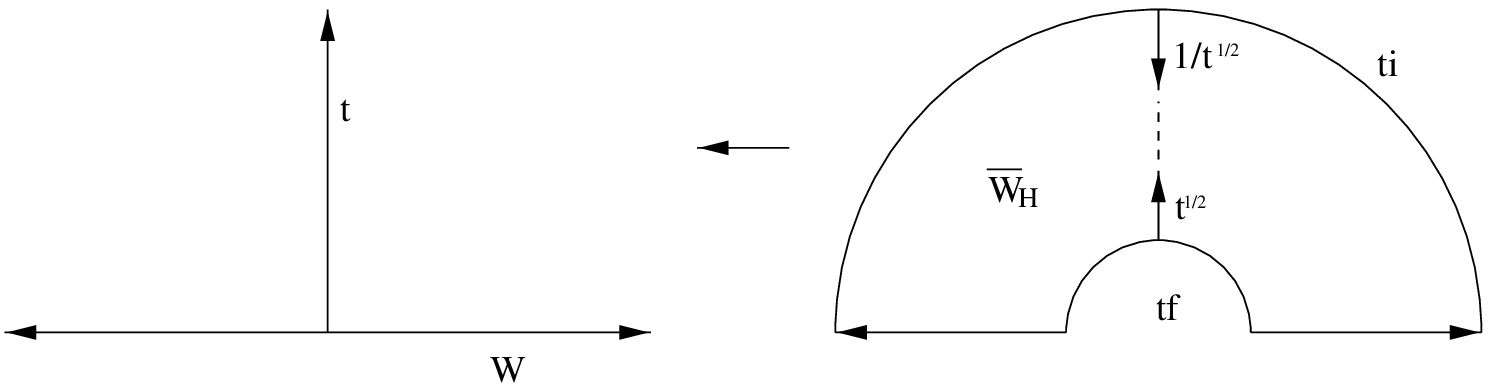}}
 \centerline{Figure 8: The Simple Heat Space $\overline{W}_H$}
 \end{table}
 Denote by
$$W_H=[W\times[0,\infty[_t;\{0,0\},\dd t]$$
 the ``simple heat
space'' associated to $W$ and by $\overline{W}_H$ its radial
compactification at spatial and temporal infinity. As before, the
blown up face at $t\rightarrow 0$ will be denoted by $\fatf$. The
face at infinity will be called ${\rm ti}$. Away from $t=0$ the
coordinate $t^{1/2}$ is a defining function for $\fatf$, and away
from spatial infinity the coordinate $t^{-1/2}$ is a defining
function for ${\rm ti}$. We can identify the interior of {\em
both} faces $\fatf$ and ${\rm ti}$ with $W$ such that for a
constant vector field $V$ in $TW$ the restriction
$$t^{1/2}V|_{\fatf}\quad \mbox{corresponds to}\quad V\in T
\fatf\quad\mbox{and}\quad t^{1/2}V|_{\rm ti}\quad
\mbox{corresponds to}\quad V\in T{\rm ti}.$$ Equivalently the
linear coordinate ${\bf Y}$ on $\fatf$ or ${\rm ti}$ corresponds
to ${\bf y}/t^{1/2}$ on the interior $W\times]0,\infty[$.
\par For
$l,p>0$ the associated ``calculus'' is defined by $$\Psi_{H,{\rm
s} ,{\rm cl}}^{l,p}(W):=\rho_{\fatf}^l\rho_{\rm
ti}^p\cdotinfty{\fatf,{\rm ti }}(\overline{W}_H, t^{-n/2}\frac{\dd
t}{t}\Omega(W)). $$
 This
will be a recipient space for  solutions of the heat equation.
However these solutions cannot necessarily be composed with one
another -- at least we will not say how. Nevertheless  we have
normal maps $$0\longrightarrow \Psi_{H,{\rm s},{\rm
cl}}^{l-1,p}(W) \longrightarrow \Psi_{H,{\rm s},{\rm cl}}^{l,p}(W)
\stackrel{N'_{H,\fatf,l}}{\longrightarrow}{\cal S}(W,\Omega(W)\dd
t)\rightarrow 0,$$ $$0\longrightarrow \Psi_{H,{\rm s},{\rm
cl}}^{l,p-1}(W) \longrightarrow \Psi_{H,{\rm s},{\rm cl}}^{l,p}(W)
\stackrel{N'_{H,{\rm ti},p}}{\longrightarrow}{\cal
S}(W,\Omega(W)\dd t)\rightarrow 0,$$ given by
$N'_{H,\fatf,l}(A):=t^{1-l/2}A|_{\fatf}$ and
$N'_{H,\fatf,p}(A):=t^{p/2}A|_{\rm ti}$.
\par
\medskip Also, a rescaled version of the simple heat ``calculus''
with coefficients in an algebra ${\cal A}$ can be defined.
Assuming that the above algebra is of the form ${\cal
A}=\Lambda^{ev} U^*$ for some vector space $U$
 we can introduce the heat calculus $\Psi_{G,{\rm s},{\rm cl}}(W)$ with coefficients in
$\Lambda U^*$ and rescaled w.r.t the corresponding filtration
(which already is a grading) at $\fatf$. Using $\dd t^{1/2}$ to
trivialize the conormal bundle $N^*\fatf$ we get the rescaled
normal operator
 $$0\longrightarrow \Psi_{G,{\rm s},{\rm cl}}^{l-1,p}(W)
\longrightarrow \Psi_{G,{\rm s},{\rm cl}}^{l,p}(W)
\stackrel{N'_{G,\fatf,l}}{\longrightarrow}{\cal
S}(W,\Lambda^*U\otimes\Omega(W)\dd t)\longrightarrow 0.$$ We then
have
\begin{proposition}\label{mehlerresult}
\begin{enumerate}
\item The operator $H$ maps
 $$\Psi_{G,{\rm s},{\rm cl}}^{l-2,p}(W)\stackrel{H}{\longrightarrow}
\Psi_{G,{\rm s},{\rm cl}}^{l,p}(W)$$ with rescaled normal operator
$N_{G,\fatf,l-2}(H\circ A)=H N_{G,\fatf,l}(A)$.
\item The heat kernel for $H$ is
$${\sf K}_H({\bf y},t):=p(R,{\bf y},t)\exp(q({\bf y})/2)\dd
t\dvol_W \in \Psi^{2,0}_{G,{\rm s },{\rm cl}}(W)$$ with normal
operators at $\fatf$ and ${\rm ti}$ given by $$ N_{G,\fatf,2}({\sf
K}_H)({\bf Y})=p(R,{\bf Y},1)\exp(q({\bf Y})/2)\dd t\dvol_W,$$ $$
N_{G,{\rm ti},0}({\sf K}_H)({\bf Y})=
\sqrt{4\pi}^{-n/2}\exp(-\|{\bf Y}\|^2/4) \dd t\dvol_W.$$
\end{enumerate}
\end{proposition}
\begin{proof}
This follows from Mehler's formula once one has shown
that $N_{G,\fatf,l-2}(H\circ A)=H N_{G,\fatf,l}(A)$.
This is done as in Section \ref{secrescnormal} by
 calculating the normal action of the expressions appearing in
 (\ref{operatorh}),
 $g_W^{-1}$, $\dd$,
$g_W(R)$ and $q$, using their mapping properties as described at
the beginning of this Appendix.
\qed
\end{proof}
\par
\smallskip
We will also use the obvious generalization of these results to
families, where $W$ is really the fibre of some vector bundle, and
apply it to the operators $H={\cal H}_X, {\cal H}_B$ in Section
\ref{rescmodprob}.

\subsection{Heat Kernels of Vertical Families}\label{vertfam}
In this Appendix we introduce the heat calculus for vertical
families of operators over the boundary fibration $M\rightarrow
Y$. Define the heat space for vertical families as the parabolic
blow up $$M^2_{H,{\rm
fib}}:=[M\times_YM\times[0,\infty[_{\tau};\bigtriangleup_M\times\{0\},\dd
\tau ].$$ The interior of the blown up space
$\fatf\sowie\fatf(M^2_{H,{\rm fib}})$ can be identified with the
vector bundle $VM\rightarrow M$ in the usual manner. Also, we
introduce the space $\overline{M}^2_{H,{\rm fib}}$, compactified
at temporal infinity by introducing $1/\tau^{1/2}$ as a
coordinate. The corresponding face will be denoted by ${\rm ti}$.
Define the density bundle $$\KD_{M/Y,H}:=\tau^{-v/2}\frac{\dd
\tau}{\tau}\beta_L^*\dvol_{M/Y}.$$ For $l>0$ the associated family
heat calculus is then defined by $$\Psi_{H,{\rm fib },{\rm cl
}}^{l,p}(M/Y):=\rho_{\fatf}^l\rho_{\rm
ti}^{p-v}\cdotinfty{\fatf,{\rm ti}}(\overline{M}^2_{H,{\rm fib
}},\KD_{M/Y,H}).$$
 In the case $l=0$ we have to add the mean value
condition as in Section \ref{dheatbasics}: $$\Psi_{H,{\rm fib
},{\rm cl }}^{0,p}(M/Y):=\cdotinfty{\fatf}(\overline{M}^2_{H,{\rm
fib }},\KD_{M/Y,H})^{\fatf}\otimes \CC.$$ Leaving aside the
behavior at temporal infinity for the moment, it can be shown as
before that the composition formula holds
\begin{theorem}   \label{fibrecomposition}
Composition gives a map $ \Psi^{l}_{H,{\rm fib},{\rm cl}}(M)
\times \Psi^{r}_{H,{\rm fib},{\rm cl}}(M)
\stackrel{\circ}{\longrightarrow} \Psi^{l+r}_{H,{\rm fib},{\rm
cl}}(M).$ \qed
\end{theorem}
\par
It is also straightforward to describe the rescaled version of
these constructions. The coefficient bundle from Section
\ref{basicresc} is $$\Lambda(\fatff)
=\Lambda\phi^*{}^bT^*B\otimes\END_{\Cl(\phi^*{}^bT^*B)}(E),$$
where we have used $\dd x$ to trivialize $N^*\fatff$.
 Writing $\nu_{\fatff}$ for the normal vector $P_{\fatff}/\dd x$ the
 connection on $\Lambda(\fatff)$  for $T\in \beta_L^*TM/Y$ is then given  by
$$\nnabla_T^{\Lambda(\fatff)}=\nnabla_T-{\mathbb
K}(T,\nu_{\fatff})-\frac{1}{2}\nabla_{\nu_{\fatff}}{\mathbb K}
(T,\nu_{\fatff}):\Gamma(M^2_{H,{\rm
fib}},\Lambda(\fatff))\rightarrow \Gamma(M^2_{H,{\rm fib
}},\Lambda(\fatff)).$$ It is calculated more explicitly in
Corollary \ref{explicitconnection}(b). As before, the filtration
of $\Lambda(\fatff)$ over the front face $\fatf$ is given by the
total Clifford filtration. Using this data, we can define the
rescaled bundle ${\rm Gr}_{\fatf}(\Lambda(\fatff))$ as in Section
\ref{basicresc}. Restriction to $\fatf$ is defined using the
connection $\nabla^{\Lambda(\fatff)}$: $$|_{{\rm Gr},\fatf}:
\Gamma({\rm Gr}_{\fatf}(\Lambda(\fatff)))\longrightarrow
C^{\infty} (VM,{\cal N}\Lambda(\fatf,\fatff)).$$
 For
lack of (obvious) better alternatives we have used the notation on
the RHS to denote
 the bundle ${\cal N}{\cal N }\Lambda(\fatf,\fatff)$ in
(\ref{cornerbundle}) with factors $(N^*\fatff)^l$ trivialized by
$\dd x$, i.e. we have written $${\cal
N}\Lambda(\fatf,\fatff)=\mbox{$\bigoplus$}_{k}N^*\fatf(M^2_{H,{
\rm fib }})^k
\otimes\Lambda^k({}^dT^*X)\otimes\END_{\Cl({}^dT^*X)}(E).$$ The
normal operator
\begin{eqnarray*}
N'_{G,\fatf,l>0}(A)&:=&
  \tau^{1-l/2}A|_{\fatf,{\rm Gr}}\in
{\cal S}_{\rm fib} (VM,\Omega(VM/M)\dd \tau\otimes {\cal
N}\Lambda(\fatf,\fatff))
\\ N'_{G,\fatf,0}(A\oplus c)&:=&
  \tau A|_{\fatf,{\rm Gr}}\oplus c \in{\cal S}_{\rm fib}
  (VM,\Omega(VM/M)\dd \tau\otimes {\cal N}
  \Lambda(\fatf,\fatff))^{\fatf} \oplus \CC
\end{eqnarray*}
fits into the short exact sequence for $l\geq 1$:
  $$ 0\rightarrow \Psi_{G,{\rm fib},{\rm cl }}^{l-1}
  (M/ Y)\rightarrow\Psi_{G,{\rm fib },{\rm cl}}^l(M/Y)
  \stackrel{N_{H,\fatf ,l}}{\longrightarrow}{\cal S}_{\rm
fib} (VM,\Omega(VM/M)\dd \tau\otimes {\cal
N}\Lambda(\fatf,\fatff)) \rightarrow 0, $$
 and similarly for $l=1$. The notation $N'_G$ is used to remind
 the reader that these normal operators are defined using the
 temporal variable $\tau$ instead of $t$.
\par
\medskip Our aim in the rest of this Appendix will be to construct
the heat kernel of the operator
\begin{equation}\label{bsquare}
{\bf B}^2:=
g_{VM/M}^{-1}(\nnabla^{\Lambda(\fatff)},\nnabla^{\Lambda(\fatff)})+\frac{\kappa_{M/Y}}{4}
+{\sf m}_d(\delta_h(x^2F^{E/S,d}(\cdot,\cdot))
\end{equation} as an element in $\Psi_{G,{\rm
fib},{\rm cl}}^{2,0}(M/Y)$.
 The following Proposition calculates the
normal operator of $\nnabla^{\Lambda(\fatff)}$ at $\fatf$.
\begin{proposition}\label{normcalcrescapp}
Let $T\in \beta_L^*T\,M/Y$ as before. Then the maps $$
\Psi_{G,{\rm cl}}^{l}(M)\stackrel{
\nnabla^{\Lambda(\fatff)}_{\tau^{1/2}T}}{\longrightarrow}
   \Psi_{G,{\rm cl}}^{l}(X),\qquad \Psi_{G,{\rm cl}}^{l}(M)
   \stackrel{\frac{\partial}{\partial \tau}}{\longrightarrow}
   \Psi_{G,{\rm cl}}^{l-2}(M)$$
   have the following general form when restricted to the face $\fatf$
\begin{enumerate}
\item $N'_{G,\fatf,l-1}(\nnabla^{\Lambda(\fatff)}_{\tau^{1/2}T} A)=
(\nnabla^{\Lambda(\fatff)}_{\tau^{1/2}T}-\frac{1}{2}P_{\fatf}
\nabla^{\Lambda(\fatff)}_{\cdot} {\mathbb
K}^{\Lambda(\fatff)}(\tau^{1/2}T,\cdot))
 N'_{G,\fatf,l}(A)$
 \item  $N'_{G,\fatf, l-2}(\frac{\partial}{\partial \tau}A)
 =\left\{
 \begin{array}{ll}
     -\frac{1}{2}( L_{R^{VM}}+2-l-{\sf N} )N'_{G,\fatf,l}(A) & \mbox{for}\quad l>2 \\
     -\frac{1}{2}( L_{R^{VM}}+2-l -{\sf N})N'_{G,\fatf,2}(A)\oplus
    -\int_{VM/M }N'_{H,\fatf,2}(A) & \mbox{for}\quad l=2\\
 \end{array}   \right.$
\end{enumerate}
Here ${\sf N}$ denotes the number operator in $\Lambda(\fatf)$.
Note also that $\nnabla^{N^*\fatf}_{R^{\phindx}}\dd
\tau^{1/2}=0$.\qed
\end{proposition}
As in Section \ref{secrescnormal} the form of the normal operator
of the connection in (a) can be calculated more explicitly:
\begin{proposition} Write $\nu_{\fatf}:=P_{\fatf}/\dd \tau^{1/2}$.
Then $$P_{\fatf} \nabla^{\Lambda(\fatff)}_{\cdot} {\mathbb
K}^{\Lambda(\fatff)}(\tau^{1/2}T,\cdot)|_{U\in VM}=\frac{1}{2}(\dd
\tau^{1/2})^2\ext(\delta_h  \langle
x^2R^{d,V}(\cdot,\cdot)[\nu_\fatf+\beta_L^*U],
\beta_L^*T\rangle_{\phi}).$$
\end{proposition}
\begin{proof}
The proof proceeds just as in Proposition \ref{curvterms}.
\begin{eqnarray*}
&&P_{\fatf} \nabla^{\Lambda(\fatff)}_{\cdot} {\mathbb
K}^{\Lambda(\fatff)}(\tau^{1/2}T,\cdot)|_{U\in VM}=(\dd
\tau^{1/2})^2{\mathbb
K}^{\Lambda(\fatff)}(\beta_L^*T,\beta_L^*U+\nu_{\fatf})|_{\fatf}
\\ &&\qquad\qquad =(\dd
\tau^{1/2})^2({\mathbb
K}(\beta_L^*T,\beta_L^*U+\nu_{\fatf}))|_{\fatf,{\rm Gr}} =(\dd
\tau^{1/2})^2\delta_h({\mathbb
K}(\beta_L^*T,\beta_L^*U+\nu_{\fatf}))|_{\fatf}
\\ &&\qquad\qquad=\frac{1}{2}(\dd
\tau^{1/2})^2\ext (\delta_h  \langle
x^2R^{d,V}(\cdot,\cdot)[\nu_\fatf+\beta_L^*U],
\beta_L^*T\rangle_{\phi} )|_{\fatf}
\end{eqnarray*}
 Here we have used  Corollary \ref{easyrestrict} and Proposition \ref{inducedtff}(c).\qed
\end{proof}
Multiplication with the function $\tau\kappa_{M/Y}$ does not
contribute anything at $\fatf$, since the function vanishes there.
Finally, the partial Clifford action $${\sf m}_d(\alpha)=
\ext_d({\sf n}\alpha)+\cl_d({\sf v}\alpha),\quad \alpha\in
{}^dT^*X,$$ becomes exterior multiplication at $\fatf$:
\begin{proposition}\label{normcalcrescappcliff}
The map $$ \Psi_{G,{\rm fib},{\rm cl}}^{l}(M/Y)\stackrel{ {\sf m
}_d(\alpha)}{\longrightarrow}
   \Psi_{G,{\rm fib},{\rm cl}}^{l-1}(M/Y)$$
has normal operator $N_{G,\fatf,l-1}({\sf
m}_d(\alpha)A)=\ext_d(\alpha)N_{G,\fatf,l}(A)$.\qed
\end{proposition}
Using the methods described in the body of the text it is now easy
to show the following
\begin{theorem}
The heat kernel for ${\bf B}^2$ is an element $${\sf K}_{{\bf
B}^2}\in \Psi_{G,{\rm fib}, {\rm cl}}^{2,0}(M/Y).$$ It has normal
operators at $\fatf$ and ${\rm ti}$ given by
\begin{enumerate}
\item $N'_{G,\fatf,2}({\sf K}_{{\bf B}^2})({\bf z})=$
\par\quad $\delta_h\left[p({\bf z},x^2R^{d,V},1)\exp(\frac{1}{2}
\langle x^2R^{d,V}\nu_{\fatf},{\bf z}
\rangle_{VM/M})\exp(-x^2F^{E/S,d})\right]\dvol_{VM/M}\dd\tau$
\item $N'_{G,{\rm ti},0}({\sf K}_{{\bf B}^2})=[\Pi_{\circ}]\dd \tau$, the
projection onto the null space of $\Dir^{\phi,V}$
\end{enumerate}
We emphasize once more that all these normal operators commute
with $\ext(\frac{\dd x}{x})$.
\end{theorem}
\begin{proof}
The first part is obtained using the specialization of the methods
in Chapters \ref{heatchapter} and \ref{rescaledchapter} based on
the series of Propositions above. For part (b) we just mention
that $${\bf B}^2=(\Dir^{\phi,V})^2+\quad\mbox{a vertical operator
of order $\geq$ 1 in $\Lambda\,\phi^*{}^bT^*B$ }.$$ This is just
of the type used in the Berline-Vergne theorem (Theorem 9.19 in
\cite{BGV}).
 \qed
\end{proof}
\par
\smallskip
Denote by $\str_{E/B}$ the relative supertrace in
$\End_{\Cl({}^bT^*B)}(E)$ as defined in (\ref{tracetrick}),
and by ${\sf N}_B$ the number operator
on $\Lambda\,\phi^*{}^bT^*B$. Also, write ${\sf T}$ for the time
axis $[0,\infty]_{\tau^{1/2}}$, compactified by $\tau^{-1/2}$ at
infinity. Then
\begin{corollary}
\begin{enumerate}
\item $\tau^{-{\sf N}_B/2}\,\str_{E/B}({\sf K}_{{\bf B}^2}(\tau)|_{\bigtriangleup_M})\in
C^{\infty}(M\times {\sf T},\Lambda\,\phi^*{}^bT^*B\otimes
\Omega(VM/M)\dd \tau).$
\item $\Int(x\frac{\partial}{\partial x})\tau^{-{\sf N}_B/2}
\str_{E/B}({\sf K}_{{\bf
B}^2}(\tau)|_{\bigtriangleup_M})$
 vanishes in $\tau^{1/2}=0$ and $\tau^{-1/2}=0$.\qed
\end{enumerate}
\end{corollary}
This Corollary allows us to define the (``twisted'') family eta
invariant as the form on the base $Y$ given by
\begin{equation}\label{deffibredeta}
\widehat{\eta}(\Dir^{\phi,V},E):=\int_0^\infty\int_{M/Y}\Int(x\frac{\partial}{\partial
x})\tau^{-{\sf N}_B/2} \str_{E/B}({\sf
K}_{{\bf B}^2}(\tau)|_{\bigtriangleup_M}).
\end{equation}
This is an even (resp. odd) form on $Y$, whenever $v$, $h+1$ are
odd (resp. even).
\begin{remark}\rm  By Remark \ref{oldfibredetanew} this
is just the classical family eta invariant, whenever $B_{\phi}=0$.
Note however that in the literature
(\cite{BC}, \cite{BC2}, \cite{Daithesis}, (\ref{tracetrick}))
authors hide of $2$, $\pi$, and $i$ in
their
definitions of the family eta invariant in many different ways.
\end{remark}

\subsection{Conormal Functions, Pullbacks and Pushforwards}
\label{appzeronotation}
In this Appendix we explain some of the essential notions
introduced in \cite{Meaomwc}, see also \cite{Melaps} or
\cite{Paul}, and adapt them to our needs.
\par
\smallskip
 Let $N$ be a manifold
with corners. Let ${\cal M}(N)=\{H_1,\ldots, H_f\}$ be the set of
boundary faces of $N$ and choose
 a boundary defining function $\rho_{H_j}\equiv\rho_j$ for
each face $H_j$. As usual, denote by ${}^b{\cal V}(N)$ the space
of vector fields on $N$ tangent to the boundaries.
\par
\smallskip
The space of functions conormal to the faces of $N$ is defined as
$$ {\cal A}(N):= L^{\infty}H^{\infty}_b(N)
   = \{ u\in L^{\infty}(N) \quad | \quad {}^b{\cal V}(N)\cdot u
     \subset L^{\infty}(N) \}.
$$ Elements in ${\cal A}(N)$ are $C^{\infty}$ in the interior and
bounded, but control at the boundary is only very moderate. One
might think of such functions as being ``$C^{\infty}$ along the
boundary''.
 Choosing
a f-tuple ${\bf a}:=((a_1,m_1),\ldots,(a_f,m_f))$ with a pair
$(a_j,m_j)\in\CC\times\NN$ for each face $H_j$ of $N$, we can also
define the space of conormal functions on $N$ with {\em conormal
bounds} ${\bf a}$ as
\begin{equation}\label{conbound}
{\cal A}^{\bf a}(N):=\rho_1^{a_1}\log(\rho_1)^{m_1}\cdots
\rho_f^{a_f}\log(\rho_f)^{m_f}\cdot {\cal A}(N).
\end{equation}
\par
A special class of conormal functions is given by functions which
possess  ``generalized Taylor expansions''
into powers of the type
\begin{equation}\label{powers}
 \rho_j^z\log(\rho_j)^{m}, \qquad (z,m)\in \CC\times\NN
\end{equation}
at the boundary faces $H_j$ of $N$. To describe these expansions,
we list the admissible powers (\ref{powers}) in an ``index set''.
To explain this concept, introduce a partial ordering on
$\CC\times\NN$, describing the comparative vanishing order at
$H_j$ of expressions of type (\ref{powers}), by declaring $$
(z_1,m_1) \leq (z_2,m_2) \quad \Longleftrightarrow \quad
\left[\Re(z_1)< \Re(z_2) \quad\mbox{or}\quad \Re(z_1)= \Re(z_2)
\quad\mbox{and}\quad m_1 \geq m_2\right] .$$ Now, a
($C^{\infty}$-){\em index set} is a countable subset $I\subset
\CC\times \NN$ such that
\begin{itemize}
\item $I$ is bounded from below in the above ordering, i.e. there
is a ``worst'' power in the expansion,
\item $I\cap B_N(0)\times \{0,\ldots,N \}$ is finite for any $N\in \NN$,
i.e. powers in the expansion ``improve'' steadily,
\item $(z,m)\in I$ implies $(z+n,m-l)\in I$ for any $l,n\in \NN$ and $l\leq m$. This
is the $C^{\infty}$-condition.
\end{itemize}
Often, a complex number $z$ will be interpreted as the index
element $(z,0)$.
\par\smallskip
Now given  a collection ${\cal I}=(I_{H_1},\ldots,I_{H_f})$
of ($C^{\infty}$-) index sets  associated to the faces of $N$,
we can recursively define the space of (classical) conormal functions
associated to ${\cal I}$ as follows: A section  $f\in{\cal A}^{\cal
I}(N)$ has an expansion at each face $H$ of the form
\begin{equation}\label{expansion}
f \sim \sum_{(z,m)\in I_H}\rho_{H}^z\log(\rho_H)^{m}
f_{H,z,m},\qquad f_{H,z,m}\in {\cal A}^{{\cal I}'}(N), \quad
I'_F=\left\{\begin{array}{ll}E_F & F\neq H  \\ \NN & F=H\end{array}.
\right.
\end{equation}
Obviously, for ${\bf a}=((a_1,m_1),\ldots,(a_f,m_f))$  as above,
we have $$ {\cal A}^{\cal I}(N)\subset {\cal A}^{\bf
a}(N),\quad\mbox{if} \quad (a_j,m_j)\leq I_{H_j}.$$ We leave it to
the reader to define spaces of conormal functions when one or
several index sets are finite. The largest or ``best'' index is
then to be interpreted as a conormal bound in the sense of
(\ref{conbound}), and the $C^{\infty}$-condition is assumed to
hold up to that best index.
\par
It will be important to keep track of the behavior of the index
sets of conormal functions w.r.t. different operations.
\begin{definition}
Let $I$, $J$ be two $C^{\infty}$-index sets. Besides the simple
set-theoretic union $I\cup J$ we will use the following other
operations on $I$ and $J$:
\begin{enumerate}
\item $I+J:=(\, (z+z',m+m')\,| (z,m)\in I, (z',m')\in J\,)$
\item $I \overline{\cup} J:=(\, (z,j)\,|\, j\leq m+m'+1,\, (z,m)\in I, (z,m')\in J\,)$
\item $\widehat{I}:=\overline{\bigcup}_{l\in\NN}(l+I)$
\end{enumerate}
\end{definition}
It should be clear that    $$ f+g \in {\cal A}^{{\cal I}\cup {\cal
J }}(N) \quad\mbox{and}\quad f\cdot g\in {\cal A}^{{\cal I}+{\cal
J }}(N)\quad\mbox{for}\quad f\in{\cal A}^{\cal I }(N),\quad
g\in{\cal A}^{\cal J}(N).$$ Especially, the space of conormal
functions ${\cal A}^{\cal I}(N)$ is  $C^{\infty}(N)$-module for
any $C^{\infty}$-index set ${\cal I}$.
\par
The reader should verify
that the coefficients $f_{H,z,m}$ in the expansion
(\ref{expansion}) of $f\in{\cal A}^{\cal I}(N)$ are not unique.
However, if we have $(z,m)\in I_H$ but $(z-n,m)\notin I_H$ for any
$n\in\NN_+$, then the coefficient $$f_{H_j,z,m}|_{H_j}\in {\cal
A}^{(I_{H_1},\ldots,I_{H_{j-1}},I_{H_{j+1}}, \ldots)}(H_j)$$ {\em
is} well defined. We therefore define the ``leading part'' of an
index set $I$ as $${\rm lead}( I):=\{(z,m)\in I\quad |\quad
(z-n,m)\notin I \quad\mbox{for any}\quad n\in\NN_+\}. $$ Note the
rule
 $$I-{\rm lead}(I)=1+I \quad !$$
Information associated to the leading coefficients of a conormal
function  is more persistent, in a sense that we describe now.
\par
\medskip
In the case of the Dirac operator $\Dir^{\phi}$ over the
 ``manifold with corners'' $N=X$ an important piece
of information is given by the decomposition of the space of
sections over the boundary into zero modes and nonzero modes: By
our main assumption (\ref{mainassumption})
 the null space of the vertical Dirac family $\Dir^{\phi,V}$ over the
boundary $\partial X$ is the space of sections of a vector bundle
${\cal K}\rightarrow Y$ and the projections $\Pi_{\circ}$,
$\Pi_{\perp}$ onto this null space and its orthogonal complement
give a decomposition $$ C^{\infty}(\partial
X,E)=\Pi_{\circ}C^{\infty}(\partial X,E)\oplus
\Pi_{\perp}C^{\infty}(\partial X,E)\sowie \Gamma(Y,{\cal
K})\oplus\Gamma(Y,{\cal K})^{\perp}.$$
\par
\smallskip
For a section $\xi \in {\cal A}^I(X,E)$ it now clearly makes sense
to ask, whether the {\em leading} coefficients $\xi_{(z,n)}$ in
its asymptotic expansion of type  (\ref{expansion}) restrict to
the zero modes over the boundary or not. To keep track of this
information, we agree to let an index $(z,m)^{\circ}$ denote a
coefficient in the zero modes and $(z,n)^{\perp}$ denote a
coefficient perpendicular to the zero modes: $$
 \xi_{(z,n)^{\circ}}|_{\partial X}\in{\cal K},\qquad
 \xi_{(z,n)^{\perp}}|_{\partial X}\perp{\cal K}.$$
More generally, given an index set $I$, we we define $[\circ]I$ to
be the same index set, but carrying  the additional information
that the leading coefficients lie in the zero modes: $$[\circ] I:=
{\rm lead}(I)^{\circ} \cup (I-{\rm lead}(I)),$$ and $[\perp]I$ is
defined analogously.
\par
Of course, these notations also make sense for the manifold
$N=X^2_{\phi}$. Here the null space of the vertical family
$\Dir^{\phi,V}$ defines vector bundles over the bases of $\falf$,
$\faphibf$, $\faff$ and we can use the above notation at each of
these faces. We can now describe the mapping properties of the
Dirac operator $\Dir^d$ (lifted from the left to $X^2_{\phi}$):
\begin{lemma}\label{mapping}
\begin{enumerate}
\item $\Dir^d:C^{\infty}(X,E)^{\circ}\rightarrow
C^{\infty}(X,E)$
\item $\Dir^d: {\cal A}^{[\circ] I}(X,E)\rightarrow
      {\cal A}^{I}(X,E)$.
\item $\Dir^d:C^{\infty}(X^2_{\phi},\END)^{\circ}
\rightarrow \rho_{\faff}^{-1}C^{\infty}(X^2_{\phi},\END)$
\item $ \Dir^d : {\cal A}^{ (I_{\farf} ,[\circ]I_{\falf},
[\circ]I_{\faphibf} ,I_{\faff}) } (X^2_{\phi},\END)
\rightarrow {\cal A}^{(I_{\farf},I_{\falf},
I_{\faphibf},I_{\faff}-1)}(X^2_{\phi},\END)$
\qed
\end{enumerate}
\end{lemma}
\par
\bigskip
\noindent{\bf b-Fibrations, Pullback and Pushforward}
\par
\medskip
Let $N$, $N'$ be two manifolds with corners. The corresponding
sets of boundary faces are ${\cal M}(N)$, ${\cal M}(N')$. A
$C^{\infty}$-map $F:N_1\rightarrow N_2$ will be called a {\em
$b$-map}, if it maps boundary faces to boundary faces in the sense
that for every face $H\in {\cal M}(N')$ we require
$$F^*\rho'_H=a_H\prod_{G\in {\cal
M}(N)}\rho_G^{e(G,H)},\quad\mbox{or}\quad F^*\rho'_H\equiv 0,$$
with $e(G,H)\in \NN$, $0< a_H\in C^{\infty}(N)$. If no
$F^*\rho'_H$ is identically $0$, $F$ will be called an {\em
interior} $b$-map. The matrix ${\bf e}=(e(G,H))$ will be called
the {\em coefficient matrix} of $F$. The readers can convince
themselves, that composition of two $b$-maps is again a $b$-map,
and that the coefficient matrix of the composition is the product
of the coefficient matrices of the components.
\par
The following result about the pullback of conormal functions
under $b$-maps will be used several times in the text
\begin{theorem}[Pullback]  \label{pullback}
Let $F:N\rightarrow N'$ be an interior $b$-map, with coefficient
matrix ${\bf e}$. Let ${\cal I}'$ a (possibly finite) index set
for $N'$. Then $$F^*:{\cal A}(N')^{{\cal I}'}\longrightarrow {\cal
A }(N)^{F^{\sharp}{\cal I}'},$$ where the index set
$F^{\sharp}{\cal I}'$ is defined as $$(F^{\sharp}{\cal I}')_G=(\,
(\sum_{H\in {\cal M}(N')}e(G,H)z'_H, \sum_{H\in {\cal M}(N')}m'_H)
\, | \, (z'_H,m'_H) \in I'_H \, ),$$ which we read as $\NN$, if
$G$ is mapped to the interior of $N'$ by $F$.\qed
\end{theorem}
\par
The corresponding result for the pushforward of functions is not
quite as easy to formulate. Also we will have to restrict the
class of admissible maps $F$ further. First, note that an interior
$b$-map $F:N\rightarrow N'$ maps $$F_*:{}^bT_pN\rightarrow
{}^bT_{F(p)}N'\quad\mbox{and}\quad F_*:{}^bN_p N\rightarrow
{}^bN_{F(p)}N'$$ for all $p\in N$. If the first of these ($b$-)
differentials is always surjective, $F$ will be called a {\em
$b$-submersion}. It will be called $b$-normal, if the second
($b$-) differential is always surjective. An interior $b$-map $F$,
which is both, a $b$-submersion and $b$-normal, will be called a
$b$-fibration. A more intuitive description of $b$-fibrations is
stated in \cite{Meaomwc}:
\begin{proposition}
Let $F:N\rightarrow N'$ be a $b$-submersion. Then $F$ is a $b$-
fibration {\em if and only if} for each face in the preimage $G\in
{\cal M}(N)$ there is {\em at most one} face in the image $H\in
{\cal M}(N')$, such that $e(G,H)\neq 0$.\qed
\end{proposition}
We can now formulate the pushforward theorem in the version given
 in \cite{Paul}:
\begin{theorem}[Pushforward]  \label{pushforward}
Let $F:N\rightarrow N'$ be a $b$-fibration, and let ${\cal I}$ be
an index set for $N$ such that $I_G>0$ for all faces $G$ of $N$
which are mapped to the interior of $N'$. Then $$F_*:{\cal
A}^{{\cal I}}(N,\Omega(N))\longrightarrow {\cal
A}^{F_{\sharp}{\cal I}}(N',\Omega(N')),$$ $$\mbox{with}\quad
(F_{\sharp}{\cal I})_H:=\overline{\bigcup} _{G\in{\cal
M}(N),F(G)\subset H}e(G,H)^{-1}I_G,$$ for all $H\in {\cal
M}(N')$.\qed
\end{theorem}

\end{document}